\documentstyle[11pt,amssymb,amsmath,amstext,amsfonts,xypic]{article}

\DeclareMathOperator{\lstl}{%
\hbox{\:\,\vbox to 0pt{\vss{\vbox to 0pt
{\vss\hbox to 0 pt {\Huge\hss$\circ$\hss}\vss}%
\vskip -14.2pt%
\vbox to 0pt {\vss\hbox to 0 pt {\hss$\ltimes$\hss}\vss}%
\vskip 3pt}}}\:\,}

\newcommand{\qed}{$\quad\square$}

\newtheorem{theorem}[equation]{Theorem}
\newtheorem{proposition}[equation]{Proposition}
\newtheorem{corollary}[equation]{Corollary}

\newtheorem{definition}[equation]{Definition}
\newtheorem{question}[equation]{Question}
\newtheorem{conjecture}[equation]{Conjecture}
\newtheorem{lemma}[equation]{Lemma}
\newtheorem{claim}[equation]{Claim}

\newcommand{\dis}{{\displaystyle}}

\newcommand{\beq}{\begin{equation}\label}

\newcommand{\iso}{{\;\;\stackrel{_\sim}{\longrightarrow}\;\;}}
\newcommand{\cd}{\!\cdot\!}
\newcommand{\bo}{{\mathbf{o}}}
\newcommand{\sminus}{\smallsetminus}
\newcommand{\bgm}{{\mathfrak{M}}}
\newcommand{\bbm}{{\mathbb{M}}}

\newcommand{\triv}{{\mathsf{triv}}}

\newcommand{\lst}{\#\Gamma}

\newcommand{\ve}{\varepsilon}
\newcommand{\vi}{${\sf {(i)}}\;$}
\newcommand{\vii}{${\sf {(ii)}}\;$}
\newcommand{\viii}{${\sf {(iii)}}\;$}
\newcommand{\vvi}{${\sf {(iv)}}\;$}
\newcommand{\CM}{{\mathsf{L}}}
\newcommand{\irreps}{{\mathsf{Irrep}}}
\newcommand{\sset}{\subset}
\newcommand{\GG}{{\mathbb{G}}}
\newcommand{\ga}{{{\boldsymbol{\Gamma}}_{\!n}}}
\newcommand{\gb}{{{\boldsymbol{\Gamma}}_{\!n-1}}}
\newcommand{\btimes}{{\,\mbox{$\bigotimes$}\,}}
\newcommand{\DD}{{\mathbb{D}}}
\newcommand{\Cat}{\,{\mathcal{C}}\!\mbox{\it at}\,}
\newcommand{\bn}{{\mathbf{n}}}
\newcommand{\BV}{{\mathbf{V}}}
\newcommand{\ann}{{\mathtt{Ann}}V}
\newcommand{\Ann}{{\mathtt{Ann}^{\,}}}
\newcommand{\uga}{{\mathcal{U}}\bigl(\pgg\bigr)}
\newcommand{\RQ}{{{\mathsf{R}}Q}}
\newcommand{\G}{\Gamma}
\newcommand{\Sym}{{\mathtt{Sym}^{\,}}}
\newcommand{\hb}{{\sf{H}}_{0,0}}
\newcommand{\SS}{{\sf{S}}}
\newcommand{\OO}{{\mathbb{O}}}

\newcommand{\so}{{\Huge{\mathfrak{s}\mathfrak{o}}}_{_{^3}}}
\newcommand{\id}{{{\mathtt {Id}}}}
\newcommand{\sq}{{$\enspace\square$}}
\newcommand{\into}{\,\,\hookrightarrow\,\,}
\newcommand{\too}{\,\,\longrightarrow\,\,}
\newcommand{\onto}{\,\,\twoheadrightarrow\,\,}

\newcommand{\ad}{{\mathtt{{ad}}}}
\newcommand{\Ad}{{\mathtt{{Ad}}^{\,}}}
\newcommand{\Lie}{{\mathtt{Lie}^{\,}}}
\newcommand{\Spec}{{\mathtt{Spec}^{\,}}}
\newcommand{\End}{{\mathtt{End}}}
\newcommand{\rk}{{\mathtt{rk}}}

\newcommand{\Hom}{{\mathtt{Hom}}}
\newcommand{\Aut}{{\mathtt{Aut}}}
\newcommand{\PAut}{{\mathtt{PAut}}}
\newcommand{\Ext}{{\mathtt{Ext}}}
\newcommand{\Ind}{{\mathtt{Ind}}}
\newcommand{\bplus}{\,\mbox{$\bigoplus$}\,}

\newcommand{\Tr}{{{\mathtt {Tr}}}}

\newcommand{\hr}{{\mathfrak{h}^{^{\mathsf{reg}}}}}
\newcommand{\gr}{{\mathfrak{g}^{^{_{\mathsf{rs}}}}}}
\newcommand{\dd}{{\mathcal{D}}}
\newcommand{\I}{{\mathcal{I}}}

\newcommand{\rep}{{{\mathsf {Rep}}}}

\newcommand{\rr}{{\mathcal{R}}}
\newcommand{\Vx}{{\mathcal{V}}}
\newcommand{\GL}{\operatorname{GL}}
\newcommand{\SL}{\operatorname{SL}}
\newcommand{\PGL}{\operatorname{PGL}}
\newcommand{\PG}{{\operatorname{PG}}}
\newcommand{\GGG}{{{\operatorname{G}}_{_{\Gamma\!,n}}}}
\newcommand{\GGp}{{{\operatorname{G}}'_{_{\Gamma\!,n}}}}
\newcommand{\pgg}{{\mathfrak{p}\mathfrak{g}}_{_{\Gamma\!,n}}}
\newcommand{\sgg}{{{\mathfrak{s}\mathfrak{g}}_{_{\Gamma\!,n}}}}
\newcommand{\ggn}{{{\mathfrak{g}}_{_{\Gamma\!,n}}}}
\newcommand{\pggn}{{\mathfrak{p}\mathfrak{g}}_{_{\Gamma\!,n}}}
\newcommand{\hh}{{\mathsf{H}}}
\newcommand{\fff}{{\mathcal{F}}}
\newcommand{\thh}{{\mathsf{\widetilde{H}}}}
\newcommand{\bhh}{{\mathbb{H}}}

\newcommand{\e}{{\mathbf{e}}}
\newcommand{\BA}{{\mathcal{A}}}
\newcommand{\BBA}{{\sf{H}_{_{\KK}}}}
\newcommand{\vv}{{\mathcal{V}}}
\newcommand{\JJ}{{{\mathcal{D}}{\mathtt{Ann}}}}
\newcommand{\FI}{{\mathfrak{I}}}
\newcommand{\CC}{{\sf{C}}}
\newcommand{\Weyl}{{\sf{A}}}
\newcommand{\ff}{{\mathtt{F}}}
\newcommand{\ka}{\kappa}
\newcommand{\ro}{{\sf{r}}}
\newcommand{\bcal}{{\mathcal{B}}}
\newcommand{\KK}{{{\mathcal{K}}}}
\newcommand{\II}{J}
\newcommand{\prh}{{\mathtt{pr}}_{\!_\h}}
\newcommand{\zz}{{\mathcal{Z}}}
\newcommand{\ZZ}{{\sf{Z}}}
\newcommand{\ZB}{{\sf{Z}}_{0,0}}
\newcommand{\mm}{{\mathcal{M}}}

\newcommand{\derv}{{\mathtt{Der}}_{\omega}}
\newcommand{\derw}{{\mathtt{Der}}_{_W}}

\newcommand{\Ug}{{{\cal{U}}\g}}

\newcommand{\Ker}{{\mathtt {Ker}}}
\newcommand{\Image}{{\mathtt {Im}}}

\newcommand{\pr}{{\mathtt {pr}}}
\newcommand{\lL}{\underround}
\newcommand{\rR}{{}}
\newcommand{\h}{{\mathfrak{h}}}
\newcommand{\bfrak}{{\mathfrak{b}}}
\newcommand{\afrak}{{\mathfrak{a}}}

\newcommand{\om}{\omega}
\newcommand{\A}{{\mathcal{A}}}
\newcommand{\Om}{\Omega}

\newcommand{\hreg}{{\mathfrak{h}^{^{_{\mathsf{reg}}}}}}

\newcommand{\mreg}{{\mathcal{M}^{^{_{\mathsf{reg}}}}}}

\newcommand{\eps}{{\epsilon}}

\def\C{{\mathbb{C}}}
\def\grd{{\mathtt{gr}}}

\def\rr{{\mathcal{R}}}
\def\Z{{\mathbb{Z}}}

\def\oo{{\mathcal O}}

\def\pp{{\mathfrak{p}}}

\def\LL{{\mathbb{L}}}
\def\HH{{{H\!H}}}
\def\HC{{{H\!C}_{\!_{^{\sf{ev}}}}}}

\def\LB{{\mathbf{L}}}

\def\hp{\hphantom{x}}
\def\Pf{{\it Proof}}
\def\CP{{\mathbb{P}}}
\def\O{{\sf O}}

\def\eu{{\mathsf{eu}}}
\def\reg{{\!}^{^{\mathsf{reg}}}}

\def\frakn{{\mathfrak{n}}}

\def\g{{\mathfrak{g}}}
\def\gln{{\mathfrak{g}\mathfrak{l}}_n}
\def\gl{{\mathfrak{g}\mathfrak{l}}}
\def\sln{{\mathfrak{s}\mathfrak{l}}_n}
\def\slt{{\mathfrak{s}\mathfrak{l}}_2}

\def\adg{{\ad^{\,}\g}}
\def\cc{{\mathcal C}}

\def\H{{\mathbf{H}}}

\def\FF{{\mathsf{F}}}

\def\bp{{\mathsf{p}}}
\def\ccirc{{{}_{^{\,^\circ}}}}

\makeatletter

\def\downroundfill{$\m@th \setbox\z@\hbox{$\braceld$}%
  \braceld\leaders\vrule height\ht\z@ depth\z@\hfill\bracerd$}
\def\uproundfill{$\m@th \setbox\z@\hbox{$\braceld$}%
  \bracelu\leaders\vrule height\ht\z@ depth\z@\hfill\braceru$}
\def\overround#1{\mathop{\vbox{\m@th\ialign{##\crcr\noalign{\kern3\p@}
      \downroundfill\crcr\noalign{\kern3\p@\nointerlineskip}
      $\hfil\displaystyle{#1}\hfil$\crcr}}}\limits}
\def\underround#1{\mathop{\vtop{\m@th\ialign{##\crcr
      $\hfil\displaystyle{#1}\hfil$\crcr\noalign{\kern3\p@\nointerlineskip}
      \uproundfill\crcr\noalign{\kern3\p@}}}}\limits}
\makeatother
\begin{document}
\setlength{\parindent}{6mm}

\thispagestyle{empty}

\centerline{\Large {\textbf{Symplectic reflection
algebras,\,  Calogero-Moser space,}}}
\vskip 3pt
\centerline{\Large {\textbf{and deformed Harish-Chandra
homomorphism}}}
\vskip 6pt
\centerline{\large {\sc {Pavel Etingof and Victor Ginzburg}}}
\setlength{\parskip}{-3pt plus 0pt minus 1pt}

{\nopagebreak {\tiny
\begin{abstract}
{\footnotesize {To any
finite group $\G\subset Sp(V)$ of automorphisms of
a  symplectic vector space $V$ we associate a new
multi-parameter deformation, $\hh_\ka,$ of the algebra
$\C[V]\#\G$, smash product of $\G$
with the polynomial algebra on $V$. The parameter $\ka$  runs over
 points of $\CP^r$, where $r=\,$number
of conjugacy classes of {\it symplectic reflections} in $\G$.
The algebra
$\hh_\ka$, called a {\it symplectic reflection algebra},
is  related to the
coordinate ring of a  Poisson deformation of the quotient singularity
$V/\G$. This leads to a symplectic analogue of McKay correspondence,
which is most complete in case of wreath-products. 
If $\G$ is the Weyl group of a root system in a vector space
$\h$ and
$V=\h\oplus\h^*$, then  the algebras $\hh_\ka$ are certain
`rational' degenerations of the  double affine Hecke algebra
introduced earlier by Cherednik.

Let $\G=S_n,$ the Weyl group of $\g=\gln.$
We construct a 1-parameter deformation of the Harish-Chandra
homomorphism
from $\dd(\g)^\g$, the algebra of invariant polynomial
differential operators on ${\frak {gl}}_n$, to
the algebra of $S_n$-invariant differential operators 
with { rational} coefficients on
the space
$\C^n$ of diagonal matrices.
The second order Laplacian on $\g$ goes, under the
 deformed homomorphism,
to the Calogero-Moser differential operator on $\C^n$,
 with  rational potential. 
Our crucial idea is  to
reinterpret  the deformed
Harish-Chandra homomorphism
as a homomorphism: $\dd(\g)^{\g}
\twoheadrightarrow \mbox{\it  spherical
subalgebra in }\hh_\ka,$
where  $\hh_\ka$ is the symplectic reflection algebra
associated to the group
$\G=S_n.$ 
This way, the deformed
Harish-Chandra homomorphism becomes nothing but a
description of the spherical subalgebra in terms 
of `quantum' Hamiltonian reduction.
  
In the  `classical' limit
$\ka\to \infty$,
our construction
 gives an 
 isomorphism between the  spherical subalgebra in $\hh_\infty$
and the coordinate ring of the Calogero-Moser space. We prove that
all simple $\hh_\infty$-modules have dimension $n!$, and are parametrised
by points of the Calogero-Moser space.
The family of these modules forms
 a distinguished vector bundle  on the  Calogero-Moser space,
 whose fibers carry 
the regular representation
of $S_n$. Moreover,
 we prove that
the algebra $\hh_\infty$ is isomorphic to the endomorphism
algebra
of that vector  bundle.}}
\end{abstract}}
{\centerline{\bf Table of Contents}
\vskip -5mm
$\hspace{20mm}$ {\footnotesize \parbox[t]{115mm}{\,

\hp${}_{}$\hp1.{ $\;\,\,$} {\tt Introduction} \newline
{\bf PART 1. Symplectic reflection  algebras}\newline
\hp2.{ $\;\,\,$} {\tt Koszul patterns and deformations}\newline
\hp3.{ $\;\,\,$} {\tt Representation theory
of the algebra $\hh_{0,c}$}\newline
\hp4.{ $\;\,\,$} {\tt The rational Cherednik algebra} \newline
\hp5.{ $\;\,\,$} {\tt Automorphisms and derivations
of the algebra $\hh_\ka$} \newline
{\bf PART 2. Harish-Chandra homomorphism}\newline
\hp6.{ $\;\,\,$} {\tt The radial part construction} \newline
\hp7.{ $\;\,\,$} {\tt Deformation of the Harish-Chandra homomorphism} 
\newline
\hp8.{ $\;\,\,$} {\tt Example: $\g={\mathfrak{s}\mathfrak{l}}_2$} \newline
\hp9.{ $\;\,\,$} {\tt The kernel of the Harish-Chandra homomorphism} \newline
10.{ $\;\,\,$} {\tt Proof of Theorem \ref{inject}} \newline
11.{ $\;\,\,$} {\tt Calogero-Moser space for wreath-products} \newline
{\bf APPENDICES.}\newline
12.{ $\;\,\,$} {\tt Appendix A: Almost commuting  matrices} \newline
13.{ $\;\,\,$} {\tt Appendix B: Geometric construction of $V_k$} \newline
14.{ $\;\,\,$} {\tt Appendix C: Small representations}\newline
15.{ $\;\,\,$} {\tt Appendix D: Deformations, Poisson brackets and 
                                               cohomology}\newline
16.{ $\;\,\,$} {\tt Appendix E: Some examples}\newline
17.{ $\;\,\,$} {\tt Appendix F: Open questions}\newline
18.{ $\;\,\,$} {\tt Appendix G: Results of Kostant and shift operators}\newline
}}}

}
\setlength{\parskip}{3pt plus 5pt minus 0pt}

\section{Introduction}\label{Quant}
\setcounter{equation}{0}

Given an action of a finite group $\G$ on a smooth
affine algebraic variety $X$, one can study the algebraic
structure of the orbifold
$X/\G$ through its coordinate ring $\C[X]^\G$, the ring of
$\G$-invariant regular functions on $X$. It is known, however,
that if the $\G$-action on $X$ is not free the (possibly singular)
space $X/\G$ does not necessarily reflect the features
of the  $\G$-action quite adequately.
For instance, one is led to introduce various "stringy"
topological invariants, cf.  [Ba], which do not coincide
with the ordinary topological invariants of $X/\G$,
but are rather related to certain "resolutions" of $X/\G$.
A standard approach to this problem, see e.g. [Co],
 is to replace
the commutative  coordinate ring $\C[X]^\G$ by the smash-product algebra
$\C[X]\#\G$. The two algebras are {\it Morita equivalent},
hence lead to equivalent theories, provided the  $\G$-action on
$X$ is free, but are not equivalent in general. Although the
algebra $\C[X]\#\G$ is {\it non-commutative}, it is believed
that the "right" geometry of the $\G$-action on $X$ can be read off
from the
"non-commutative algebraic geometry" of $\C[X]\#\G$. This is
the approach that we adopt in the present paper.

Specifically, we will study deformations of some interesting orbifolds
$X/\G$ through non-commutative deformations of the corresponding algebra
$\C[X]\#\G$. We will see that a systematic use of this ideology
leads, via Morita equivalence,
 to various  new results in conventional
"commutative" algebraic geometry which have been unaccessible
by purely "commutative" methods. 

Another  fascinating feature
of our approach is that genuine 
 deformations of $X/\G$ coming from 
"commutative" algebraic varieties
stand on an equal footing with   purely "quantum"
deformations that are not associated with any
 algebraic varieties whatsoever.
"Commutative" and  "quantum"
deformations only differ by the values of  
deformation parameters. Thus, "quantum" deformations may be
viewed as quantizations
of "commutative" ones, and
Kontsevich's quantization plays a role here.
 On the other hand, the algebras arising in such
"quantum" deformations are closely related to
the algebras of invariant polynomial differential operators
on semisimple Lie algebras and Quiver varieties.
Thus, on the "quantum" side, we enter the realm of
representation theory, which was our original motivation.
\medskip

\noindent {\bf {Symplectic reflection algebras.}}\quad
Let $V$ be a finite dimensional
vector space over $\C$, and $TV$ its tensor algebra.
Let $\G \subset \GL(V)$ be a finite group,
and $(TV)\#\G$ the smash-product of $TV$ with $\C\G$,
the group algebra of $\G$.
Given a
skew-symmetric $\C$-bilinear pairing $\kappa: V\times V \to \C\G$,
put
\begin{equation}\label{rel}
\hh_\kappa\; :=\; (TV\#\G)/I\langle x\otimes y - y\otimes x - \kappa(x,y) \in
T^2V\, \oplus\, \C\G\rangle_{\,x,y \in V}\,\,,
\end{equation}
where $I\langle\ldots\rangle $ stands for the two-sided ideal in
$(TV)\#\G$ generated by the indicated set. Thus, $\hh_\kappa$
is an associative algebra.

In the special  case $\kappa=0$, the situation
simplifies
 drastically since all the generators $x\in V$ commute with each other.
The resulting algebra  is therefore isomorphic
to the smash-product algebra:
$\hh_{\kappa=0} \,=\, (SV)\,\#\,\G.$
Note that this algebra has a natural grading obtained by placing
$\C\G$ in grade degree zero, and $V$
 in grade degree 1.

In general, if $\kappa\neq 0$, the defining relations, see (\ref{rel}),
become inhomogeneous, hence assigning
$\C\G$ grade degree zero, and $V$
 grade degree 1, one only gets an increasing {\it filtration}
$F_\bullet(\hh_\kappa)$ on the corresponding algebra $\hh_\kappa$,
such that, for any $i,j\ge 0,$
we have:\linebreak
 $F_i(\hh_\kappa)\cdot F_j(\hh_\kappa)\subset F_{i+j}(\hh_\kappa)$.
Let $\grd(\hh_\kappa)= \bigoplus_i\;F_i(\hh_\kappa)/F_{i-1}(\hh_\kappa)$
denote the associated graded algebra.
It is immediate from (\ref{rel}) that, for any $x,y\in
V\subset\grd(\hh_\kappa)$, in
 $\grd(\hh_\kappa)$ we have: $xy-yx=0$.
Therefore the tautological imbedding: $V\into \grd(\hh_\kappa)$
extends to a well-defined and surjective graded algebra
homomorphism: $\hh_{\kappa=0}=(SV)\#\G \too \grd(\hh_\kappa)$. We
say that {\it Poincar\'e-Birkhoff-Witt (PBW-)property} holds for
$\hh_\kappa$ if this morphism is an isomorphism.

     From now on we assume, in addition, that
$V$ is a
symplectic vector space with symplectic 2-form $\om$, and that
$\G \subset Sp(V)$.
The triple $(V,\om,\G)$ is said to be
{\bf indecomposable} if there is no $\om$-orthogonal
direct sum decomposition: $V=V_1 \oplus V_2$,
where $V_i$ are $\G$-stable proper symplectic vector subspaces
in $V$. 

\begin{definition}\label{symplectic_refl} 
An element $s\in \G$ is called a {\sl symplectic
reflection} if $\;\;\rk(\id-s)=2.$
\end{definition}

Let $S$ denote the set of symplectic reflections in $\G.$
The group $\G$ acts on $S$ by conjugation.
For each $s \in S$, write $\om_s$ for the (possibly degenerate)
skew-symmetric form on $V$ which coincides with $\om$
on $\Image(\id-s)$, and has $\Ker(\id-s)$ as the
 radical
(recall that there is an $\om$-orthogonal direct sum decomposition
$V= \Image(\id-s) \oplus \Ker(\id-s)$).
\medskip

\begin{theorem}\label{H_PBW} Assume that  $(V,\om,\G)$ is an
indecomposable
triple.
 Then, PBW-property holds for $\hh_\kappa$ if and only if
there exists a constant
$t \in \C$, and an $\Ad\G$-invariant function $c: S \to \C\,,\, s\mapsto
c_s,$
such that the pairing $\kappa$ has the form:
$$\kappa(x,y) = t \cdot \om(x,y)\cdot 1 + \sum\nolimits_{s\in S}\; c_s\cdot
\om_s(x,y)\cdot s\quad,\quad\forall x,y \in V\,.$$
\end{theorem}

\noindent
{\bf Remarks } \vi Our proof of the Theorem shows that
the  algebra  $\hh_\kappa$ associated to a
pairing  $\kappa$ of the type described by the formula
above always satisfies the
PBW-property, even if the triple $(V,\om,\G)$ is
decomposable.

\vii An analogue of Theorem \ref{H_PBW} holds in a more
general
setup, where $V$ is any (not necessarily symplectic) vector space,
and $\G\subset \GL(V)$ is a finite subgroup containing no 
 elements $s\in \G$ such that $\rk(\id-s)=1$. 
The proof of  Theorem \ref{H_PBW} given
in \S2 shows that, in this generality, 
PBW-property holds for $\hh_\kappa$ if and only if
$\kappa(x,y) =$\linebreak
$ \om(x,y)\cdot 1 + \sum\nolimits_{s\in S}\; 
\om_s(x,y)\cdot s,\,$
$\forall x,y \in V,\,$
where $\om$ is a (possibly {\it degenerate}) 
skew-symmetric $\G$-invariant form
on $V$, symplectic reflections are defined as in 
Definition \ref{symplectic_refl} (involving no reference to any symplectic
structure), and $\,\{\om_s\}_{s\in S}\,$ stands for an arbitrary
$\G$-equivariant collection
of skew-symmetric forms such that:
$\Ker(\id-s)\subset {\mathtt{Radical}}(\om_s)\,,\,\forall s\in S$.

Theorems \ref{Coh_Mac_intro}--\ref{proj_intro}
 below also hold in this generality.$\quad\lozenge$
\medskip

\noindent
Let $\CC$ denote the vector space of $\Ad\G$-invariant functions on $S$.
We have
\begin{equation}\label{number}
\dim\CC = \;\mbox{\it  number of $\Ad\G$-conjugacy classes in $S$}\,.
\end{equation}
     From now on, we assume that the pairing $\kappa$ 
in (\ref{rel}) has the form
described in Theorem \ref{H_PBW}, for a certain $t\in \C$ and $c\in\CC$.
The corresponding algebra
$\hh_{t,c}:=\hh_\kappa$ will be referred to as a {\it symplectic reflection
algebra}.
\smallskip

Let $\e= \frac{1}{|\G|}\sum_{g\in \G} g \in \C\G,$ be the
`symmetriser' idempotent, viewed as an element of  $\hh_{t,c}$.
The algebra $\e \hh_{t,c} \e\subset  \hh_{t,c}$ will be called the
{\it spherical subalgebra} in $\hh_{t,c}$.
For $t=0$ and $c=0$, we have:
$\,\e \hh_{_{0,0}}\e=\e\bigl(SV\#\G\bigr)\e\simeq
(SV)^\G,\,$ the algebra of $\G$-invariants. Furthermore, 
for any $(t,c)\in \C\oplus \CC,$ Theorem \ref{H_PBW} yields:
$\grd\bigl(\e \hh_{t,c} \e\bigr) \simeq (SV)^\G,$ see 
(\ref{grad}).
\medskip

The space $\hh_{t,c} \e$ has an obvious {\it left}
$\hh_{t,c}$-module structure, and {\it  right}
$\e\hh_{t,c} \e$-module structure.
Given an algebra $A$, and a left, resp. right,
$A$-module $M$, we set: $M^\vee=\Hom_{A\mbox{\tiny{-left}}}(M,A),\,$
resp. ${}^{\vee\!}M=\Hom_{A\mbox{\tiny{-right}}}(M,A)$.

\begin{theorem}\label{Coh_Mac_intro}\vi
 $\;\;\e\hh_{t,c} \e$ is a finitely generated Gorenstein
 algebra without zero~divisors.

\vii $\;\hh_{t,c} \e$ is a finitely generated,
Cohen-Macaulay right
$\e\hh_{t,c} \e$-module.

\viii
We have:
${}^{\vee\!}(\hh_{t,c} \e)= \e\hh_{t,c}$, and
$(\e\hh_{t,c})^\vee=\hh_{t,c} \e$.
In particular, $\hh_{t,c}\e$ is reflexive.

\vvi The left $\hh_{t,c}$-action induces an algebra isomorphism:
$\,\hh_{t,c} \iso \End_{_{\e\hh_{t,c} \e}}(\hh_{t,c} \e)$.
\end{theorem}

\begin{remark} For a definition of a (not necessarily commutative)
Gorenstein ring, resp. Cohen-Macaulay module, which,
in the  commutative case  reduces to the
standard definition as given, e.g. in [BBG], the reader
is referred to \S3 (proof of Thm. \ref{Coh_Mac_intro}).
\end{remark}\medskip

For any $t\neq 0$, let $\,\Weyl_t=TV/I\langle x\cdot y - y\cdot x - t\cdot
\om(x,y)\rangle_{ x,y\in V}\,$ be
the {\it Weyl algebra} of the symplectic vector space $(V,t\cdot\om)$;
the algebras $\Weyl_t$ are isomorphic to each other, for all  $t\neq 0$.
We observe that, for  $c=0$ (and any $t\neq 0$), we have:
 $\hh_{t,0} = \Weyl_t\#\G$.
 We will show (see Theorem \ref{deform_t}) that,
for any fixed $t\neq 0$, the
family $\,\{\hh_{t,c}\}_{c\in\CC}\,$ gives a  universal
deformation of the algebra $\Weyl_t\#\G$.
Furthermore, the spherical subalgebra $\e \hh_{t,0} \e$
is isomorphic to $\Weyl_t^{^{\!\G}}$, the subalgebra of $\G$-invariants in
the Weyl algebra, and
(for any fixed $t\neq 0$) the
family $\,\{\e \hh_{t,c} \e\}_{c\in\CC}\,$ gives a {\sf universal
deformation} of the algebra $\Weyl_t^{^{\!\G}}$.

\begin{theorem}\label{comm_conj}
For any $c\in\CC$, the algebra $\e \hh_{t,c} \e$
is commutative if and only if $t=0$.
\end{theorem}
\medskip

\noindent
{\bf Example: Kleinian singularity.}\quad Let
$\G\subset SL_2(\C)$ be 
a finite subgroup acting on $V=\C^2$.
The action on
$V\smallsetminus\{0\}$ is free, and
any element of $\G\smallsetminus\{1\}$
 is a symplectic reflection. This case has been
studied by Crawley-Boevey and
Holland [CBH].$\quad\lozenge$\medskip

\noindent
{\bf{Quasi-classical case: ${\mathbf{\ka=(0,c)}.}$}}\quad
Assume next that  $t=0$  and $c\neq 0$.
The algebra $\hh_{0,c}\,,\, c\neq 0,\,$ has some very interesting special features
making it quite different from other members of the family
$\,\{\hh_\ka\}_{\ka\in {\C\oplus {\sf{C}} }}\,$. First of all, the algebra
 $\hh_{0,c}$ has a large center, $
\ZZ_{0,c}=\ZZ(\hh_{0,c}),$
as opposed to the algebras $\hh_\ka$ for generic $\ka$.
In more detail, the standard filtration on $\hh_{0,c}$
induces a filtration on $\ZZ(\hh_{0,c})$, and  we will show that:
$$
\grd^{\,}\ZZ(\hh_{0,c})
\;\;=\;\;\ZZ(\hh_{0,0})\;\;=\;\;\ZZ(SV\#\G)
\;\;=\;\;(SV)^\G\,.
$$
Further, for any $c\in \CC$, we use the construction of Hayashi [Ha],
to make both the spherical subalgebra $\e\hh_{0,c}\e$
and the center $\ZZ_{0,c}$  Poisson algebras, which reduce to the
standard Poisson structure on $(SV)^\G$, if $c=0$.
Unlike the case $c=0$,  the center of
$\hh_{0,c}$ is harder to describe explicitly. To this end, we
generalize the classical Satake isomorphism for an affine Hecke
algebra (as proposed by Lusztig [Lu]), and show that the map:
$z\mapsto z\cdot\e$ establishes a  Poisson  algebra isomorphism:
$\ZZ_{0,c}\iso\e\hh_{0,c}\e$, which we call the {\it Satake
isomorphism}.

Let $\rr$ denote the coherent
sheaf on  $\Spec{\ZZ_{0,c}}$ corresponding to $\hh_{0,c} \e$,
viewed as a (finitely generated) $\ZZ_{0,c}$-module, that is:
$\Gamma(\Spec{\ZZ_{0,c}}\,,\,\rr)=\hh_{0,c} \e$.
We will view points of
$\Spec \ZZ_{0,c}$ as algebra homomorphisms $\chi: \ZZ_{0,c}\too\C$.
Then, for any point $\chi\in\Spec \ZZ_{0,c}$,
the geometric fiber of the sheaf $\rr$ at $\chi$
equals: $\hh_{0,c}\e\bigotimes_{_{\ZZ_{0,c}}}\chi$.
The left action of $\hh_{0,c}$ on $\hh_{0,c} \e$
gives an $\hh_{0,c}$-module structure on each geometric
fiber of $\rr$.

In the theorem below, we 
fix $U\subset \Spec\ZZ_{0,c}$, a Zariski-open affine subset,
and use the subscript `$U$' to denote `restriction to $U$'.
In particular we write $\ZZ_{_U} :=\C[U]\supset \ZZ_{0,c}$, and
$\hh_{_U}:= \ZZ_{_U} \bigotimes_{_{\ZZ_{0,c}}}\,\hh_{0,c},$
and also $\rr_{_U}:= \rr\big|_U$.

\begin{theorem}\label{proj_intro}\label{proj_prop_intro}
\label{struc} Let $U$ be an  affine Zariski-open 
subset contained in the smooth locus of $\Spec{\ZZ_{0,c}}$. Then we have:

 \vi The sheaf $\rr_{_U}$ is locally free, and the
algebra $\hh_{_U}$ is isomorphic to $\End_{_{\ZZ_{_U}}}\!(\rr_{_U})$.

 \vii
The algebra $\hh_{_U}$ is Morita equivalent to $\ZZ_{_U}$,
i.e. the categories of finitely generated $\hh_{_U}$-modules
and $\ZZ_{_U}$-modules
are equivalent.

\viii
 Any simple $\hh_{_U}$-module is isomorphic to one
 of the form: $\,\hh_{_U} \e\bigotimes_{_{\ZZ_{_U}}}{\Huge{\chi}},\,$
for a certain character $\chi: \ZZ_{_U}\to\C$. In particular,
the action of the center defines a bijection between the isomorphism classes
of simple $\hh_{_U}$-modules and points of $U$.

\vvi Any simple $\hh_{_U}$-module has dimension $|\G|$ and is isomorphic,
 as a $\G$-module,
to the regular representation of $\G$.

If $\Spec{\ZZ_{0,c}}$ is itself
smooth then claims \vi-\vvi hold for $U=\Spec{\ZZ_{0,c}}$, and 
$\hh_{_U}=\hh_{0,c}$.
\end{theorem}
\medskip

Following [Alv], on the group algebra $\C\G$,
define an increasing filtration $F_\bullet(\C\G)$
by letting $F_k(\C\G)\,,\, k\geq 0,$ be the $\C$-linear span
of the elements $g\in\G$ such that $\rk(\id-g)\leq k$.
This filtration is obviously
compatible with the algebra structure
on $\C\G$. Let $F_\bullet(\ZZ\G)$ denote the induced
filtration on $\ZZ\G$, the center of $\C\G$.
Write $\grd^F_\bullet(\ZZ\G)$ for the corresponding
associated graded algebra.

Further, write $\HH^\bullet(\A)$ for the Hochschild cohomology
of an associative algebra $\A$, and $H^\bullet(-)$ for the
ordinary singular
cohomology of a topological space (with complex
coefficients). Part (i) of the theorem below
is a generalization of [AFLS] and [Alv], while part (ii),
may be viewed as a `symplectic generalization' of 
{\bf McKay correspondence}, cf.~[Re].

\begin{theorem}\label{mckay} For any $c\in\CC$ we have:

\vi For all $t\in\C^*$, except possibly a countable set,
there is a natural graded algebra isomorphism:
$\HH^{\bullet}(\hh_{t,c})\simeq \grd^F_\bullet(\ZZ\G)$.

\vii If
 the variety
$\Spec{\ZZ_{0,c}}$ is smooth, then the Poisson structure on
$\ZZ_{0,c}$ makes\linebreak $\Spec{\ZZ_{0,c}}\,$  a symplectic
manifold; Furthermore, there is a
graded algebra isomorphism:\linebreak
$H^{\bullet}(\Spec\ZZ_{0,c}) \simeq \grd^F_\bullet(\ZZ\G)$.
\end{theorem}

\begin{corollary}\label{betti}\label{length}
 If $\Spec{\ZZ_{0,c}}$ is smooth, then
$\,\dim H^{^{{\sf{odd}}}}(\Spec\ZZ_{0,c}) =0,\,$
and for all $i\ge 0,$ we have: 
$\dim H^{2i}(\Spec\ZZ_{0,c}) =\bn(i),\,$ where
\[
\bn(i)\,:=
\;\mbox{\textsl number of conjugacy classes of
$\,g\in \G\,$ such that}\quad\rk(\id-g)= 2i\,.\qquad\square
\]
\end{corollary}

Recall the 
result of Batyrev [Ba, Thm. 4.8]
saying that, for any {\it crepant} (hence,
 symplectic, see e.g. [Ve]) resolution: $\widehat{V/\Gamma}\onto V/\Gamma$,
one has: $\dim H^{2i}(\widehat{V/\Gamma})=\bn(i)$.
It has been
conjectured by E. Vasserot [Va] (and also by the second author)
that, for a resolution $\widehat{V/\Gamma}$ as above,
there is a graded algebra isomorphism:
$H^{\bullet}(\widehat{V/\G}) \simeq \grd^F_\bullet(\ZZ\G)$.
This conjecture is related to  Theorem \ref{mckay} 
as follows. According to [Ka2], 
 if the symplectic  orbifold $V/\Gamma$
has  a smooth
 symplectic resolution of singularities,
$\widehat{V/\Gamma}$, then 
this resolution can be deformed to  a smooth affine
 symplectic manifold.
We expect that the latter coincides with the 
variety  $\Spec{\ZZ_{0,c}}$, for an appropriate
 value of the parameter $c\in\CC$,
in particular, $\Spec{\ZZ_{0,c}}$ is smooth for generic $c\in\CC$.
Our expectation is based on Proposition 17.4 (of Appendix F)
that would conjecturally imply that {\it any}
Poisson deformation of $V/\Gamma$,
in particular the one constructed by Kaledin,
is given by a certain variety $\Spec{\ZZ_{0,c}}$.
If our expectations are correct, then
any smooth symplectic resolution, $\widehat{V/\G},$
has the same cohomology algebra as $\Spec{\ZZ_{0,c}}$.
Thus,  Theorem 
\ref{mckay}  would imply  the conjecture in [Va].
In many concrete examples the space $\widehat{V/\G}$  has a natural
hyper-K\"ahler structure
such that $\Spec{\ZZ_{0,c}}$ is obtained from
 $\widehat{V/\G}$ by `rotating' the complex structure.
Hence the  two spaces are even diffeomorphic, and Theorem 
\ref{mckay} applies, cf. Question 17.2.
This is the case, for instance, for
$\G=S_n$ acting diagonally on
$V=\C^n\oplus\C^n$.
In that case, $\widehat{V/\G}=
{\mathsf{Hilb}}^n(\C^2)$ is the Hilbert scheme of
$n$ points on the plane, and 
the corresponding result on the algebra structure of 
$H^{\bullet}\bigl({\mathsf{Hilb}}^n(\C^2)\bigr)$
has been recently proved in [LSo] and [Va] (independently).
Their approach is totally different from ours.

Further, 
write $K(X)$ for the algebraic $K$-theory of an algebraic variety
$X$, that is, for the Grothendieck group
formed by
 algebraic vector bundles on $X$.
Now let $X=\Spec\ZZ_{0,c}$,
and assume that it is smooth
so that the sheaf $\rr$ is locally free, by
Theorem \ref{proj_intro}(i).
We view $\rr$ as an   algebraic vector bundle on $\Spec\ZZ_{0,c}$.
 Given an
irreducible $\G$-representation $E$,
we let $\rr^E:= E\otimes_{_\G}\rr=
\Hom_{_\G}(E^*, \rr)$ denote the $E$-isotypic component
of $\rr$, which is again
an   algebraic vector bundle  on $\Spec\ZZ_{0,c}$.

\begin{proposition}\label{k_theory} If
 the variety
$\Spec{\ZZ_{0,c}}$ is smooth, then 
$K(\Spec{\ZZ_{0,c}})$ is a free abelian group
with the basis $\,\{\rr^E\}_{_{E\in{\sf{Irreps}}(\G)}}.\,$
Moreover,
the Chern character map gives an isomorphism:
$\C\otimes_{_\Z}K(\Spec{\ZZ_{0,c}})
\iso H^{^{{\sf{ev}}}}(\Spec\ZZ_{0,c})$.
\end{proposition}
\medskip

\noindent 
{\bf Wreath-product case.}\quad An infinite series of
interesting examples of groups generated by symplectic 
reflections\footnote{
see \cite{Coh} for a classification of finite groups
generated by symplectic reflections.}
is provided by a so-called {\it wreath-product} construction.
Fix a finite subgroup $\G\subset SL_2(\C)$ and an integer $n>1$.
We write  $L=\C^2$ for the tautological
2-dimensional $\G$-module, form the vector space
$V= L^{\oplus n}$ with the symplectic structure
$\om_{_V}$ induced from the standard one on   $L=\C^2$, and let $S_n$,
the Symmetric group, act
on $V$  by permutations of the
direct summands.
The group
$\,\ga:=S_n\ltimes (\G\times \G\times\ldots\times \G)\subset Sp(V)$
acts naturally on $V= L^{\oplus n}$ and is called the
wreath-product of
$S_n$ with $\G$.

It is known that  $\CC(\ga)$,
the space of class functions on the
set of symplectic reflections in $\ga$, can be identified,
cf. (\ref{conj_class}),
with $\ZZ\G$, the center of the group algebra of $\G$.
Thus there is a symplectic reflection
algebra, $\hh_{0,c}(\ga),$ attached to the
triple $(V,\om_{_V},\ga),$ and any parameter
 $c\in \ZZ\G\simeq\CC(\ga)$.

On the other hand, let $\BV=\C\G^{\oplus n}$
be the direct sum of $n$ copies of the regular representation of $\G$.
Following Nakajima [Na1,Na2], to any $c\in \ZZ\G$ one associates an affine
algebraic variety of {\em quiver data}:
\begin{equation}\label{bbm}
\bbm^c_{_\G}(\BV) \;:=\;
\big\{(\nabla,I,J)\in \bigl(\Hom_{_\C}(\BV,\BV\otimes_{_\C}L)
\;\bplus\;\BV\;\bplus\;\BV^*\bigr)^\G\enspace
\Big|\enspace
[\nabla,\nabla] + I\otimes J = c|_{_\BV}\big\}\,,
\end{equation}
where $I\otimes J\in \End_{_\G}\BV,$ and
 $[\nabla,\nabla]\in (\End_{_\G}\!\BV)\,\,\btimes\bigwedge^2\!L\simeq
\End_{_\G}\BV$.
The natural action on
$\,\bigl(\Hom_{_\C}(\BV,\BV\otimes_\C L)
\,\bplus\,\BV\,\bplus\,\BV^*\bigr)^\G$ of the
group $G_{_\G}(\BV)$ of $\Gamma$-equivariant
automorphisms of $\BV$
preserves the subvariety
$\bbm^c_{_\G}(\BV)$.
The quotient: $\bgm_{_{\G\!,n,c}}= \bbm^c_{_\G}(\BV)/G_{_\G}(\BV)$
is a particular case of the Nakajima {\it quiver variety}.
In general, the quotient is understood in the sense of 
Geometric Invariant
Theory, but in our case the GIT quotient coincides
with the naive one because of the following  well-known
 result (which  is essentially Proposition
\ref{G_Wilson}), see also [Na1]
\begin{lemma}\label{Wilson2}
If $c\in \ZZ\G$ is generic enough, then
the induced $G_{_\G}(\BV)$-action on $\bbm^c_{_\G}(\BV)$
is  free.
Furthermore, $\bgm_{_{\G\!,n,c}}$ 
is a smooth  symplectic affine algebraic variety
of dimension  $2n$. 
\end{lemma}

Write $\ZZ_{0,c}(\ga)$ for the center of the algebra
$\hh_{0,c}(\ga)$.
One of our main results is
\begin{theorem}\label{wreath_thm}
If $c\in \ZZ\G$ is generic enough, then $\Spec\ZZ_{0,c}(\ga)$ is 
isomorphic to $\bgm_{_{\G\!,n,c}}$,
as a
Poisson  algebraic variety.
\end{theorem}
\begin{corollary}\label{wreath_cor}
If $c\in \ZZ\G$  is generic enough, then we have:

\vi $\;\;\Spec\ZZ_{0,c}(\ga)$ is a
smooth  symplectic variety.

\vii Simple $\hh_{0,c}(\ga)$-modules all have dimension
$\, n!\cdot|\G|^n,\,$ and  are parametrized
by points of the quiver variety $\bgm_{_{\G\!,n,c}}$.

\viii There is a canonical algebraic 
vector bundle $\rr$ on $\bgm_{_{\G\!,n,c}}$
whose fibers carry the regular representation of the group $\ga$.
The classes: $\,\{\rr^E\}_{_{E\in {\sf{Irrep}}(\ga)}}\,$
form a basis of $K(\bgm_{_{\G\!,n,c}})$.

\vvi There is an algebra isomorphism:
$\,H^\bullet(\bgm_{_{\G\!,n,c}})
\simeq \grd^F_\bullet\ZZ\ga.$\quad\qed
\end{corollary}\medskip

\noindent {\bf {Rational Cherednik algebra.}}\quad 
All indecomposable triples
$(V,\om,\G)$ fall into two large groups:
{\bf{(1)}} those where $V$ is an irreducible $\G$-module,
and {\bf{(2)}} those where $V=\h\oplus\h^*$ is a direct sum
of two  $\G$-stable irreducible Lagrangian
subspaces. In the second case, the action on $\h$ 
by any symplectic reflection $g\in\G$ has a pointwise fixed
codimension 1 hyperplane in $\h$, hence is a
{\it complex reflection} in $\h$. Thus, 
in case (2), saying that
 $\G$ is generated by symplectic reflections amounts
 to saying that $\G=W$ is 
a finite complex reflection group  in a vector space
$\h$, acting on $V= T^*\h \simeq\h + \h^*$ by
 induced symplectic automorphisms.
Any (complex) reflection $s\in W$ induces the corresponding
symplectic reflection of $\h + \h^*$.

In this paper, we will be mostly interested in the special case
where $W$ is the Weyl  group of  a  finite  reduced root system
$R\subset\h^*$. The set $S$ of reflections in $W$ forms either one
or two conjugacy classes, depending on whether all roots
$\alpha \in R$ are of the same length or not. Thus, giving
a $W$-invariant function $c: S\to \C$ amounts to
giving a map: $\alpha \mapsto c_\alpha\in \C$, where
$c_\alpha$ depends only on the length of $\alpha \in R$.

Let $\hh_{t,c}$ be the symplectic reflection algebra
associated to
a complex number $t\in \C$, a map $c$ as above,
and the vector space $V=\h + \h^*$.
Explicitly, write $\alpha^\vee\in \h$ for the coroot
 corresponding to a root $\alpha \in R$, and
$s_\alpha\in W$ for the reflection relative to $\alpha$.
Then, formula (\ref{rel}) shows that the algebra $\hh_{t,c}$
is generated by the spaces $\h,\h^*,$ and the group $W$,
subject to the following defining relations:
\begin{equation}\label{P-bracket}
\begin{array}{lll}\displaystyle
&{}_{_{\vphantom{x}}}w\cdot x\cdot w^{-1}= w(x)\;\;,\;\;
w\cdot y\cdot w^{-1}= w(y)\,,&
\forall y\in \h\,,\,x\in \h^*\,,\,w\in W\break\medskip\\
&{}^{^{\vphantom{x}}}{}_{_{\vphantom{x}}}[x_1,x_2] = 0 = [y_1,y_2]\,, &
\forall y_1,y_2\in \h,\;x_1\,,\,x_2 \in \h^*\,\break\medskip\\
&{}^{^{\vphantom{x}}}[y,x] = t \cdot \langle y,x\rangle
-\frac{1}{2}\cdot \sum_{\alpha\in R}\;
c_\alpha\cdot\langle y,\alpha\rangle
\langle\alpha^\vee,x\rangle \cdot s_\alpha\,,& \forall y\in
\h\,,\,x\in \h^*\,,
\end{array}
\end{equation}
(the factor $\frac{1}{2}$ on the RHS of the last
equation accounts to the fact that roots
$\alpha$ and $-\alpha$ give rise to the same reflection).

We will refer to the algebra  $\hh_{t,c}$ above as
 a {\it rational Cherednik algebra} because it is a
 certain `rational'
degeneration of the double affine Hecke algebra studied by Cherednik
\cite{Ch}.
In the case of type ${\mathbf{A_1}}$,
the algebra $\hh_\ka$ was discussed in detail in \cite{ChM}.

Olshanetsky-Perelomov [OP], have associated
a  Calogero-Moser type  integrable system to any
Weyl group $W$. We show that the space
 $\,\Spec \ZZ(\hh_{0,c})\,$ for the
rational Cherednik algebra attached to $W$
is a natural  completion of the phase space of that
 integrable system. The  affine algebraic
variety  $\,\Spec \ZZ(\hh_{0,c})\,$
 is always an irreducible normal Gorenstein  
 variety
with a Poisson structure.
We propose to call it
the {\it Calogero-Moser space associated to} the Weyl
group $W$, with parameter $c$.
For Weyl groups of type: 
${\mathbf{A_{n}, B_n,C_n}}$ and sufficiently general
parameters $c\in\CC$, the variety
$\Spec \ZZ(\hh_{0,c})$ is a {\it smooth} symplectic manifold,
by Theorem \ref{wreath_thm}.
\medskip

\noindent
{\bf{Relation to the double-affine Hecke algebra.}}\quad
The Cherednik algebra $\hh_{t,c}$ may be thought of as the following
two-step degeneration of the double affine Hecke algebra
$\H$ defined in [Ch]:
$$\H\;=\;\H^{^{\,_{\text{elliptic}}}}
\enspace\rightsquigarrow \enspace
\H^{^{\,_{\text{trigonometric}}}}\enspace\rightsquigarrow \enspace
\H^{^{\,_{\text{rational}}}}\;=\;\hh_{t,c}\,.
$$
In more detail, let $G$ be the Lie
 group corresponding to the Lie algebra $\g$, and
$T$  a maximal torus of $G$.
The algebra $\H$ contains the group algebras of two
(dual) lattices: $X=\Hom(\C^*,T)$ and $Y=\Hom(T,\C^*)$.
Degeneration from $\H$ to $\H^{^{_{\text{trigonometric}}}}$
amounts, effectively, to replacing the group algebra of the
lattice $X=\Hom(\C^*,T)$ by the Symmetric algebra of the
vector space $\Hom_{_\C}(\C,\Lie{T})=\h$,
 while
 the group algebra of the second
 lattice $Y$ remains unaffected.
Degeneration from $\H^{^{\,_{\text{trigonometric}}}}$
to $\H^{^{\,_{\text{rational}}}}=\hh_{t,c}$
is obtained similarly,
by further replacing the group algebra of the
lattice $Y=\Hom(T,\C^*)$
by the Symmetric algebra of the
vector space $\Hom_{_\C}(\Lie{T},\C)=\h^*$.

An analogous  degeneration pattern  applies to
other objects considered in this paper, so that
many of the constructions of the paper,
for $\G$ being a Weyl group,
have their `trigonometric' and `elliptic'
analogues, see (\ref{trig}).
\medskip

\noindent
{\bf{Harish-Chandra homomorphism.}}\quad
The goal of Part II is to describe the spherical
subalgebra, $\e\hh_{t,c}\e,$ of the symplectic reflection algebra
associated to the root system of type ${\bf{A}}$ as a `quantum Hamiltonian
reduction'
of the algebra $\dd(\gln)$ of polynomial differential operators on 
the Lie algebra $\gln$. A similar description conjecturally
exists for the  spherical
subalgebra of the symplectic reflection algebra
associated to any wreath-product $\ga$, in which case
$\gln$ gets replaced by a (vector) space of representations
of a Dynkin quiver ($\gln$ is the space of $n$-dimensional
representations of the quiver with one vertex and one edge-loop).

We begin in a greater generality of  an arbitrary complex
reductive Lie algebra $\g$. Let $\h$ be a Cartan subalgebra of
$\g$, and $W$ the Weyl group of $(\g,\h)$. Let $\dd(\g)^\g$
be the algebra of $\adg$-invariant polynomial
differential operators on $\g$,
 and $\dd(\h)^W$ the algebra of $W$-invariant  polynomial
differential
 operators on $\h$. Write $\Delta_\g$,
resp. $\Delta_\h$, for the second order
Laplacian on $\g$, resp. on $\h$, associated to a
nondegenerate invariant bilinear form on $\g$.
In 1964, Harish-Chandra \cite{HC} defined an
algebra homomorphism
$\Phi: \dd(\g)^\g\to \dd(\h)^W$,
that reduces to the restriction map:
$\C[\g]^\g\to\C[\h]^W\,,\,f\mapsto f\big|_\h,\,$
on zero order differential operators, and
 such that $\Phi(\Delta_\g)=\Delta_\h$.

In this paper, we construct a 1-parameter
deformation
of the Harish-Chandra homomorphism $\Phi$ in the special case
 $\g={\frak {gl}}_n$,
and show that this deformed  homomorphism
may (and should) be viewed as an isomorphism between the
quantum Hamiltonian
reduction
of the algebra $\dd(\gln)$ and the spherical subalgebra
$\e\hh_{t,c}\e$.
 To this end, in \S6 we work out
the formalism of {\it radial parts} of invariant differential
operators in a slightly more general framework than treated in
classical texts.

Now let $\h=\C^n$ be the Cartan subalgebra of diagonal matrices
in $\g=\gln$. Write $x_1,\ldots,x_n$ for coordinates in $\C^n$,
 and $\hreg\subset\C^n $ for the open subset of
points with pairwise distinct coordinates.
In \S7 we apply the formalism of radial parts
to construct a family of algebra homomorphisms
$\,\Phi_k: \dd(\g)^{\g}\to \dd(\hreg)^W\,$
depending on a parameter $k\in\C$,
which specializes to the  Harish-Chandra homomorphism $\Phi$
at $k=0$. For each $k\in\C$, this homomorphism $\Phi_k$
reduces to the restriction map:
$\C[\g]^\g\to\C[\h]^W\,,\,f\mapsto f\big|_\h,\,$
on zero order differential operators. Further,
we have: $\Phi_k(\Delta_\g)=\CM_k$, where
\begin{equation}\label{H}
\CM_k=\sum\nolimits_j\;\frac{\partial^2}{\partial x_j^2}\;-\;
\sum\nolimits_{i\ne j}\;\frac{k(k+1)}{(x_i-x_j)^2}\,
\end{equation}
is the Calogero-Moser operator with rational potential,
corresponding to the parameter~$k$.

We now describe the image of the map  $\Phi_k$.
Let $\dd(\hreg)^W$ denote the algebra
of $W$-invariant differential operators on $\hreg$. 
Let $\dd(\hreg)^W_-$ denote the algebra spanned by
homogeneous elements $D$ (with respect to the $\C^*$-action on $\hreg$)
 of $\dd(\hreg)^W$ such that 
$\text{\it order}(D)+\text{\it degree}(D)\le 0$. Clearly, 
$\CM_k$ belongs to ${\dd(\hreg)_-^W}$.  
Let $\cc_k$ denote the centraliser of $\CM_k$ in the
algebra ${\dd(\hreg)_-^W}$.
If $k=0$, then it follows from Lemma \ref{symbols} below that $\,
\cc_0= S\h^W,$
is the algebra of $W$-invariant
differential operators with
 constant coefficients. In general, follows from 
\cite{O} and Lemma \ref{symbols} 
that for each $k$, there exists an algebra isomorphism
$\,\sigma_k: S\h^W\iso \cc_k,\,$ such that, for any homogeneous element
$u\in S\h^W$ of degree $\ell\geq 0$,
the differential operator $\sigma_k(u) $ has order $\ell$, and its
principal symbol  equals $u$; moreover,
$\sigma_k(\Delta_\h)=\CM_k$.
It turns out that the image of the  homomorphism  $\Phi_k$ equals
$\bcal_k$, the associative 
subalgebra in $\dd(\hreg)^W$
 generated by $\cc_k$ and by
the polynomial subalgebra $\C[\h]^W\subset \dd(\h)^W,$ of zero order
 $W$-invariant operators on $\h$.

Thus, we summarize:
\begin{theorem}\label{maintheorem} There exists a (flat) family of
 surjective algebra homomorphisms $\Phi_k:$\linebreak
$\dd(\g)^{\g}\onto \bcal_k,\,$ all of which
reduce to the restriction map:
$f\mapsto f\big|_\h,\,$
on zero order differential operators.
This family
specializes,  at $k=0$, to the Harish-Chandra homomorphism $\Phi$;
furthermore, 
 $\,\Phi_k(\Delta_\g)=\CM_k$, for any $k\in \C$.
\end{theorem}

A drawback of Theorem \ref{maintheorem} is that
elements of
the algebra
$\bcal_k\subset \dd(\hreg)^W$ may have quite complicated singularites
at the divisor $\h\smallsetminus\hreg$.
Therefore, it is tricky to give an explicit description
of $\bcal_k$, as a subalgebra inside $\dd(\hreg)^W$.
Remarkably, the algebra $\bcal_k$
admits an imbedding into $\hh_\ka$,
the rational Cherednik algebra of type 
${\mathbf{A_{n-1}}}$. The image of $\bcal_k$
has a very simple  description that we now explain.

If $R={\mathbf{A_{n-1}}},$ is the root system of $\g=\gln$ then
all roots are of the same length, and $\CC=\C$. We
treat the parameter $\ka=(t,c)$ as a point on the Riemann sphere $\CP^1=
(\C\oplus\CC)/\C^*$
with homogeneous coordinates $(c:t)$, i.e.,
as a complex
number $\ka=c/t\in \C\cup\{\infty\}$ (as will be explained after
(\ref{grad}), rescaling: $(t,c)\mapsto
(r\cdot t,r\cdot c)\,,\, r\in \C^*,$ has no essential effect on the
algebra $\hh_{t,c}$).
Thus, in the ${\mathbf{A_{n-1}}}$-case we make no distinction
between the parameter: $\ka=c/t$, and the complex parameter `$k$'
entering the deformed Harish-Chandra homomorphism $\Phi_k$.
Note that the limit $\ka\to \infty$ corresponds to $t\to 0$.

For any Weyl group,
Cherednik defined a {\it faithful}
 representation of the double affine Hecke algebra, given by  the so-called
Demazure-Lusztig-Dunkl operators.
A version of Cherednik's construction
yields  a faithful representation of  $\hh_\ka$
in the vector space $\C[\hreg]$.
It is easy to see that the restriction of this representation to the
spherical subalgebra $\,\e\hh_\ka\e \subset \hh_\ka$
keeps the subspace $\C[\hreg]^W\subset \C[\hreg]$ stable,
and the resulting $\e\hh_\ka\e$-action on $\C[\hreg]^W$ is given
by  $W$-invariant regular
 differential operators on $\hreg$.
This gives an algebra homomorphism
\begin{equation}\label{HC_Spherical}
\Theta^{^{\tt{spher}}}_\ka:\, \e\hh_\ka\e\too\dd(\hreg)^W.
\end{equation}
Thus, for $\g=\gln$ and $W=S_n$, we get the following diagram:
\begin{equation}\label{diag}
{
\diagram
\dd(\g)^{\g} \rrto^<>(.5){\Phi_\ka}
&& \dd(\hreg)^{W}&& \e\hh_\ka\e\llto_<>(.5){\Theta^{^{\tt{spher}}}_\ka}
\enddiagram}
\end{equation}
The remarkable fact is, that the  homomorphisms $\Phi_\ka$
and $\Theta^{^{\tt{spher}}}_\ka$  above have the same image, the subalgebra
$\bcal_\ka\subset\dd(\hreg)^{S_n}$. Further,
the map $\Theta^{^{\tt{spher}}}_\ka$ in (\ref{diag}) is a bijection
onto its image.
Inverting this map,
we thus obtain an algebra homomorphism:
\begin{equation}\label{main_map}
\Phi_\ka^{^{\tt{spher}}}\,=
\,(\Theta^{^{\tt{spher}}}_\ka)^{\,-1}\ccirc\Phi_\ka\;:
\quad
\dd(\gln)^{\gln} \too \e\hh_\ka\e\,.
\end{equation}
This homomorphism provides a link between
invariant differential operators on $\gln$
and the Cherednik algebra $\hh_\ka$; it plays a crucial role in this
paper. Furthermore,
 replacing the Lie algebra $\g=\gln$ by the group 
$G=GL_n$, one can define a `trigonometric' version of
the deformed
Harish-Chandra homomorphism (\ref{main_map}):
\begin{equation}\label{trig}
\Phi_\ka^{^{\,_{\text{trigonometric}}}}:\;\; \dd(G)^G\too
\e\cdot\H^{^{\,_{\text{trigonometric}}}}\cdot\e\;.
\end{equation}
One might expect that replacing further $\dd(G)$ by the
so-called {\it Heisenberg double} of the corresponding quantum
group, see e.g. \cite{STS}, one gets an elliptic analogue of the
sperical Harish-Chandra homomorphism with image
 $\e\H\e$.
This agrees well with the results of \cite{EK},
although we do not pursue this line here.
\medskip

\noindent
{\bf{Calogero-Moser space.}}\quad
 It is a
well-known principle
 that  the `classical' limit of the algebra
$\dd(\g)$ is the commutative algebra of regular functions on
$T^*\g$, the cotangent bundle on $\g$.
Since $\,T^*\g=\g\oplus\g^*\simeq \g\oplus\g,\,$
the  `classical' limit of $\dd(\g)^\g$ is the
algebra $\C[\g\oplus\g]^\g$ of $\ad^{\,}\g$-invariant polynomials on
$\g\oplus\g$. Now, let $\g=\gln$.
There are reasons to think of the
limit:
$\ka=(t,c)\too\infty= (0,1)\in \CP^1,$ as a  `classical' limit.
Accordingly, we show  that the family of Harish-Chandra
homomorphisms
$\,\Phi_\ka^{^{\tt{spher}}}:\,
\dd(\g)^{\g} \onto \e\hh_\ka\e,\,$
see (\ref{main_map}), has a `classical'
 limit, a surjective algebra homomorphism
$\,\Phi_\infty^{^{\tt{spher}}}:\,
\C[\g\oplus\g]^\g \onto \e\hh_\infty\e.\,$
 It turns out that  this homomorphism is
intimately
related to
 the so-called {\it Calogero-Moser space} ${\mathcal M}_n$, a symplectic
algebraic manifold introduced by Kazhdan-Kostant-Sternberg [KKS],
and studied further by G. Wilson [Wi]. This space,
which is a particular case of the quiver variety
${\frak M}_{_{\G\!,n,c}},$ for  $\G=\{1\}$, is defined as
a categorical quotient:
 ${\mathcal M}_n\;=\;M/\!/\text{Ad}\,\GL_n,\,$ where
\begin{equation}\label{CMM}
M \;:= \;
\lbrace{(X,Y)\in \g\times\g\quad\big|\quad
[X,Y]+\id=\;\text{rank } 1 \text{ matrix}\rbrace}\,,
\end{equation}
is an affine subvariety in $\g\oplus\g,$
a special case of the quiver data variety $\bbm^c_{\G}(\BV)$.
Thus, by definition, the coordinate ring of ${\mathcal M}_n$
equals $\,\C[{\mathcal M}_n]=\C[\g\oplus\g]^\g/(\I_M)^\g,\,$
where $\I_M\subset \C[\g\oplus\g]$ denotes the defining ideal of $M$,
and the superscript `$\,\g\,$' stands for `$\ad^{\,}\g$-invariants'.
A `commutative' counterpart of homomorphism (\ref{main_map}) reads:
\begin{theorem}\label{key_intro}
The kernel of the homomorphism $\Phi_\infty^{^{\tt{spher}}}$ equals
$(\I_M)^\g$; hence one has the following chain of Poisson
algebra isomorphisms:
$$
\C[{\mathcal M}_n]\;\;=\;\;
\C[\g\oplus\g]^\g/\Ker\bigl(\Phi_\infty^{^{\tt{spher}}}\bigr)
\;\;\overset{{}^{\Phi_\infty^{_{\tt{spher}}}}}{\stackrel{_\sim}{\too}}\;\;
\e\hh_\infty\e\;\;\overset{^{_{{\sf{Satake}}^{-1}}}}{\stackrel{_\sim}{\too}}\;\;
\ZZ(\hh_\infty)\,.
$$
\end{theorem}
This result provides, in the special case $\G=\{1\}$,
a transparent conceptual interpretation of the
isomorphism of Theorem \ref{wreath_thm}.
Thus, for $\ga=S_n$ and  $\ka=\infty$,
Corollary \ref{wreath_cor} yields
\begin{theorem}\label{Morita} If $W=S_n$ is the Symmetric group
and $\ka=\infty=(0,c)\,,\,c\neq 0,$ then

\vi  The algebra $\hh_\infty=\hh_{0,c}(S_n)$
is Morita equivalent to the algebra:
$\ZZ_\infty\simeq \C[{\mathcal M}_n]$.
\vskip 1pt

\vii Any simple $\hh_\infty$-module is isomorphic to one
 of the form: $\,\hh_\infty \e\otimes_{_{\ZZ_\infty}}{\Huge{\Huge\chi}},$
for a certain character $\chi: \ZZ_\infty\to\C$. Thus, the
 isomorphism classes
of simple $\hh_\infty$-modules are parametrized
by points of the Calogero-Moser space ${\mathcal M}_n.$
\vskip 1pt

\viii
Any simple
 $\hh_\infty$-module has dimension $n!$ and is isomorphic,
 as an $S_n$-module,
to the regular representation of $S_n$.
\end{theorem}

 Berest and Wilson have introduced
an (infinite-dimensional) group ${\mathbb{G}}$ of the
automorphisms of  $\C\langle x,y\rangle$, the free associative
algebra on two generators, that keep the element
$\,xy-yx\in\C\langle x,y\rangle\,$ fixed. We show that
the group ${\mathbb{G}}$
acts naturally on $\dd(\gln)$ and on $\hh_\ka$ by algebra automorphisms,
preserving subalgebras $\dd(\gln)^{\gln}$ and $\e\hh_\ka\e,\,$
respectively. Moreover, the map
$\Phi_\ka^{^{\tt{spher}}}:
\dd(\gln)^{\gln} \too \e\hh_\ka\e\,$
commutes with the ${\mathbb{G}}$-action.
Further, in [BW], Berest-Wilson  defined a transitive action of
${\mathbb{G}}$ on the Calogero-Moser space
$\mm_n$. It follows from our results
 that the
${\mathbb{G}}$-action on $\mm_n$ can be canonically lifted to a
${\mathbb{G}}$-action on the vector bundle $\rr$,
see Theorem \ref{proj_intro}. It seems to us
 quite important to obtain a better
 geometric understanding of this  vector bundle.
\medskip

\noindent {\bf Guide for the reader.}\quad Theorem \ref{H_PBW},
Theorem \ref{comm_conj}, and part (i) of  Theorem
\ref{mckay} are proved in \S2. Theorem \ref{Coh_Mac_intro},
Theorem \ref{proj_prop_intro}, part (ii) of   Theorem
\ref{mckay}, and  Proposition \ref{k_theory} are proved in \S3.
For the proof of
the Satake isomorphism see Theorem \ref{ZeHe},
and normality of $\Spec\ZZ_{0,c}$ is Lemma \ref{normal}.
Rational Cherednik algebras are discussed in \S4.
Theorem \ref{wreath_thm} about the structure of
 $\Spec\ZZ_{0,c}$ for wreath-products is a combination of
Theorem \ref{irreps},   and Theorem \ref{wreath00} (Appendix E).
The image and the kernel of the deformed Harish-Chandra
homomorphism $\Phi_k$ are described in Theorem \ref{inject},
which contains in particular Theorem \ref{maintheorem} as part of the
statement. The proof of  Theorem \ref{inject} is given in \S10,
but it depends in an essential way on 
Theorem \ref{maintech}, and on  Theorem \ref{WLS}, proved in \S9 and in
references [LS1],[LS2],[Wa], respectively.
Finally, 
Theorem \ref{key_intro} is Theorem
\ref{class_isom}(ii).

\medskip
\noindent {\bf Acknowledgments.} {\footnotesize We are very
grateful to
 A. Braverman for many discussions and helping us
 with a number of proofs. We thank I. Cherednik, V. Drinfeld,
D. Gaitsgory, and V. Vologodsky for useful discussions. We are
especially indebted to Bert Kostant for communicating to us his
(partly unpublished) results, which we reproduced in Appendix G.
The second author
 thanks T. Stafford for consultations on non-commutative algebra.
 The work of the first author was partly conducted for 
the Clay Mathematics Institute, and was partially
supported by the NSF grants DMS-9700477 and DMS-9988796}.
\bigskip\smallskip

\centerline{\bf PART 1. SYMPLECTIC REFLECTION ALGEBRAS}
\medskip

\section{Koszul patterns and deformations}
\setcounter{equation}{0}

Recall the setup outlined in the Introduction.
Thus $(V, \om)$ is a finite dimensional symplectic
vector space, and
$\G \subset Sp(V)$ is a finite group,
such that the triple $(V,\om,\G)$ is indecomposable, i.e.,
there is no $\om$-orthogonal
direct sum decomposition: $V=V_1 \oplus V_2$,
where $V_i$ are $\G$-stable proper symplectic vector subspaces
in $V$. Schur lemma implies that there are exactly two
situations in which $(V,\om,\G)$ is indecomposable:
\vskip 2pt

\noindent
$\bullet\enspace$  \parbox[t]{140mm}{
$V$ is an irreducible $\G$-module of quaternionic type,
i.e. $(\Lambda^2 V^*)^\G\neq 0$.}

\noindent
$\bullet\enspace$  \parbox[t]{140mm}{ $V=\h\oplus \h^*$, 
where $\h$ is an irreducible
$\G$-module of either real or complex type.}
\vskip 2pt

\noindent
An application of Schur lemma shows that a triple
 $(V,\om,\G)$ is indecomposable if and only if the space
$(\Lambda^2 V^*)^\G$ is 1-dimensional, in which case
this space  is spanned by $\om$.

Given a
skew-symmetric $\C$-bilinear pairing $\kappa: V\times V \to \C\G$,
we introduce an associative algebra
$\hh_\ka$ as defined in (\ref{rel}).\medskip

\noindent
{\bf Proof of Theorem \ref{H_PBW}:}\quad
 First of all, it is clear that a necessary condition
for the PBW-property is that the form $\kappa$ is
$\G$-invariant (where $\G$ acts on itself by conjugation);
otherwise the PBW property fails already in degree two of the filtration.
Thus from now on we will assume that $\kappa$ is $\G$-invariant.


Write: $v\mapsto v^g$ for the action of $g\in \G$ on $V$.
Set $K=\C\G$, and let $E=V\otimes_{_\C}\C\G$.
The space $E$ has a natural $K$-bimodule structure, with left
$\G$-action given by $g: v\otimes a \mapsto v^g\otimes (ga),$
and right $\G$-action given by $g: v\otimes a \mapsto v\otimes(ag)$.
Let $\,T_KE= \bigoplus_i\; T^i_KE=
K \oplus E \oplus (E\otimes_K E) \oplus\ldots\, $
denote the tensor algebra of the $K$-bimodule $E$.
We observe that, for any $i\geq 0$, there is a natural
isomorphism: $T^i_KE \simeq (T^i_{_\C}V)\otimes_{_\C}\C\G$.
This gives a canonical graded algebra isomorphism:
$T_KE \simeq (TV)\#\G$. Thus we can write:
$\hh_\kappa = T_KE/I\langle P\rangle$, where $I\langle P\rangle$ is the two-sided ideal
in $T_KE$ generated by a certain $K$-bimodule
$P\subset K \oplus (E\otimes_K E).$
Explicitly, write
$E\wedge E$ for the $K$-subbimodule in $E\otimes_K E
=V\otimes_{_\C} V \otimes_{_\C}\C\G$ spanned by elements of the
form: $x\otimes y \otimes g - y\otimes x\otimes g\,,\, x,y \in V, g\in
\G,\,$ and
let $\ka_{_K}$ be the (unique) $K$-bimodule map:
$E\wedge E\to K$ that extends the pairing $\kappa: V\otimes V \to K$.
Then
we have:
$$P= \{-\ka_{_K}(p)+p \in K \,\oplus\, (E\otimes_{_K} E)\;\;\;\big|\;\;\; p\in E\wedge
E\}\,.$$

We see that $\hh_\kappa = T_KE/I\langle P\rangle
$ is a {\it nonhomogeneous quadratic
$K$-algebra}, in the terminology of [BG] (Braverman and Gaitsgory
only consider
algebras over a field but, as explained e.g. in [BGS, n.2.7],
everything works for quadratic algebras over any ground ring $R$,
provided $R$ is a finite dimensional semisimple $\C$-algebra,
e.g., $R=\C\G$). Similarly, we have:
$SV\#\G \simeq  T_KE/I\langle E\wedge E\rangle$, is a
{\it homogeneous}  quadratic
$K$-algebra. Moreover, since $SV$ is a Koszul $\C$-algebra, cf. [BGS],
it is easy to deduce that $SV\#\G$ is  a Koszul $K$-algebra.
Thus, to prove that the PBW-property holds for $\hh_\kappa =
T_KE/I\langle P\rangle$
we may apply (a `$K$-version' of) the criterion of [BG];
similar  criteria have been obtained earlier by
Drinfeld [Dr1], and Polishchuk-Positselsky [PP]. According to
[BG, Thm. 4.1] we must verify conditions (i)-(iii) of [BG, Lemma 3.3].
In our situation, conditions (i) and (iii) of that lemma become vacuous, and
condition (ii) reads:
\begin{equation}\label{BG_cond}
\ka_{_K} \otimes_{_K} \id_E - \id_E \otimes_{_K}\ka_{_K} \enspace
\mbox{vanishes on}\enspace (E\wedge E)\otimes_{_K} E \;\bigcap\;
E\otimes_{_K} (E\wedge E)\,\subset\, T^3_KE\,.
\end{equation}

Now, it is straightforward to see that
$(E\wedge E)\otimes_{_K} E \,\bigcap\,
E\otimes_{_K} (E\wedge E)= \Lambda^3V\otimes_{_\C}\C\G,\,$
where $\Lambda^3V \subset T^3V$ denotes the subspace of totally
skew-symmetric tensors. Furthermore, it is clear that condition (\ref{BG_cond})
amounts to the identity
${\mathtt{Alt}}([xy]z-x[yz])=0$, 
where $[a,b]$ stands for $a\cdot b - b\cdot a.\,$
The latter   identity is  the standard
 Jacobi identity: $\,[z,[x,y]]= [[z,x],y] +[x,[z,y]].\,$
Thus, we find that \ref{BG_cond} is equivalent to
that 
\begin{equation}\label{jacobi}
[z,\ka(x,y)]\;=\;[\ka(z,x),y]\,+\,[x,\ka(z,y)]
\quad,\quad
\forall x,y,z\in V.
\end{equation}

To study (\ref{jacobi}) in more detail,
write $\kappa(x,y)= \sum_{g\in \G} b(g,x,y)\cdot g,$
for some $b(g,x,y)\in \C$.
We claim $b(g,x,y)=0$ unless $g=1$ or $g\in S.$
Using an identity: $\,[z,g] = z\cdot g- g\cdot z$
$=
z\cdot g- g\cdot z\cdot g^{-1}\cdot g = (z-z^g)\cdot g,\,$ we find:
$$[z,[x,y]]= \sum\nolimits_{g\in\G}\;
 b(g,x,y)\cdot  [z,g] = \sum\nolimits_{g\in\G}\; b(g,x,y)\cdot  (z-z^g)\cdot g$$
Writing similar expressions for
$[[z,x],y]$ and $[x,[z,y]]$, from (\ref{jacobi})
we deduce
$$\sum\nolimits_g\; b(g,x,y)\cdot  (z-z^g)\cdot g =\sum\nolimits_g
 \; b(g,z,x)\cdot  (y^g-y)\cdot g
 +\sum\nolimits_g \; b(g,z,y)\cdot  (x-x^g)\cdot g.$$
Therefore, for each $g\in \G$ we must have:
\begin{equation}\label{main_eq}
b(g,x,y)\cdot  (z-z^g) = b(g,z,x)\cdot  (y^g-y) + b(g,z,y)\cdot  (x-x^g)
\quad,\quad x,y,z\in V
\end{equation}
Fix $g\ne 1$ and assume $b(g,x,y)$ is not identically zero.
Choosing $x,y$ generic enough so that $b(g,x,y) \neq 0$, we see
that $(y-y^g)$ and $(x-x^g)$ span $\Image(\id-g)$. Hence $g\in S.$

Fix $s\in S$. The assignment: $x,y \mapsto b(s,x,y)$
gives a skew-symmetric  bilinear form:
$V\times V \to \C$. Assuming this form is non-zero, we find
$x,y \in V$ such that $b(s,x,y)\neq 0$. Then, for
$g=s$ and
$z\in \Ker(\id-s)$, equation (\ref{main_eq}) yields:
$b(s,x,z)=0=b(s,y,z)$. This proves that the
space  $\,\Ker(\id-s)\,$ is contained in the radical
of the bilinear form:
$x,y \mapsto b(s,x,y)$. Hence, this form is proportional to $\om_s$,
since any two 2-forms on the two-dimensional
vector space $\Image(\id-s)$ must be proportional.

Thus, $\G$-equivariance of the pairing $\kappa$ implies the
existence  of an $\Ad\G$-invariant function $c: S \to \C$
and an $\Ad\G$-invariant skew-symmetric bilinear form
$f: V\times V \to \C$ such that we have
$$\kappa(x,y) = f(x,y)\cdot 1 + \sum\nolimits_{s\in S}\; c_s\cdot
\om_s(x,y)\cdot s\quad,\quad\forall x,y \in V\,.$$
 The form $f$ belongs to $(\Lambda^2 V^*)^\G$ and therefore
is proportional to $\om$, since this space is
one dimensional, by  indecomposability of the triple
$(V,\om,\G)$.\sq
\medskip

Let $\CC$ denote the vector space of $\Ad\G$-invariant functions on $S$,
so that:

\centerline{$
\dim\CC = \;\mbox{\it  number of $\Ad\G$-conjugacy classes in $S$}\,.
$}

     From now on, we assume that the pairing $\kappa$ in (\ref{rel})
has the form described in Theorem \ref{H_PBW}, for a certain
$t\in \C$ and $c\in\CC$. The corresponding algebra
$\hh_{t,c}:=\hh_\kappa$ will be referred to as a {\it symplectic
reflection algebra}. Thus, Theorem \ref{H_PBW} says:
\begin{equation}\label{grad}
\grd(\hh_{t,c})\;\simeq\;SV\,\#\,\G\,,\quad\mbox{and}\quad
\hh_{_{0,0}}\;=\;SV\,\#\,\G\,.
\end{equation}

We note
that, for any $r\in\C^*$, the rescaling map:
$V\to V\,,\, v\mapsto \sqrt{r}\cdot v,\,$ gives  an algebra  isomorphism:
$\hh_{r\cdot t\,,\,r\cdot c}\iso\hh_{t,c}$.
Thus, the  family of algebras
$\,\dis \big\lbrace\hh_{t,c}\,,\,(t,c)\in (\C\oplus\CC)
\smallsetminus\{(0,0)\}\big\rbrace
\,$ is parametrized, effectively, by the projective space
$\overline{\CC}:= (\C\oplus\CC)/\C^*$ of dimension
$\dim \CC$. Points of the form $(0,c) \in \overline{\CC}$
will play an especially important role.\medskip

\begin{remark} Define  an oriented quiver ${\mathbb{Q}}$ whose vertex set $I$
is the set of (isomorphism classes of) simple $\G$-modules.
Given two vertices, i.e. two simple $\G$-modules $\rho,\sigma\in I,$
we let $\,\dim \Hom_{_\G}(\rho\otimes V, \sigma)\,$
be
the number of arrows in ${\mathbb{Q}}$ with head $\rho$ and tail
$\sigma$. Note that since $V$ is a self-dual $\G$-module,
the  number of arrows from  $\rho$ to $\sigma$ equals
the  number of arrows from  $\sigma$ to $\rho$.
One can show, mimicking the argument in [CBH, \S3],
that the algebra $\hh_{t,c}$ is Morita equivalent to a 
{\sf proper quotient}
of a suitable {\it deformed preprojective algebra} associated to
${\mathbb{Q}}$ (viewed as the {\it double} of some other quiver), see [CBH].
For a general group $\G\subset Sp(V)$, the quiver
${\mathbb{Q}}$ is quite complicated, as illustrated by the following example.

Let $\G=S_n$ be
the Symmetric group acting on $V=\C^{n-1}\oplus\C^{n-1}$,
a sum of two copies of the  permutation
representation in  $\,\C^{n-1}=\{(x_1,\ldots,x_n)\in \C^n
\;\big|\; \sum x_i=0\}\,$. Then one has:
$$
\dim \Hom_{_\G}(\rho\otimes \C^{n-1}, \sigma)=
\begin{cases}
N-1 & \mbox{if}\enspace\rho=\sigma\\
1   & \mbox{if}\enspace\rho\enspace\mbox{and}\enspace\sigma\enspace\mbox{are {\it
adjacent}}\\
0   & \mbox{otherwise},
\end{cases}
$$
where $N$ is the number of (outer) corners in the
Young diagram corresponding to $\rho$, and
two representations are said to be {\it adjacent} if the
corresponding Young diagrams are obtained from each other
by moving one corner box from one place to another.
$\quad\enspace\lozenge$
\end{remark} \medskip

\noindent {\bf Proof of part (i) of Theorem \ref{mckay}.}$\;$
 We follow the notation of the proof of Theorem \ref{H_PBW},
in particular we have: $\hh_\kappa = T_KE/I\langle P\rangle,$  
a  nonhomogeneous Koszul
$K$-algebra. Let $\Lambda^pV\subset V^{\otimes p}$ denote
the $\G$-submodule of totally skew-symmetric $p$-tensors,
and write $\,\Lambda^pE := \Lambda^pV \otimes_{_{\C}}\C\G
\subset T_K^pE$ for the corresponding $\G$-bimodule.

To compute Hochschild cohomology we will use
a Koszul type complex of the form
\begin{equation}\label{koszul}
\ldots
\stackrel{\partial}{\too}
\hh_\ka\otimes_{_{\C\G}}\Lambda^2E\otimes_{_{\C\G}}\hh_\ka
\stackrel{\partial}{\too}
\hh_\ka\otimes_{_{\C\G}}E\otimes_{_{\C\G}}\hh_\ka
\stackrel{\partial}{\too}
\hh_\ka\otimes_{_{\C\G}}\hh_\ka\,\stackrel{\epsilon}{\onto}\,\hh_\ka
\end{equation}
where the augmentation $\epsilon$ is induced by
multiplication in the algebra $\hh_\ka$, and the 
differential $\partial\,$ is defined by the formula
\begin{align}\label{koszul_d}
\partial:\;\;& h \otimes v_1\wedge\ldots \wedge v_p \otimes h'
\;\;\;\mapsto\;\;\;\sum\nolimits_{j=1}^p\;(-1)^{j-1}
\Bigl(h v_j\otimes
v_1\wedge\ldots \wedge\widehat{v}_j\wedge\ldots \wedge v_p 
\otimes h'\nonumber\\
&-\;h \otimes 
v_1\wedge\ldots \wedge\widehat{v}_j\wedge\ldots \wedge v_p
\otimes v_j h' \Bigr)\quad,\quad
h,h'\in\hh_\ka\,,\,v_1,\ldots v_p\in V.
\end{align}
We claim that although, in general,
$v_i\cdot v_j\neq  v_j\cdot v_i$ in $\hh_\ka$, 
 we have $\partial\ccirc\partial=0,$
that is (\ref{koszul}) is indeed a complex. To verify the claim,
for any $p\geq 1$ and a subset $D\subset
\{1,2,\ldots,p+1\},$ put:
${\boldsymbol{\Delta}}_D=
\bigwedge_{_{\ell\in \{1,2,\ldots,p+1\}\smallsetminus D}}\, v_{_\ell}
\in \Lambda^\bullet V\,;$
in particular, for
 $D=\{i,j\},$ where $1\leq i<j\leq
p+1,$ we have:
$${\boldsymbol{\Delta}}_{\{i,j\}}=
v_1\wedge\ldots\wedge\widehat{v}_i\wedge\ldots \wedge
\widehat{v}_j\wedge\ldots \wedge v_{p+1}\,.
$$
We
compute
\begin{align}\label{koszul_long}
&\partial\ccirc\partial\,(h \otimes v_1\wedge\ldots 
\wedge v_{p+1}\otimes h')\;=\;\sum_{1\leq i<j\leq p+1}\;(-1)^{i+j-1}
\Bigl(h v_i  v_j\otimes{\boldsymbol{\Delta}}_{\{i,j\}}
\otimes h'-\\
&\qquad\qquad\qquad -\;h \otimes{\boldsymbol{\Delta}}_{\{i,j\}}
\otimes v_i  v_j h'\;-\;h v_j  v_i\otimes
 {\boldsymbol{\Delta}}_{\{i,j\}}\otimes h'\;-\;h \otimes
 {\boldsymbol{\Delta}}_{\{i,j\}}
\otimes v_j v_i h'\Bigr)\nonumber \\
&=\sum_{1\leq i<j\leq
p+1}\;(-1)^{i+j-1}\Bigl(\bigl(h\cd \ka(v_i,v_j)\bigr)\otimes
{\boldsymbol{\Delta}}_{\{i,j\}}
\otimes h'\;-\;h
\otimes {\boldsymbol{\Delta}}_{\{i,j\}}
\otimes \bigl(\ka(v_i,v_j)\cd h'\bigr)\Bigr)\,.\nonumber 
\end{align}
Now, a Leibniz type formula says that, for 
each $n\geq 0$ and any $u_1,\ldots,u_n\in V$
and $a\in \C\G$,
in $\hh_\ka\otimes_{_{\C\G}}\Lambda^nE\otimes_{_{\C\G}}\hh_\ka$ one has:
$$
a\otimes  u_1\wedge\ldots \wedge u_n\otimes 1 \;-\;
1\otimes u_1\wedge\ldots \wedge u_n\otimes a=
\sum_{j=1}^n\;
u_1\wedge\ldots u_{j-1}\wedge [a, u_j]\wedge 
u_{j+1}\wedge\ldots u_n,
$$
where $\,[a, u]\,$ stands for $a\cdot u-
u\cdot a.$ Hence, the last sum in
(\ref{koszul_long}) reads
$$
\underset{\mbox{\tiny{${\Big\lbrace i,j,k\Big|
{\begin{array}{l}1\leq i,j,k\leq
p+1\\
\;i<j\;\,\&\,\; k\neq i,j
\end{array}}\!\!\!\Big\rbrace}$}}}{\sum}\!\!\!\!\!\!
(-1)^{i+j+k}\cdot{\mathsf{sign}}\bigl((k-i)(i-j)(j-k)\bigr)\cdot
h\otimes [[v_i,v_j],{v}_k]\wedge
{\boldsymbol{\Delta}}_{\{i,j,k\}}\otimes h'.
$$
This sum vanishes due to 
Jacobi identity (\ref{jacobi}). Thus, (\ref{koszul}) is a complex.

We remark that each term in (\ref{koszul}) is a
{\it projective} (but not free) $\hh_\ka$-bimodule, and that
the equation $\partial\ccirc\partial=0$ would fail
if we had taken in (\ref{koszul}) tensor products
over $\C$ instead of tensor products
over $\C\G$. Further, it is straightforward to check
that the complex obtained from (\ref{koszul}) by passing
to the associated graded algebras, $\grd \hh_\ka$,
is the standard Koszul complex, hence, is exact. It follows
that (\ref{koszul}) is a projective resolution of
$\hh_\ka$, viewed as an $\hh_\ka$-bimodule. Thus, the groups
$\Ext^\bullet_{_{\hh_\ka{\mathtt{-bimod}}}}(\hh_\ka,\hh_\ka),\,$
may be computed as
cohomology groups of the complex with the terms
$$\Hom_{_{\hh_\ka-{\mathtt{bimod}}}}(\hh_\ka
\otimes_{_{\C\G}}\Lambda^pE\otimes_{_{\C\G}}\hh_\ka
\,,\,\hh_\ka)=
\Hom_{_{\G}}(\Lambda^pV\,,\,\hh_\ka),
$$
where elements $g\in \G$ act on $\hh_\ka$ by conjugation:
 $ h\mapsto g\cd h\cd g^{-1}.$ 
We conclude that  Hochschild cohomology of $\hh_\ka$ may
be computed using the  complex:
\begin{equation}\label{koszul_fin}
\ldots
\stackrel{\partial^\star}{\too}
\Hom_{_{\G}}(\Lambda^{p}V\,,\,\hh_\ka)
\stackrel{\partial^\star}{\too}
\Hom_{_{\G}}(\Lambda^{p+1}V\,,\,\hh_\ka)
\stackrel{\partial^\star}{\too}\ldots\;.
\end{equation}
with the following differential, $\partial^\star$, induced
by (\ref{koszul_d}):
\begin{equation}\label{koszul_d*}
\partial^\star\beta:\,v_1\wedge\ldots \wedge v_{p+1}
\;\mapsto\;
\sum\nolimits_{j=1}^p\;(-1)^{j-1}
[v_j\,,\,\beta(
v_1,\ldots,\widehat{v}_j,\ldots, v_p)]
\enspace,\enspace \beta\in \Hom_{_{\G}}(\Lambda^{p}V\,,\,\hh_\ka)
\end{equation}
The algebra structure on $\hh_\ka$ 
makes the graded vector space:
$$\Hom_{_{\G}}(\Lambda^\bullet V\,,\,\hh_\ka)\,:=\,
\bplus_{\!_{n\geq 0}}\;\Hom_{_{\G}}(\Lambda^nV\,,\,\hh_\ka)
\simeq \bigl(\oplus_{n\geq 0}\;\Lambda^nV^*\bigr)
\mbox{$\bigotimes_{_{\C\G}}$}
\hh_\ka\,$$
into a differential graded algebra. 
Moreover, one veriefies in a standard way that the induced
cup-product on the cohomology of this algebra coincides
with the Yoneda (equivalently, cup-) product on the 
$\Ext$-algebra
$\Ext^\bullet_{_{\hh_\ka-{\mathtt{bimod}}}}(\hh_\ka,\hh_\ka).$

Now given $g\in \G$,  set: $V_{\!_{{\sf{Ker}}}}=
\Ker(\id-g)\subset V$, resp. $V_{\!_{{\sf{Im}}}}
=\Image(\id-g)$. Thus,
 $V=V_{\!_{{\sf{Im}}}}\oplus
V_{\!_{{\sf{Ker}}}}$, and if $\dim V_{\!_{{\sf{Im}}}}=2 d(g)$,
then we
have the corresponding decomposition:
$\dis
\Lambda^{2 d(g)}V= \bplus_{\!p+q={2 d(g)}}\;\Lambda^pV_{\!_{{\sf{Im}}}}
\otimes \Lambda^qV_{\!_{{\sf{Ker}}}}. $ 
We let $\om_g$ be the $2 d(g)$-form on $V$ 
whose restriction to the line $\Lambda^{2 d(g)}V_{\!_{{\sf{Im}}}}$
coincides with  $\wedge^{d(g)}\om\big|_{\Lambda^{2
d(g)}V_{\!_{{\sf{Im}}}}}$, and that vanishes on
any other direct summand $\Lambda^pV_{\!_{{\sf{Im}}}}
\otimes \Lambda^qV_{\!_{{\sf{Ker}}}}\,,\,q>0,$
in the decomposition above. Thus,
for $g\in S$ a symplectic reflection, $\om_g$ is
the 2-form involved in the bilinear form $\ka$, see
Theorem \ref{H_PBW}. 
Further, following an idea of Alvarez [Alv],
for each  $g\in \G$, we put
$\,\psi_g := \om_g\otimes g\; \in \;
\Lambda^{2 d(g)}V^*\bigotimes \hh_\ka=
\Hom_{_{\C}}(\Lambda^{2 d(g)}V\,,\,\hh_\ka)\,.
$

We remark 
 that although the element $\psi_g\,,\,g\in \G,$ is 
{\it not} a $\G$-equivariant map,
i.e. not an element
of
$\bigl(\Hom_{_{\C}}(\Lambda^\bullet V\,,
\,\hh_\ka)\bigr)^\G =
\Hom_{_{\G }}(\Lambda^\bullet V\,,
\,\hh_\ka),\,$ one can still formally
apply the differential $\partial^\star$ given by
(\ref{koszul_d*}) to  $\psi_g$
(the reader should be warned that the map 
$\,\partial^\star:
\Hom_{_{\C}}(\Lambda^\bullet V\,,
\,\hh_\ka) \to \Hom_{_{\C}}(\Lambda^{\bullet+1} V\,,
\,\hh_\ka)
 $
does {\it not} however square to zero).
We claim that $\psi_g$ is a `cocycle', that is:
\begin{equation}\label{psi_closed}
\partial^\star\psi_g=0\quad,\quad\forall g\in \G.
\end{equation}
To see this, it suffices to check the equation:
$\partial^\star\psi_g(v_1\wedge\ldots \wedge v_{2d(g)+1})=0,$
assuming that each vector $v_j$ belongs either
to $V_{\!_{{\sf{Im}}}}$ or to $
V_{\!_{{\sf{Ker}}}}$.
If the number of the $v_j$'s which belong to $V_{\!_{{\sf{Ker}}}}$
is not equal to 1,
then all terms $\om_g(v_1,\ldots,\widehat{v}_j,\ldots, v_{2d(g)+1})$,
in the corresponding sum (\ref{koszul_d*})
vanish by the definition of $\om_g$. 
So assume $v_j \in V_{\!_{{\sf{Ker}}}}$ is this one vector.
Then, we have:
$[v_j, g]= v_j\cdot g-g\cdot v_j=(\id-g)(v_j)\cdot g=0,$
since $v_j\in \Ker(\id-g).$ Hence, for the  corresponding  term
in (\ref{koszul_d*}) we get:
$\,[v_j\,,\,\psi_g(
v_1,\ldots,\widehat{v}_j,\ldots, v_{2d(g)+1})]=0,\,$ and
(\ref{psi_closed}) follows.

Observe next that 
the algebra structure on $\hh_\ka$ induces a natural cup-product
map:
\[\cup:\;\Hom_{_{\C}}(\Lambda^p V,\hh_\ka)
\;\mbox{$\bigotimes$}\;\Hom_{_{\C}}(\Lambda^qV,\hh_\ka)
\to \Hom_{_{\C}}(\Lambda^{p+q}V,\hh_\ka\otimes\hh_\ka)
\stackrel{{}_{\sf{mult}}}{\too}
\Hom_{_{\C}}(\Lambda^{p+q}V,\hh_\ka).\,\]
An elementary linear algebra shows
  that for any $g_1,g_2\in \G$ one has an inequality:
$d(g_1)+d(g_2)\geq d(g_1g_2)$ and, moreover,
 in $\Hom_{_{\C}}(\Lambda^\bullet V\,,\,\hh_\ka)$
we have:
\begin{equation}\label{psi_product}
\psi_{g_1}\cup \psi_{g_2}=
\begin{cases}
\psi_{g_1g_2}\quad\mbox{if}\quad d(g_1)+d(g_2)=d(g_1g_2)\\
\;\;0\qquad\mbox{if}\quad d(g_1)+d(g_2)>d(g_1g_2)\,.
\end{cases}
\end{equation}
Therefore, the assignment: $g\mapsto \psi_g$ extends by $\C$-linearity
to an
$\Ad\G$-equivariant
 graded algebra homomorphism $\psi^\ka: \grd^F_\bullet\C\G \to 
\Hom_{_{\C}}(\Lambda^\bullet V\,,\,\hh_\ka),\,$
where `$F$' is the filtration on $\C\G$ introduced before
 Theorem \ref{mckay}.
Equation (\ref{psi_closed}) implies that the
homomorphism $\psi^\ka$
induces, when restricted to $\G$-invariants,
a well-defined  graded algebra homomorphism into Hochschild
cohomology
$\psi^\ka_{_{^\ZZ}}: \grd^F_\bullet\ZZ\G\to \HH^\bullet(\hh_\ka).$

To complete the proof of part (i) of Theorem \ref{mckay}
we  use an explicit calculation carried out in [AFLS]
of Hochschild cohomology
for $\Weyl_t\#\G$, the smash-product of
 the Weyl algebra, $\,\Weyl_t=TV/I\langle
x\cdot y - y\cdot x - t\cdot
\om(x,y)\rangle_{ x,y\in V},\,$ with $\G$.
Specifically, 
it has been shown in [AFLS] that for Hochschild cohomology
one has: 
\begin{equation}\label{alev}
\dim H\!H^j(\Weyl_t\#\G)= \begin{cases} 
\bn(i)\quad\mbox{if}\quad j=2i\\
\,0\qquad\,\mbox{if\quad$j$ is odd}\,,
\end{cases}
\end{equation}
 where
the numbers
$\bn(i)$ have been defined in (\ref{length}). 
 We note that, for  $c=0$ (and any $t\neq 0$), we have:
 $\hh_{t,0} = \Weyl_t\#\G$.  
Furthermore,
it is easy to verify by going through the calculation in [AFLS]
that, in the special case: $\ka=(t,0)$, the map
$\psi^\ka_{_{^\ZZ}}$ gives  in effect a bijection
\begin{equation}\label{alev_bij}
\psi^\ka_{_{^\ZZ}}
:\; \grd^F_\bullet\ZZ\G\iso \HH^\bullet(\Weyl_t\#\G)\,.
\end{equation}

Now fix $\ka=(t,c)\in \C^*\times\CC$, 
introduce an auxiliary complex parameter $r\in \C$,
and consider the 
family
$\,\{\hh_{t,r\cdot c}\}_{r\in \C}\,$
as a single $\C[r]$-algebra, $\bhh\,,$ such that: $\bhh/r\cd\bhh=$
$
\Weyl_t\#\G$.
Write $\HH^\bullet(\bhh)$ for 
the  {\it relative} Hochschild cohomology of $\bhh$ with respect
to the subalgebra
$\C[r]$.
The morphisms $\psi^\ka_{_{^\ZZ}}\,,\,
\ka=(t,r\cdot c),\,$ assemble together to give
a graded $\C[r]$-algebra morphism:
\begin{equation}\label{psi_bhh}
\psi_{_{^\ZZ}}
:\; (\grd^F_\bullet\ZZ\G)[r] \too\HH^\bullet(\bhh)\,.
\end{equation}

The algebra  $\bhh$ is flat over $\C[r]$, due to the PBW-property.
The canonical filtration on the algebra
$\,\hh_{t,c}\,$ makes $\bhh$ a filtered $\C[r]$-algebra
whose degree zero component is
$\C\G[r]$.
This filtration induces a filtration on a relative
version of
Koszul  complex (\ref{koszul_fin}) used  for computing 
$\HH^\bullet(\bhh)$,
the  {\it relative} Hochschild cohomology.
This way
the  relative Koszul  complex 
is filtered by an increasing  sequence of subcomplexes formed by
free $\C[r]$-modules of finite rank.
Therefore, since
$\,H\!H^{^{{\sf{odd}}}}(\Weyl_t\#\G)=0,\,$
the Euler-Poincar\'e principle  and 
semi-continuity 
imply that, for any $i$, the $\C[r]$-module
 $\,H\!H^{2i+1}(\bhh)$ is supported at most on
a countable subset of  $\C$ and, similarly,
$H\!H^{2i}(\bhh)$ is, up to torsion
on 
a countable subset of  $\C$,  a free $\C[r]$-module
of rank:
$\dim H\!H^{2i}(\Weyl_t\#\G)= \bn(i).\,$
Further, since the Koszul complex for
$\bhh$ has finite length, there are only finitely
many non-zero cohomology groups involved. 
Thus, the map (\ref{psi_bhh}) is a map
between $\C[r]$-modules of the same generic rank.
Moreover, both modules have no torsion at $r=0$,
and at this point the map specializes to a bijection.
It follows that there is  at most a countable set
$\mbox{BAD} \subset \C$ such that  for 
any  $r\in \C\smallsetminus\mbox{BAD},$
 the map (\ref{psi_bhh}) specializes
to a bijection, hence an algebra isomorphism:
$\psi_{_{^\ZZ}}^{t,r\cdot c}
:\; (\grd^F_\bullet\ZZ\G) \iso\HH^\bullet(\hh_{t,r\cdot c})$.
This proves part (i) of Theorem \ref{mckay},
since $\,\hh_{t,r\cdot c}\simeq
\hh_{t/r, c}$.
\quad\qed\medskip

Recall the `averaging' idempotent
 $\e= \frac{1}{|\G|}\sum_{g\in \G} g \in \C\G$.
If the group $\G$ acts on an associative algebra $\A$ by algebra 
automorphisms, then one has a canonical algebra
isomorphism:
\begin{equation}\label{smash_iso}
\A^{\G} \iso \e\cd(\A\#\G)\cd\e\quad,\quad a\,\mapsto\, a\cd\e = \e\cd a
\end{equation}

Recall now
 the 
{\it spherical subalgebra}: $\,\e \hh_{t,c} \e\subset  \hh_{t,c}\,$
(we warn the reader that this subalgebra does
{\it not} contain the
unit of $\hh_\ka$).
For $t=0$ and $c=0$, from (\ref{grad}) and (\ref{smash_iso})
one finds:
$\,\e \hh_{_{0,0}}\e=\e\bigl(SV\#\G\bigr)\e\simeq
(SV)^\G,\,$ the algebra of $\G$-invariants. Furthermore, isomorphism
(\ref{grad}) yields readily:
$\grd\bigl(\e \hh_{t,c} \e\bigr) \simeq (SV)^\G.$
Similarly, as we have noticed in the proof of Theorem
\ref{mckay}, for  $c=0$ and any $t\neq 0$, we have:
 $\hh_{t,0} = \Weyl_t\#\G$, where
$\,\Weyl_t=TV/I\langle
x\cdot y - y\cdot x - t\cdot
\om(x,y)\rangle_{ x,y\in V},\,$
is the  {\it Weyl algebra} 
 of the symplectic vector space $(V,t\cdot\om)$.
By (\ref{smash_iso}),
the spherical subalgebra $\e \hh_{t,0} \e$
is isomorphic
to $\Weyl_t^{^{\!\G}}$, the subalgebra of $\G$-invariants in
the Weyl algebra.

\begin{theorem}\label{deform_t} For any fixed $t\neq 0$, the
family $\,\{\hh_{t,c}\}_{c\in\CC}\,$ gives a universal
deformation of the algebra $\Weyl_t\#\G$,
and  the
family $\,\{\e \hh_{t,c} \e\}_{c\in\CC}\,$ gives a universal
deformation of the algebra~$\Weyl_t^{^{\!\G}}$.
\end{theorem}

To prove the theorem, we fix $t\neq 0$ and given $c\in \CC$,
 treat the family $\,\{\hh_{t,\ve\cdot c}\}_{\ve\in\C}\,$
as a single $\C[\ve]$-algebra, still denoted $\hh_{t,\ve c}$,
 where $\ve$ is now regarded
as an indeterminate.

\begin{lemma}\label{pasha2} For any  nonzero  $c\in \CC$,
the $\C[\ve]/(\ve^2)$-algebras: 
$\hh_{t,{\ve}c}/{\ve}^2\cd\hh_{t,{\ve}c}$ and
$\hh_{t,0}\bigotimes_{_{\C}}\, (\C[{\ve}]/{\ve}^2)$ are not isomorphic. 
\end {lemma}

\Pf . Observe first that  we have
a natural isomorphism of $\C$-vector spaces $\,\phi: \hh_{t,0} \iso
\ve\cd\hh_{t,{\ve}c}/{\ve}^2\cd \hh_{t,{\ve}c}$.
Assume,  contrary to the claim, that there is an
isomorphism $\,\hh_{t,{\ve}c}/{\ve}^2\cd \hh_{t,{\ve}c}
\simeq \hh_{t,0}\bigotimes_{_{\C}}\, (\C[{\ve}]/{\ve}^2),\,$
as algebras over $\C[\ve]/(\ve^2)$, the ring of `dual numbers'.
Then the natural imbedding:
$\Weyl_t\into \Weyl_t\#\G=\hh_{t,0}$ lifts to the dual numbers, that is, one
can find a $\C$-linear map $f:V\to \hh_{t,0}$
such that in $\hh_{t,{\ve}c}/{\ve}^2\cd \hh_{t,{\ve}c}$, the following equation holds:
\[
\bigl(x+\phi\ccirc f(x)\bigr)\cdot \bigl(y+\phi\ccirc f(y)\bigr) -
\bigl(y+\phi\ccirc f(y)\bigr)\cdot
\bigl(x+\phi\ccirc f(x)\bigr)\;=\;t
\cdot\omega(x,y)\cdot 1\quad,\quad \forall x,y\in
V\,.
\]
Using the defining relations in $\hh_{t,{\ve}c}$
and applying the isomorphism $\phi^{-1}$,  this gives
the following equation in $\hh_{t,0}$:
\begin{equation}\label{p3}
[f(x),y]+[x,f(y)]+\sum\nolimits_s c_s\cdot \om_s(x,y)\cdot s=0 
\end{equation}
Write $f(x)=\sum_{g\in \G}\; f_g(x)\cdot g$, where
$f_g(x)\in \Weyl_t$. Fix  $s\in S$ such that $c_s\ne 0$.
Taking the $s$-term 
of the equation (\ref{p3}), we get 
\begin{equation}\label{p4}
[f_s(x),y]+f_s(x)\cdot (y^s-y)+[x,f_s(y)]+f_s(y)\cdot (x-x^s)=-c_s\cdot\om_s(x,y).
\end{equation}
One may choose
 $x,y$, a basis of the 2-dimensional vector space
$\Image(s-\id)$  such that $\om_s(x,y)=1$,
and such that $s(x)={\lambda}x$, $\,s(y)={\lambda}^{-1}y$,
where  ${\lambda}\ne 1$.
Put
$a:=\frac{1}{c_s}f_s(x)\neq 0$, and $b:= \frac{1}{c_s}f_s(y)\neq 0$.
    From equation (\ref{p4}) we find
\begin{equation}\label{p7}
[a,y]-(1-{\lambda}^{-1})a\cdot y+[x,b]+(1-{\lambda})b\cdot x=-1\,.
\end{equation}

Let $\Weyl^{^{{\sf{Ker}}}}$, resp. $\Weyl^{^{{\sf{Im}}}}$,
denote the subalgebra in  $ \Weyl_t$ generated by the subspace
$\Ker(\id-s) \subset V$, resp. $\Image(\id-s)$.
Thus, $\Weyl_t = \Weyl^{^{{\sf{Ker}}}} \bigotimes \Weyl^{^{{\sf{Im}}}}$,
and
 the monomials: $\,\{x^my^n\}_{m,n\geq 0}\,$ form a 
$\C$-basis of $\Weyl^{^{{\sf{Im}}}}$.
Write $a,b$ as a $\Weyl^{^{{\sf{Ker}}}}$-linear combination of the monomials $x^my^n$.
Then begin  treating $x$ and $y$ as commuting variables, so that
the expressions for $a,b$ become certain polynomials: $a,b \in 
\Weyl^{^{{\sf{Ker}}}}[x,y]$.
Using the identity:
$b x = -[x,b]+ xb,$ and setting $t=1$, from equation (\ref{p7}) we deduce
\[
\Bigl(\frac{\partial}{\partial x}-(1-{\lambda}^{-1})\cdot y\Bigr)\,a+
\Bigl({\lambda}\frac{\partial}{\partial y} +(1-{\lambda})\cdot
x\Bigr)\,b=-1
\quad\mbox{in}\quad \Weyl^{^{{\sf{Ker}}}}[x,y]\,.
\]
Set $A:=-a\,,\,B:=-{\lambda}b$,
and $\nu:= 1-{\lambda}^{-1}$. Our equation reads:
\begin{equation}\label{p5}
\bigl(\frac{\partial}{\partial x}-\nu\cdot y\bigr)A +
\bigl(\frac{\partial}{\partial y} -\nu\cdot x\bigr)B=1
\quad,\quad A,B\in \Weyl^{^{{\sf{Ker}}}}[x,y]\,.
\end{equation}
The proof of  the Lemma is now completed by the following
\medskip

\noindent
{\bf Claim:}\quad
For any $\nu\in\C^*$,
equation (\ref{p5})
has no polynomial solutions $A,B\in \Weyl^{^{{\sf{Ker}}}}[x,y]$.
\medskip

\noindent
{\it Proof of Claim.}$\;$
Clearly, it suffices to prove that equation (\ref{p5})
has no  solutions $A,B\in \C[x,y]$. Further,
replacing the function $A$, by the function:
$(x,y) \mapsto \nu^{-1}\cdot A(\nu\cdot x, \nu\cdot y)$,
and making a similar substitution for $B$,
one reduces the problem to the case: $\nu=1$.
Assuming this, we further specialize: $x=z,y=\bar z$,
 where $z$ is a complex variable,
and $ \bar z$ its complex conjugate.
Integrating both sides of  equation (\ref{p5})
for  $\nu=1$ against
the Gaussian measure $e^{- z{\bar z}}dzd{\bar z}$, we get
\begin{equation}\label{int}
\int_\C\Bigl(\bigl(\frac{\partial}{\partial z}- {\bar z}\bigr)
A(z,  {\bar z}) +
\bigl(\frac{\partial}{\partial {\bar z}} -z\bigr)B(z,  {\bar z})\Bigr)
e^{- z{\bar z}}dzd{\bar z}\;
=\;\int_\C e^{- z{\bar z}}dzd{\bar z}\;,
\end{equation}
where the integral on the LHS is absolutely convergent, for any
polynomials $A,B$.
But, integration by parts shows that the LHS(\ref{int})$\,=0$, while clearly
RHS(\ref{int})$\,> 0$. The contradiction completes the proof.  \quad\sq
\medskip

\noindent
{\bf Proof of Theorem \ref{deform_t}:}$\;$
Deformations of the algebra $\Weyl_t\#\G$ are controlled
by $H\!H^2\bigl(\Weyl_t\#\G\bigr),$ the second Hochschild cohomology
group
of $\Weyl_t\#\G$, and obstructions to deformations are controlled
by the third  Hochschild cohomology
group
of $\Weyl_t\#\G$. It was shown in  [AFLS] that
$\dim H\!H^2\bigl(\Weyl_t\#\G\bigr) =  \;\mbox{\it
number of $\Ad\G$-conjugacy classes in $S$},\,$ and
$H\!H^3\bigl(\Weyl_t\#\G\bigr) = 0$. It follows that the deformations are
unobstructed, the point $\Weyl_t\#\G$ is a smooth point of
the moduli space of all  deformations,
and the tangent space at the point $\Weyl_t\#\G$ to the moduli
space has dimension equal to $\dim\CC$. 
The cohomology classes 
corresponding to deformations from the
family $\,\{\hh_{t,c}\}_{c\in\CC}\,$ 
form a certain vector subspace $H\subset 
H\!H^2\bigl(\Weyl_t\#\G\bigr)$.
Clearly, $\dim H= \dim \CC$ if and only if, for any
non-zero $c\in \CC$, the algebra $\hh_{t,\ve c}$ gives a non-trivial
deformation of $\Weyl_t\#\G$ over $\C[\ve]/(\ve^2)$.
Hence,  Lemma \ref{pasha2} implies:
$\dim H= \dim \CC= \dim H\!H^2\bigl(\Weyl_t\#\G\bigr)$,
and therefore $H= H\!H^2\bigl(\Weyl_t\#\G\bigr)$.
Thus,
our family is
 a universal family. 

To prove the second claim 
observe that by (\ref{smash_iso}) we have: $\e(\Weyl_t\#\G)\e\simeq
\Weyl_t^{\!^{\G}}$.
Hence the family: $\,\{\e\hh_{t,c}\e\}_{c\in\CC}\,$ gives a flat
deformation of the algebra $\Weyl_t^{\!^{\G}}$. The rest
of the proof is entirely similar:
 Morita equivalence of the algebras
$\Weyl_t\#\G$ and $\Weyl_t^{\!^{\G}}$ implies that these algebras have 
isomorphic Hochschild cohomology, {see [AFLS].~\sq}\medskip

Observe that
 the family $\{\A_t\}_{t\in \C},$
where $\,\A_t:=\Weyl_t^{^{\!\G}}$ if $t\neq 0,$
and $\A_{t=0}:=(SV)^\G,$ is flat,
hence, gives a deformation of the algebra
$\A_{t=0}=(SV)^\G$. Since the latter algebra
is commutative, the general construction of \S15 yields
 a Poisson structure on
the algebra $(SV)^\G=\C[V^*]^\G$. It is well-known and easy to verify that
the corresponding  Poisson
bracket $\{-,-\}_\om$ is obtained by restricting the  Poisson
bracket on $\C[V^*]$ induced by the symplectic structure on $V^*$.
Note that, for the standard grading
$(SV)^\G= \bigoplus_{i\geq 0}\;(S^iV)^\G$, we
have $\{-,-\}_\om: (S^iV)^\G \times (S^jV)^\G\too (S^{i+j-2}V)^\G.\,$
More generally, a Poisson bracket $B$ on $(SV)^\G$ is said to
be {\it of degree $\ell$} if, for any $i,j\geq 0$,
we have $B: (S^iV)^\G \times (S^jV)^\G\too (S^{i+j+\ell}V)^\G.\,$

\begin{lemma}\label{pois}
\vi Any Poisson bracket on $(SV)^\G$ of degree
$(-2)$ is proportional to  $\{-,-\}_\om$.

\vii Any Poisson bracket on $(SV)^\G$ of degree $\ell<-2$ is zero.
\end{lemma}

\Pf .
We identify $V$ with $V^*$ using the symplectic form.
Let $Y$ be the set of points of $V$ with a nontrivial isotropy group
in $\G$.
Then $(V\smallsetminus Y)/\G$ is a smooth manifold. A Poisson bracket on
$(SV)^\G$ in particular defines a Poisson bracket
on $(V\smallsetminus Y)/\G$, which is the same as
a $\G$-invariant Poisson bracket, that is a bivector field, on
$V\smallsetminus Y$.

Since $\G\subset Sp(V)$ is a finite group,
any nontrivial element of $\G$ has at least two eigenvalues
$\neq 1$. This implies that the codimension of
$Y$ is at least $2$. Thus, by Hartogs theorem, any
regular bivector on $V\smallsetminus Y$ extends by continuity to
a regular algebraic  bivector on the whole of $V$.
But a  regular bivector on $V$ of homogeneity degree $\ell$ must vanish
if $\ell<-2$, and comes from an element of
$\Lambda^2V$ if $\ell=-2$. In the latter case, since our bracket 
is $\G$-invariant, by indecomposability of $(V,\om,\G)$, this bracket 
must be a multiple of $\{-,-\}_\omega$.
\sq\medskip

\begin{lemma} \label{genchar}
Let $t,c$ be such that the algebra $\hh_{t,c}$ is commutative. Let
$\chi: \e\hh_{t,c}\e \to \C$ be a generic character, and let
$T_\chi=\hh_{t,c}\e\otimes_{\e\hh_{t,c}\e}\chi $ be the
corresponding induced $\hh_{t,c}$-module. Then $T_\chi$ is isomorphic, as a
$\Gamma$-module, to the regular
representation of $\Gamma$.
\end{lemma}  

\Pf . Consider the family of algebras $\{\hh_{rt,rc}\}_{r\in \C}$.
These algebras 
are isomorphic to each
other for all $r\ne 0$, so it is sufficient to prove the lemma
for a generic member of the family. To this end, we treat
the family $\{\hh_{rt,rc}\}_{r\in \C}$
as a single $\C[r]$-algebra $\thh$. This algebra is {\it flat} over
 $\C[r]$, due to PBW-property. 
It follows that $\e \thh\e$ is a commutative flat  $\C[r]$-algebra.

Further, $\thh\e$ is  a finitely generated $\e \thh\e$-module,
due to Theorem \ref{Coh_Mac_intro}(ii) to be proved below
(which is independent of the material in between).
 Fix a simple $\Gamma$-module $E$, and set
$F=\Hom_{_\G}(E, \thh\e)$. By standard results in invariant theory,
see proof of Theorem  \ref{Coh_Mac_intro} below,
$F$ is a finitely generated $\e \thh\e$-module
which is, moreover, flat over  over $\C[r]$, by the PBW-property.
We write $F_r=F/r\cd F= \Hom_{_\G}(E, \hh_{rt,rc}\e)$
for the corresponding $\e\hh_{rt,rc}\e$-module.

Given a  commutative $\C$-algebra $\A$ without zero-divisors,
 and a $\A$-module $M$, let $\rk_{_\A}(M)$
 denote
the generic rank of $M$, i.e., $\rk_{_\A}(M):=\dim_{_{Q(\A)}}\bigl(
Q(\A)\otimes_\A M\bigr),\,$
where $Q(\A)$ stands for the field of fractions of $\A$.
 It is clear that
$\rk_{_{\e\hh_{0,0}\e}}(F_0) =\dim E$. On the other hand, proving the
lemma amounts to showing that $\rk_{_{\e\thh\e}}F =\dim E$.

To this end, let $\pp=(\e\thh\e)\cd r\subset \e\thh\e$ be the prime ideal 
generated by 
$r$, and set $\A:=(\e\thh\e)_{_{(\pp)}}$,  the localization of $\e\thh\e$ with respect
to $\pp$. Thus, $\A$  is a local ring with the residue
class field $\A_0:=\A/\A\cd r\simeq Q(\e\hh_{_{0,0}}\e),\,$
the field  of fractions of the algebra $\e\hh_{_{0,0}}\e$.
Furthermore, $\A$  is a local ring of dimension $1$ with
local parameter $r$, due to the
$\C[r]$-flatness of $\e\thh\e$.
Let
$\fff= \A\otimes_{_{\e\thh\e}} F$ denote the $\A$-module
 obtained by localizing
at $\pp$ the $\e\thh\e$-module $F$, and set $\fff_0:= \fff/\fff\cd r$,
an $\A_0$-module.
Observe that since $F$ is $\C[r]$-flat, and $r$ is a local parameter
on $\A$, the module $\fff$ is free over $\A$. Thus, we obtain
\[
\rk_{_{\e\thh\e}}F\; =\;  \rk_{_{\A}}\fff\;=\;\rk_{_{\A_0}}(\fff_0)\;=\;
\rk_{_{\e\hh_{_{0,0}}\e}}(F_0)\;.\qquad\square
\]
\smallskip

 \noindent {\bf Proof of Theorem \ref{comm_conj}.}$\;$
 Introduce
an auxiliary variable `$\hbar$', and set $TV[\hbar]= TV\otimes
\C[\hbar]$. We treat $TV[\hbar]$ as a graded algebra with $\deg
\hbar=2$, and assume that the group $\G$ acts trivially on
$\hbar$. Let
$${\sf{\widetilde{H}}}:=\;
(TV[\hbar]\#\G)/I\langle x\otimes y - y\otimes x - \kappa(x,y)\cdot \hbar \in
T^2V\,\, \oplus\,\, \C\G[\hbar]\rangle_{x,y \in V}\,,
$$
be a `homogenized' version of the algebra $\hh_\ka$. We set
$\A= \e{\sf{\widetilde{H}}}\e,$ the corresponding
spherical subalgebra. Clearly, $\A$ is a flat
$\C[\hbar]$-algebra such that:
$\A/\hbar\A= \grd(\e\hh_{t,c}\e)$.
Moreover, since $\A/(\hbar-1)\A=
\e\hh_{t,c}\e,\,$ one may view the algebra $\e\hh_{t,c}\e$ as a
flat $\hbar$-deformation of the commutative algebra
$\grd(\e\hh_{t,c}\e)= (SV)^\Gamma$. Thus, for any $(t,c)$,
by the general
construction of \S15,
 this deformation gives rise to  a  Poisson
bracket $B_{t,c}$ on $(SV)^\G$. 
Let $\,{\mathbf{m}}_{t,c}\in \{1,2,\ldots,\infty\},$
 be  the number involved in  the
construction of the Poisson bracket, see Lemma \ref{commut}.
We have:

\begin{claim}\label{claim_pf} \vi There are only two alternatives:
either $\;{\mathbf{m}}_{t,c}=1$, or 
$\;{\mathbf{m}}_{t,c}=\infty$. 
 The algebra $\e\hh_{t,c}\e$ is non-commutative if 
$\;{\mathbf{m}}_{t,c}=1,$ and commutative if 
$\;{\mathbf{m}}_{t,c}=\infty$. 

\vii If $\,{\mathbf{m}}_{t,c}=1,$ then $B_{t,c}=
f(t,c)\cdot\{-,-\}_\omega,$
where $f: \C\oplus \CC \to \C$ is a non-zero
{\sf linear} function.
\end{claim}

To prove part (i) of the Claim, note that since  $\deg \hbar=2$, it is clear that
 the Poisson bracket  $B_{t,c}$ has
degree $(-2{\mathbf{m}}_{t,c})$. Thus, by Lemma \ref{pois}(ii),
in order for $B_{t,c}$ to be non-zero
 we must have ${\mathbf{m}}_{t,c}=1$. But, 
by Lemma \ref{commut}, vanishing of $B_{t,c}$ implies
commutativity of the algebra $\A$, and (i)
follows.

To prove part (ii) of Claim \ref{claim_pf},
we  change our point of view, and regard $t$ and $c$ as
variables, and assign them grade degree 2, while set $\deg\hbar=0$.
Thus, $\A$ becomes now a graded $\C[\C\oplus\CC]$-algebra,
depending on $\hbar$, the deformation parameter.
Applying the Poisson bracket construction of \S15 in this
new setting we get, for ${\mathbf{m}}=1$,
a bracket $B$ on $(SV)^\G\bigotimes \C[\C\oplus\CC]$
of degree $(-2)$.
By Lemma \ref{pois}(i), we have:
$B=f(t,c)\{-,-\}_\omega,$ for some $ f(t,c)\in \C$.
Since the relations in
$\hh_{t,c}$  become homogeneous in our new grading: 
$\deg(t)=\deg(c)=2$, it follows
that: $t,c\mapsto f(t,c)$ is a certain linear function on $\C\oplus\CC$.
Moreover, for any  $(t,c)\in \C\oplus\CC$ such that $B_{t,c}\neq 0$,
the Poisson bracket $B_{t,c}$ is  clearly a specialisation
of $B$ at the point $(t,c)$, that is: $B_{t,c}
= f(t,c)\{-,-\}_\omega$. This completes the proof of 
Claim \ref{claim_pf}(ii).\medskip

Claim \ref{claim_pf} implies that the algebra
 $\e\hh_{t,c}\e$ is commutative
if and only if  $f(t,c)=0$.
We see that the parameters $(t,c)$ such that $\e\hh_{t,c}\e$ is
commutative form a hyper-plane in $\C\oplus\CC$ given by
the equation:  $f(t,c)=0$.
Therefore, to complete
 the proof of the theorem, it  suffices to show that $f(t,c)$
is a multiple of $t$,
i.e. that the above hyper-plane is the one given
by the equation $t=0$. To prove this assume that, for some $t,c$,
one has $f(t,c)=0$, and hence the algebra $\e\hh_{t,c}\e$ is
commutative. In this case, choose a generic  character $\chi:
\e\hh_{t,c}\e \to \C$, and let
$T_\chi=\hh_{t,c}\otimes_{\e\hh_{t,c}\e}\chi $ be the
corresponding induced $\hh_{t,c}$-module. By lemma
\ref{genchar}, this module is isomorphic, as a $\G$ module, to the
regular representation of $\G$.

Now, for any $g\in \G$ such that $g\neq 1$,
in the regular $\G$-representation
one has ${\mathtt{tr}}(g)=0$.
Taking the traces in $T_\chi$ of both sides of
the main commutation relation: $x\otimes y - y\otimes x=
\ka(x,y)$ in the algebra
 $\hh_{t,c}$ we deduce: $t=0$. Thus, the two hyperplanes:
$f(t,c)=0$, and $t=0$ coincide, and we are done.
\sq\medskip

According to Theorem \ref{comm_conj},
the algebra $\e\hh_{0,c}\e$ is commutative 
for any  $c\in \CC$, and the family 
$\{\e\hh_{t,c}\e\}_{t\in\C}$ gives, by the PBW-theorem,
a flat deformation of $\e\hh_{0,c}\e$.
Hence, by the general
construction
of \S15 this deformation induces a Poisson bracket on $\e\hh_{0,c}\e$,
to be denoted $\{-,-\}$. Recall further that 
the algebra $\e\hh_{0,c}\e$  has a
natural increasing filtration
$F_\bullet(\e\hh_{0,c}\e)$ induced from that on $\hh_{0,c}$.

Let $B$ be any  Poisson bracket on  $\e\hh_{0,c}\e$.
We say that $B$ has {\it filtration degree} $\ell\in\Z$,
if the following two conditions hold:

\noindent
$\bullet\quad\; B\bigl(F_i(\e\hh_{0,c}\e)\,,\,F_j(\e\hh_{0,c}\e)\bigr)\,\subset
F_{i+j+\ell}(\e\hh_{0,c}\e)\quad,\quad\forall i,j\geq 0;$

\noindent
$\bullet\quad$ There exist 
$i,j\geq 0,$ such that:  $
B\bigl(F_i(\e\hh_{0,c}\e)\,,\,F_j(\e\hh_{0,c}\e)\bigr)\,\not
\subset\,F_{i+j+\ell-1}(\e\hh_{0,c}\e).$

\begin{lemma}\label{pois_filt}
Let $B$ be a non-zero Poisson bracket on  $\e\hh_{0,c}\e$ of filtration degree
$\ell\leq -2$. Then $\ell =-2$ and, moreover,
$B=\lambda\cdot  \{-,-\}$, for some non-zero constant $\lambda\in\C$.
\end{lemma}

\Pf . 
The Poisson bracket $B$ induces naturally a {\it non-zero}
 Poisson bracket $\grd^{\,}B$ on $\grd(\e\hh_{0,c}\e)$ of degree $\ell$.
Recall that $\grd(\e\hh_{0,c}\e)\simeq (SV)^\G$.
Thus, we may view  $\grd^{\,}B$ as a  {\it non-zero} Poisson bracket on
$(SV)^\G$. Hence, Lemma \ref{pois}(ii) yields:
$\ell=-2$. Moreover, part (i) of the same Lemma implies
that there exists a constant  $\lambda\in\C$
such that $\grd^{\,}B=\lambda\cdot\{-,-\}_{_\om}$.
It follows that $B-\lambda\cdot  \{-,-\}$
is a
Poisson bracket on  $\e\hh_{0,c}\e$ of filtration degree $<\ell=-2$.
The argument above shows that this bracket must vanish, and we are done.
\qed

\section{Representation theory of the algebra $\hh_{0,c}$}
\setcounter{equation}{0}

\noindent
{\bf Proof of Theorem \ref{Coh_Mac_intro}:}\quad
 It is known from invariant theory
 that $(SV)^\G$ is a finitely generated
algebra and, moreover, $SV$ is a finitely generated
$(SV)^\G$-module:
by Hilbert's {\it Basis theorem}, each isotypic component of
$SV$ is finitely generated over $(SV)^\G$, see e.g. [PV].
 Hence, the
isomorphisms: $\grd\bigl(\e \hh_{t,c} \e\bigr) \simeq (SV)^\G$
and $\grd\bigl(\hh_{t,c} \e\bigr) \simeq SV\,$ insured by
Theorem \ref{H_PBW}
imply that $\e \hh_{t,c} \e$ is a finitely generated
 algebra without zero divisors,
and $\hh_{t,c} \e$  is a finitely generated $\e\hh_{t,c} \e$-module.

Next we need to recall some definitions. Given a not necessarily
commutative associative $\C$-algebra $\A$, one defines
the notion of a {\it rigid Auslander dualizing complex}
for $\A$, see [VB], [YZ]. This is an object
$\DD_\A \in D_f(\A\mbox{-{\tt{bimod}}})= $ bounded derived category
of complexes of $\A$-bimodules, whose cohomology groups are
finitely generated both as left and as right $\A$-modules.
Such a dualizing complex might or might not exist, in general,
but if it exists it is unique [VB], and enjoys nice functorial
properties, as explained in  [YZ]. Furthermore, it has
been shown in [VB], [YZ] that $\DD_\A$ does exist  if $\A$ is
a positively filtered $\C$-algebra, such that $\grd(\A)$ is
a finitely generated commutative algebra with $\grd_0(\A)=\C$.
In particular, for a commutative finitely generated algebra $\A$,
the dualizing complex, $\DD_\A$,
exists and coincides, due to the
uniqueness property, with the one
known in Commutative Algebra.

Fix an algebra $\A$ that has the rigid Auslander dualizing complex,
$\DD_\A.$ We say that $\A$ is {\it Cohen-Macaulay} if all the cohomology
groups $H^i(\DD_\A)$ vanish, except for $i=0$.
In the situation of a filtered algebra $\A$, as above, our definition
agrees, by [YZ, Prop. 6.18], with other definitions of
(non-commutative) Cohen-Macaulay rings that involve
 Gelfand-Kirillov dimensions,
e.g., those used in [Bj] and [CBH].
For a Cohen-Macaulay algebra $\A$, we will abuse the notation, and write
$\DD_\A$ for the bimodule $H^0(\DD_\A)$, quasi-isomorphic to it.
A Cohen-Macaulay algebra $\A$ is called {\it Gorenstein}
if $\DD_\A=\A$.
Further, a finitely generated (left) module $M$ over
a Cohen-Macaulay algebra $\A$ is said to be a
{\it Cohen-Macaulay module of degree} $d$ if,
for all $i\neq d$ we have:
$\Ext^i_{\A\mbox{\tiny{-left}}}(M,\DD_\A)=0$.
For commutative algebras, all these definitions agree with the standard
ones.

The commutative algebra $\grd\bigl(\e \hh_{t,c} \e\bigr) \simeq (SV)^\G$
is known to be Gorenstein [Wat].
Part (i) of the Theorem
now follows from a result proved by Bjork [Bj] saying that:
{\it  for a positively filtered
associative $\C$-algebra $\A$, with $\grd_0(\A)=\C$, one has}
$$\;\grd \A\quad \mbox{\it is}\enspace
\Big\{
\begin{array}{l}
{\mbox{\it Cohen-Macaulay,}} \\ 
{\mbox{\it resp. Gorenstein}}
\end{array}
\Big\}
\quad\Longrightarrow\quad \A
\quad \mbox{\it is}\enspace\Big\{
\begin{array}{l}
{\mbox{\it Cohen-Macaulay,}}\\
{\mbox{\it resp. Gorenstein}}
\end{array}
\Big\}\,.
$$
(For the Gorenstein property, Bjork actually proves the implication
above with $\DD_\A$ being replaced by $\A$. But, if
$\grd(\A)$ is Gorenstein, then an explicit description
of $\DD_\A$ given in [VB, Prop. 8.2, Prop. 8.4] shows that
$\DD_\A=\A$ in this case). This proves part (i).

To prove part (ii) of the Theorem,
 recall the following standard result from commutative algebra:
\vskip 2pt

\noindent
{\sl Claim ([BBG, Theorem 2.1]).}\quad {\it Let $A\subset B$
be finitely generated commutative $\C$-algebras with unit,
and $M$ a $B$-module which is
 finitely generated over $A$.
Then $M$ is Cohen-Macaulay over $B$ if
and only if it is Cohen-Macaulay over $A$.}
\smallskip

\noindent
This implies that $SV$ is a Cohen-Macaulay $(SV)^\G$-module
of degree $0$ (i.e. its support has full dimension).
Further, Theorem \ref{H_PBW} shows that $\grd(\hh_\ka\e)\simeq SV$.
We can now apply [Bj] to $\A= \e\hh_\ka\e$ and
$M=\hh_\ka\e$, to conclude that for any $i\geq 0$:
$\,\Ext^i_{\grd(\A)}\bigl(\grd(M)\,,\,\grd(\A)\bigr)=0\,$  implies:
$\,\Ext^i_\A(M,\A)=0.\,$ It follows, due to the
equality $\DD_\A=\A$, that $\hh_\ka\e$ is a Cohen-Macaulay
$\e\hh_\ka\e$-module.

We now prove (vi). Equip $M=\hh_\ka\e$ with the induced
increasing filtration
$F_0M \subset F_1M \subset \ldots,$
and choose a finite set, $u_1,\ldots,u_\ell,$ of homogeneous
generators of the $\grd(\e\hh_\ka\e)$-module
$\grd(M)=\grd(\hh_\ka\e)$. Let $d_i=\deg(u_i)$,
and let $\widetilde{u}_i\in F_{d_i}M$ be a representative
of $u_i$. Then $\,\widetilde{u}_1,\ldots,\widetilde{u}_\ell\,$
generate $M$ as a $\e\hh_\ka\e$-module.
Given an $\e\hh_\ka\e$-linear map
$f: M\to M$, we may find an integer $m\geq 0$ such that
for all $i=1,\ldots,\ell,$ one has: $f(u_i)\in F_{m+d_i}M$.
It follows that there exists an integer $m$ such that
for all $n\in\Z$ one has: $f(F_nM)\subset F_{n+m}M$.
Thus, we may define an increasing $\Z$-filtration $\ff_\bullet\End$
on
$\End_{\e\hh_\ka\e}(\hh_\ka\e)$ by setting:
$\,\ff_m\End = \{f\in \End_{\e\hh_\ka\e}(\hh_\ka\e)
\;\;\big|\;\;f(F_nM)\subset F_{n+m}M\,,\,\forall n\in\Z\}.$
This makes $\End_{\e\hh_\ka\e}(\hh_\ka\e)$ a filtered ring.
Observe further that $\dim_{_\C}(\ff_m\End) <\infty$,
for any $m\in \Z$, and
if $d={\mathtt{Max}}(\deg u_1,
\ldots, \deg u_\ell),\,$ then we necessarily have:
$\ff_n\End=0$, for all $n<-d$.

Now, the canonical action-map
$\eta: \hh_\ka \to \End_{_{\e\hh_\ka \e}}(\hh_\ka \e)
$ is a filtration preserving
map of filtered vector spaces.
Thus, to show that $\eta$ is an isomorphism, it  suffices
to show that so is $\grd(\eta)$.
To this end, we consider the composite map:
$$ SV\#\G=
\grd(\hh_\ka) \;\stackrel{\grd(\eta)}{\too}\;
\grd\bigl(\End_{_{\e\hh_\ka \e}}(\hh_\ka \e)\bigr)
\;\stackrel{j}{\too}\;
\End_{_{\grd(\e\hh_\ka \e)}}\grd(\hh_\ka \e)\,=\,
\End_{_{(SV)^\G}}(SV)\,,$$ 
where $j$ is the canonical map, which is clearly {\it injective}.
It suffices to  prove  that the composite map
$j\ccirc \grd(\eta)$ is bijective.
The injectivity is obvious, since $j\ccirc \grd(\eta)$ becomes an isomorphism
after tensoring with $\C(V^*)^\G$,
the field of $\G$-invariant rational functions on $V^*$.
To prove the surjectivity, we need to show that any
endomorphism of $SV$ over $(SV)^\G$ comes from
an element of $SV\#\G$.

Let $a:SV\to SV$ be an $(SV)^\G$-linear map.
Then it defines a $\C(V^*)^\G$-linear map
$\C(V^*)\to \C(V^*)$. So $a=\sum a_g\cdot  g$, where
$a_g\in \C(V^*),g\in \Gamma$. It is clear
that the functions
$a_g$ must be regular on $V\smallsetminus Y$, where $Y$ is the set of points
with a nontrivial isotropy group in $\Gamma$
(since on $V\smallsetminus Y$, the surjectivity is obvious).
But $Y$ has codimension $\ge 2$, so by Hartogs theorem
$a_g$ are regular, as desired.

To prove (iii), we consider the multiplication pairing:
$\,\e\hh_\ka\otimes \hh_\ka\e \too\e\hh_\ka\e.$ This pairing
induces a right ${\e\hh_\ka\e}$-linear map $\psi:\e\hh_\ka\to
\Hom_{{\e\hh_\ka\e}}(\hh_\ka\e\, ,\,{\e\hh_\ka\e}),\,$
given  by  $\psi(x)(y)=xy$. This map is injective by part (iv),
since $\psi$ is the restriction of $\eta$ to $\e \hh_\ka$. It is
also surjective. Indeed, by part (iv), for any $a\in
(\hh_\ka\e)^\vee$, there exists $x\in \hh_\ka$ such that
$xy=a(y)$ for all $y\in \hh_\ka\e$. But then we also have $\e
xy=a(y)$, since $a(y)\in \e\hh_\ka\e$. So $a=\psi(\e x)$.
Hence,
$\psi$ is an isomorphism, and
 the first statement is proved. The second statement of part (iii)
follows similarly.
The theorem is proved.
\sq\medskip

\begin{remark}
There is an alternative proof of the isomorphisms in
parts (ii) and (iii) of  Theorem \ref{Coh_Mac_intro}, 
along the same lines as the proof
of Lemmas 1.3, 1.4 in [CBH]. That proof is less elementary and
involves a lot of non-commutative algebra machinery.
$\enspace\lozenge$
\end{remark}\medskip

  Note that the algebra $\hh_{t,0}$
is isomorphic to $\Weyl_t\lst$, provided $t\neq 0$. The latter is
known, see e.g. [Mo], to be
a simple algebra with trivial center. It follows by semi-continuity
that, for generic values  $(t,c)\in {\sf{\overline{C}}}$,
the algebra  $\hh_{t,c}$  also has the trivial center.
We will see below that the points $(0,c)$ are  very special:
the algebra $\hh_{0,c}$ has a large center, $\ZZ_{0,c}:=\ZZ(\hh_{0,c})$.

Further, for any  $c\in \CC$,
the construction of Hayashi [Ha], cf.  \S15,
applied to the  flat family $\,\{\hh_{t,c}\}_{t\in\C}\,,$ 
 gives a Poisson algebra structure
 on $\ZZ_{0,c}$.
\begin{theorem}[Satake isomorphism]\label{ZeHe}\label{ZtoeHe}
For any $c\in \CC$,
 the map: $\,
\ZZ_{0,c}\to \e\hh_{0,c} \e\,,\,z\mapsto z\cdot\e,\,$ is a Poisson
algebra isomorphism.
\end{theorem}

\Pf. It is clear that the map: $z\mapsto z\cdot\e$ gives an
algebra homomorphism $\zeta: \ZZ(\hh_{0,c})\to \e\hh_{0,c} \e$.
Since the algebra $\e\hh_{0,c} \e$ is commutative, we
can also construct  an
algebra homomorphism $\xi: \e\hh_{0,c} \e\to \ZZ(\hh_{0,c})$
as follows. For any, $a \in\e\hh_{0,c} \e$, right multiplication
by $a$ gives, due to commutativity
of $\e\hh_{0,c} \e$, a right $\e\hh_{0,c} \e$-linear endomorphism of $\hh_{0,c} \e$.
This endomorphism must arise from the left
action of an element $\xi(a)\in \hh_{0,c}$, since
$\End_{_{\e\hh_{0,c} \e}}(\hh_{0,c} \e)=\hh_{0,c} ,$
by Theorem \ref{Coh_Mac_intro}(ii). Moreover,
the element $\xi(a)$ belongs to the
center of $\hh_{0,c}$ since the left
$\xi(a)$-action on $\hh_{0,c} \e$, being induced by multiplication
on the right, commutes with
the left action on $\hh_{0,c} \e$ of any other element
of $\hh_{0,c}$.
It is clear that the assignment: $a\mapsto \xi(a)$ gives an
algebra homomorphism $\xi: \e\hh_{0,c} \e\to \ZZ(\hh_{0,c})$,
and that this homomorphism is inverse to $\zeta$.

The compatibility of the map: $z\mapsto z\cdot\e$ with the Poisson
brackets follows directly from the "deformation-construction"
of the brackets on $\e\hh_{0,c} \e$ and on $\ZZ(\hh_{0,c})$
given in \S15, once we know that the integer ${\mathbf{m}}_{0,c}$
involved in the construction of the bracket on  $\e\hh_{0,c}\e$
equals 1,
cf. Lemma \ref{commut}. 
But this has been established in the course of the proof of
Theorem \ref{comm_conj}. \quad\sq\medskip

\noindent
{\bf Remark.} The isomorphism  of Theorem
\ref{ZeHe} is analogous to the Satake isomorphism
between the center and the spherical
subalgebra of the affine Hecke algebra,
due to Lusztig~\cite{Lu}.$\quad\lozenge$

\begin{lemma} For $t=c=0,$ we have: $\ZZ(\hh_{_{0,0}})=(SV)^\G$.
\end{lemma}

\Pf . It is clear that any element
of $SV\,\#\,\G $ that commutes with
the subalgebra $SV$
belongs to $SV$. Hence, $\ZB\subset
SV$. Now, an element of
$SV$ commutes with $\G$ if and only if
it is $\G$-invariant.
\sq\medskip

The imbedding: $\,\ZZ_{0,c}\into \hh_{0,c}$ induces a canonical
 map: $\,\grd(\ZZ_{0,c})\to\ZZ(\hh_{0,0})=\ZB$.
\begin{theorem}\label{center} The canonical
map is an algebra isomorphism:
$\grd(\ZZ_{0,c})\simeq\ZB=(SV)^\G.$
\end{theorem}

\Pf . Write $\varphi: z\mapsto z\cdot\e$ for the map
of Theorem \ref{ZtoeHe},
and $\overline{\varphi}: \ZB\to\e\hb\e$ for
a similar map for the algebra $\hb=\grd(\hh_{0,c})$.
It is clear that the maps $\varphi$ and $\overline{\varphi}$
are algebra
 homomorphisms. Observe that
 the associated graded homomorphism:
$\grd(\varphi): \grd(\ZZ_{0,c})\too \grd(\e\hh_{0,c}\e)$
can be identified, in view of the PBW-theorem,
with the composite map:
$$
\grd(\ZZ_{0,c})\into \ZB=SV^\G\;
\stackrel{\overline{\varphi}}{\too}
\e\hb\e\;=\e\cdot\bigl(SV\lst \bigr)
\cdot\e\,.
$$
It is immediate that the second map, $\overline{\varphi}$, is an
isomorphism. It follows that the map $\grd(\varphi)$ is
injective. Now we use the
following obvious claim:

\noindent
{\it
Let $\,E=\cup_{i\geq 0}\,E_i\,$ and $\,E'=\cup_{i\geq 0}\,E'_i\,$
be two filtered vector spaces such that, for any $i\geq 0$, we
have $\dim E_i=\dim E'_i <\infty.\,$ Assume $\,\varphi: E\to E'$
is a filtration preserving linear map  such that $\grd(\varphi)$
is injective. Then both $\grd(\varphi)$ and $\varphi$ are
bijective maps.}

Applying this  observation in our situation, we
conclude that the map $\grd(\varphi)$ gives an isomorphism
 $\grd(\ZZ_{0,c})\iso\ZB$, since
$\grd(\ZZ_{0,c})\subset \ZB$ and the two sides
have equal Poincar\'e series (by Theorem \ref{ZeHe}).
\hfill\sq\break

The argument in the proof Theorem \ref{center} above  yields the following
strengthening of Theorem \ref{ZtoeHe}

\begin{proposition}\label{satake_gr}
The associated graded map: $\,\grd(\varphi):
\grd(\ZZ_{0,c})\to \grd(\e\hh_{0,c} \e)\,$ is an
algebra isomorphism.\hfill\sq
\end{proposition}
\smallskip

\noindent
{\bf Proof of theorem \ref{proj_intro}:}\quad
Recall that $U$ is a smooth Zariski open dense 
subset in $\Spec\ZZ_{0,c}$, and
$\rr_{_U}$ denotes the coherent
sheaf on  $U$ corresponding to $\hh_{_U} \e$,
viewed as a (finitely generated) $\ZZ_{_U}$-module, so that
$\Gamma(U, \rr_{_U})=\hh_{_U} \e$.
The left $\hh_{_U}$-action  on $\hh_{_U} \e$
gives an $\hh_{_U}$-module structure on each geometric
fiber of $\rr_{_U}$.

Proving that the sheaf $\rr_{_U}$ is locally free
amounts to showing that $\hh_{_U} \e$ is
a finite rank {\sf projective} $\ZZ_{_U}$-module.
To this end, recall the following standard result 
of commutative algebra. Let $A$
be a finitely generated commutative $\C$-algebra with unit,
such that $\Spec{A}$ is a smooth algebraic
variety. Then, {\it a finitely-generated $A$-module $M$
is Cohen-Macaulay if and only if it is projective},
see  \cite[Ch.4, Cor.2]{Se}.

Since  $U$ is smooth by our assumptions,
it follows from Theorem \ref{Coh_Mac_intro}(ii) and Satake isomorphism
that
 $\hh_{_U}\e$ is a finite rank projective $\ZZ_{_U}$-module.
 Part (i) now follows from
Theorem \ref{Coh_Mac_intro}(vi). 

Further,
the endomorphism algebra of an algebraic
 vector bundle on an affine algebraic
variety is always Morita equivalent to the coordinate ring
on this variety. Explicitly, in our case,
the equivalence is provided
by the following mutually inverse functors
$$
\hh_{_U}\text{-mod}\;
\underset{{}^F}{\overset{{}_G}{\rightleftarrows}}\;\ZZ_{_U}\text{-mod}\quad,
\quad F(M)= \hh_{_U} \e\otimes_{_{\ZZ_{_U}}}M\quad,
\quad G(N)=\e\hh_{_U}\otimes_{_{\hh_{_U}}}N =\e\cd N\,.
$$

To prove (vi), notice that since $\rr_{_U}$
is a vector bundle it is sufficient, due to
 rigidity of $\G$-modules, to establish
that the generic fiber  of $\rr_{_U}$ is isomorphic
to the regular representation of $\G$. But this 
follows from Lemma \ref{genchar}.\quad\sq\medskip

\noindent
{\bf Remarks.} \,\vi We have established 
a special case of the following well-known
  general result: {\it If  $A$
 is an associative algebra and
$\e\in A$  an idempotent, then $A$ is
 Morita
equivalent to $\e A\e$ if and only if we have: $A\e A = A$.}

\vii We have shown that each  geometric fiber of $\rr$ over 
the regular locus of $\Spec\ZZ_{0,c}$ is isomorphic, as a $\C\G$-module,
to the regular representation of $\G$.

\noindent
{\bf Proof of Proposition \ref{k_theory}.}\quad
Given an associative
algebra $A$, one has, see [Lo], a 
{\it Chern character} map $\,ch:
K(A)\to\HC(A)$, where $K(A)$ stands for the Grothendieck group
of finitely generated projective  $A$-modules,
and $\HC(A)$ denotes the even degree part of the
 periodic cyclic homology groups of $A$, see [Lo].
The Chern character map is compatible with Morita
equivalences. 
It is also known, see [Lo], 
 that if $A=\C[X]$ is the coordinate ring
of a smooth affine algebraic variety,
then: $\,\HC(A)\simeq \bplus_{\!i}\,H^{2i}(X),$
the direct sum of even cohomology of $X$, as a topological space.
Moreover, the map $ch$ reduces in this case to the
classical Chern character map: $K(X)\to \bplus_{\!i}\,H^{2i}(X)$.

Recall now that the algebra $\hh_{0,c}$ has a canonical increasing
filtration, $\,F_\bullet(\hh_{0,c}),$ such that $F_0(\hh_{0,c})=\C\G$.
Since $\hh_{0,c}$ is a projective $\C\G$-module, and the algebra  
$\grd^F_\bullet\hh_{0,c}$ has
finite homological dimension, a general result of Quillen [Q]
says that the functor: $V\mapsto \hh_{0,c}\otimes_{_{\C\G}}V$
induces an isomorphism of $K$-groups: $K(\G)\iso K(\hh_{0,c}).$
Similarly, a result of Block [Bl] says that
the imbedding: $\C\G\into \hh_{0,c}$ induces an isomorphism:
 $\HC(\G)\iso \HC(\hh_{0,c}).$
Thus, we obtain the following
commutative diagram: 
{\small{$$
{\small
\diagram
K(\C\G)\enspace\dto_{ch}\rdouble^<>(.5){\mbox{\footnotesize{Quillen}}}
&\enspace K(\hh_{0,c})\enspace\dto_{ch}
\rdouble^<>(.5){{_{\mbox{\footnotesize{Morita}}}}}_<>(.5)
{^{^{\mbox{\footnotesize{equiv.}}}}}
&\enspace K(\e\hh_{0,c}\e)\enspace \dto_{ch}
\rdouble^<>(.5){\mbox{\footnotesize{Satake}}}
&\enspace K(\ZZ_{0,c})\dto_{ch}\rdouble&
K(\Spec\ZZ_{0,c})\dto_{ch}
\\
\HC(\C\G)\enspace\rdouble^<>(.5){\mbox{\footnotesize{Block}}}
&\enspace\HC(\hh_{0,c})\enspace
\rdouble^<>(.5){{_{\mbox{\footnotesize{Morita}}}}}_<>(.5)
{^{^{\mbox{\footnotesize{equiv.}}}}}
&\enspace\HC(\e\hh_{0,c}\e)\enspace 
\rdouble^<>(.5){\mbox{\footnotesize{Satake}}}
&\enspace \HC(\ZZ_{0,c})\rdouble&
\bplus_{\!_i}\,H^{2i}(\Spec\ZZ_{0,c})
\enddiagram
}
$$}}
where Morita equivalence is guaranteed
by Theorem \ref{proj_intro}(ii).

To complete the proof, it suffices to show that the leftmost vertical
map in the diagram induces an isomorphism
$\,ch_{_\G}: \C\otimes_{_\Z}K(\C\G)\too \HC(\C\G)$. But the algebra $\C\G$
is Morita equivalent to a direct sum of several copies
of the field $\C$. Hence
one may replace  $\C\G$ in the above map $ch_{_\G}$ by 
the
algebra $\C$. The Chern character map: $\C\otimes_{_\Z}K(\C)
\too\HC(\C)$ is clearly an isomorphism, and
we are done.
\quad\qed\medskip

\begin{lemma}\label{normal}\vi For any $c\in\CC$, the variety
$\Spec \ZZ_{0,c}$ is normal.

\vii If the  variety
$\Spec \ZZ_{0,c}$ is smooth, then the Poisson structure on
$\ZZ_{0,c}$ is non-degenerate, i.e., 
$\Spec \ZZ_{0,c}$ is a symplectic manifold.
\end{lemma}

\Pf . Note that, for $c=0$, the singular locus of 
$\Spec \ZZ_{0,0}=V/\G$ has codimension $\ge 2$ in
$\Spec \ZZ_{0,0}$. The codimension of the  singular locus
being semi-continuous under going to the
asymptotic cone it follows that,
for any $c\in \CC$,  the singular locus of 
$\Spec \ZZ_{0,c}$ has codimension $\ge 2$ in
$\Spec \ZZ_{0,c}$. Further, we know that 
$\Spec \ZZ_{0,c}$ is a Cohen-Macaulay variety, by
Theorem \ref{Coh_Mac_intro}(i). But any
Cohen-Macaulay variety which is non-singular in codimension
one is normal, see e.g. proof in [CG, \S2.2].
Part (i) follows.

To prove (ii), assume
 $\Spec \ZZ_{0,c}$ is smooth, 
write $\beta_c$ for the bivector 
on $\Spec \ZZ_{0,c}$ giving the  Poisson structure,
and set $2d=\dim V= \dim(\Spec \ZZ_{0,c})$.
The  set of points where the Poisson structure degenerates 
equals the set  of points where the $2d$-vector
$\wedge^d\beta_c$ vanishes. This set is clearly a divisor
(a codimension one closed subvariety) in $\Spec \ZZ_{0,c}$.
Let $J\subset \ZZ_{0,c}$ denote the defining ideal of
that subvariety, and $\grd(J) \subset
\grd(\ZZ_{0,c}) = \C[V]^\G$, the corresponding associated
graded ideal with respect to the canonical filtration
on $\hh_{0,c}$. Then $\grd(J)$ is a proper non-zero ideal
in $\C[V]^\G$, 
and the zero-variety of $\grd(J)$ is a 
codimension one subvariety in $V/\G$, the
`asymptotic cone' to the zero-set  of $\wedge^d\beta_c$.
It is clear that $\wedge^d\beta_0$ is the 
leading homogeneous component of $\wedge^d\beta_c$,
hence, vanishes on this
asymptotic cone. Hence, the Poisson structure
on $V/\G$ must degenerate at a divisor.
On the other hand,
 the Poisson structure on $V/\G$ is non-degenerate
at any smooth point, and the singular locus of 
 $V/\G$ has codimension $\ge 2$ in $V/\G$.
The contradiction completes the proof.\quad\qed\medskip

For any  $c\in \CC$, let $\rep_{_{\C\G}}(\hh_{0,c})$  be the
variety of all algebra homomorphisms: $\hh_{0,c}\too \End_{_\C}(\C{\G})$,
whose restriction to $\C{\G}$ is the $\C{\G}$-action 
by left multiplication. This is
an affine algebraic variety.
Given an $\hh_{0,c}$-module $M\in \rep_{_{\C\G}}(\hh_{0,c})$, 
the space $\e M$ is
clearly a line which is stable under the action of $\ZZ_{0,c}$.
 The action on this line gives an
algebra homomorphism $\chi_{_M}:\ZZ_{0,c}\to \C$. 
 It is clear that the
assignment: $M\mapsto \chi_{_M}$ gives
a morphism of algebraic varieties,  
\begin{equation}\label{con_hom}
\pi:\; 
\rep_{_{\C\G}}(\hh_{0,c}) \too \Spec \ZZ_{0,c}\,.
\end{equation}

Let 
$\Aut_{_\G}(\C{\G})$ be the group of $\C$-linear automorphisms of
the vector space $\C{\G}$ that commute with left multiplication by
$\G$.
 The
 group $\Aut_{_\G}(\C{\G})$ is reductive,
and it acts naturally on $\rep_{_{\C\G}}(\hh_{0,c}),$
preserving the fibers of the map (\ref{con_hom}).

\begin{theorem}\label{components}
\vi There is a unique irreducible component,
$\rep_{_{\C\G}}^\circ(\hh_{0,c}),$ of the variety
$\rep_{_{\C\G}}(\hh_{0,c})$ 
whose image under (\ref{con_hom}) is dense in
$\Spec \ZZ_{0,c}$. The generic point
of $\rep_{_{\C\G}}^\circ(\hh_{0,c})$ is a
simple $\hh_{0,c}$-module.

\vii The  map  (\ref{con_hom}) induces an algebra isomorphism:
$\ZZ_{0,c} \iso \C[\rep_{_{\C\G}}^\circ(\hh_{0,c})]^{
\Aut_{_\G}(\C{\G})}.$
\end{theorem}
\Pf .
It follows from theorem \ref{proj_intro}
that any point of $\rep_{_{\C\G}}(\hh_{0,c})$
containing, as a subquotient, a simple
$\hh_{0,c}$-module on which $\ZZ_{0,c}$ acts by
a character that belongs to the smooth locus of
$\Spec \ZZ_{0,c}$ is itself simple and, moreover,
such points form a Zariski open  subset in
 $\rep_{_{\C\G}}(\hh_{0,c})$. This set is irreducible,
since so is $\Spec \ZZ_{0,c}$, hence its closure,
$\rep_{_{\C\G}}^\circ(\hh_{0,c})$, is an irreducible component
of  $\rep_{_{\C\G}}(\hh_{0,c})$. This proves (i).

To prove (ii), notice that the isotropy group 
in $\Aut_{_\G}(\C{\G})$ of a simple $\hh_{0,c}$-module is reductive.
Hence the corresponding $\Aut_{_\G}(\C{\G})$-orbit in
$\rep_{_{\C\G}}(\hh_{0,c})$ is closed. It follows that
generic $\Aut_{_\G}(\C{\G})$-orbits in $\rep_{_{\C\G}}^\circ(\hh_{0,c})$
are closed, and the induced map
$\pi: 
\dis \rep_{_{\C\G}}^\circ(\hh_{0,c})/\Aut_{_\G}(\C{\G})=
\Spec \C[\rep_{_{\C\G}}^\circ(\hh_{0,c})]^{\Aut_{_\G}(\C{\G})}
\too \Spec \ZZ_{0,c}$ is a birational isomorphism.
Thus, any element
$f\in \C[\rep_{_{\C\G}}^\circ(\hh_{0,c})]^{\Aut_{_\G}(\C{\G})}$
descends to a rational function
$\psi$ on $\Spec \ZZ_{0,c}$. Moreover, Theorem \ref{proj_intro}
implies that $\psi$ is regular at any smooth point of 
$\Spec \ZZ_{0,c}$ (and the complement of smooth points has codimension
$\geq 2$). Thus, $\psi$ extends to a regular
function on the whole of $\Spec \ZZ_{0,c}$, due to normality
of the latter, see Lemma \ref{normal}(i).
Hence, $\psi$ is represented by an element
of $\ZZ_{0,c}$, i.e. there is $z\in \ZZ_{0,c}$
such that the functions $f$ and $\pi^*(z)$ coincide
on an open dense subset of $\rep_{_{\C\G}}^\circ(\hh_{0,c})$.
Thus, $f=\pi^*(z)$ on  $\rep_{_{\C\G}}^\circ(\hh_{0,c})$,
and we are done.\quad\qed\medskip

Given an associative algebra $\A$, we write $[\A,\A]$ for the
$\C$-vector subspace in $\A$ spanned by the commutators:
$\,[a,b]_{a,b\in \A}\,$. Here is a useful criterion
for the smoothness of $\Spec \ZZ_{0,c}$.

\begin{proposition}\label{criterion}
For any $c\in \CC$, and a Zariski open affine subset 
$U\subset  \Spec \ZZ_{0,c}$,
the following conditions are equivalent:

\vi The dimension of any simple $\hh_{_U}$-module is
$\geq |\G|$.

\vii $U$ is a smooth algebraic variety.

\viii For any $g\in \G\,,\, g\neq 1,\,$ we have: $\;g\in
[\hh_{_U}, \hh_{_U}]$.

\vvi For any simple $\hh_{_U}$-module $M$, and
 any $g\in \G\,,\, g\neq 1,\,$ we have: $\Tr(g,M)=0$.

\noindent
In particular, if any of the conditions 
{\sf{(i), (iii)}} or {\sf{(vi)}} holds for $U=\Spec \ZZ_{0,c}$,
then $\Spec \ZZ_{0,c}$ is smooth.
\end{proposition}

We refer the reader to \S16 for 
some applications of Proposition \ref{criterion}.

\begin{remark}  It will be clear from 
the proof of the Proposition that conditions (i)-(vi) 
imply automatically that the
restriction to $\G$ of any simple $\hh_{_U}$-module is
isomorphic
 to the regular representation
of $\G$.$\quad\lozenge$
\end{remark}

Our proof of the Proposition will
 be based on the following important construction
 that will also play a role later, in \S11.

Let $\widehat{\rr}: \C\G\times\rep_{_{\C\G}}^\circ(\hh_{0,c})
\too\rep_{_{\C\G}}^\circ(\hh_{0,c})$ be the trivial vector bundle
on the variety $\rep_{_{\C\G}}^\circ(\hh_{0,c})$,
viewed as an $\Aut_{_\G}(\C{\G})$-equivariant
 $\oo$-sheaf with geometric fiber $\C\G$
(so that the group $\Aut_{_\G}(\C{\G})$ acts 
simultaneously on both $\rep_{_{\C\G}}^\circ(\hh_{0,c})$
and $\C\G$).
There is a canonical $\hh_{0,c}$-module structure
on $\widehat{\rr}$ defined as follows. A point
$\xi\in \rep_{_{\C\G}}^\circ(\hh_{0,c})$ is by definition
an algebra homomorphism $\xi: \hh_{0,c} \to \End_{_\C}(\C\G).$
We let $\hh_{0,c}$ act on $\widehat{\rr}_\xi\simeq\C\G$,
the fiber of $\widehat{\rr}$ at $\xi$, via the 
homomorphism $\xi$. The  $\hh_{0,c}$-module structure
on $\widehat{\rr}$ thus defined clearly commutes with
the $\Aut_{_\G}(\C{\G})$-action.

Let $\PAut_{_\G}(\C{\G})=\Aut_{_\G}(\C{\G})\big/\C^*$ be
the quotient modulo the scalars.
We have the natural projection: $\Aut_{_\G}(\C{\G}) \onto \PAut_{_\G}(\C{\G})$,
and the $\Aut_{_\G}(\C{\G})$-action on $\rep_{_{\C{\G}}}(\hh_{0,c})$
(but not on $\C{\G}$) factors throught the group $\PAut_{_\G}(\C{\G})$.
We define a canonical
splitting: $\PAut_{_\G}(\C{\G})\into \Aut_{_\G}(\C{\G})$
of the natural projection
by the requirement, that the image of $\PAut_{_\G}(\C{\G})$ in
$\Aut_{_\G}(\C{\G})$ acts trivially on the element $\e\in \C{\G}$.
This splitting makes it possible to view
the bundle: $\C{\G}\times\rep_{_{\C{\G}}}^\circ(\hh_{0,c})
\to
\rep_{_{\C{\G}}}^\circ(\hh_{0,c})$
as a $\PAut_{_\G}(\C{\G})$-equivariant 
$\hh_{0,c}$-vector bundle.
\medskip

\noindent
{\bf Proof of Proposition \ref{criterion}:}$\;$
\vi$\;\Longrightarrow\;$\vii 
If all simple $\hh_{_U}$-modules have
dimension $\geq |\G|$ then each point
$\xi\in \rep_{_{\C\G}}^\circ(\hh_{_U})$ must represent a simple
$\hh_{_U}$-module. Hence, the isotropy group
of $\xi$ in $\Aut_{_\G}(\C{\G})$ reduces to scalars
$\C^*\subset\Aut_{_\G}(\C{\G})$,
by Schur lemma. Therefore,
$\,\PAut_{_\G}(\C{\G})$,
a  reductive group,
acts freely on $\rep_{_{\C\G}}^\circ(\hh_{_U})$.
Hence, all $\PAut_{_\G}(\C{\G})$-orbits in
$\rep_{_{\C\G}}^\circ(\hh_{_U})$
are closed, and 
are exactly the fibers of the projection
$\,\pi_{_U}
:\Spec \rep_{_{\C\G}}^\circ(\hh_{_U})\to
\Spec
\bigl(\C[\rep_{_{\C\G}}^\circ(\hh_{_U})]^{\PAut_{_\G}(\C{\G})}\bigr).\,$
It follows that $\Spec \rep_{_{\C\G}}^\circ(\hh_{_U})$ is
a principal $\PAut_{_\G}(\C{\G})$-bundle over
$\,\Spec
\bigl(\C[\rep_{_{\C\G}}^\circ(\hh_{_U})]^{\PAut_{_\G}(\C{\G})}\bigr)\,$
(in \'etale topology).
Therefore,
 the $\PAut_{_\G}(\C{\G})$-equivariant
vector bundle $\widehat{\rr}_{_U}$ 
 on $\rep_{_{\C\G}}^\circ(\hh_{_U})$ 
descends to a well-defined
algebraic vector bundle $\,\widetilde{\rr}_{_U}$ on
$
\Spec
\bigl(\C[\rep_{_{\C\G}}^\circ(\hh_{_U})]^{\Aut_{_\G}(\C{\G})}\bigr).\,$
Moreover, the $\hh_{_U}$-module structure
on $\widehat{\rr}_{_U}$ descends to
an $\hh_{_U}$-module structure
on $\widetilde{\rr}_{_U}$.

Further, by Theorem \ref{components}(ii), we
may (and will) identify $\,\Spec
\bigl(\C[\rep_{_{\C\G}}^\circ(\hh_{_U})]^{\Aut_{_\G}(\C{\G})}\bigr)$
with $U$, hence
view  $\widetilde{\rr}_{_U}$ as a locally free coherent
sheaf on  $U$. The (constant) section
of $\widehat{\rr}_{_U}$ corresponding to the element
$\e \in \C\G$ gives rise to a non-vanishing regular section
of $\widetilde{\rr}_{_U}$, to be denoted $\widetilde{\e}$.
Using the $\hh_{_U}$-action on $\widetilde{\rr}_{_U}$,
we define an $\hh_{_U}$-module map:
$\hh_{_U}\e \too\G(U, \widetilde{\rr})$
by the assignment: $u\cd\e \mapsto u\cd\widetilde{\e}$.
This gives an $\hh_{_U}$-equivariant morphism
$\eps: \rr_{_U} \to \widetilde{\rr}_{_U}$
of coherent sheaves on $U$.
The morphism $\eps$ is clearly an isomorphism over the
generic point of $U$, hence is injective
since the sheaf $\rr_{_U}$, being reflexive, is torsion free.
On the other hand, for any $\chi \in U$,
the induced map on the fibers:
$\rr_\chi \to \widetilde{\rr}_\chi$
is a non-zero map of $\hh_{_U}$-modules,
hence is surjective, since $\widetilde{\rr}_\chi$ is a simple 
$\hh_{_U}$-module by our hypothesis (i).
Therefore, $\eps$ is  an isomorphism, by  Nakayama lemma.  
Thus, the sheaf $\rr_{_U}$ is also  locally free.

Recall that, according to a well known criterion, 
a finitely generated commutative algebra $\A$ is smooth
(i.e. the scheme $\Spec\A$ is smooth)
if and only if $\A$ has
finite (global) homological dimension.
Thus, 
$U$ is smooth if and only if $\ZZ_{_U}$ has
finite (global) homological dimension.
On the other hand, since $SV$ is smooth,
the algebra $\grd(\hh_{0,c})\simeq
SV\#\G$, clearly has finite homological dimension.
It follows that both the algebra $\hh_{0,c}$ and
$\hh_{_U}$ have  finite homological dimension.
Further, by Theorem
\ref{proj_intro} we have: $\hh_{_U}\simeq \End_{_{\ZZ_{_U}}}\rr_{_U}$.
It follows that, if $\rr_{_U}$ is a vector bundle, then
the algebra $\hh_{_U}$ is Morita equivalent to
$\ZZ_{_U}$. Hence the algebras $\ZZ_{_U}$ and $\hh_{_U}$
have equal
 global homological dimensions, and \vii follows.

\vii$\;\Longrightarrow\;$\viii Set $\ZZ=\ZZ_{_U}.\,$
We know that
$\hh_{_U}\simeq \End_{_{\ZZ}}\rr_{_U}$, by  Theorem \ref{proj_intro}.
Note that since
any element of $[\End_{_{\ZZ}}\rr_{_U}, 
\End_{_{\ZZ}}\rr_{_U}]$ acts by a trace zero
operator in each fiber of $\rr_{_U}$, we have: $\ZZ
\cap [\End_{_{\ZZ}}\rr_{_U}, \End_{_{\ZZ}}\rr_{_U}]=0.\,$ Further,
the algebra $\End_{_{\ZZ}}\rr_{_U}$ being Morita equivalent
to $\ZZ$, the imbedding:
$\ZZ\into \hh_{_U}=\End_{_{\ZZ}}\rr_{_U}$ induces
a bijection: $\ZZ=\ZZ/[\ZZ,\ZZ]
\iso  \hh_{_U}/[\hh_{_U},\hh_{_U}].\,$
It follows that there is a direct sum
decomposition:
$$ \hh_{_U} = \ZZ_{_U}\;\bigoplus\;[\hh_{_U},\hh_{_U}].\,$$

Now fix $g\in \G\,,\, g\neq1.$
Using the decomposition above we can write: $g= z+ u,$
where $z\in \ZZ_{_U},$ and $u \in [\hh_{_U},\hh_{_U}]$.
Since the trace of $g$ 
in the regular representation vanishes,
the element $g= z+ u$ acts by a trace zero
operator in each fiber of $\rr_{_U}$, by
Theorem  \ref{proj_intro}(vi). It follows that $z$=0.
Thus, $g=u \in [\hh_{_U},\hh_{_U}]$.

The implication \viii$\;\Longrightarrow\;$\vvi is clear.

\vvi$\;\Longrightarrow\;$\vi 

Assume that for any $g\in \G\,,\, g\neq 1,$
 the action of $g$ in
any simple (in particular, finite-dimensional)
$\hh_{_U}$-module $M$ has vanishing trace: $\Tr(g,M)=0$.
It follows from representation theory of finite groups
that the restriction
of any simple $\hh_{_U}$-module $M$ to $\G$ is isomorphic to
a multiple of the regular representation of $\G$. 
\qed\medskip

\begin{remark} There is an alternative proof of the Morita equivalence
in the implication \vi$\;\Longrightarrow\;$\vii of the Proposition
above, based on a general result of M. Artin
 (1968) which in the special situation we need reads:

\noindent
{\it Let $A$ be a $\C$-algebra with center $\ZZ$, finitely generated as 
a $\ZZ$-module, and such that:\newline
\noindent
$\bullet\enspace$  $\ZZ$ has no zero divisors; 
\newline
\noindent
$\bullet\enspace$  $A$ injects into $A\otimes_{_{\ZZ}}Q(\ZZ)$,
where $Q(\ZZ)$ is the field of fractions of $\ZZ$;
\newline
\noindent
$\bullet\enspace$  Any irreducible $A$-module has dimension $N$; 
\newline
\noindent
$\bullet\enspace$  $A\otimes_{_{\ZZ}}Q(\ZZ)=
{\mathtt{Mat}}_{\!_{N\!}}\bigl(Q(\ZZ)\bigr).$

Then $A$ is Morita equivalent to $\,\ZZ.$ }\quad\qed

\noindent
The proof of this result makes a  nontrivial
use of the theory of polynomial identities.$\quad\lozenge$
\end{remark}
\medskip

\noindent
{\bf Proof of part (ii) of
theorem \ref{mckay}:}\quad The first claim of part (ii)
is just Lemma \ref{normal}(i).

To proceed further, 
fix  $c\in \CC$ such that
$\Spec(\e\hh_{0,c}\e)$ is smooth.
The family of 
non-commutative associative
algebras $\,\{\hh_{t,c}\}_{t\in \C}\,$
may be viewed as a single $\C[t]$-algebra, which is
a polynomial deformation
 of the algebra $\hh_{0,c}$.
Let ${\sf{\hat{H}}}_{t,c}$ 
denote the $t$-adic completion of that  $\C[t]$-algebra,
a formal flat deformation of  $\hh_{0,c}$.
Thus, $\e{\sf{\hat{H}}}_{t,c}\e$ is 
a formal flat deformation of  $\e\hh_{0,c}\e$.
Recall   that this deformation gives a Poisson algebra structure
on
$\e\hh_{0,c}\e$. Since  $\Spec(\e\hh_{0,c}\e)$ is smooth,
hence a symplectic manifold,
we may apply the results of Appendix D
to that  deformation.

Write $\KK=\C((t))$ for the field of formal
Laurent  series,
and let $\BBA:= \KK\widehat{\otimes}_{_{\C[[t]]}}{\sf{\hat{H}}}_{t,c}$
be the  $\KK$-algebra obtained from
${\sf{\hat{H}}}_{t,c}$ by localizing at $t$. 
Thus, Lemma \ref{split} and
Proposition \ref{konts} applied
to the algebra $\,\A=\e{\sf{\hat{H}}}_{t,c}\e\,$  yield
a graded $\KK$-algebra isomorphism:
\begin{equation}\label{tangent3}
H\!H^\bullet(\e\BBA\e)\; \simeq\; \KK\otimes_{_\C}
H^\bullet\bigl(\Spec(\e\hh_{0,c}\e)\bigr)
\,.
\end{equation}

Theorem 
\ref{proj_intro} and Proposition \ref{morita_e} imply
that the $\KK$-algebra
$\e\BBA\e$ is Morita equivalent to $\BBA$.
Hence these two algebras have canonically isomorphic
Hochschild cohomology. Thus,  part (i) of
Theorem \ref{mckay} and isomorphism (\ref{tangent3})
yield the following chain of $\KK$-algebra isomorphisms,
that completes the proof of part (ii) of Theorem \ref{mckay}
$$\KK\otimes_{_\C}
\grd^F_\bullet\ZZ\G \iso H\!H^\bullet(\BBA)
\;\simeq\; \HH^\bullet(\e\BBA\e)\;\simeq \;
\KK\otimes_{_\C}H^\bullet\bigl(\Spec(\e\hh_{0,c}\e)\bigr)\,.
\qquad\square
$$

Let $\Gamma\subset Sp(V)$ be any finite group, and
 $\G^{\circ}\subset \G$   the subgroup generated by 
all symplectic reflections in $\G$. Write $\hh_{t,c}^{^{_\circ}}$
 for the symplectic reflection
algebra constructed out of $V$ and $\G^{\circ}$
(instead of $\G$), i.e.,
$\hh_{t,c}^{^{_\circ}}
 \subset \hh_{t,c}$ is 
the subalgebra generated by $V$ and $\G^{\circ}$.
Then $\G^{\circ}$ is a normal subgroup of $\G$, and 
the algebra $\hh_{t,c}$ is an
`extension' of $\hh_{t,c}^{^{_\circ}}$ by $\G/\G^{\circ}$.
Specifically,  we have: $\hh_{t,c} \,\simeq\, 
\hh_{t,c}^{^{_\circ}}\,\#_{_{\G^{\circ}}} \G\,.$
Here, for any normal subgroup $H \lhd \G$ and any
$\G$-algebra $A$ containing $\C{H}$, we write
$A\#_{_H} \G$ for the vector space $A\bigotimes_{_{\C{H}}}\,\C{\G}$ 
equipped with the natural 
algebra structure. Thus,  understanding general symplectic reflection 
algebras reduces to understanding those associated with groups $\G$
 generated by 
symplectic reflections.

Verbitsky has shown in [Ve]
that  a necessary condition for the existence of a smooth
 symplectic resolution of singularities
$p :\widehat{V/\Gamma}\onto V/\Gamma$
is that the group $\Gamma$ is
generated by symplectic reflections.
We have an analogue of that result for 
varieties  $\Spec \ZZ_{0,c}$.
Specifically, given any finite group
 $\Gamma\subset Sp(V)$, fix $c\in \CC$
and let $\hh_{0,c}$ be the 
corresponding symplectic reflection algebra.
Let $\G^\circ
\subset \Gamma$ be the subgroup generated by those symplectic 
reflections $s\in S$ for which $c_s\ne 0$. 

\begin{proposition}\label{verbit}
 If $\Spec \ZZ_{0,c}$ is smooth then $\G^\circ=\Gamma$.
\end{proposition}

\Pf . It is clear that $\G^\circ$ is a normal subgroup in $\Gamma$.
Set $K=\Gamma/\G^\circ$.
Let $\hh_{0,c}^{^{_\circ}}$ be the symplectic 
reflection 
algebra corresponding to the group $\G^\circ$ (as opposed 
to $\Gamma$). As we have mentioned above
 the algebra $\hh_{0,c}$ is an extension of 
$\hh_{0,c}^{^{_\circ}}$ by $K$, and in particular it is a free 
$\hh_{0,c}^{^{_\circ}}$-module (both left and right) of rank $|K|$. 

Let $\ZZ_{0,c}^{^{_\circ}}$ be the center of $\hh_{0,c}^{^{_\circ}}$. Then
$K$ acts on $\ZZ_{0,c}^{^{_\circ}}$ by conjugation, and  
$\ZZ_{0,c}=\bigl(\ZZ_{0,c}^{^{_\circ}}\bigr)^K$. 
So if $\mm_c$ and $\mm_c^\circ$ denote the spectra 
of these algebras then $\mm_c=\mm_c^\circ/K$. 

We have assumed that $\mm_c=
\Spec \ZZ_{0,c}$ is smooth. We claim that in this case 
$\mm_c^\circ$ is smooth also. Indeed, by Proposition 
\ref{criterion}, every finite dimensional $\hh_{0,c}$-module
 is of dimension at least $|\Gamma|$.
If there were a
 finite dimensional  $\hh_{0,c}^{^{_\circ}}$-module
of 
dimension $<|\G^\circ|$, one could induce from
$\hh_{0,c}^{^{_\circ}}$ to $\hh_{0,c}$, and 
get a contradiction. Thus, all
simple  $\hh_{0,c}^{^{_\circ}}$-modules
are of dimension $\ge|\G^\circ|$, hence by the same 
proposition we have $\mm_c^\circ$ is smooth. 
We claim that $K$ acts freely on $\mm_c^\circ$.
Indeed, assume  for some
$\gamma\in K\,,\,\gamma\neq 1,$ the fixed point set
$(\mm_c^\circ)^\gamma$ is not empty.
Let $m\in (\mm_c^\circ)^\gamma$
be a smooth point, and write $K'\subset K$ for the
finite group generated by $\gamma$.
The group $K'$ acts linearly
on $T_m(\mm_c^\circ)$, the tangent space at $m$.
We observe that, by Chevalley theorem,
the orbifold $T_m(\mm_c^\circ)\big/K'$
can {\it not} be  smooth, since the group $K'$
acts on $T_m(\mm_c^\circ)$ by symplectic automorphisms,
hence is not generated by (ordinary) reflections.
It follows that $m$ projects to a singular point
of $\mm_c=\mm_c^\circ/K$.
But since $\mm_c$ has been shown to be smooth, it follows
that $K$ acts freely on $\mm_c^\circ$.

Let $\ZZ\G^\circ$ be the center of the group algebra $\C\G^\circ$.
By Theorem \ref{mckay}, the cohomology of $\mm_c^\circ$ is 
isomorphic, as a graded space, to $\grd(\ZZ\G^\circ)$. 
Moreover, it is clear from the 
proof of this theorem that this isomorphism can be chosen to
be $K$-equivariant.

Now let $g$ be a nontrivial element of 
$K$. By the Lefschetz fixed point formula for $g$,
an automorphism of finite order, we  have
$\,
\Tr\bigl(g,\,H^*(\mm_c^\circ)\bigr)=0.
\,$
On the other hand, 
$$
\Tr\bigl(g,\,H^*(\mm_c^\circ)\bigr)=
\Tr\bigl(g,\,\grd(\ZZ\G^\circ)\bigr)=\Tr\bigl(g,\,\ZZ\G^\circ\bigr),
$$
where the rightmost term  equals
the number of conjugacy classes in
$\G^\circ,$ hence,
 is $>0$. The contradiction implies that $g$ does not exist.
Thus, the group $K$ is trivial,
and {$\G^\circ=\Gamma.$~\qed}

\section{The rational Cherednik algebra}
\setcounter{equation}{0}

Let $\h$ be a
finite dimensional complex vector space
with a nondegenerate inner product,
and $W$ a finite group of orthogonal transformations of $\h$
generated by reflections.
For any reflection $s\in W$, fix a nonzero linear function
$\alpha_s\in \h^*$ uniquely determined, up to a nonzero scalar factor,
by the condition: $s(\alpha_s)=-\alpha_s$.
We choose the scalar factors in such a way that the function:
$s\mapsto \alpha_s^2$ is $W$-invariant. Write $R$ for
the collection of linear functions
$\lbrace{\pm \alpha_s\rbrace}_{s=\text{reflection in }W}$.
 For $\alpha\in R\subset \h^*$, we
let $s_\alpha$ denote the corresponding reflection $s$,
and set
$\,\alpha^\vee :=2\frac{(\alpha,-)}{(\alpha,\alpha)}\in \h$
(we note that $(\alpha,\alpha)\ne 0$ since $\C\cdot\alpha$ and
$\alpha^\perp$ are eigenspaces of $s$ with different eigenvalues, and
therefore have zero intersection).
A point $x\in \h$ is called regular if, for all $\alpha\in R$,
we have: $\alpha(x)\neq 0$. 

     From now on, we  fix a complex-valued
$W$-invariant
function $c: R\to \C$,
$\alpha\mapsto c_\alpha$, and
 $t\in \C$.
Define a {\it rational Cherednik algebra} with parameters
$(t,c)$ to be the symplectic reflection
algebra corresponding to the $W$-diagonal action on
$\h\oplus\h^*$. Explicitly, it is an associative $\C$-algebra
$\hh_{t,c}$
generated by the spaces $\h,\h^*,$ and the group $W$,
with  defining relations (\ref{P-bracket}).
\medskip

\noindent
{\bf Example: ${\mathbf{S_n}}$-case.}\quad
In the special case of the root system $R=
{\mathbf{A_{n-1}}}$ in the vector space
$\h=\C^n$ (not in $\C^{n-1}$)
we have: $W=S_n$.
Use the standard coordinates on $\C^n$ to 
write $\C[\h]=\C[x_1,\ldots,x_n]$ and
$\C[\h^*]=\C[y_1,\ldots,y_n]$.
Recall that, in the ${\mathbf{A_{n-1}}}$-case,
all roots are $W$-conjugate and the function $c: R\to\C$
reduces to a constant;  we will assume that
$c=1$, and  will view the parameter $\ka=(t,1)$ as a point of
$\CP^1$. Thus, the
case $t=0$ corresponds to the point $\ka=\infty$.

Write $s_{ij}\in S_n$ for the transposition: $i\longleftrightarrow j$.
The algebra $\hh_t=\hh_{t,1}(S_n)$ has generators
$x_1,...,x_n,y_1,...,y_n$ and the group $S_n$, with
the following defining relations, which are a specialization of 
\ref{P-bracket}:
\begin{equation}\label{a_n}
\begin{array}{lll}\displaystyle
&{}_{_{\vphantom{x}}}s_{ij}\cdot x_i=x_j\cdot s_{ij}\quad,\quad
 s_{ij}\cdot y_i=y_j\cdot s_{ij}\,,&
\forall i,j\in\{1,2,\ldots,n\}\;,\;i\neq j\break\medskip\\
&{}_{_{\vphantom{x}}}{}^{^{\vphantom{x}}}
[y_i,x_j]=  s_{ij}\quad,\quad[x_i,x_j]=0=[y_i,y_j]
\,,&
\forall i,j\in\{1,2,\ldots,n\}\;,\;i\neq j\break\medskip\\
&{}^{^{\vphantom{x}}} [y_k,x_k]= t\cdot 1-\sum_{i\neq k}\;s_{ik}\,.
\end{array}
\end{equation}
The relations (\ref{a_n}) have  an intriguing `hidden' symmetry.
Specifically, for each $i=1,2,\ldots,n,$ put:
$z_i=x_i\cdot y_i-\sum_{j<i}\;s_{ij}\in\hh_t.$
Let ${\mathcal{H}}$ be the {\it degenerate affine Hecke algebra}
 of type ${\mathbf{A_{n-1}}}$.
Recall that the latter algebra 
contains $\C{S_n}$ and $\C[u_1,\ldots,u_n]$ as subalgebras,
subject to certain commutation relations, defined e.g. in [Dr2].
Given any non-zero triple $(a,b,c)\in \C^3,$ 
we consider an assignment:
\begin{equation}\label{trig_ass}
\varpi_{_{a,b,c}}:\;u_i\; \mapsto\; ax_i+by_i+cz_i\;\;,
\;\,i=1,2,\ldots,n\;,\quad\mbox{and}\quad
w\mapsto w\,\,,\;\;\forall w\in S_n\,.
\end{equation}

The Proposition below
is verified by a straightforward calculation
(Alternatively, one may observe that 
in the faithful representation ${\widetilde\Theta}_{t,c}$, see
Proposition \ref{Dunkl} below,
the $z_i$'s go, after a change of variable, to `trigonometric'
Dunkl operators. The latter are known to commute
with each other by Cherednik).
We do not know if this Proposition has  a generalization
to root systems other than  ${\mathbf{A_{n-1}}}$.
\begin{proposition}\label{trig_dunkl}
\vi The elements: $\,\{ax_i+by_i+cz_i\}_{i=1,2,\ldots,n}\,$
commute in $\hh_t$, for any triple $(a,b,c)\in \C^3\sminus
\{(0,0,0)\}$.

\vii For  any $c\neq 0$, the assignment (\ref{trig_ass})
extends to an algebra imbedding
$\varpi_{_{a,b,c}}: {\mathcal{H}}\into\hh_t$. 
For $c=0$ 
(and any $(a,b)\neq (0,0)$), the map $\varpi_{_{a,b,c}}$
 degenerates to an algebra imbedding
$\varpi_{_{a,b,0}}:
\C[u_1,\ldots,u_n]\#S_n\into \hh_t$.
\qed
\end{proposition}
\medskip
 
We now return to Cherednik algebras associated
to a general root system.
Theorem \ref{H_PBW} yields another useful version of the PBW-theorem
for Cherednik algebras:
\begin{corollary}
\label{PBW3} For any $\ka\in {\sf{\overline{C}}}$,
multiplication in $\hh_\ka$ induces a vector space isomorphism:
$\,\C[\h]\otimes \C[\h^*]\otimes \C{W}\iso \hh_\ka$.\quad\sq
\end{corollary}\medskip

 \noindent
{\bf{Representation via Dunkl operators.}}\quad Regular points
form a Zariski open subset $\hreg\subset \h$,
and we write $\dd(\hreg)$ for the algebra of algebraic differential
operators on $\hreg$. Given $\ka=(t,c)\in
\C\oplus
{\sf{C}}$,
to any $y\in \h$ one associates the
following {\it Dunkl operator} ([D1,D2]):
$$
D_y:=t\frac{\partial}{\partial y}+\frac{1}{2} \sum\nolimits_{\alpha\in R}\;c_\alpha
\cdot\frac{\langle\alpha, y\rangle}{\alpha}\cdot(s_{\alpha}-1)
\enspace\in\enspace \dd(\hreg)\#W \,.
$$

We recall the
following well known result, due to Cherednik
in the trigonometric case.

\begin{proposition} \label{Dunkl} For any $\ka=(t,c)\in \C\oplus
{\sf{C}}\,,\,
t\neq 0,\,$ the assignment:
$\,w \mapsto w\,,\,
x \mapsto x\,,\,
y \mapsto D_y\,,\,w\in W,$
$x\in \h^*\,,\,y\in \h,\,$
extends to an {\sf injective} algebra homomorphism
$\,{\widetilde{\Theta}}_{t,c}:\, \hh_{t,c}\,\into\, \dd(\hreg)\#{W}\,.$
\end{proposition}

\Pf . 
The construction of the homomorphism
${\widetilde{\Theta}}_{t,c}$ has the following conceptual interpretation.
Consider the $\hh_{t,c}$-module
$\,{\mathtt{Ind}}_{S\h\#W}^{\hh_{t,c}}{\mathbf{1}},\,$
induced from the trivial representation of
the subalgebra $S\h\#W\subset \hh_{t,c}$. The underlying vector space
of this module can be identified, by  PBW-isomorphism
\ref{PBW3}, with $\C[\h]$, the polynomial algebra.
 Writing explicitly the action on polynomials
of the generators of the algebra $\hh_{t,c}$ yields the
formulas for ${\widetilde{\Theta}}_{t,c}$.

To prove injectivity of  ${\widetilde{\Theta}}_{t,c}$,
equip the algebra $\hh_{t,c}$ with an increasing filtration
$F_\bullet(\hh_{t,c})$,
 by placing
both $\C{W}$ and $\h^*$  in  filtration degree zero,
and $\h$ in  filtration degree 1 (this filtration
is different from the canonical filtration on $\hh_{t,c}$ introduced earlier).
Clearly, for any $k\geq 0$,
the homomorphism $\,{\widetilde{\Theta}}_{t,c}: \hh_{t,c}\too \dd(\hreg)\#{W}\,$
takes elements of $F_k(\hh_{t,c})$ into $\C{W}$-valued
differential operators on $\hreg$ of order $\leq k$.
Thus, there is a well-defined associated graded map
$$\grd_{_{\!F}}\!({\widetilde{\Theta}}_{t,c}):\;\;\grd_{_{\!F}}\!(\hh_{t,c})
\too \grd\bigl(\dd(\hreg)\#{W}\bigr)\;\simeq\;
\C[\hreg\times\h^*]\,\#\,W\,.
$$
Note that the explicit formula for the Dunkl operator shows
that the image of $\,\grd_{_{\!F}}\!(\hh_{t,c})\,$ equals
the subalgebra: $\C[\h\times\h^*]\,\#\,W\subset
\C[\hreg\times\h^*]\,\#\,W\,.$
We consider the composition:
$$
\C[\h\oplus \h^*]\,\#\,W\;=\;\hb\;\iso\;{\grd}(\hh_{t,c})\;
\stackrel{\grd_{_{\!F}}\!({\widetilde{\Theta}}_{t,c})}{\too}\;
\C[\h\times\h^*]\,\#\,W\;\subset\;
\C[\hreg\times\h^*]\,\#\,W\,.
$$
 It is clear that the composition above
gives the identity map on $\C[\h\times\h^*]\,\#\,W$.
 It follows that the map
$\grd_{_{\!F}}\!({\widetilde{\Theta}}_{t,c})$, hence ${\widetilde{\Theta}}_{t,c}$, is injective.
Furthermore, we deduce that $\grd_{_{\!F}}\!(\hh_{t,c})\simeq
\C[\h\times\h^*]\,\#\,W$. 
\sq\medskip

The definition of Dunkl operators as well as
the 
Proposition above can be extended easily to the case of any
complex reflection group $W$, see  [BMR]. In particular, we
get the following result, whose  verification
by (quite non-trivial) direct computations is first
due to Dunkl [D1],[D2] and
 Cherednik [Ch] in the Weyl group case, and in [BMR] in general
(for complex reflection groups).

\begin{corollary}\label{dunk_comm}
The Dunkl operators $\,D_y\,,\,y\in\h,\,$ commute with each other.
$\;\;\square$
\end{corollary}

Following Olshanetsky-Perelomov [OP],
define the rational Calogero-Moser operator with parameters
$\ka=(t,c)\in \overline{\CC}$, to be
the following $W$-invariant {differential operator~on~$\hreg$:}
$$
\CM_\ka (=\CM_{t,c})\;\;=\;\; t^2\cdot
\Delta_\h- \frac{1}{2}\sum\nolimits_{\alpha\in
R}\;\,c_\alpha(c_\alpha+t) \cdot
\frac{(\alpha,\alpha)}{\alpha^2} \quad \text{\bf{(Calogero-Moser
operator)}}\,.
$$

     From now on we fix $c\in \CC$, and set $\delta_c=\prod_{\alpha\in
R}\alpha^{c_\alpha/2}$. For $t\in \C$, we consider the modified
homomorphism $\,\Theta_{t,c}:=\frac{1}{\delta_{c/t}}\ccirc
{\widetilde{\Theta}}_{t,c}\ccirc \delta_{c/t}:
 \hh_{t,c}\too\dis \dd(\hr)\#W,\,$
sending $u \in \hh_{t,c}$ to the operator:
$f \mapsto
\frac{1}{\delta_{c/t}}\cdot {\widetilde{\Theta}}_{t,c}(u)(\delta_{c/t}\cdot f).\,$ 
Although $\delta_{c/t}$ is a  multivalued
function on $\h$, the homomorphism $\frac{1}{\delta_{c/t}}
\ccirc {\widetilde{\Theta}}_{c/t}
 \ccirc \delta_{c/t}$
is  well-defined on generators, hence, on the whole of $\hh_\ka$.
Thus, Proposition \ref{Dunkl} implies that 
$\,\Theta_{t,c}: \hh_{t,c}\too\dd(\hr)\#W \,$
is an injective algebra homomorphism.
The restriction of  $\,\Theta_{t,c}\,$
to the spherical subalgebra gives an algebra homomorphism
\begin{equation}\label{theta_spher}
\Theta_\ka^{{\tt{spher}}}:\;\; \e\hh_\ka \e\;\too\;
\e\bigl(\dd(\hr)\#W\bigr)\e
\;=\;\dd(\hr)^W\,,
\end{equation}
see (\ref{smash_iso}). The map $\Theta_\ka^{{\tt{spher}}}$
clearly respects the filtrations on both sides.
Furthermore, 
a straighforward calculation yields:
$\Theta_\ka^{{\tt{spher}}}(\e\cd\Delta_\h)=
{\mathsf{L}}_\ka,$ where
$\Delta_\h\in S\h^W=\C[h^*]^W$ is the quadratic Casimir,
viewed as an element of $\C[h^*]^W\subset \hh_\ka\,$
(note that $\e\cd\Delta_\h=\Delta_\h\cd\e)$.

Let $\cc_{t,c}$ be the centraliser of $\CM_{t,c}$ in 
${\dd(\hreg)_-^W}$ (= the algebra spanned by
homogeneous elements $D\in\dd(\hreg)^W$ such that 
$\text{\it order}(D)+\text{\it degree}(D)\le 0$).
Let $\bcal_{t,c}$ be the subalgebra in $\dd(\hreg)^W$
generated by $\cc_{t,c}$ and
$\C[\h]^W\subset\dd(\h)^W$, the subalgebra of $W$-invariant polynomials.

\begin{theorem}\label{generic} 
\vi For any fixed $c\in \CC$, and generic $t\neq 0$
(i.e. outside of a countable set), we have:
$\Image(\Theta_{t,c}^{{\tt{spher}}})=\bcal_{t,c}$. 

(ii) If, for some $(t,c) \in \C\oplus\CC$, the image of
 $\Theta_{t,c}^{{\tt{spher}}}$ is contained in
$\bcal_{t,c}$ then
the following maps are algebra isomorphisms
$$\Theta_{t,c}^{{\tt{spher}}}:\;\;\e\hh_{t,c}\e\iso \bcal_{t,c}\,,
\qquad\mbox{$\operatorname{and}$}\qquad
\grd(\Theta_{t,c}^{{\tt{spher}}}):\;\;
\grd(\e\hh_{t,c}\e)\iso \grd(\bcal_{t,c})\,.
$$
\end{theorem}

\Pf : \vi
For generic $t$, 
the algebra $\e\hh_\ka\e$ is generated, as an
associative algebra, by the subspaces: $\C[\h]^W\!\cdot\!\e$ and
$\C[\h^*]^W\!\cdot\!\e$. Indeed, 
this is known \cite{LS1} to be true in the Weyl
algebra case: $t=1, c=0$, hence it is true for
all $t$ except a countable set. 

On the other hand, the image of $\C[\h]^W\!\cdot\!\e$ under 
the map $\Theta_{t,c}^{{\tt{spher}}}$ equals $\C[\h]^W$. 
Also, the
results of Opdam (see \cite{O},\cite{Ki}) 
imply that the image of $\C[\h^*]^W\!\cdot\!\e$ 
under this map is the algebra $\cc_\ka$. This implies (i).  

\vii In the proof of (i) we showed that $\bcal_\ka$ is always 
contained in the image of $\Theta_{t,c}^{{\tt{spher}}}$
Thus, under our assumptions,
the image of $\Theta_\ka^{{\tt{spher}}}$ equals $\bcal_\ka$.
We know that the map $\grd(\Theta_\ka^{{\tt{spher}}})$ is injective
(it is independent on $\ka$). Therefore, 
the map $\Theta_\ka^{{\tt{spher}}}$ is also injective. 
We conclude that 
$\Theta_\ka^{{\tt{spher}}}:\e \hh_\ka \e\to \bcal_\ka$ is an isomorphism. 

It remains to show that $\grd(\Theta_\ka^{{\tt{spher}}}):
\grd(\e \hh_\ka \e)\to \grd(\bcal_\ka)$ is an isomorphism. 
But this is clear, since if $f$ is any filtered map 
of $\Z_+$-graded spaces 
such that $\grd(f)$ is injective and $f$ is an isomorphism then
$\grd(f)$ is also an isomorphism. We are done.  
\sq\medskip

\begin{proposition}\label{typeEF} Suppose that $W$ is a Weyl group with no
factors of type ${\mathbf{E}}$ or
 ${\mathbf{F}}$. Then  part \vi of Theorem \ref{generic},
in particular,
the conclusion of part \vii\!, holds for all
values of $t$.
\end{proposition}

\Pf . We need to show that the image of $\Theta_\ka^{{\tt{spher}}}$
is contained in $\bcal_{t,c}$.
For this, 
it is sufficient to prove that $\e\hh_{t,c}\e$ is generated, as an algebra, 
by the subalgebras $\C[\h]^W$ and $\C[\h^*]^W$. This follows from a stronger 
statement that  
$\grd(\e\hh_{t,c} \e)= \C[\h\oplus\h^*]^W$ is generated, as a
Poisson algebra, 
by the subalgebras $\C[\h]^W$ and $\C[\h^*]^W$. But this
is shown in \cite{Wa} for Weyl groups  containing no factors
of type ${\mathbf{E}}$ or
 ${\mathbf{F}}$.\sq\medskip

We expect, but cannot show, 
that  the same result holds for all finite
reflection groups.

\medskip   \noindent
{\bf{Relation of $\hh_\ka$ to the double-affine Hecke algebra.}}\quad
The algebra $\hh_\ka$ is a rational version of the double-affine Hecke algebra.
By this we mean that if
$W$ is the Weyl group associated to a finite root system $R$,
 then $\hh_\ka$ is obtained
from the double affine Hecke algebra attached to $W$
by a certain limiting procedure.

To
illustrate this, we will consider the special
case of the root system of type ${\mathbf{A_1}}$.
In this case $W=\lbrace{1,s\rbrace}$.
Thus, if $t=1$ and $c=c_\alpha$, for the unique
positive root $\alpha$, then,
the algebra $\hh_{1,c}$ is generated
by 3 elements: $s,x,y$, with defining relations
$$
s^2=1\quad,\quad
sx=-xs\quad,\quad sy=-ys\quad,\quad [y,x]=1-2c\cdot s.
$$

Let us show how the algebra $\hh_{1,c}$
is obtained from the double-affine Hecke algebra by a  limiting procedure.
Let us consider the version of $\hh_{1,c}$ over $\C[c]$, which we will denote
in the same way.
Recall that the double-affine Hecke algebra
${\mathbf{H}}$ is the $\C[q^{\pm 1},\tau^{\pm 1}]$-algebra with 3
generators: $T,X,Y$ (where
$X,Y$ are invertible), and defining relations:
$$
(T-\tau)(T+\tau^{-1})=0\;\;,\;\;
TXT=X^{-1}\;\;,\;\; T^{-1}YT^{-1}=Y^{-1}\;\;,\;\;
Y^{-1}X^{-1}YXT^2=q\,.
$$
Define a
 completed double-affine Hecke algebra ${\mathbf{\widehat{H}}}$ 
 to be a $\C[c][[h]]$-algebra topologically generated
(in the $h$-adic topology) by $s,x,y$, with the same defining
relations, for $X=e^{h^{}x}$,
$Y=e^{h^{}y}$, $T=se^{h^2 cs}$,
$q=e^{h^2}$, $\tau=e^{h^2 c}$.
It is clear that ${\mathbf{\widehat{H}}}$ is a certain completion of
${\mathbf{H}}$.

The proposition below says that the algebra ${\mathbf{\widehat{H}}}$ 
is a flat deformation
of $\hh_{1,c}$.

\begin{proposition} The algebra ${\mathbf{\widehat{H}}}$ is
flat over $\C[[h]]$, and ${\mathbf{\widehat{H}}}/h\cdot{\mathbf{\widehat{H}}}=
\hh_{1,c}$.
\end{proposition}

\Pf . The last statement is obtained directly from the relations:
one shows that the relations of ${\mathbf{\widehat{H}}}$ are deformations
of the relations of $\hh_{1,c}$. The flatness  statement
follows from the fact, due to Cherednik, that the representation
of $\hh_{1,c}$ by Dunkl operators can be deformed to a 
faithful representation
of ${\mathbf{\widehat{H}}}$
by q-difference analogs of Dunkl operators.~\sq

A similar relation between $\hh_\ka$ and the double-affine Hecke algebra
exists for any  Weyl group.
The statement and proof of this fact is analogous
to the ${\mathfrak{s}\mathfrak{l}}_2$ case.

\medskip   \noindent
{\bf{Fourier transform.}}\quad
There is an important algebra automorphism $\FF:\hh_\ka\to \hh_\ka,$
to be called the {\it Fourier transform}.
The automorphism $\FF$ is a rational version of
Cherednik's ``difference Fourier transform'' for double affine Hecke
algebras. In particular, in the
${\mathfrak{s}\mathfrak{l}}_2$-case, it is the classical
Fourier-Hankel transform, see \cite{ChM}.

Let  $(-,-)$ be the  $W$-invariant bilinear form on $\h$ that has been fixed
throughout.
For $x\in \h^*$, let $\check{x}$ be the
element of
$\h$ such that $x=(\check{x},-),$
 and for $y\in \h$, let $\check{y}=(y,-)$ be the corresponding
element of
$\h^*$. The Fourier transform on $\hh_\ka$ is defined to be
 an algebra homomorphism
 $\FF:\hh_\ka\to
\hh_\ka$  given on generators by the following assignment:
$$
\FF:\;x \mapsto \check{x}\quad,\quad
y \mapsto -\check{y}\quad,\quad w \mapsto w, \quad\quad
\forall  x\in \h^*\,,\, y\in \h\,,\, w\in W.$$
It is clear that $\FF$ extends to an
automorphism of $\hh_\ka$  preserving the subalgebra
$\e\hh_\ka \e$. It follows  that (in cases when
$\e \hh_\ka \e$ is isomorphic to $\bcal_\ka$, see Theorem \ref{generic}
and Proposition \ref{typeEF})
it also acts on $\bcal_\ka$, permuting
$\C[\h]^W$ and $\cc_\ka$.
In particular, this is so for type ${\mathbf{A_n}}$.
\medskip

\noindent
{\bf{Homomorphism ${\mathbf{\Theta_{0,c}}.}$}}\quad
We now fix $c\in \CC\,,\, c\neq0,\,$ and consider the limit
of the map $\Theta_{t,c}$, see above (\ref{theta_spher}),
 when $t\to 0$.
Let $\dd_t(\h)$ be the Weyl algebra of
 the  vector space
$\h\oplus\h^*$ equipped with the standard  symplectic form multiplied
by the parameter $t$. Thus, the algebra  $\dd_t(\h)$ is generated by
$\h$ and $\h^*$, so that elements of $\h$ commute,
elements of $\h^*$ commute, and we have:
$[y,x]=t\cdot\langle x,y\rangle$, $\forall x\in \h^*,y\in \h$. Let
$\dd_t(\hr)$ be the localized algebra
obtained by inverting the Weyl denominator $\delta_x\in \C[\h]\subset
\dd_t(\h).\,$ For  $t\ne 0$, the assignment:
$y\mapsto t \partial_y\,,\, x\mapsto x$ gives an algebra
isomorphism $\,\phi_t: \dd_t(\hr)\iso\dd(\hr)$.
For $t\ne 0$,  the composite map
$\Theta_{t,c}^\sharp=(\phi_t)^{-1}\ccirc \Theta_{t,c}:
\,\hh_{t,c}\too \dd_t(\hr)\#W\,$ is the identity map on $W$, and is given
 on generators: $x\in \h^*\,,\,y\in\h,\,$ 
by the following formulas:
$$
x\mapsto x\quad,\quad
y\,\mapsto\,\frac{1}{\delta_{c/t}}\cdot y\cdot\delta_{c/t}+\frac{1}{2}
\sum_{\alpha\in R}\;c_\alpha\cdot
\frac{\langle\alpha, y\rangle}{\alpha}\,(s_{\alpha}-1)
\;= \;y+\frac{1}{2}
\sum_{\alpha\in R}\;c_\alpha\cdot
\frac{\langle\alpha, y\rangle}{\alpha}\,s_{\alpha}
\;.
$$
Since the formulas above do not involve the parameter $t$,
the homomorphism  $\Theta_{t,c}^\sharp$ makes sense for $t=0$,
and we set $\,{\Theta}_{0,c} :=\Theta^\sharp_{0,c}$.

In the case  $t=0$ we have:
$\,\dd_{t}(\hr)=\C[\hr\times \h^*]$,  a commutative
algebra. Copying the proof of
theorem \ref{generic}  one obtains:

\begin{proposition}\label{theta2}
The map $\,{\Theta}_{0,c}:\,
 \hh_{0,c}\too \C[\hr\times \h^*]\,\#\,W
$ is an injective  algebra
homomorphism. \quad\sq
\end{proposition}

Write: $p\mapsto |p|^2$ for the quadratic
polynomial on $\h^*$ corresponding to the 
squared norm relative to the invariant
form. View it as an element of $\C[\hr\times \h^*]$
constant along the first factor,
 and view $\alpha \in R$ as a function on $\hreg\times \h^*$
constant  along the second  factor.

\begin{definition}[\cite{OP}]\label{CM_ham}
The function $\,\dis {\sf{L}}_{0,c}:\, (x,p) \mapsto |p|^2- \frac{1}{2}
\sum\nolimits_{\alpha\in R}\;\frac{c_\alpha(c_\alpha+1)}{\alpha(x)^2},\,$
on $\hr\times \h^*$, is called the (classical) Calogero-Moser
hamiltonian.
\end{definition}

A straightforward calculation yields:
$\Theta_{0,c}^{{\tt{spher}}}(\e\cd|p|^2)={\sf{L}}_{0,c},$
where $|p|^2$ denotes the
element of $\hh_{0,c}$ corresponding to the polynomial:
$p\mapsto |p|^2$, viewed as an element of
$\C[\h^*]^W\subset \hh_{0,c}$.

Further, 
the standard symplectic structure on $\,T^*(\hr)=\hr\times \h^*\,$ makes
$\C[\hr\times \h^*]$
a Poisson algebra. We have: ${\sf{L}}_{0,c}\in
\C[\hr\times \h^*]^W.$ Let $\C[\hr\times \h^*]^W_-$
denote the algebra spanned by 
bi-homogeneous functions $P\in\C[\hr\times \h^*]^W$ such that 
$\text{\it total degree}(P)+\text{\it $p$-degree}(P)\le 0$.
Let
$\cc_{0,c}$ denote the centralizer in $\C[\hr\times
  \h^*]^W_-$ of the
Calogero-Moser
hamiltonian $\,{\sf{L}}_{0,c}$
 with respect to the
Poisson bracket, and let ${\mathcal{B}}_{0,c}$ denote the
Poisson subalgebra of $\C[\hr\times \h^*]^W$ generated by
$\cc_{0,c}$ and the  algebra $\C[\h^*]^W$.
On the other hand, recall that the spherical subalgebra
$\e\hh_{0,c}\e$ has a natural
Poisson algebra structure defined using the deformation construction
of \S15.

\begin{proposition}\label{theta_infty}
  If $W$ is a Weyl group which has  no
factors of type ${\mathbf{E}}$ or
 ${\mathbf{F}}$, then the map ${\Theta}_{0,c}$ induces
a Poisson algebra isomorphism
$\,
\Theta_{0,c}^{{\tt{spher}}}: \e\hh_{0,c}\e \iso \bcal_{0,c}$.
Moreover, the associated graded map
$\,\grd(\Theta_{0,c}^{{\tt{spher}}}): \grd(\e\hh_{0,c}\e) \iso
\grd(\bcal_{0,c})$ is also bijective.
\end{proposition}

\Pf .
Since ${\Theta}_{0,c}^{{\tt{spher}}}$
is the specialisation  at $t=0$ of the family of associative algebra homomorphisms
$\,\{\Theta_{t,c}^{{\tt{spher}}}\}_{t\neq 0},\,$
general `deformation formalism' implies
that ${\Theta}_{0,c}^{{\tt{spher}}}$ is a homomorphism of
Poisson algebras.
The rest of the argument is identical to that in the proof of
Theorem \ref{generic} and Proposition \ref{typeEF}.\quad\sq\medskip

Since
$\grd(\ZZ_{0,c})=\ZB$
is a   finitely generated algebra without  zero divisors,
it follows that $\ZZ_{0,c}$ is a   
finitely generated algebra without  zero divisors
as well. Therefore, it defines an irreducible affine
algebraic variety $\Spec{\ZZ_{0,c}}$ of dimension equal to $\dim V$.
This variety has a canonical  Poisson structure.
Thus,  $\Spec{\ZZ_{0,c}}$ is a Poisson deformation of the
Poisson variety $V^*/W$, cf. also \S17.
We propose the following

\begin{definition} The Poisson
variety $\Spec{\ZZ(\hh_{0,c})}$ is called the Calogero-Moser
space associated to~$W$, and the parameter $c\in\CC$.
\end{definition}

 The variety $\Spec{\ZZ(\hh_{0,c})}$ may be regarded as 
the ``correct" completion of the phase space for the classical 
Calogero-Moser system 
acssociated to the Weyl group~${W.\,}$ 
We have shown that $\Spec{\ZZ(\hh_{0,c})}$ is an irreducible
Gorenstein normal variety. The problem of {\it smoothness}
of  $\Spec{\ZZ(\hh_{0,c})}$ will be considered, in a few 
special cases, in Appendix E.
\medskip

We now study the center of the
algebra $\hh_{0,c}$ in some detail.

\begin{proposition}\label{finite} 
For an arbitrary finite complex reflection group $W$, we have:

\vi The subalgebras $\C[\h^*]^W\subset \hh_{0,c}$ and
$\,\C[\h]^W\subset \hh_{0,c}$  are both contained in  $\ZZ_{0,c}$.

\vii The algebra $\ZZ_{0,c}$
is a free  $\C[\h]^W\otimes\C[\h^*]^W$-module
of rank $|W|$.
\end{proposition}

\Pf . The  inclusion: $\,\C[\h^*]^W\subset \ZZ_{0,c}$
 follows from the injectivity
of $\Theta_{0,c}$ and the fact that,
for $f\in \C[\h^*]^W$, we have $\Theta_{0,c}(f)\in \C[\hr]^W$,
which is obviously in the center of $\C[\hr\times \h^*]\#W $.
The inclusion: $\,\C[\h]^W\subset \ZZ_{0,c}$ follows from the
previous one
by applying the Fourier  automorphism $\FF$.
Thus, we have shown that: $\C[\h]^W\otimes\C[\h^*]^W\,\subset \,
\ZZ_{0,c}$.

Further,  $\,\C[\h\oplus\h^*]$ is
a free module over its  subalgebra
$\C[\h]^W\otimes\C[\h^*]^W\subset \C[\h\oplus\h^*]$,
by the Chevalley theorem. Moreover, one has a $W\times W$-module isomorphism:
$$\C[\h\oplus\h^*]\;\simeq\;
\C[W\times W]\;\bigotimes\; \bigl(\C[\h]^W\otimes\C[\h^*]^W\bigr)\,.
$$
Taking the invariants with respect to the diagonal
subgroup $W\subset W\times W$ shows
that $\C[\h\oplus\h^*]^W$ is a
free  $\C[\h]^W\otimes\C[\h^*]^W$-module
of rank $|W|$. Part (ii) now follows from
the isomorphism $\grd(\ZZ_{0,c})\iso\C[\h\oplus\h^*]^W$
of  Theorem \ref{center}.\hfill\sq
\medskip

The imbedding: $\C[\h]^W\otimes\C[\h^*]^W
\into \ZZ_{0,c}$ of Proposition
\ref{finite} gives rise to a {\it finite} surjective map
$\Upsilon: \Spec\ZZ_{0,c}\onto \h/W\times \h^*/W,$
of degree $|W|$. Let $|W/\Ad|$ denote the number of conjugacy
classes in the finite
complex reflection group $W$.

\begin{proposition}\label{nilp_rep} If $\Spec\ZZ_{0,c}$ is
smooth, then the fiber $\Upsilon^{-1}(0)\subset \Spec\ZZ_{0,c}$ consists of 
 $|W/\Ad|$ points.
\end{proposition}

\Pf . The $\C^*$-action on $\h\oplus\h^*$ given by
the assignment $\,\C^*\ni \lambda: 
\,y\oplus x \mapsto \lambda\cd y\oplus \lambda^{-1}\cd x$
extends to a $\C^*$-action on $\hh_{t,c}$ by
algebra automorphisms keeping the group $W$ pointwise fixed
(this is part of the $SL_2(\C)$-action on 
$\hh_{t,c}$ considered in Corollary \ref{rotation}).
The $\C^*$-action on $\hh_{0,c}$ induces an action
on $\ZZ_{0,c}$,
hence on  $\Spec\ZZ_{0,c}$. 

Furthermore, the  $\C^*$-action on $\h\oplus\h^*$
induces one on $\h/W\times \h^*/W$, so that the map
$\Upsilon$ becomes $\C^*$-equivariant. Observe
that the point $(0,0) \in \h/W\times \h^*/W$
is the only fixed point of  the  $\C^*$-action on 
 $\h/W\times \h^*/W$.
It follows, since $\C^*$ is connected and the fiber
$\Upsilon^{-1}(0)$ is finite, that each point
of $\Upsilon^{-1}(0)$ is a $\C^*$-fixed point,
and these are the only fixed points on  $\Spec\ZZ_{0,c}$. 

Now, the Lefschetz fixed point formula for the
 $\C^*$-action on  $\Spec\ZZ_{0,c}$ yields:
$|\Upsilon^{-1}(0)| = \boldsymbol{\chi}(\Spec\ZZ_{0,c})$,
where $\boldsymbol{\chi}(-)$ stands for the 
 Euler characteristic. On the other hand, if 
$\Spec\ZZ_{0,c}$ is smooth, then from Theorem \ref{mckay}
one finds: $\boldsymbol{\chi}(\Spec\ZZ_{0,c})
=|W/\Ad|.$ The result follows.\quad\qed\medskip

Simple $\hh_{0,c}$-modules correspoding to points of $\Upsilon^{-1}(0)$
may be called 
nilpotent, because  they are characterized by the property that all
 elements of both
$\h\subset \hh_{0,c}$ and $\h^*\subset \hh_{0,c}$ 
act by nilpotent operators in them. 
\medskip

\noindent
{\bf Remark.}
For any finite Coxeter group $W$,
results similar to  Theorem \ref{Coh_Mac_intro} are valid if
the symmetrizer $\e$ is replaced with the antisymmetrizer
$\,{\mathbf{e_-}}=\frac{1}{|W|}\sum_{w\in W}\;(-1)^{\ell(w)}
\cdot w,\,$ since there is an anti-involution of $\hh_{0,c}$
that acts trivially on $\h\oplus \h^*$ and acts on $W$ via:
 $w\mapsto (-1)^{\ell(w)}\cdot w$; this
anti-involution interchanges $\e$ and ${\mathbf{e_-}}.\quad
\lozenge$\medskip

\noindent
{\bf Partial brackets.}\quad
A slight generalization
of the  construction of the
Poisson brackets on  commutative
algebras $\e\hh_{0,c}\e$ and $\ZZ_{0,c}$ (see \S15)
can be also applied to construct
 bilinear pairings:
\begin{equation}\label{poisson}
\C[\h]^W \times \hh_{0,c}\too\hh_{0,c}
\quad\text{and}\quad
\C[\h^*]^W\times \hh_{0,c}\too\hh_{0,c}
\quad,\quad
(P,u)\mapsto \{P,u\}\,.
\end{equation}
To define these pairings, observe that, for any $(t,c)\in\C\oplus\CC$,
one has a {\it canonical} algebra imbedding
$\,\epsilon_{_{t,c}}: \C[\h]^W \into \hh_{t,c}\,$ and, moreover,
the algebra $\epsilon_{_{0,c}}(\C[\h]^W)$ is central
in $\hh_{0,c}$. Now, given an element $u\in \hh_{0,c}$
choose an algebraic family $u_t\in \hh_{t,c}$ such that $u_{0,c}=u$.
For any $z\in \C[\h]^W$ define an element
$\{z,u\}\in \hh_{0,c}$ to be the linear term
in the expansion of $[\epsilon_{_{t,c}}(z),u_t]=
\epsilon_{_{t,c}}(z)\cdot u_t - u_t\cdot \epsilon_{_{t,c}}(z)$
at $t=0$. This term is independent of the
choice of the family $u_t\in \hh_{t,c}$ such that $u_t=u$,
since $\epsilon_{_{0,c}}(z)$ is a central element.
For each $z\in \C[\h]^W$, the assignment:
$u\mapsto \{z,u\}$ gives a derivation of the algebra $\hh_{0,c}$,
which reduces, on $\ZZ_{0,c}$, to the map given by taking the
Poisson bracket (in $\ZZ_{0,c}$) with
$\epsilon_{_{0,c}}(z)$.

A similar construction applies to the algebra
$\C[\h^*]^W$ instead of $\C[\h]^W$. The reader should be warned
 that we have not constructed a  `Poisson bracket':
$\ZZ_{0,c}\times \hh_{0,c}\to\hh_{0,c}$ since, for $t\neq 0$,  there is no
obvious family of vector space
imbeddings $\,\epsilon_{_{t,c}}: \ZZ_{0,c} \into \hh_{t,c}.$

\section{ Automorphisms and derivations
of the  algebra $\hh_\ka$.}
\setcounter{equation}{0}

Let $(\h,W,R)$ be as in \S4. In \cite{Ch2}, Cherednik
made an observation, equivalent to 

\begin{proposition}\label{aut_ka} Let $f\in \C[\h]^W$.
Then, for any $t\neq 0$, the assignment
$$w\mapsto w\quad,\quad
x\mapsto x\quad,\quad
y\mapsto y-\frac{\partial f(x)}{\partial x}\quad,\quad
x\in \h^*\,,\,w\in W\,,\, y\in \h
$$
uniquely extends to an automorphism $a_{_f}:\hh_{t,c}\to\hh_{t,c}, $
acting trivially on $\C{W}\subset \hh_{t,c}$. A similar result holds
for $f\in \C[\h^*]^W$.
\end{proposition}

\Pf . It is easy to see that, for any $u\in \dd(\hreg)\#W$,
the operator $e^{-f}\ccirc u\ccirc e^f$
is again an element of $\dd(\hreg)\#W$.
Thus, the assignment: $u\mapsto
e^{-f}\ccirc u\ccirc e^f$ gives an automorphism of the algebra
 $\dd(\hreg)\#W$. We transport this automorphism to
$\hh_{t,c}$ via the (faithful)  representation of $\,\Theta_{t,c}:
\hh_{t,c}\to \dd(\hreg)\#W$ by Dunkl operators. A
straightforward calculation shows that the
transported automorphism is given, on generators, by the formulas
of the Proposition, and in particular, maps
$\Image(\Theta_{t,c})$ to itself. The proposition is proved.\quad\sq\medskip

\noindent
{\bf Remark.} The automorphism
$a_{_f}$ can be  informally written as ``$\text{Ad}(e^f)$'', i.e. 
$a_{_f}$  is, essentially, an "inner" automorphism induced by
the element $e^f$ (which belongs to an appropriate completion of
$\hh_{t,c}$).$\quad\lozenge$

\begin{corollary}\label{aut_infty} The assignment of Proposition
\ref{aut_ka} uniquely extends to an automorphism $a_{_f}:\hh_{0,c}
\to\hh_{0,c}$.
\end{corollary}

\Pf . Since the relations (\ref{P-bracket}) among the generators
of the algebras $\hh_{t,c}$ depend on $t$ in a continuous fashion, and the
assignment  of Proposition
\ref{aut_ka} is independent of the parameter $t$,
the result follows from Proposition
\ref{aut_ka} `by continuity'.\quad\sq\medskip

Write $\check{x}$, resp. $\check{y}$,
for the element of $\h$ corresponding to $x\in \h^*$,
resp.  the element of $\h^*$ corresponding to $y\in \h$,
under the invariant form $(-,-)$ on $\h$. The following
result is immediate.

\begin{corollary}\label{rotation} The assignment:
$\;
x\to ax+b\check{y}\quad,\quad y\to c\check{x}+dy\;,\;$
$\,(ad-bc=1),\,$
extends to an $SL_2(\C)$-action on each algebra
 $\hh_\ka\,,\,\ka\in{\sf{\overline{C}}},$ by
algebra automorphisms.\quad\sq
\end{corollary}
\medskip

The above results can be applied to
 study derivations of the algebra $\hh_\ka$,
viewed as `infinitesimal' automorphisms of  $\hh_\ka$.
Fix $f\in \C[\h]^W,$ and consider the  one-parameter group:
$\,\tau\mapsto a_{_{\tau\cdot f}}\,,\,{\tau\in\C},\,$
of automorphisms
of the algebra $\hh_\ka$. Taking the derivative of this family
 at $\tau=0$ we conclude that the map
$\;\theta_{_f}:=\frac{da_{_{\tau\cdot f}}}{d \tau}\big|_{\tau=0}:
\,\hh_\ka\to\hh_\ka\,$ is a derivation. Formulas
of Proposition \ref{aut_ka} yield:
\begin{equation}\label{der_ka}
\theta_{_f}:\quad w\mapsto 0\quad,\quad
x\mapsto 0
\quad,\quad y\mapsto
-\frac{\partial f(x)}{\partial x}
\quad,\quad
x\in \h^*\,,\,w\in W\,,\, y\in \h\,.
\end{equation}
Note that, if
$\ka=(t,c)$ where $t\neq 0$, we have: $\theta_{_f}=\ad f$,
i.e., $\theta_{_f}$ is an inner derivation.
For $\ka=(0,c)$, the derivation $\theta_{_f}$ is not inner.
It is clearly obtained from a family of inner derivations
on the family of algebras $\hh_{t,c}\,,\,t\neq 0,$ by a "limiting procedure"
as $t\to 0$. This argument yields the following expression
for $\theta_{_f}$ in terms of the `partial bracket'
$\,\lbrace{-,-\rbrace}: \C[\h]^W\times\hh_{0,c}\to \hh_{0,c},$
see (\ref{poisson}).
\begin{proposition}\label{two_brackets} For any $\;f\in \C[\h]^W,\,$ resp.
$\;f\in \C[\h^*]^W,\,$
in $\hh_{0,c}$ one has:
$\;
\theta_{_f}(u)=
\lbrace{f,u\rbrace}\,,\,\forall u\in\hh_{0,c}
\,.$\sq
\end{proposition}
\medskip

\noindent 
{\bf{The $\gln$-case.}}\quad From now until the end of this section,
let $R$ be the root system of type
${\mathbf{A_{n-1}}}$ in the vector space
$\h=\C^n$ (not in $\C^{n-1}$).
Write $\,x_1,\ldots,x_n,$ for the standard
coordinates on $\h=\C^n$, and $\,y_1,\ldots,y_n,$
for the dual coordinates on $\h^*$.
We introduce the following  "power sum" elements:
\begin{equation}\label{power_sum}
\psi_\ell=\sum\nolimits_{i=1}^n\;x_i^\ell\;\in\; \C[\h^*]^W\quad
\mbox{and}\quad
\phi_\ell=\sum\nolimits_{i=1}^n\;y_i^\ell\;\in\; \C[\h]^W\quad,
\quad\ell\geq 0.
\end{equation}
Let  $\hh_\infty=\hh_{0,c}(S_n)\,,\,c\neq 0,$
be  the 
rational Cherednik algebra of type ${\mathbf{A_{n-1}}}$, see (\ref{a_n}).
We may  view $\psi_\ell, \phi_\ell$ as elements
of $\ZZ_\infty$, the center of $\hh_\infty$.
Then, an easy calculation yields:
$$
\theta_{_{\psi_k}}(\phi_l)=
-k\cdot\!\!\sum\nolimits_{p+q=l-1}\;\sum\nolimits_{i=1}^n\;\; 
y_i^px_i^{k-1}y_i^q.
$$
Using skew-symmetry of the Poisson bracket on $\ZZ_\infty$,
we obtain, as a corollary, the following  a priori  non-obvious result

\begin{corollary}\label{cen} The element $\lbrace{\psi_k,\phi_l\rbrace}
\in \ZZ_\infty$
is given by either of the following two expressions:
$$\lbrace{\phi_l,
\psi_k\rbrace}\;=\;k\cdot\!\!\sum\nolimits_{p+q=l-1}\;\;
\sum\nolimits_{i=1}^n\;\;
 y_i^px_i^{k-1}y_i^q\;
=\;l\cdot\!\!\sum\nolimits_{p+q=k-1}\;\;\sum\nolimits_{i=1}^n\;\;
 x_i^py_i^{l-1}x_i^q.\qquad\square
$$
\end{corollary}

 We do not know how to check Corollary
\ref{cen} by
a direct computation.\medskip

Let $\C\langle x,y\rangle$ be the free associative algebra on
two generators. Let $F$ be either an automorphism or a derivation of
the algebra $\C\langle x,y\rangle$. It is clear that $F$ is uniquely  determined by
two non-commutative polynomials: $F(x)=: P_{_F}(x,y)$ and $F(y)=:
Q_{_F}(x,y)$. In particular,
for each integer $\ell\geq 0$, we define two one-parameter
groups of automorphisms
$a_{\tau,\ell}\,,\,b_{\tau,\ell},\, (\tau\in\C),\,$
resp.  two derivations $\alpha_\ell\,,\,\beta_\ell$,  of
the algebra $\C\langle x,y\rangle$ by
the following assignments:
\begin{equation}\label{ab}
\begin{array}{lll}\displaystyle
&{}_{_{\vphantom{x}}}\!\!\!\!
a_{\tau,\ell}(x)= x+\tau y^\ell\;\;,\;\;a_{\tau,\ell}(y)=y,\enspace\text{and}\enspace
b_{\tau,\ell}(x)=x\;\;,\;\;b_{\tau,\ell}(y)=y+\tau x^\ell,
\enspace\text{resp.}\break\medskip\\
&{}_{_{\vphantom{x}}}{}^{^{\vphantom{x}}}\!\!\!\!
\alpha_\ell(x)=y^\ell\;\;,\;\; \alpha_\ell(y)= 0,\quad\text{and}\quad
\beta_\ell(x)= 0\;\;,\;\;\beta_\ell(y)= x^\ell\quad,\quad \ell= 0,1,\ldots\,.
\end{array}\end{equation}

\begin{definition}\label{LLL}
Let $\GG$ denote the group of the automorphisms of  $\C\langle x,y\rangle$
generated by the set $\{a_{\tau,\ell},b_{\tau,\ell}\}_{\ell\geq 0},$
and let $\LL$ denote
the Lie algebra of all derivations
 of  $\C\langle x,y\rangle$
generated by the set $\,\{\alpha_\ell,\beta_\ell\}_{\ell\geq 0}.\,$
\end{definition}

\begin{remark}
Write $\Aut_\om\C\langle x,y\rangle$, resp.
$\derv\C\langle x,y\rangle$,
 for the group of all automorphisms $F: \C\langle x,y\rangle\to\C\langle x,y\rangle$ such that
$\,F(xy-yx)=xy-yx\,$, resp.
 for the Lie algebra of all
derivations $f:\C\langle x,y\rangle\to\C\langle x,y\rangle$ such that $\,f(xy-yx)=0\,$.
It is clear that $\GG\subset \Aut_\om\C\langle x,y\rangle$, and it is known,
see \cite{ML},
that in effect $\GG=\Aut_\om\C\langle x,y\rangle$.
Similarly, one has: $\LL\subset \derv\C\langle x,y\rangle$.
Yet, this latter inclusion is {\it not} an equality, cf.
Question 17.10.~$\lozenge$
\end{remark}
\medskip

Write
$\Aut_{_W}(\hh_\ka)$ for the group of  automorphisms
$\,a: \hh_\ka\to \hh_\ka$ that restrict to the identity on $\C{W}\subset \hh_\ka.$
Similarly, let
$\derw(\hh_\ka)$ denote the Lie algebra of all derivations
$\gamma: \hh_\ka\to\hh_\ka$ such that $\gamma(\C{W})=0$.
\begin{theorem}\label{main_der} Let $\ka\in\CP^1$
and $\hh_\ka=\hh_\ka(S_n).$ For any $F\in \GG$,
resp. $F\in \LL$,
the assignment:
$$x_i \mapsto P_{_F}(x_i,y_i)\quad\text{and}\quad
y_i\mapsto Q_{_F}(x_i,y_i)
\quad,\quad i=1,\ldots,n,$$
extends uniquely to an algebra 
automorphism $a_{_F}: \hh_\ka\to\hh_\ka,\,$ resp. a
 derivation $\gamma_{_F}: \hh_\ka\to\hh_\ka.\,$ The map:
$F\mapsto a_{_F},\,$ resp. $F\mapsto \gamma_{_F},\,$ thus defined
gives a group homomorphism: $\GG\to \Aut_{_W}(\hh_\ka),\,$ resp.,
 a Lie algebra homomorphism:
$\LL\to\derw(\hh_\ka)$.
\end{theorem}

\Pf . The proofs of the `group' and `Lie algebra'
parts are completely analogous, so we prove only the Lie algebra part.

Formula (\ref{der_ka}) applied to the element $f=\psi_{\ell+1}$,
resp. $f=\phi_{\ell+1}$, see
(\ref{power_sum}), yields readily
$$\gamma_{_{\psi_{\ell+1}}}(y_i)=-(\ell+1)\cdot
x_i^\ell\quad\text{and}\quad
 \gamma_{_{\phi_{\ell+1}}}(x_i)=0\quad,\quad\forall i=1,\ldots,n\,.
$$
It follows that, for $F=\alpha_\ell\in \LL,$
the assignment of the Theorem gives indeed  a well-defined
derivation $\gamma_{_F}\in \derw(\hh_\ka)$.
Similarly, applying formula  (\ref{der_ka}) to $f=-\phi_{\ell+1}$
shows that, for $F=\beta_\ell\in \LL$,
the assignment of the Theorem  gives indeed  a well-defined
derivation $\gamma_{_F}\in \derw(\hh_\ka)$.

It is clear now that, if   $F_1,F_2 \in \LL$ are such that
 the assignments $\gamma_{_{F_1}}$ and $\gamma_{_{F_2}}$
of the Theorem extend to derivations of $\hh_\ka$,
then $\,[\gamma_{_{F_1}},\gamma_{_{F_2}}]\,$
 is a derivation  of $\hh_\ka$ again, and
one has:
$$[\gamma_{_{F_1}}, \gamma_{_{F_2}}]\,:\;\;
x_i \mapsto P_{_{[F_1,F_2]}}(x_i,y_i)\quad\text{and}\quad
y_i\mapsto Q_{_{[F_1,F_2]}}(x_i,y_i)
\quad,\quad i=1,\ldots,n.
$$
Thus, we get: $[\gamma_{_{F_1}}, \gamma_{_{F_2}}]=
\gamma_{_{[F_1,F_2]}}$.
This completes the proof.\quad\sq
\bigskip\medskip

\newpage
\centerline{\large{{\bf PART 2. HARISH-CHANDRA HOMOMORPHISMS}}}

\section{The radial part construction}\label{Harish}
\setcounter{equation}{0}

We begin with the following generalization of the classical
 construction of  `radial part' of an invariant differential operator.

Given a smooth complex algebraic variety $X$,
let $\dd(X)$ denote the algebra of regular differential operators on
$X$. 
Further, given a (possibly
infinite dimensional) $\C$-vector space
$A$, we write $\dd(X, A):= A\bigotimes \dd(X)$ for the
vector space of $A$-valued regular differential operators.
If $A$ is an  associative $\C$-algebra, the space $\dd(X, A)$ acquires a
natural associative algebra structure, the tensor product
of that on $A$ and on $\dd(X)$.

Let $G$ be a  connected reductive algebraic group
over $\C$ with Lie algebra $\g$,
and $T \subset G$ a maximal torus. Let
$\gr\subset \g$ be the subset of semisimple regular elements.
Let $\h\subset \g$ be the Cartan subalgebra corresponding to $T$,  and write
 $\hr=\h\cap\gr$ for the Zariski open subset formed by the
 regular elements of $\h$,
i.e. for the complement of the root hyperplanes. Let $\Ug$ denote
the enveloping algebra of the Lie algebra $\g$, and
$(\Ug)^{\ad^{\,}\h}$ denote the centralizer of $\h$ in $\Ug$.
Then $(\Ug)^{\ad^{\,}\h}\!\cdot\!\h$ is a two-sided ideal in the
algebra $(\Ug)^{\ad^{\,}\h}$, and we set $(\Ug)_\h =
(\Ug)^{\ad^{\,}\h} /(\Ug)^{\ad^{\,}\h}\!\cdot\!\h$. The action of
the Weyl group $W$ on $\h$ preserves the subset $\hr$, hence
induces a $W$-action on $\dd(\hr)$ and on $(\Ug)_\h$. Let
$\,\dd\bigl(\hr,(\Ug)_\h\bigr)^W\,$ denote the subalgebra of
$W$-invariants in the algebra ${\dis\dd\bigl(\hr,(\Ug)_\h\bigr)}
= (\Ug)_\h\,\bigotimes_{_\C}\, \dd(\hr)$.

Given a  $\g$-module $V$, let $V\!\langle 0\rangle
:=V^\h\,$ denote
its zero weight subspace.
It is clear  that $V\!\langle 0\rangle$ is an $(\Ug)^{\ad^{\,}\h}$-stable
subspace, and the
action of $(\Ug)^{\ad^{\,}\h}$ on $V\!\langle 0\rangle$ factors through
$(\Ug)_\h$. Thus, the space
$V\!\langle 0\rangle$ has a natural
$(\Ug)_\h$-module structure, and
the space
$V\!\langle 0\rangle\otimes \C[\hreg]$
 has a natural $\dd\bigl(\hr,(\Ug)_\h\bigr)$-module structure.

\begin{proposition}\label{psi}
There exists a canonical  algebra isomorphism
\[\,\widetilde{\Psi}:\;\; \dd(\gr)^G \iso
\dd\bigl(\hr,\, (\Ug)_\h\bigr)^W\,,\]
such that, for any $f\in \bigl(V\otimes
\C[\gr]\bigr)^\g$ and
$L\in\dd(\gr)^\g $, one has: $\;\dis
(Lf)|_{\hreg}=\widetilde{\Psi}(L)(f|_{\hreg}).
$
\end{proposition}

\Pf . The map $p: G/T \times \hreg \stackrel{W}{\too} \gr\,,
\, (g,h)\mapsto \Ad g(h),\,$ gives the standard Galois covering of
$\gr$ with the Galois group $W$, the Weyl group, see e.g. [CG].
By the properties of Galois covering, the
pull-back via $p$ identifies the algebra
$\dd(\gr)$ with $\dd(G/T \times \hreg)^W$, the algebra of $W$-invariant
differential
operators on $G/T \times \hreg$. Now, the group $G$ acts on
$G/T \times \hreg$ by means of left translations on the first factor
$G/T$. The map $p$ intertwines this $G$-action with the $\Ad G$-action
on $\gr$. It follows that pull-back via $p$ induces an algebra isomorphism:
$\,\dd(\gr)^G \iso \dd(G/T \times \hreg)^{G\times W}\,$
(note that the actions of $G$ and $W$ on  $G/T \times \hreg$ commute).

In general, one has: $\,\dd(G/T \times \hreg) =\dd(G/T)\otimes
\dd(\hreg),\,$
hence, $\,\dd(G/T \times \hreg)^G = \dd(G/T)^G\otimes
\dd(\hreg),\,$ since the group $G$ acts trivially on $\hreg$.
We now recall the well-known isomorphism:
$\dd(G/T)^G \simeq (\Ug)_\h,\,$ see e.g. [He].
Thus, we obtain:
$$ 
\dd(\gr)^G \iso \dd(G/T \times \hreg)^{G\times W}
\simeq \bigl((\Ug)_\h\otimes
\dd(\hreg)\bigr)^W= \dd\bigl(\hreg\,,\, (\Ug)_\h\bigr)^W.\quad\square
$$
\medskip

\noindent
{\bf{The
`universal' Harish-Chandra homomorphism.}}\quad
 Following the classical construction of
 Harish-Chandra, we will use
not the radial part homomorphism $\widetilde{\Psi}$ but its
modification, obtained by conjugating $\widetilde{\Psi}$ by the
Weyl denominator. Namely, choose the set $R_+\subset R$ of
positive roots, and let $\,\delta=\prod_{_{\alpha\in
R_+}}\alpha\,$ denote the product of positive roots. For any
$L\in \dd(\gr)^G$, put $\,\Psi(L):= \delta\ccirc
{\widetilde{\Psi}}(L)\ccirc \delta^{-1} \in
\dd\bigl(\hreg,(\Ug)_\h\bigr)^W,\,$ that is,
for any $f\in \C[\hreg, (\Ug)_\h],$ we have
$\Psi(L): f \mapsto \delta^{-1}\cdot {\widetilde{\Psi}}(L)(\delta\cdot f).\,$
This definition does not
depend on the choice of the set
 of positive roots.

Let  $\,\{e_\alpha\}_{\alpha\in R}\,,$ denote 
 root vectors from
 a Weyl basis
of the Lie algebra $\g$, normalized so that $(e_\alpha,e_{-\alpha})=1$.

\begin{proposition}\label{rapart} The image of the Laplacian $\Delta_{\g} \in \dd(\g)^G$ under
the universal Harish-Chandra homomorphism $\Psi$ is given by the
formula:
\begin{equation}\label{CM}
\Psi(\Delta_{_\g})\;=\;
\Delta_{_\h} -\sum\nolimits_{\alpha\in R}\;
\frac{e_\alpha\cdot e_{-\alpha}}{\alpha^2} 
\quad\; \text{\bf{(Spin-Calogero-Moser
operator)}}\;.
\end{equation}
\end{proposition}

\Pf .  For $x\in
\g$, write $\frac{\partial}{\partial x}$ for the operator of 
directional derivative along $x$,
abusing the notation put $e^x=\exp(x)\in G$. We have:
$\;\dis
\Delta_\g=\Delta_\h+$ $\sum\nolimits_{\alpha}\;\frac{\partial^2}{\partial
{e_{-\alpha}}\cdot 
\partial  {e_{\alpha}}}
.\,$

Let $V$ be a finite dimensional $\g$-module, and
$f$ a $\g$-equivariant $V$-valued regular function on a neighborhood
of $\hreg$ in $\gr$.
Fix $h\in \hreg$, and $s\in \C$. 
For $t\in \C$, we
have
$$
e^{t\cdot  e_\alpha}\bigl(f(h+s\cdot  e_{-\alpha})\bigr)
= f\bigl(\Ad(e^{t\cdot  e_\alpha})(h+s\cdot
e_{-\alpha})\bigr)= f\bigl(h+s\cdot e_{-\alpha}-t\alpha(h)\cdot
e_\alpha+ st\cdot h_\alpha
+O(t^2)\bigr).
$$
Differentiating this equation with respect to $t$, at $t=0$ we find
\[
e_\alpha \bigl(f(h+s\cdot  e_{-\alpha})\bigr)= -\alpha(h)\cdot
\frac{\partial f}{\partial  {e_\alpha}}(h+s
\cdot e_{-\alpha})+s\cdot
 \frac{\partial f}{\partial  { h_\alpha}}(h+s\cdot  e_{-\alpha}).
\]
We rewrite the latter equation as follows
\begin{equation}\label{seqn1}
\frac{\partial f}{\partial  {e_\alpha}}(h+s\cdot  e_{-\alpha})= 
\frac{s}{\alpha(h)}\cdot \frac{\partial f}{\partial  
h_\alpha}(h+s\cdot  e_{-\alpha})-\frac{1}{\alpha(h)}\cdot e_\alpha\bigl( f(h+s\cdot 
e_{-\alpha})\bigr).
\end{equation}
Differentiating equation (\ref{seqn1}) with respect to $s$ at $s=0$,
and using equation (\ref{seqn1}) for the root $-\alpha$, we obtain
\[
\frac{\partial^2 f}{\partial  {e_{-\alpha}}\partial  {e_\alpha}}=
\frac{1}{\alpha}\cdot\frac{\partial f}{\partial  {h_\alpha}}
-\frac{1}{\alpha^2}\cdot
e_\alpha e_{-\alpha} f\;.
\]
Thus, we find that
\begin{equation}\label{ppp}
\widetilde\Psi(\Delta_{_\g})\;=\; \Delta_{_\h}
+\sum\nolimits_{\alpha\in R}\;\frac{1}{\alpha}\cdot\frac{\partial}{\partial h_\alpha}
-\sum\nolimits_{\alpha\in R}\; \frac{e_\alpha\cdot
e_{-\alpha}}{\alpha^2}\;,
\end{equation}
Conjugating this equation by $\delta$  kills the first order term 
in (\ref{ppp}), and we are done.\sq\medskip

 Now fix a $\Ug$-module $V$. As has been explained
above, the zero-weight space $V\!\langle 0\rangle=V^\h$ is
$(\Ug)^{\ad^{\,}\h}$-stable, and is in effect an  $(\Ug)_\h=
(\Ug)^{\ad^{\,}\h}/(\Ug)^{\ad^{\,}\h}\cdot\h$-module.
 This way one gets an algebra
homomorphism $\chi: (\Ug)_\h\to \End_{_\C}V\!\langle 0\rangle$.
The composite homomorphism:
\begin{equation}
\label{psiv}
\Psi_V:\quad\dd(\g)^\g\;
\stackrel{\Psi}{\too}\;
\dd\bigl(\hreg, (\Ug)_\h\bigr)^W\;\stackrel{\chi}{\too}\;
\dd\bigl(\hreg,\,\End_{_\C}V\!\langle 0\rangle\bigr)^W.
\end{equation}
will be referred to as the {\it Harish-Chandra homomorphism
associated to} $V$. \medskip

In the case of a {\it finite-dimensional} representation $V$, the
construction of the homomorphism $\Psi_V$  can be equivalently
described as follows. We start with the natural
$\dd(\g)^\g$-action on $\C[\gr]$. The $\dd(\g)^\g$-action
clearly  commutes with  the $\adg$-action on $\C[\gr]$, hence
preserves each $\adg$-isotypic component of  $\C[\gr]$. This
gives, for any finite-dimensional $\g$-module $V$, a
$\dd(\g)^\g$-module structure on the vector space
$\Hom_\g\bigl(V^*\,,\,\C[\gr]\bigl)$. The space  $\C[\gr]$ has,
by Frobenius reciprocity, see Kostant [Ko1], the following
decomposition into $\adg$-isotypic components:
\begin{equation*}\label{kostant}
\C[\gr]\;\;=\;\;
\bigoplus_{\text{finite-dimensional simple $\g$-modules $V$}}\;\;
V^*\otimes V\!\langle 0\rangle\otimes \C[\gr]^\g
\end{equation*}
Therefore, we find:
$$\Hom_\g\bigl(V^*\,,\,\C[\gr]\bigl)\;\simeq\;
V\!\langle 0\rangle\otimes \C[\gr]^\g\;\simeq\;
V\!\langle 0\rangle\otimes \C[\hreg]^W\,.
$$
Thus, one obtains a $\dd(\g)^\g$-action on $V\!\langle
0\rangle\otimes \C[\hreg]^W.\,$ It is easy to see that, for any
$L\in \dd(\g)^\g$, the action of $L$ arising in this way is given
by a certain $\End_{_\C}V\!\langle 0\rangle$-valued differential
operator on $\Spec\bigl(\C[\hreg]^W\bigr)=\hreg/W$. But since the
projection $\pi:\hreg\onto\hreg/W$ is \'etale, the pull-back via
$\pi$ gives a differential operator $\tilde\Psi_V(L)\in \dd(\hreg,
\End_{_\C}V\!\langle 0\rangle)$. Then one can define $\Psi_V(L)$
by $\Psi_V(L)=\delta\circ \tilde\Psi_V(L)\circ \delta^{-1}$. We
leave it to the reader to check that the map: $L\mapsto\Psi_V(L)$
thus defined is nothing but the one given by formula
(\ref{psiv}).\medskip

\noindent
{\bf {`The' Harish-Chandra homomorphism.}}\quad
Let $\Phi: \dd(\g)^\g\to \dd(\hr)^W$
denote the Harish-Chandra homomorphism
associated to the trivial 1-dimensional $\g$-module $\C$.
It is easy to see that the following definition
agrees with the classical notion of the Harish-Chandra homomorphism:

\begin{definition} The map $\Phi$ is called `the'
Harish-Chandra homomorphism.
\end{definition}

By a fundamental  result of
Harish-Chandra, \cite{HC}, the image of $\Phi$ consists
of {\it regular} differential operators  on the whole of $\h$.
This is not true for Harish-Chandra homomorphisms
associated to nontrivial $\g$-modules $V$.
\medskip

Recently, thanks to Wallach, and Levasseur-Stafford, the
results of Harish-Chandra have been strengthened as follows.
Let $\ad: \g \into \dd(\g)$ denote the Lie algebra map
sending $x\in \g$ to the vector field $\ad^{\,}x$ on $\g$.
Form the left ideal $\dd(\g)\cdot\adg \subset \dd(\g)$,
generated by the image of the map $\ad$, and set
$\,I_{\ad}:= \bigl(\dd(\g){}^{\!}\cdot{}^{\!}
\adg\bigr)\,\cap\,\dd(\g)^\g\,.$
It is clear that $I_{\ad}$ is a two-sided ideal in $\dd(\g)^\g$;
moreover, the Harish-Chandra homomomorphism $\Phi$
vanishes on $I_{\ad}$, hence descends to a map $\Phi:
\dd(\g)^\g/I_{\ad}\to \dd(\h)^W$.
The theorem below summarizes the results of Wallach \cite{Wa}, and
Levasseur-Stafford
\cite{LS1,LS2}.

\begin{theorem}\label{WLS}
The map
$\Phi: \dd(\g)^\g/I_{\ad}\iso \dd(\h)^W$ is an algebra isomorphism.\sq
\end{theorem}

The algebras $\dd(\g)^\g$ and $\dd(\h)^W$ come equipped with
canonical filtrations by the order of differential operators,
and there are canonical algebra isomorphisms:
$\grd\bigl(\dd(\g)^\g\bigr)\simeq \C[\g\oplus\g]^\g,$
and $\grd\bigl(\dd(\h)^W\bigr)\simeq \C[\h\oplus\h]^W.$
In particular, $\grd(I_{\ad})\subset \grd\bigl(\dd(\g)^\g\bigr)$
may be viewed as an ideal in $\C[\g\oplus\g]^\g$.
One can describe this ideal as follows.

Let $J\subset \C[\g\oplus\g]$ denote the ideal
generated by the  functions:
$\g\oplus\g\to \C$ of the form:
$(X,Y)\mapsto \lambda([X,Y])\;,\,\lambda\in\g^*.$
Write $\sqrt{J}$ for the
{\it radical} of $J$, and let:
$\bigl(\sqrt{J}\bigr)^\g= \sqrt{J}\,\cap\, \C[\g\oplus\g]^\g$.

\begin{theorem}[\cite{Le}]\label{LJ}
One has: $\grd(I_{\ad})=\bigl(\sqrt{J}\bigr)^\g$.\quad\sq
\end{theorem}
This theorem has the following geometric meaning.
Let
$\,{\mathcal Z}=\{(X,Y)\in \g\oplus \g
\;\;\big|\;\;[X,Y]=0\}\,$ be
  the commuting variety of $\g$.
By Hilbert Nullstellensatz,
 $\,\C[{\mathcal Z}]=\C[\g\oplus\g]/\sqrt{J}$.
Theorem \ref{LJ} says:
$\grd\bigl(\dd(\g)^\g\bigr)/\grd(I_{\ad}) =
\C[{\mathcal Z}].\,$
Observe now that the natural inclusion
$\,i:\h\oplus\h\into {\mathcal Z}$ gives the
restriction map $\,i^*: \C[{\mathcal Z}]\too\C[\h\oplus\h].$

\begin{theorem}\label{Joseph} \vi
The map $i^*$ induces an algebra  isomorphism:
$\C[{\mathcal Z}]^G\iso \C[\h\oplus\h]^W$.

\vii
The associated graded
map $\,\grd(\Phi):\, \grd\bigl(\dd(\g)^\g\bigr)/\grd(I_{\ad})
\too \grd\bigl(\dd(\h)^W\bigr)\,$
coincides with the map in \vi, hence is an isomorphism.
\sq
\end{theorem}

 The
map $\,i^*: \C[{\mathcal Z}]^G\to \C[\h\oplus\h]^W$ is a
double-analog of the Chevalley restriction isomorphism.
Injectivity of $i^*$ is easy for $\g=\gln$, and is due to Richardson, in
general.
Surjectivity of $i^*$ follows,
for $\g=\gln$, from classical invariant theory
[We], and is due to Joseph [J] in the
general case.\medskip

\noindent
{\bf Remark.}
The ideal $I_{\ad}
\subset
\dd(\g)^\g$ can be thought of as a non-commutative analogue of the ideal
$J^\g\subset\C[\g\oplus \g]^\g$. In view of Theorem \ref{LJ}
we can regard  Theorem
\ref{WLS}, saying that $I_{\ad}={\Ker}(\Phi)$,
as a non-commutative analog of the
equality $J^\g=\bigl({\sqrt{J}}\bigr)^{\,\g}$.
We cannot prove this equality;
it is a special case of a well known
conjecture saying that $J=\sqrt{J}$.

\section{Deformation of the Harish-Chandra homomorphism}\label{mainr}
\setcounter{equation}{0}

 We would like to construct a 1-parameter deformation of the
Harish-Chandra homomorphism. To do this, we  construct
a 1-parameter
family of representations $\{V_k\}_{\,k\in\C}$ of $\g$,
such that, for any $k$, the space $V_k\langle 0\rangle$ is
1-dimensional, and $V_0\langle 0\rangle$ is the trivial representation
of $(\Ug)_\h$. Then we may define
a deformation of $\Phi$ by setting
$\Phi_\ka=\Psi_{_{V_k}}$.

We now produce a particular 1-parameter
family of representations $V_k$ in the case $\g=\gln$
(it turns out that such a family does
not exist for Lie algebras which are not
of type ${\mathbf{A}}$). We think of $\gln$ as the Lie algebra of
linear vector fields on the vector space
$\C^n$, with coordinates  $x_1,...,x_n$,
so that
$E_{ij}=x_i\frac{\partial}{\partial x_j}$
are the  matrix units forming the standard basis of $\gln$.
Define $V_k$ to be the $\C$-vector space formed by all the
expressions
$(x_1\cdot\ldots\cdot x_n)^k\cdot P\,,$ where $P$ is
a Laurent polynomial in $x_1,...,x_n$ of total degree zero,
with the Lie algebra $\gln$ acting by formal
differentiation. Specifically, 
for each $i,j$, there is a well-defined action on the  expressions
above 
of the differential operator
$x_i\frac{\partial}{\partial x_j}$.
The assignment: $
E_{ij}\mapsto x_i\frac{\partial}{\partial x_j}\,$
makes $V_k$ an $\gln$-module.
The $\gln$-module that we need is obtained from this one by
redefining the action of the center of $\gln$. Specifically,
we first consider  $V_k$ as an $\sln$-module, obtained by restricting
the action from $\gln$ to $\sln$. Now, define a {\it new} $\gln$-action on
 $V_k$ by pulling-back the $\sln$-action via the 
natural Lie algebra projection: $\gln\onto\sln$,
so that the center of $\gln$ acts trivially.
This is the $\gln$-module $V_k$ that we wanted to produce.
The representation $V_k$ has an interpretation in terms of twisted
differential operators explained in Appendix B.

It is clear that
the $\g$-module $V_k$ thus obtained
 has a 1-dimensional zero weight space
$V_k\langle 0\rangle$, spanned by the element
$(x_1\cdot\ldots\cdot x_n)^k\cdot 1$, and
$V_0\langle 0\rangle$ is the trivial representation of  $(\Ug)_\h$.

\begin{definition} For $\g=\gln$, or $\g=\sln$, we call 
$\,\Phi_k=\Psi_{_{V_k}}\,:\; \dd(\g)^\g\to \dd(\hr)^W$
the {\it deformed Harish-Chandra homomorphism} with parameter $k\in \C$.
\end{definition}

It is clear that if $k=0$, then the map $\Phi_k$ coincides with
$\Phi$ (this follows from the fact that any equivariant
function on $\gr$ with values in $V_0$ is, in effect, an invariant
function
taking values in the trivial 1-dimensional subrepresentation of $V_0$).

Let $\cc_k$, denote the centralizer in $\widetilde{\dd(\hr)^W}$
of the Calogero-Moser differential operator $\CM_k$, see (\ref{H}).
Identify $(S\g)^\g$ with the algebra of
$\adg$-invariant constant
coefficient differential operators on $\g$.

\begin{proposition}\label{phi(cc)} 
\vi For any $k\in \C$, we have:
$\;\Phi_{k}(\Delta_\g)=\CM_k$.

\vii The map $\Phi_k$ gives an algebra isomorphism:
$(S\g)^\g\iso \cc_k$.
\end{proposition}

\Pf .  For any root $\alpha$, one finds that:
$e_\alpha e_{-\alpha}\big|_{_{V_k\langle
0\rangle}}=k(k+1)\cdot\id_{_{V_k\langle
0\rangle}}$.
Part (i) now follows from Proposition \ref{rapart}.
To prove (ii), observe first that any element of
$\Phi_k\bigl((S\g)^\g\bigr)$ commutes with
$\Phi_k(\Delta_\g)=\CM_k,$ since $(S\g)^\g$ is a commutative
algebra, and belongs to $\widetilde{\dd(\hr)^W}$. It 
follows that $\Phi_k\bigl((S\g)^\g\bigr)\subset
\cc_k$. Next, equip $(S\g)^\g$ and $\cc_k$ with filtrations given
by the order of  differential operator. The map $\Phi_k$ is
clearly filtration preserving, hence, induces a map
$\grd(\Phi_k): \grd\bigl((S\g)^\g\bigr) \to \grd(\cc_k)$. We have
an obvious isomorphism: $\grd\bigl((S\g)^\g\bigr)\simeq (S\g)^\g$.
Further, according to Opdam \cite{O}, cf. also \cite{VSC}, there
exists an algebra isomorphism $\,\sigma_k: S\h^W\iso \cc_k,\,$
such that, for any homogeneous $u\in S\h^W$ of degree $\ell$, the 
the operator $\sigma_k(u) $ has order $\ell$, and its principal symbol
 equals $u$. It follows that $\grd(\cc_k)=S\h^W$.
Therefore, we may identify the map $\grd(\Phi_k)$ with the Chevalley
isomorphism: $(S\g)^\g\iso S\h^W$. Hence,
$\grd(\Phi_k)$ is bijective. Thus, $\Phi_k$ is a bijection as
well.~\sq\medskip

\noindent
{\bf Remark.} Part (i) also follows from \cite{E,EFK}. Specifically,
 in \cite{EFK} the authors consider radial parts of
differential operators on the group $G$, rather than on its Lie
algebra, and obtain Calogero-Moser operators with trigonometric
potential $\,k(k+1)/\sin^2(x)\,$ rather than $\,k(k+1)/x^2\,$.
The rational case may be obtained by taking an appropriate
limit of the result of \cite{EFK}.
$\enspace\lozenge$\medskip

Now recall that we have defined in the introduction the algebra
$\bcal_k\subset \dd(\hreg)^W$
 generated by $\cc_k$ and by
the polynomial algebra, $\C[\h]^W\subset \dd(\h)^W$, of zero order
 $W$-invariant operators on $\h$.

Recall the Lie algebra map
$\ad : \g \to \dd(\g)\,,\, x\mapsto \ad x.\,$
This Lie algebra map extends uniquely to an associative algebra homomorphism
$\ad: \Ug \to \dd(\g)$.
Further, let
$\ann_k\subset
\Ug $ denote the annihilator of the representation $V_k$,
or,
equivalently, the kernel of (\ref{varrho}). Thus, $\ann_k$ is a two-sided ideal
in $\Ug$, and it is easy to see that $\adg$-invariants in the left ideal
$\dd(\g)\cdot\ad(\ann_k)\subset \dd(\g)$ form a  two-sided ideal:
$\bigl(\dd(\g)\cdot\ad(\ann_k)\bigr)^\g\,\subset \dd(\g)^\g$.

The following theorem is one of the main results of Part 2 of this paper.

\begin{theorem}\label{inject} \vi
For any $k\in\C$, we have: $\;\bcal_k \,=\,\Image(\Phi_k)$.

\vii $\;\Ker(\Phi_k)\,= \,\bigl(\dd(\g)\cdot\ad(\ann_k)\bigr)^\g,\,$
for all $k\in\C,$ except possibly a finite set.
\end{theorem}

\begin{remark} The definition of the homomorphism
$\Phi_k$ easily implies that, for any $k\in\C,$ one has
the  inclusion:
$\,\bigl(\dd(\g)\cdot\ad(\ann_k)\bigr)^\g\subset\Ker(\Phi_k),\,$
which is (for generic $k$) a part of claim (ii) of the Theorem.
The  opposite inclusion is  harder,
it will be proved only in \S10.
Further, notice that $\Phi_k\bigl(\C[\g]^\g\bigr)=\C[\h]^W$.
Hence,  Proposition \ref{phi(cc)}
yields an inclusion: $\;\bcal_k\subset\Image(\Phi_k)$,
which is part of claim (i) of the Theorem.
The opposite inclusion is  more difficult and will be proved
in \S10.$\quad\lozenge$
\end{remark}

\begin{corollary}\label{isom}
For almost all $k\in\C$, the map $\Phi_k$ induces an algebra
isomorphism
\vskip 2pt

\centerline{
$\,\Phi_k:\dd(\g)^\g/\bigl(\dd(\g)\cdot\ad(\ann_k)\bigr)^\g\iso\bcal_k.\qquad\square$}
\end{corollary}
\medskip 

\noindent
{\bf Hamiltonian interpretation:}\quad
 The $\Ad G$-action on $\g=\sln$
gives rise to a Hamiltonian $G$-action on $T^*\g$, the cotangent bundle
on $\g$, with moment map $\mu: T^*\g \to \g^*$.
The algebra homomorphism $\ad: \Ug \to \dd(\g)$ may be viewed 
as a "quantization" of
the pull-back morphism $\mu^*: \C[\g^*] \to \C[T^*\g]$.

Further, it is known, see e.g. [BJ],
 that $\ann_k$ is the primitive ideal
in $\Ug$ associated to the `minimal' coadjoint orbit 
$\,\O_k:=\{s\in \sln\simeq \g^*\;\,\big|\;\, 
\rk\bigl(s-(k+\frac{1}{2})\cd\id\bigr)=1\},$
via the {\it Geometric quantization} construction.  Put
$\FI_{\mathtt{quantum}}:=\ann_k$, and let
$\FI_{\mathtt{classic}}$ denote the defining ideal of the (closed)
orbit $\O_k$.
Thus, according to the "geometric quantization"
philosophy, the algebra $\Ug/\FI_{\mathtt{quantum}}$ should be thought
of as
a quantization of $\C[\O_k]=\C[\g^*]/\FI_{\mathtt{classic}}$.

As has been discovered by  Kazhdan-Kostant-Sternberg [KKS], the classical
Calogero-Moser
space can be obtained as the Hamiltonian reduction: $\mu^{-1}(\O_k)/G$,
cf. also \S11 below. The coordinate ring of $\mu^{-1}(\O_k)/G$, an
affine algebraic variety, by definition is:
$\C[\mu^{-1}(\O_k)/G] = \C[T^*\g]^\g\big/\bigl(\C[T^*\g]\cdot
\mu^*\FI_{\mathtt{classic}}\bigr)^\g.$ Thus,
the non-commutative algebra
$\,\dd(\g)^\g\big/\bigl(\dd(\g)\cdot\ad(\FI_{\mathtt{quantum}})\bigr)^\g\,$
is
a natural  "quantization" of the coordinate ring
of the classical
Calogero-Moser
space. According to Corollary \ref{isom}, this is nothing but
the image of the deformed Harish-Chandra homomorphism.
$\quad\lozenge$\medskip

\noindent 
{\bf Spherical Harish-Chandra isomorphism.}\quad 
Let  $\hh_{t,c}$ be the rational Cherednik
algebra associated to the
root system ${\mathbf{A_{n-1}}}$. In this case
 the function $c: R\to\C$
reduces to a constant. As has been explained in the Introduction,
we may (and will) view the parameter $\ka=(t,c)$ as a
point  $\ka=c/t\in\C\cup\{\infty\}.$
The special case: $t=0$ corresponds to the point:
$\ka=\infty \in \CP^1$. 

Let $t\neq 0$.
Then, inverting the bijection $\Theta_\ka^{^{\tt{spher}}}$ of Theorem
\ref{generic}(i),
see also Proposition \ref{typeEF},
we obtain an algebra homomorphism, called the `spherical Harish-Chandra
homomorphism', cf., (\ref{main_map}):
\begin{equation}\label{spher_map}
\Phi_\ka^{^{\tt{spher}}}\,=\,
\bigl(\Theta_\ka^{^{\tt{spher}}}\bigr)^{\,-1}\ccirc\Phi_\ka\;:
\quad
\dd(\gln)^{\gln} \too \e\hh_\ka\e\quad,\quad
 \forall \ka\in \CP^1\smallsetminus\{\infty\}\,.
\end{equation}
By Corollary \ref{isom}, the  induced map
$\,\Phi_\ka^{^{\tt{spher}}}:
\dd(\g)^\g/\bigl(\dd(\g)\cdot\ad(\ann_k)\bigr)^\g\too
\e\hh_\ka\e\,$ is 
an isomorphism for almost all $\ka\in \C$,
making $\e\hh_\ka\e$ a quantum Hamiltonian reduction of 
the algebra $\dd(\g)$.

\begin{proposition}\label{fourier_s}
For $\g=\gln,$ the spherical
 Harish-Chandra homomorphism, $\Phi_\ka^{^{\tt{spher}}}\!,$
intertwines the standard Fourier  automorphism
on $\dd(\g)^\g$ and the  automorphism $\FF$ on {$\e\hh_\ka\e.$}
\end{proposition}

\Pf . This is clear, since $\FF$ permutes $\C[\h]^W$ and $\cc_\ka$.
\sq\medskip

Let  $x_{ij}\,,\,i,j=1,\ldots,n,$ be the coordinates on
$\g=\gln$ corresponding to the standard
basis formed by the matrix units. Let $X=\|x_{ij}\|,Y=\|\frac{\partial}
{\partial x_{ji}}\|$.
For each $\ell=0,1,2,\ldots,$ the assignments, cf. \ref{ab} 
$$
a_{\tau,\ell}(X)= X+\tau
Y^\ell\;\;,\;\;a_{\tau,\ell}(Y)=Y\,,\enspace\text{resp.}\,,
\enspace
b_{\tau,\ell}(X)=X\;\;,\;\;b_{\tau,\ell}(Y)=Y+\tau X^\ell\;,\enspace\tau\in\C,
$$
give rise to well-defined algebra automorphisms
$\,a_{\tau,\ell},b_{\tau,\ell}:\dd(\g)\to\dd(\g).\,$
This way, one defines
 an action of the group $\GG$, see Definition \ref{LLL},
on $\dd(\g)$ by algebra automorphisms.
Furthermore, it is straightforward to verify that
the $\GG$-action commutes with the $\ad^{\,}\g$-action on
$\dd(\g)$. In particular, the subalgebra
$\dd(\g)^\g$ is stable under the
$\GG$-action.

Similarly, one defines an action of the Lie algebra
$\LL$ on $\dd(\g)^\g$ by algebra derivations.
The following result incorporates Proposition
\ref{fourier_s}   as a very special case.

\begin{proposition}\label{GG_phi}
For any $\ka\in\CP^1\smallsetminus\{\infty\}$,
the spherical Harish-Chandra homomorphism (\ref{spher_map})
intertwines the $\GG$-actions, resp. the $\LL$-actions,
on $\dd(\gln)^{\gln}$ and on $\e\hh_\ka\e$.
\end{proposition}

\Pf . Let us prove the statement in the Lie algebra case
(the group case is entirely analogous).
By definition, for any $\ell$,
the action of the derivation $
\alpha_\ell$  on $\dd(\gln)$, resp. on $\hh_\ka$, 
is given by the formula:
$$
\alpha_\ell=\frac{d}{d \tau}\Big|_{\tau=0} \Ad(e^{-\tau\cdot
 {\mathtt{Tr}}(X^{\ell+1}/(\ell+1))});
\quad\mbox{resp.}\quad
\alpha_\ell=\frac{d}{d \tau}\Big|_{\tau=0} 
\Ad(e^{-\tau\sum_i x_i^{\ell+1}/(\ell+1)})\,.
$$
Since $\,\Phi_\ka^{^{\tt{spher}}}$ clearly takes
the function: $X\mapsto {\mathtt{Tr}}(X^{\ell+1})$,
viewed as an element of $\dd(\g)^\g$, to
$\sum_i x_i^{\ell+1}\in \hh_\ka$,
it follows that $\,\Phi_\ka^{^{\tt{spher}}}$
commutes with $
\alpha_\ell\,,\,\forall\ell\geq 0.$

Further, we know that $\Phi_\ka(\Delta_\g)={\mathsf{L}}_\ka$,
and $\Theta_\ka^{^{\tt{spher}}}(\e\cd\Delta_\h\cd\e)={\mathsf{L}}_\ka,$
where ${\mathsf{L}}_\ka$
is the Calogero-Moser operator,
and $\Delta_\h=
\sum_i y_i^2$ denotes the
element of $\hh_\ka$ corresponding to the quadratic Casimir, 
viewed as an element of
$\C[\h]^W\subset \hh_\ka.$ Hence,
$\Phi_\ka^{^{\tt{spher}}}(\Delta_\g)$
$=(\Theta^{\text{spher}}_\ka)^{-1}\ccirc\Phi_\ka(\Delta_\g)=
(\Theta^{\text{spher}}_\ka)^{-1}({\mathsf{L}}_\ka)=\e\cd\Delta_\h\cd\e$.

It follows that 
the map $\Phi_\ka^{^{\tt{spher}}}$ commutes with
the derivation $\beta_2$. But the derivations 
$\,\{\alpha_2,\beta_2\}\,$ generate the Lie algebra $\slt$
that corresponds to the group of standard
$\SL_2$-automorphisms, see Corollary \ref{rotation}.
Therefore, $\Phi_\ka^{^{\tt{spher}}}$ commutes
with the standard
$\SL_2$-automorphisms, in particular, with
 the Fourier transform ${\mathsf{F}}$. Thus,
for any  $\ell$,
the map $\,\Phi_\ka^{^{\tt{spher}}}$ 
commutes
with $\beta_\ell={\mathsf{F}}\ccirc\alpha_\ell\ccirc
{\mathsf{F}}^{-1},$
and the Proposition follows.
\sq

\section{Example: $\g={\frak {sl}}_2$}
\setcounter{equation}{0}

In the case $\g={\frak {sl}}_2$, we will explicitly compute the
source and the target of the map $\Phi_k$ 
of Corollary \ref{isom} and verify  directly
 that both are flat deformations and that this map
is an isomorphism.

We first give an explicit description of the algebra
$\dd(\g)^\g/\bigl(\dd(\g)\cdot\ad(\ann_k)\bigr)^\g$
 in terms of generators and relations.
To this end, it is convenient to identify
$\g$ with the Lie algebra ${\so}$,
so that the adjoint representation of $\g$ becomes the
tautological representation of
${\so}$ on $\C^3$. Thus,
the algebra $\dd(\g)^\g$ gets identified with
 the algebra  of rotation-invariant differential operators
on $\C^3$. We have the following classical
result (see e.g. \cite{Sch}).

\begin{proposition}\label{so3} The algebra $\dd(\C^3)^{{\so}}$ is generated by
the squared radius $r^2$, the second order Laplacian $\Delta$, and the
Euler field $\eu$.
\end{proposition}

\Pf . It is the
standard result in elementary invariant theory
that the commutative algebra
$\C[\C^3\oplus\C^3]^{{\so}}$
is generated by 3 elements:
$\bold r^2,{\mathbf  r}{\mathbf
  p},\bold p^2$,
where $\bold r$
is the radius vector, and $\bold p$ is the momentum.
Now, equip $\dd(\C^3)^{{\so}}$
with filtration by the order of differential operators.
Then, we have: $\grd\bigl(\dd(\C^3)^{{\so}}\bigr)
=\C[\C^3\oplus\C^3]^{{\so}}$,
and our claim follows from its analogue for
$\C[\C^3\oplus\C^3]^{{\so}}$,
mentioned above.\sq\medskip

Let $e,h,f$ be the standard ${\frak {sl}}_2$-triple, $\g=\langle
e,h,f\rangle$. Proposition \ref{so3}
shows that $\dd(\g)^\g$ can be identified
with $\Ug $, via the map given by the assignment:
 $-\frac{1}{2}r^2\mapsto f$, $\frac{1}{2}\Delta\mapsto e$.
\medskip

Next, we give an explicit description of the algebra $\bcal_k$
in terms of generators and relations.
By the definition, $\bcal_k$ is
generated by the operators $x^2$ and $\partial^2-\frac{k(k+1)}{x^{2}}$
in one variable ($\partial$ denotes the differentiation
with respect to $x$).

Let
$C=ef+fe+h^2/2\in {\mathcal{U}}({\mathfrak{s}\mathfrak{l}}_2)$ 
be the quadratic Casimir, and set:
 $c_k=\frac{1}{2}(k-\frac{1}{2})(k+\frac{3}{2})$.

\begin{proposition}\label{assign}
The assignment: $f\mapsto -\frac{1}{2}x^2\,$, $\,e\mapsto \frac{1}{2}(\partial
^2-\frac{k(k+1)}{x^{2}})\,$
extends to a well-defined algebra isomorphism
$\,\eta:\Ug /(C-c_k)\!\cdot\!\Ug\iso
\bcal_k$.
\end{proposition}

\Pf .  It is straightforward
to check that $\eta$ extends to a surjective algebra homomorphism:
$\Ug\onto\bcal_k$. To show that $\eta$ is an isomorphism,
write $M(\mu)$ for the  Verma
module over $\g={\frak {sl}}_2$ with highest weight $\mu$.

The annihilator of the module $M(k-1/2)$ is
equal to the ideal $(C-c_k)\!\cdot\!\Ug\subset\Ug $, so that
the algebra $\Ug /(C-c_k)\!\cdot\!\Ug$
acts faithfully on this  Verma
module. On the other hand, consider the $\bcal_k$-module spanned by
$x^{2m-k}$, $m=0,1,2,...$.
It is immediate to check
that the pullback of this $\bcal_k$-module via  $\eta$
is an ${\frak {sl}}_2$-module isomorphic to $M(k-1/2)$.
It follows that the map $\eta$ is injective.\sq\medskip

\noindent
{\bf Remark.} It is known
(see e.g. \cite{St}) that the algebras $\Ug /(C-c)\!\cdot\!\Ug$,
$c\in \C$, are pairwise nonisomorphic.
Therefore, Proposition \ref{assign} shows that the algebras
$\bcal_k$ and $\bcal_l$ are not isomorphic unless $k=l$ or $k=-l-1.
\quad\lozenge$
\medskip

According to Propositions \ref{so3} and \ref{assign}, we may
view  the deformed Harish-Chandra homomorphism
 $\Phi_k: \dd(\g)^\g\to\bcal_k$,
as a certain map: $\Ug\too \Ug /(C-c_k)\!\cdot\!\Ug.$
We are going to check, by an explicit calculation of radial parts,
that this map coincides with the canonical projection:
 $\Ug\onto \Ug /(C-c_k)\!\cdot\!\Ug.$

For simplicity consider $k\in \Bbb Z^+$, the general situation being
similar.
In this case, the zero weight component
of an equivariant function: $\g\to V_k$,
in spherical coordinates has the form:
$g(r)\cdot P_k(\text{cos}\,\theta)$,
where $g(r)$ is a radial function, and $P_k$ is the $k$-th Legendre
polynomial. Thus, we get: $\Phi_k(r^2)=r^2$,
$\Phi_k(\eu)=r\partial_r-1$ (because of conjugation by $\delta=r$).
We compute $\Phi_k(\Delta)$. In spherical coordinates,
the Laplacian $\Delta$, acting on
longitude-independent functions reads
\[
\Delta=\frac{\partial^2}{\partial r^2}
+2\frac{1}{r}\frac{\partial}{\partial r}
+\frac{1}{r^2}\left(\frac{\partial^2}{\partial \theta^2}
+\text{cotan}\,\theta\cdot \frac{\partial}{\partial \theta}\right).
\]
Legendre polynomials satisfy the differential equation:
$$
\left(\frac{\partial^2}{\partial \theta^2}+
\text{cotan}\,\theta\cdot\frac{\partial}{\partial \theta}\right)
P_k(\text{cos}\,\theta)=
-k(k+1)\cdot P_k(\text{cos}\,\theta).
$$
Therefore, for $f=g(r)\cdot P_k(\text{cos}\,\theta)$, we get
$$
\Delta(f)=
\left[\left(\frac{1}{r}\ccirc\bigl(\frac{\partial^2}{\partial r^2}-
\frac{k(k+1)}{r^2}\bigr)\ccirc r\right)g\right]\cdot
P_k(\text{cos}\,\theta).
$$
Thus, $\Phi_k(\Delta)=\frac{\partial^2}{\partial r^2}-\frac{k(k+1)}{r^2}$
 (= Calogero-Moser
operator in one variable).

It is obvious now that $\Phi_k$ becomes,
after the identification of $\dd(\g)^\g$ with $\Ug $,
the projection: $\Ug \onto \Ug /(C-c_k)\!\cdot\!\Ug$.
\medskip

We may now  verify
Theorem  \ref{inject} for ${\frak {sl}}_2$, using the
explicit description of the homorphism $\Phi_k$ given above.
To prove the theorem, we must show that the ideal
$\bigl(\dd(\g)\cdot\ad(\ann_k)\bigr)^\g$ contains $C-c_k$.
It is easy to see that the eigenvalue of $C$ in $V_k$
equals
$2k(k+1)$. Thus,  $C-2k(k+1)\in
\ann_k$. One can check by a direct calculation that
$\ad\bigl(C-2k(k+1)\bigr)$ goes to $4(C-c_k)$ under the
isomorphism: $\dd(\g)^\g\iso \Ug$. Thus, the kernel of
$\Phi_k$ equals $\bigl(\dd(\g)\cdot\ad(\ann_k)\bigr)^\g$.\sq

\section{The kernel of the Harish-Chandra homomorphism}
\label{proof1}
\setcounter{equation}{0}

The main technical tool in the proof of Theorem \ref{inject}
is the following result (for $\g=\gln$)

\begin{theorem}\label{maintech} In $\dd(\g)^\g$, one has an equality:
 $\;\bigl(\dd(\g)\cdot\adg\bigr)^\g=
\bigl(\dd(\g)\cdot\ad(\ann_0)\bigr)^\g$.
\end{theorem}

\noindent
{\bf Remark.} Recall that according to Theorem \ref{WLS}:
$\;\bigl(\dd(\g)\cdot\adg\bigr)^\g=\Ker(\Phi)$,
is the kernel of the homomorphism $\,\Phi_{k=0}$.
Thus, Theorem \ref{maintech} insures compatibility of
Theorem \ref{WLS} with Theorem \ref{inject}.
Notice further that the statement of Theorem \ref{maintech}
might look
 somewhat surprising,
since in $\Ug$ one has a {\it strict} inclusion of ideals:
$\,\ann_0\;
\underset{{}^{\neq}}{\subset}\; \Ug\cdot\adg\,.\quad\lozenge$
\medskip

The remainder of this section is devoted to the proof of Theorem \ref{maintech}.
The argument below  is a refined version of the proof
of a commutative analog of this theorem, given in Appendix A.
Understanding this section, and its geometric meaning,
 should be easier for the reader
after having first read the proof of Theorem \ref{schemeiso}, and 
Lemma \ref{Rud}
of the Appendix.

We introduce some notation. Let $X=\|x_{ij}\|$ be an
$n\times n$-matrix of indeterminates, and $Y=
\|\partial_{ji}\|$ be the corresponding matrix of partial derivatives.
For any noncommutative polynomial $Q(x,y)$,
let $\lL{ Q(X,Y)}\rR $ denote the normal ordering of $Q(X,Y)$.
That is, each entry of $\lL{ Q(X,Y)}\rR $ is obtained by first
computing the corresponding entry of $Q(X,Y)$
(which is a sum of products of $x_{ij}$ and $\partial_{kl}$), and then
moving all the $\partial_{kl}$ to the right from $x_{ij}$.
In particular, we write: $\lL{[X,Y]}\rR =
\lL{ XY}\rR - \lL{ YX}\rR,$  where
$\lL{ XY}\rR =XY$,  $\,\lL{ YX}\rR =(X^tY^t)^t=YX-n\cdot\id$,
and "$t$" stands for "transposed".\medskip

It is known, see  e.g. \cite{BJ}, that the annihilator of $V_0$
in $\Ug$ is generated by the elements
$$
E_{ij}\cdot E_{kl}-E_{il}\cdot E_{kj}-\delta_{kj}\cdot E_{il}+\delta_{kl}\cdot E_{ij}
\,.
$$
Therefore, the left ideal $\JJ:=\dd(\g)\cdot \ad(\ann_0)\subset \dd(\g)$
is  generated
by the elements
\begin{equation}\label{tilI}
R_{jk}^{il}=\lL{ [X,Y]}\rR _{ji}\cdot \lL{ [X,Y]}\rR _{lk}-
\lL{ [X,Y]}\rR _{li}\cdot \lL{ [X,Y]}\rR _{jk}-\delta_{jk}\cdot \lL{ [X,Y]}\rR _{li}
+\delta_{lk}\cdot \lL{ [X,Y]}\rR _{ji}.
\end{equation}

 Let $I_{\ad}=\bigl(\dd(\g)\cdot\adg\bigr)^\g,\;$
and $I_0=\bigl(\dd(\g)\cdot\ad(\ann_0)\bigr)^\g =\JJ\cap \dd(\g)^\g.\,$
 It is clear that $I_0\subset
I_{\ad}\,.$ Hence, to prove the Theorem it remains to show
that any element $D\in I_{\ad}$
vanishes in $\dd(\g)/\JJ$.

Below, we view $\g=\gln$ as an associative (matrix) algebra.
We also consider the associative algebra $\dd(\g,\g):=
\g\otimes \dd(\g)$ of $\g$-valued differential operators on $\g$,
equipped with the standard filtration by the order of
differential operators. Then, $\grd\bigl(\dd(\g,\g)\bigr)
=\C[\g\oplus\g\,,\,\g]$,
is the (non-commutative) associative algebra
of polynomial maps: $\g\oplus\g\to\g$.
The filtration on $\dd(\g,\g)$ induces a filtration on the
subalgebra
$\dd(\g,\g)^{\adg}\subset \dd(\g,\g)$ of
$\adg$-invariant differential operators
(with respect to the simultaneous $\adg$-action on both
the source and the target spaces).

\begin{lemma} \label{Weyl}
The vector space $\,\dd(\g,\g)^{\adg}\,$ is spanned by operators of the form
$$
\prod\nolimits_{i=1}^m \;\Tr\bigl(P_i(X,Y)\bigr)\cdot Q(X,Y)\;,
\quad\mbox{where}\quad P_1,...,P_m,Q\quad\mbox{are noncommutative polynomials.}
$$
\end{lemma}

\Pf .  We prove the
lemma for an $\adg$-invariant $\g$-valued operator $D$
by induction on $d$, the order of $D$.
The principal symbol $D_d$ of $D$ with
respect to the above filtration is a $\g$-invariant
polynomial: $\g\oplus \g\to \g\,,\,
(X,Y)\mapsto D_d(X,Y).\,$ By Weyl's fundamental
theorem of invariant theory \cite{We},
this polynomial must be a linear combination of
expressions of the form $\prod_{i=1}^m \Tr(P_i(X,Y))\cdot Q(X,Y)$.
Therefore, subtracting from $D$ expressions of such form
(but now with noncommuting matrix entries),
we  reduce the order of $D$. \sq

\begin{lemma} \label{iprime} \quad
Any element of $I_{\ad}$ is a linear combination of
elements of the form: \quad

\noindent
$\dis
\prod\nolimits_i\;\Tr\bigl(P_i(X,Y)\bigr)\cdot
\Tr\bigl(Q(X,Y)\cdot \lL{ [X,Y]}\rR \bigr)
\,,\,$
where $P_i,Q$ are noncommutative
polynomials.
\end{lemma}

\Pf . Let $E_{ij}$ be an elementary matrix.
It is easy to show that the linear vector field
$x\to [E_{ij},x]$ on $\g$ equals
 $\lL{ [X,Y]}\rR _{ji}$. Therefore, the ideal $I_{\ad}$ is the set of
 $\g$-invariants in the left ideal in $\dd(\g)$ generated by the entries
of $\lL{ [X,Y]}\rR $. Thus, any element of $I_{\ad}$ is of the form
$\Tr\bigl(Q(X,Y)\cdot \lL{ [X,Y]}\rR \bigr)$, where $Q$ is an invariant differential
 operator on $\g$ with values in $\g$. Thus the claim follows
 from Lemma \ref{Weyl}.
\sq \medskip

We say that
$D\in I_{\ad}$ has level $\le d$ if
it is a linear combination of
elements of the form
$\prod_i\Tr\bigl(P_i(X,Y)\bigr)\cdot
\Tr\bigl(Q(X,Y)\cdot \lL{ [X,Y]}\rR \bigr)$
with  degree of $Q$ being $\le d$.
By Lemma \ref{iprime}, the level is well-defined for all
elements $D\in I_{\ad}$. Write: $\,\,Q:= Q(X,Y)$, for short.

\begin{lemma}\label{twocom} Let $Q_1$, $Q_2$ be noncommutative
  polynomials of $X,Y$ of degrees $d_1,d_2$
such that $d_1+d_2\le d-4$. Then, in $\dd(\g)/\JJ$,
cf. (\ref{tilI}), we have:
$
{\Tr\bigl(Q_1\cdot  \lL{ [X,Y]}\rR\cdot   Q_2\cdot  \lL{ [X,Y]}\rR\bigr)\!=\!0.}
$
\end{lemma}

\Pf . We start with the identity: 
\[
\Tr\bigl(Q_1\cd \lL{ [X,Y]}\rR\cd  Q_2\cd \lL{ [X,Y]}\rR \bigr) =
\sum\nolimits_{pqrs}\;\,(Q_1)_{pq}\cdot \lL{ [X,Y]}\rR _{qr}\cdot (Q_2)_{rs}\cdot
\lL{ [X,Y]}\rR _{sp}.\,
\]
Therefore, we have in $\dd(\g)/\JJ$:
$$
\Tr\bigl(Q_1\cdot \lL{ [X,Y]}\rR\cdot  Q_2\cdot \lL{ [X,Y]}\rR \bigr)=
\sum\nolimits_{pqrs}\,\;(Q_1)_{pq}\cdot (Q_2)_{rs}\cdot \lL{ [X,Y]}\rR _{qr}\cdot
\lL{ [X,Y]}\rR _{sp}
$$
Indeed, this identity holds
up to terms arising from reordering  the second and the third
 factors, which have lower
level and therefore are zero modulo $\JJ$
(note that all these terms belong to $I_{\ad}$ since they have
$\lL{ [X,Y]}\rR $ on the right end).

Looking at the generators (\ref{tilI}) 
of $\JJ$, we conclude that in $\dd(\g)/\JJ$ one has:
{\small
{$$
\Tr\bigl(Q_1\cdot \lL{ [X,Y]}\rR\cdot  Q_2\cdot \lL{ [X,Y]}\rR \bigr)=
\sum_{pqrs}\;(Q_1)_{pq}\cdot (Q_2)_{rs}\cdot
\bigl(\lL{ [X,Y]}\rR_{sr}\cdot \lL{ [X,Y]}\rR_{qp}+\delta_{qp}\cdot
\lL{ [X,Y]}\rR_{sr}
-\delta_{sp}\cdot \lL{ [X,Y]}\rR_{qr}\bigr),
$$}}
where $\delta_{ab}= $Kronecker delta.
The latter sum thus splits into three parts.
The first part, up to lower level terms, equals
$\;\dis
\Tr\bigl(Q_1\cdot \lL{ [X,Y]}\rR \bigr)\cdot\Tr\bigl(Q_2\cdot \lL{ [X,Y]}\rR \bigr).
\,$
Therefore it is zero in $\dd(\g)/\JJ$,
by the induction assumption.
The other two parts have level lower than $d$, so they  also
vanish in  $\dd(\g)/\JJ$.
\sq

\begin{lemma}\label{twocom2} Let $Q_1$, $Q_2$ be noncommutative
  polynomials of $X,Y$ of degrees $d_1,d_2$
such that $d_1+d_2\le d-4$. Then,  in $\dd(\g)/\JJ$ we have
an equality:
$$
\Tr\bigl(Q_1\cdot X\cdot Y\cdot Q_2\cdot \lL{ [X,Y]}\rR \bigr)=
\Tr\bigl(Q_1\cdot Y\cdot X\cdot Q_2\cdot \lL{ [X,Y]}\rR \bigr)\;.
$$
\end{lemma}

\Pf .
We have
$$
\Tr\bigl(Q_1\cdot X\cdot Y\cdot Q_2\cdot \lL{ [X,Y]}\rR \bigr)=
\Tr\bigl(Q_1\cdot Y\cdot X\cdot Q_2\cdot \lL{ [X,Y]}\rR \bigr)
+\Tr\bigl(Q_1\cdot \lL{ [X,Y]}\rR\cdot  Q_2\cdot \lL{ [X,Y]}\rR \bigr)+\;...\;,
$$
where "$...$" denotes a $\g$-invariant operator
of level $\le d-1$.
By the induction assumption, such an operator is zero in
$\dd(\g)/\JJ$. By Lemma \ref{twocom},
in $\dd(\g)/\JJ$ one has the equation:\quad
$\;\dis
\Tr\bigl(Q_1\cdot \lL{ [X,Y]}\rR\cdot  Q_2\cdot \lL{ [X,Y]}\rR \bigr)=0.
\,$ The Lemma follows.
\sq

\begin{lemma}\label{last}
For any $a,b,c\geq 0$, in $\dd(\g)$, one has the
identity:
$$
\Tr\bigl(\lL{ X^aY^bX^c}\rR \lL{ [X,Y]}\rR \bigr)-
\Tr\bigl(\lL{ X^aY^bX^c[X,Y]}\rR \bigr)=
\Tr(X^c)\cdot\Tr(X^aY^b)-
\Tr(X^a)\cdot\Tr(X^cY^b).
$$
\end{lemma}

\Pf . First observe that:
$\;\dis
\partial_{ji}\cdot\lL{ [X,Y]}\rR -\lL{ \partial_{ji}\cdot[X,Y]}\rR =
[E_{ji},Y].
\;$
Therefore, we find
\begin{align*}
&\Tr\bigl(\lL{ X^aY^bX^c}\rR \lL{ [X,Y]}\rR \bigr)-
\Tr\bigl(\lL{ X^aY^bX^c[X,Y]}\rR \bigr)=
\\
&\sum_{b_1+b_2=b-1}\enspace\sum_{pqrs}\;
\lL{ (X^aY^{b_1})_{pq}\cdot (Y^{b_2}X^c)_{rs}\cdot [E_{rq},Y]_{sp}}\rR =
\\
&\sum_{b_1+b_2=b-1}\enspace\sum_{pqrs}\;
\lL{ (X^aY^{b_1})_{pq}\cdot (Y^{b_2}X^c)_{rs}\cdot (\delta_{rs}\cdot y_{qp}-
\delta_{qp}\cdot y_{sr})}\rR =
\\
&\sum_{b_1+b_2=b-1}\;
\lL{ \Tr(X^aY^{b_1+1})\cdot\Tr(Y^{b_2}X^c)-
\Tr(X^aY^{b_1})\cdot\Tr(Y^{b_2+1}X^c)}\rR\;.\quad\square
\end{align*}
\medskip

\noindent
{\bf {Proof of Theorem \ref{maintech}.}}\quad
We must show
that any element $D\in I_{\ad}$
is zero in $\dd(\g)/\JJ$.

We will prove that $D=0$ in
$\dd(\g)/\JJ$ by induction in the level of $D$.
The base of induction $(d=0)$ is obvious.
Assume that the level of $D$ is $d$ and for levels $\le d-1$
it has been proved that elements of $I_{\ad}$ vanish in
$\dd(\g)/\JJ$. By Lemma
\ref{iprime}, we may assume that
$D=\Tr\bigl(Q(X,Y)\cdot \lL{ [X,Y]}\rR \bigr)$, where
$Q$ has degree $d$.

Lemma \ref{twocom2} shows that it is enough to check the
statement of the theorem for the
elements: $D_m'=\Tr\bigl(Q_m)\cdot \lL{ [X,Y]}\rR
 \bigr)$,
$m=0,...,d$, where $Q_m=\sum_{k=0}^m X^{m-k}Y^{d-m}X^k$.
Let $D_m=\Tr\bigl(\sum_{k=0}^m \lL{ X^{m-k}Y^{d-m}X^k}\rR\cdot 
\lL{ [X,Y]}\rR \bigr)$.
Clearly, $D_m=D_m'$,
up to terms of lower level. Hence, it suffices
to show that $D_m=0$. We compute $D_m$ by means of  Lemma \ref{last}.
Observe that after summation over $\,0\leq k\leq m\,$
(under the trace sign)
the terms coming from  the RHS of the equation of Lemma \ref{last}
cancel out. Thus, in $\dd(\g)$,
we find:
\begin{align*}
D_m=&\sum_{k=0}^m\Tr\bigl(\lL{ X^{m-k}\cdot Y^{d-m}\cdot 
X^k}\rR\cdot \lL{ [X,Y]}\rR \bigr)=
\sum_{k=0}^m\Tr\bigl(\lL{ X^{m-k}\cdot Y^{d-m}\cdot X^k\cdot [X,Y]}\rR \bigr)=
\\
&\sum_{k=0}^m\Tr\bigl(\lL{ Y^{d-m}\cdot X^k\cdot [X,Y]\cdot X^{m-k}}\rR \bigr)=
\Tr\bigl(\lL{[Y^{d-m}\cdot X^m\,,\,Y]}\rR \bigr)=0\;.\quad\square
\end{align*}

\section{ Proof of  Theorem \ref{inject}}
\setcounter{equation}{0}

In this section: $\g={\frak {gl}}_n$, and we put
$I_k := \bigl(\dd(\g)\cdot\ad(\ann_k)\bigr)^\g$, and
$I_{\text{ad}}:=\bigl(\dd(\g)\cdot\ad\g\bigr)^\g$.
\medskip

\noindent
{\bf {Proof of Theorem \ref{inject}(i) for ${\mathbf{k=0}}$.}}\quad
By Theorem \ref{maintech} and Theorem \ref{WLS},
we have $I_0=I$.
Moreover, it is known that $W$-invariant polynomials on $\h$
and $W$-invariant differential operators with constant
coefficients generate $\dd(\h)^W$, see \cite{Wa},
which implies that $\bcal_0=\dd(\h)^W$.
Thus, the Harish-Chandra homomorphism $\Phi$ gives rise to
an isomorphism $\,\Phi_0:\dd(\g)^\g/I_0\iso \bcal_0$.\medskip

\noindent
{\bf {Proof of part (i): general case.}}\quad
For  any smooth algebraic variety $X$
with a $\C^*$-action, we write $\dd(X)(l,m)\subset \dd(X)$
for the space of  differential operators  of order  $\le m$
(with respect to the standard filtration on $\dd(X)$)
and  of homogeneity degree $l$
with respect to the $\C^*$-action.
Given any subspace $A\subset \dd(X)$, we put
$\,A(l,m):= A\cap \dd(X)(l,m)$.
We will apply these notations for $X=\g$
and for $X=\hreg$, with the natural 
 $\C^*$-action by dilations in either case.

Observe that the map $\Phi_k$ is compatible both with the 
filtrations and the gradings, furthermore,
 for any $m\geq 0$ and $l\in \Z$,
the space $\dd(\g)(l,m)$ is finite dimensional.
We claim  that, for generic  $k$ (all except countably many),
one has:
\begin{equation}\label{nondegeneracy}
\bigl(\Image^{\,}\Phi_k\bigr)(l,m)\;=\;\Phi_k\bigl((\dd(\g)^\g/I_k)(l,m)\bigr)
\quad,\quad\forall m\geq0\,,\,l\in\Z\,.
\end{equation}
To prove this, observe that by Theorem \ref{Joseph} we know that
the map ${\grd}(\Phi_0)$ is surjective.
This implies that (\ref{nondegeneracy}) holds for $k=0$.
Further, by semi-continuity, for generic  $k$, one has:
$\dim I_k(l,m)\ge \dim I_0(l,m)$.
It follows that if an element $D\in \dd(\g)^\g$ has
degree exactly $d$ modulo $I_0$, it has degree
$\le d$ modulo $I_k$, for generic $k$. However,
$\Phi_0(D)$ has degree exactly $d$, which implies that
$\Phi_k(D)$ has degree $\ge d$, for generic $k$.
This implies that (\ref{nondegeneracy})
must hold for generic $k$.

Recall that, as has been mentioned in the Remark 
following Theorem  \ref{inject},
it is clear that: $\bcal_k\subset \Image^{\,}\Phi_k$, for any $k\in\C$.
It follows from equation (\ref{nondegeneracy})
that, for generic $k$ and
any $m\geq 0$ and $l\in \Z$, we have:
$\bcal_k(l,m) \subset \bigl(\Image^{\,}\Phi_k\bigr)(l,m).\,$
In particular, for generic $k$, we deduce 
\begin{align}\label{inequality}
\dim \bcal_k(l,m) &\leq \dim\bigl(\Image^{\,}\Phi_k\bigr)(l,m)
= \dim\Phi_k\Bigl(\bigl(\dd(\g)^\g/I_k\bigr)(l,m)\Bigr)\\
& \leq\dim\bigl(\dd(\g)^\g/I_k\bigr)(l,m)
=\dim\dd(\g)^\g(l,m) -\dim I_k(l,m)\,.\nonumber
\end{align}

Further, a semi-continuity argument shows that,  for generic $k$,
one has: $\dim \bcal_k(l,m)\ge \dim
\bcal_0(l,m)$,
 and $\dim I_k(l,m)\ge \dim I_0(l,m)$.
Therefore, since the theorem holds for $k=0$, 
the inequalities above yield
\begin{align*}
\dim \bcal_k(l,m)\ge \dim
\bcal_0(l,m)&=\dim\dd(\g)^\g(l,m) -\dim I_0(l,m)\\
&\geq \dim\dd(\g)^\g(l,m) -\dim I_k(l,m).
\end{align*}
Comparing this with (\ref{inequality})
shows that,
for generic $k$, all the inequalities must be equalities,
in particular:
$\dim I_k(l,m)=\dim I_0(l,m)$, and
$\dim \bcal_k(l,m)= \dim \bcal_0(l,m).$ 
Thus, for all $m\geq 0\,,\, l\in\Z$, we have:
$\dim \bcal_k(l,m) 
= \dim \bigl(\dd(\g)^\g/I_k\bigr)(l,m),$ and the map $\Phi_k$ is
an isomorphism, for generic $k$.

It remains to show that the image of $\Phi_k$ is $\bcal_k$ not only
for generic but actually for all $k\in \C$.
It follows from the fundamental theorem of invariant theory (see \cite{Wa})
that $\C[\h\oplus \h]^W$ is generated, as a Poisson algebra,
by the two subalgebras $\C[\h]^W$ of
invariant polynomials on the first and the second copies of
$\h$, respectively. Therefore, for any $k$, we have
$\C[\h\oplus \h]^W\subset\grd(\bcal_k)$, where `$\grd$'
is taken with respect to the natural filtration
by the order of differential operators.
Also, as we mentioned,
one has $\bcal_0=\dd(\h)^W$, hence
 $\grd(\bcal_0)=\C[\h\oplus\h]^W$.
Therefore,
{\small
\begin{equation}\label{**}
{
\dim \bcal_k(l,m)= \dim \Bigl(\bigoplus_{i\leq m}\;\grd_i\bigl(\bcal_k(l,m)\bigr)\Bigr)
\ge \dim \Bigl(\bigoplus_{i\leq m}\;\grd_i\bigl(\bcal_0(l,m)\bigr)\Bigr)=
\dim \bcal_0(l,m).}
\end{equation}}

On the other hand, since $\bcal_k$ is described by generators,
the dimension $\dim \bcal_k(l,m)$ cannot go up
at special points $k$ (it can only go down). Therefore, (\ref{**}) yields:
$\dim \bcal_k(l,m)=\dim \bcal_0(l,m),$
for all $k$.
This together with the fact that $\Image^{\,} \Phi_k=\bcal_k$ for generic $k$,
implies that
$\Image^{\,} \Phi_k=\bcal_k$, for all $k$.\medskip

\noindent
{\bf {Proof of part (ii).}}\quad
Write $\bcal$ for the $\C[k]$-algebra
whose specialization at any particular value of the
parameter $k$ is the algebra $\bcal_k$.
In more detail, treat $k$ as a variable, and view
the Calogero-Moser operator (\ref{H})
as an element of $\dd(\hreg)^W\bigotimes \C[k]$.
Let $\cc$ denote the centraliser of this operator in
 ${\dd(\hreg)_-^W}\bigotimes \C[k]$, and
$\bcal\subset \dd(\hreg)^W\bigotimes \C[k]$
the algebra
 generated by $\cc$ and by
the polynomial subalgebra $\C[\h]^W\bigotimes \C[k]
\subset \dd(\h)^W\bigotimes \C[k]$.

By part (i),
we may think of the family of deformed Harish-Chandra homomorphisms
$\Phi_k$ as a single $\C[k]$-linear  homomorphism:
$\dd(\g)^\g\bigotimes \C[k]\to \bcal$.
Let ${\mathbf{Ker}}$ be the kernel of this homomorphism.
Similarly, let ${\mathbf{I}}$ be the two-sided ideal
of the algebra $\dd(\g)^\g\bigotimes \C[k]$
that  specialises to  $I_k$ at any particular value of the
parameter $k$. We know that ${\mathbf{I}}\sset {\mathbf{Ker}}$.
The claim of  part (ii) amounts to showing
that ${\mathbf{Ker}}/{\mathbf{I}}$ is supported,
as a $\C[k]$-module, at finitely many points.

To prove this,  equip $\dd(\g)^\g\bigotimes \C[k]$
 with the standard filtration (place $\C[k]$ in degee $0$).
Then, we have:
 $\grd\bigl(\dd(\g)^\g\bigotimes \C[k]\bigr)=
\C[\g\oplus\g]^\g\bigotimes \C[k]$.
The RHS here is a  finitely generated algebra,
since $\g$ is a reductive Lie algebra.
Hence $\dd(\g)^\g\bigotimes \C[k]$ and
$\bigl(\dd(\g)^\g\bigotimes \C[k]\bigr)/{\mathbf{I}}$
are  finitely generated algebras as well.
Therefore, the torsion ideal
in $\bigl(\dd(\g)^\g\bigotimes \C[k]\bigr)/{\mathbf{I}}$
 is finitely generated by the Hilbert basis theorem.
Thus, this ideal is supported at finitely many points of
$\Spec \C[k]$. But, by part (i), for all $k$ except possibly
a countable subset, we have: $I_k=\Ker(\Phi_k)$.
It follows that ${\mathbf{Ker}}/{\mathbf{I}}$ is
contained in  the torsion ideal
in $\bigl(\dd(\g)^\g\bigotimes \C[k]\bigr)/{\mathbf{I}}$,
and we are done.\sq\medskip

\noindent
{\bf Remarks} \vi
In the course of the proof of the theorem we have established the following:

{\it The algebra $\bcal$ is
a free $\C[k]$-module, and $\grd(\bcal)=\C[\h\oplus\h]^W\bigotimes\C[k]$.}

\noindent
This  result  has been already proved differently in
Part 1 of the paper.

\vii The conjecture that the ideal $J$,
(see the notation of
Remark at the end of \S6), coincides with its
radical $\sqrt{J}$ at the level of $\g$-invariants would easily imply
that $\Ker(\Phi_k)=I_k$, for all $k$. Indeed, $\grd(I_k)$
contains $(\sqrt{J})^\g$; so the conjecture would
imply that $J^\g\subset\grd(I_k)$ which, by Theorem
\ref{LJ}, equals  $\grd(I_0)$, thus $\dim I_k(l,m)\ge \dim
I_0(l,m)$, and hence $\Ker(\Phi_k)=I_k$.$\quad\lozenge$

\section{Calogero-Moser space for wreath-products.}
\setcounter{equation}{0}

Let $(L,\om_{_L})$ be a 2-dimensional symplectic vector
space, and $\Gamma\subset Sp(L)$ a finite subgroup.
Consider the wreath product $\ga  =S_n\ltimes \Gamma^n$, acting 
on $V:=L^{\oplus n}=\C^{n}\otimes L,$
a symplectic vector space.
Given $\gamma\in \Gamma$, we write $\gamma_i\in \ga$
for the element $\gamma$ regarded as an element of the
$i$-th factor $\Gamma$. Let
 $s_{ij}\in S_n$ denote  the transposition: $i\leftrightarrow j,$
in $S_n$. Then
the group $\ga$
is generated by the following symplectic reflections:
the transpositions $s_{ij}$, and
the elements $\gamma_i\in \ga$.  The conjugacy classes
of symplectic reflections in $\ga\subset Sp(V)$ are known to be
of the following two types:
\begin{equation}\label{conj_class}
\begin{array}{l}
{\sc{(S)}}\quad\mbox{{\it The set: $\;\{s_{ij}\cdot\gamma_i\cdot
\gamma_j^{-1}\}_{i,j\in[1,n],\,\gamma\in \Gamma},$
forms a single $\ga$-conjugacy class.}}
\\
{(\G)}\quad \mbox{{\it Elements: $\;\big\lbrace\gamma_i\,,\, i\in[1,n],\,
\gamma\in \Gamma\smallsetminus \{1\}\big\rbrace,$ form one 
$\ga$-conjugacy}}\\
\hp\quad\quad
\mbox{{\it class, for any given conjugacy class of $\gamma$. }}
\end{array}
\end{equation}

An important role below
will be played  by the conjugacy class
$\O \subset \sln$, formed by all $n\times n$-matrices of the form:
$P-\id$, where $P$ is a
 semisimple rank 1 matrix such that $\Tr(P) = \Tr(\id)=n$.
Thus, $\O$ is a closed $\GL_n$-conjugacy class in $\sln$,
containing in particular the following matrix:
\begin{equation}\label{s_matrix}
\bp=\begin{pmatrix}
0 & 1& 1&\ldots & 1& 1\\
1 & 0& 1&\ldots & 1& 1\\
1 & 1& 0&\ldots & 1& 1\\
1 & 1& 1&\ldots & 0& 1\\
1 & 1& 1&\ldots & 1& 0
\end{pmatrix}\,,
\end{equation}
since $P=\bp+\id$ is a rank one matrix.

Write ${\bf{e}}_{_{\G}} \in \End_{_\C}(\C\Gamma)$
for the  projector $\,{\bf{e}}_{_{\G}}:
u\mapsto \bigl(\frac{1}{|\G|}\sum_{\gamma\in\G}\,\gamma\bigr)
\cdot u,\,$  on the trivial
representation of $\G$.
Further, given a class function: $\G\smallsetminus\{1\}\,
\to\C\,,\,\gamma \mapsto c'_\gamma,$ 
let $c'=\sum_\gamma c'_\gamma\cdot\gamma$ denote
the corresponding central element 
of $\C\Gamma$,
 which is traceless in the regular representation. 
For each pair $c=(k,c'),$ where $k\in \C$ and
$c'$ is a class function as above,
we introduce
the following set: 
\begin{equation}\label{M}
\left.\!\!\begin{array}{ll}
\!\!\!\!M_{_{\G,n,c}}\;=\;
\Big\{\nabla\in
&\Hom_{_\G}\bigl(L\,,\,\End_{_\C}(\C^n\otimes
\C\Gamma)\bigr)
\;\;\Big|\;\;\nabla:\, x\mapsto \nabla_x,
\;{\mbox{such that $\forall x,y \in L$}}
\\&{\mbox{we have:}}\;\;
[\nabla_x,\nabla_y]\,\in\,\om_{_L}(x,y)\cdot
\bigl(k\cdot |\Gamma|\cdot  \O\otimes {\bf{e}}_{_{\G}}+
\id_{_{\C^n}}\otimes c'\bigr)
\end{array}\!\!\right\}
\end{equation}

It is convenient to
 choose and fix a symplectic basis $\,\{x,y\}\,$ in
$L$ and this way  identify $L$ with $\C^2$,
so that the group $Sp(L)$ gets identified with
$\SL_2(\C)$. Setting $\nabla_1:=\nabla_x\,,\,
\nabla_2:=\nabla_y,$ we can
rewrite the formula above more concretely as follows:

\begin{equation}\label{M2}
\!\begin{array}{ll}
\!\!\!\!M_{_{\G,n,c}}\;=\;
\Big\lbrace\nabla_1,\nabla_2\in\End_{_\C}(\C^n\otimes
\C\Gamma)
\;\;\Big|\,&[\nabla_1,\nabla_2]=
k\cdot |\Gamma|\cdot  {\mathbf{o}}\otimes {\bf{e}}_{_{\G}}+
\id_{_{\C^n}}\otimes c'\!,\\
&{\mbox{for some}}\;\;{\mathbf{o}}\in\O\end{array}\!\!\Big\rbrace\,.
\end{equation}

Set $\,\GGG:= \Aut_{_\G}(\C^n\otimes
\C\Gamma)$ and  $\,\PG_{_{\Gamma\!,n}}:=\PAut_{_\G}(\C^n\otimes
\C\Gamma)$. The reductive group $\PG_{_{\Gamma\!,n}}$ acts naturally
on the vector space $\,\Hom_{_\G}\bigl(L\,,\,\End_{_\C}(\C^n\otimes
\C\Gamma)\bigr)\,$ preserving the set $M_{_{\G,n,c}}$.
\begin{definition}
The quotient variety:
$\mm_{_{\G\!,n,c}}=M_{_{\G,n,c}}/\PG_{_{\Gamma\!,n}}$
is called the Calogero-Moser space for $\ga$.
\end{definition}
The same variety can be equivalently defined as follows.
First take the point $\bp\in \O$, see (\ref{s_matrix}),
and let $\PG_{_{\Gamma\!,n}}(\bp) \subset \PG_{_{\Gamma\!,n}}$ denote the isotropy group
of the element: $k\cdot |\Gamma|\cdot  \bp\otimes {\bf{e}}_{_{\G}}+
\id_{_{\C^n}}\otimes c'\,\in \End_{_\C}(\C^n\otimes
\C\Gamma)$. Then, clearly, one has: $\mm_{_{\G\!,n,c}}\simeq
M_{_{\G,n,c}}(\bp)/\PG_{_{\Gamma\!,n}}(\bp),$ where
\begin{equation}\label{M3}
M_{_{\G,n,c}}(\bp)\;=\;
\Big\lbrace\nabla_1,\nabla_2\in\End_{_\C}(\C^n\otimes
\C\Gamma)
\;\;\Big|\,[\nabla_1,\nabla_2]=
k\cdot |\Gamma|\cdot  \bp\otimes {\bf{e}}_{_{\G}}+
\id_{_{\C^n}}\otimes c'\Big\rbrace\,.
\end{equation}

We describe the structure of the groups 
$\PG_{_{\Gamma\!,n}}$ and $\PG_{_{\Gamma\!,n}}(\bp)$ in more detail.
There is an obvious
algebra isomorphism: $\End_{_\Gamma}(\C^n\otimes \C\Gamma)\simeq
\gln\otimes \End_{_\Gamma}\C\Gamma\simeq
\gln\otimes \C\Gamma.\,$ Therefore, decomposing $\C\G$ into irreducible
components,  
identifying $\GL(\C^n\otimes \triv)$ with $\GL_n$,
and using the notation
$\,\GGp\;:=
\!\prod_{E\in\, {\mathsf{Irrep}}(\Gamma)\smallsetminus{\mathsf{triv}}}
\,\GL(\C^n\otimes E)\,,$
we find
\begin{equation}\label{G_G}
\PG_{_{\Gamma\!,n}}\,=\,
\GGG/\C^*\,=\,\Bigl(\prod_{E\in {\mathsf{Irrep}}(\Gamma)}\;
\GL(\C^n\otimes E)\Bigr)\Big/\C^*
\;=\;\Bigl(\GL_n\,\times\,\GGp\Bigr)\Big/\C^*
\,.
\end{equation}
  From this we compute:
\begin{equation}\label{G_G2}
\quad\dim \PG_{_{\Gamma\!,n}}=n^2|\Gamma|-1
\quad\mbox{and}\quad
\dim \PG_{_{\Gamma\!,n}}(\bp)=n^2(|\G|-1)+(n-1)^2\,.
\end{equation}

Further, the decomposition:
$\C^n=\Ker(\bp+\id)\oplus\Image(\bp+\id)$
gives an imbedding $\imath:\GL_{n-1}\into \GL_n\,,\,
A\mapsto A\oplus \id_{_{\Image(\bp+\id)}}.\,$
We observe  that the composite map:
\begin{equation}\label{PG-split}
\GL_{n-1}\,\times\GGp
\quad\stackrel{\imath\times\id}{\into}\quad
\GL_n\,\times\,\GGp
\;\onto \;\Bigl(\GL_n\,\times\,\GGp
\Bigr)\Big/\C^*=\PG_{_{\Gamma\!,n}}
\end{equation}
induces an isomorphism of the group $\GL_{n-1}\times\GGp$
with  the group $\PG_{_{\Gamma\!,n}}(\bp)\subset \PG_{_{\Gamma\!,n}}.$
We  often invert this  isomorphism to obtain
an imbedding: 
$\PG_{_{\Gamma\!,n}}(\bp)\iso\GL_{n-1}\times\GGp
\hookrightarrow~\GGG$.\smallskip

\noindent
{\bf Hamiltonian interpretation}\quad
The   variety ${\mathcal M}_{_{\G,n,c}}$ can be usefully
interpreted as follows.
The symplectic form on the 2-dimensional space
$L$ and the {\it symmetric}
trace form $\Tr: a,b \mapsto \Tr(ab)$ on
$\End_{_\C}(\C^n\otimes
\C\Gamma),$ 
make $\,L^*\bigotimes\End_{_\C}(\C^n\otimes
\C\Gamma),$ a symplectic vector space with
symplectic 2-form $\om_{_L}\otimes \Tr.$ 
Clearly, 
$\Hom_{_\G}\bigl(L\,,\,\End_{_\C}(\C^n\otimes
\C\Gamma)\bigr)\,$ is nothing but
the subspace of $\G$-invariants in
$\,L^*\bigotimes\End_{_\C}(\C^n\otimes              
\C\Gamma),$ which is again a
symplectic vector space.
The natural $\PG_{_{\Gamma\!,n}}$-action on
$\Hom_{_\G}\bigl(L\,,\,\End_{_\C}(\C^n\otimes
\C\Gamma)\bigr)$ is  Hamiltonian,
and gives rise to a moment map                               
$\mu: \Hom_{_\G}\bigl(L\,,\,\End_{_\C}(\C^n\otimes
\C\Gamma)\bigr) \too (\Lie\PG_{_{\Gamma\!,n}})^*.$

It is convenient at this point to introduce
three Lie algebras:
\begin{align*}
&\ggn:=\Lie\GGG=\End_{_\G}(\C^n\otimes
\C\Gamma)=
\underset{E\in {\mathsf{Irrep}}(\Gamma)}{\bplus}
\;\gl(\C^n\otimes E)
\enspace,\enspace
\pgg:=\Lie \PG_{_{\Gamma\!,n}}=\ggn/\C,\\
&\sgg=\big\{a= (a_{_E})_{_{E\in {\mathsf{Irrep}}(\Gamma)}}\;
\in \bplus_{\!\!_{E\in {\mathsf{Irrep}}(\Gamma)}}\;\,
\gl(\C^n\otimes E)=\ggn
\quad\Big|\quad \sum\nolimits_E\;\,\Tr(a_{_E})=0\big\}\,.
\end{align*}
We may identify $\pgg^*$ with $\sgg\subset \ggn$
via the trace form.
Then the  moment map above becomes a map
$\mu: \Hom_{_\G}\bigl(L\,,\,\End_{_\C}(\C^n\otimes
\C\Gamma)\bigr) \to \sgg,$
 given
by the formula:
$\nabla\mapsto [\nabla_1,\nabla_2].$

For $c=(k,c')$, we set
$$\OO_{_{\G,n,c}}=\{u\in\gln\otimes \End_{_\C}(\C\Gamma)
\;\;\;\big|\;\;\;
u=k\cdot |\Gamma|\cdot  \bo\otimes {\bf{e}}_{_{\G}}+
\id_{_{\C^n}}\otimes c'\in\gln\otimes \End_{_\C}(\C\Gamma)\;,\;
\bo\;\in\; \O\}$$
which is a coadjoint orbit in $\pgg^*\simeq \sgg.$
By definition we
have: $M_{_{\G,n,c}}=\mu^{-1}(\OO_{_{\G,n,c}})$.
Thus, $\mm_{_{\G\!,n,c}}=\mu^{-1}(\OO_{_{\G,n,c}})/\PG_{_{\Gamma\!,n}}$,
and we get 
\begin{proposition}\label{ham_red}
The variety $\mm_{_{\G\!,n,c}}$
is the Hamiltonian reduction  over
$\OO_{_{\G,n,c}}$ of the symplectic vector space
$\Hom_{_\G}\bigl(L\,,\,\End_{_\C}
(\C^n\otimes
\C\Gamma)\bigr)$.\quad\qed
\end{proposition}

The Proposition below
and its proof are completely parallel to [Wi, Prop. 1.7],
a result proved by Wilson in the special case: $\G=\{1\}$.
\begin{proposition}\label{G_Wilson} 
For generic $c\in\CC(\ga)$,
 the $\PG_{_{\Gamma\!,n}}$-action on
 the variety $M_{\Gamma,n,c}$
is free. Furthermore,
$\mm_{\Gamma,n,c}$ is a smooth  symplectic affine algebraic variety
of dimension  $2n$. 
\end{proposition}

\Pf . We claim first that if $E\subset \C^n\otimes \C\Gamma$
is a proper $\G$-stable subspace, which is in addition invariant
under $\nabla_i\,,\,i=1,2,$ then $E=0$.
If $E\neq 0$, then we consider the operator
$[\nabla_1,\nabla_2]\big|_E$. The  trace of this operator is clearly
zero. We now use the following observation that will be
also exploited many times elsewhere, esp. in Appendix E.
Since $\G$ is finite, the character table $\|\rho_i(u_j)\|,$
where $\rho_i$ is the $i$-th irreducible character of $\G$
and $u_j$ is the $j$-th conjugacy class in 
$\G$, involves a finite number of algebraic numbers.
Therefore, vanishing of a linear combination
of the $\rho_i(u_j)$'s involving some {\it generic} parameters
implies that the coefficient in front of each of the  $\rho_i(u_j)$'s
vanishes identically. Thus, for generic $c=(k,c')$,
vanishing of the trace of $[\nabla_1,\nabla_2]\big|_E$
implies that the terms in the trace involving $k$, resp. $c'$,
vanish separately.

Now, looking at
$c'$-terms we see that $E$ is a multiple of $\C\Gamma$. 
Further, looking at
the $k$-term, we see exactly like in [Wi] that $E^\G$, the 
$\Gamma$-fixed part
of $E$, must be $n$-dimensional. It follows that
 $\dim E=n\cdot|\G|$, hence $E= \C^n\otimes \C\Gamma$.
Schur lemma implies the freeness of the $\PG_{_{\Gamma\!,n}}$-action. 

We now show that the  variety 
$M_{\Gamma,n,c}$ is smooth. To do this, it is
sufficient to show that $k|\Gamma| s\otimes {\bf{e}}_{_{\G}}
+1\otimes c'$ is a
regular value of the moment map $\mu: \nabla\mapsto
[\nabla_1,\nabla_2]$,
i.e., that the differential $d\mu$ is surjective. 
But the {\it surjectivity} of the differential of a moment map
for a group action on a symplectic manifold 
is well-known to be equivalent to the local freeness 
of the group action, established above.

It remains to calculate the dimension of 
$\mm_{\Gamma,n,c}$. First we calculate the dimension of 
$M_{\Gamma,n,c}(\bp)$. Since $k\cdot|\Gamma|\cdot \bp
\otimes {\bf{e}}_{_{\G}}+1\otimes c'$
is a regular value of the moment map, 
 we get
$$\dim M_{\Gamma,n,c}(\bp) = \dim\mu^{-1}(\bp)
= \dim \Hom_{_\G}\bigl(L\,,\,\End_{_\C}
(\C^n\otimes
\C\Gamma)\bigr) -\dim(\pgg)\,.
$$
The space
$\Hom_{_\G}\bigl(L\,,\,\End_{_\C}
(\C^n\otimes
\C\Gamma)\bigr)$ clearly has dimension $2n^2|\Gamma|$,
and, moreover,
 we have shown above that the $\PG_{_{\Gamma\!,n}}(\bp)$-action on
$M_{\Gamma,n,c}(\bp)$ is free.
Therefore, using  formulas (\ref{G_G}) and
(\ref{G_G2}),  we find
\begin{align*}
\dim \mm_{\Gamma,n,c}&=
\dim M_{\Gamma,n,c}(\bp) -\dim \PG_{_{\Gamma\!,n}}(\bp)
\;=\;(2n^2|\Gamma|- \dim \PG_{_{\Gamma\!,n}})-\dim \PG_{_{\Gamma\!,n}}(\bp)\\
&=\bigl(2n^2|\Gamma|- n^2|\Gamma|+1\bigr)- 
\bigl(n^2(|\Gamma|-1)+(n-1)^2\bigr)\;=\;2n.
\quad\square
\end{align*}

Next, recall the quiver variety: $
\bbm^c_{_\G}(\BV)/G_{_\G}(\BV)=\bgm_{_{\G\!,n,c}},$ defined above 
Lemma \ref{Wilson2}. 
A straighforward rewriting using
McKay correspondence yields the following isomorphisms:
\begin{equation}
\bbm^c_{_\G}(\BV)\simeq M_{\Gamma,n,c}
\quad,\quad
G_{_\G}(\BV)\simeq \PG_{_{\Gamma\!,n}}
\quad,\quad
\bgm_{_{\G\!,n,c}}
\,\simeq\,\mm_{_{\G\!,n,c}}\;\;,\enspace\forall c\in\ZZ\G\,,\, n\geq 1.
\end{equation}
\begin{corollary}\label{connected}
The variety $\mm_{_{\G\!,n,c}}$ is connected.
\end{corollary}

\Pf . This is a general property of quiver varieties proved by
Crawley-Boevey [CB] (proof of
a similar result  in [Na2] contains a gap).
\quad\qed\medskip

\noindent
{\bf The algebra $\hh_\ka(\ga)$.}\quad
According to (\ref{conj_class}), the space
$\CC(\ga)$ of class functions on
the set of symplectic reflections in $\ga$ can be identified with
$\ZZ\G$ via the map assigning to
a class function $f\in \CC(\ga)$ the
element
$k\cdot1+\sum_{c'\in \G\smallsetminus\{1\}}\;c'_\gamma\cdot\gamma\in
\ZZ\G$,
where $k\in \C$ is the value of $f$ on the
conjugacy class of type ({\sc{S}}), see (\ref{conj_class}),
and $c'_\gamma$ is
the value of $f$ on the corresponding conjugacy class of type $(\G)$.
Given $k$ and $c'=\sum_{c'\in
\G\smallsetminus\{1\}}\;c'_\gamma\cdot\gamma$
as above, we will write $c=(k,c')\in\CC(\ga)$
instead of the corresponding class function~$f$.

Let $\hh_\ka(\ga)
\,,\, \ka=(t,k,c'),$ denote the
symplectic reflection algebra attached to
$\ka=(t,c)\in \C\times\CC(\ga)$. 
The algebra $\hh_{t,c}(\ga)$ is generated by
elements of the vector space $\C^n\otimes L=L^{\oplus n},$
and the group $\ga$. It will be convenient to use
the following notation:
given $x\in L$, for any $1\leq k\leq n$, we put:
$x_k= (0,\ldots,0,x,0,\ldots,0)\in L^{\oplus n}$
($x$ placed on the $k$-th spot). Thus, $\C\ga$ and
the elements $\,\{x_1,\ldots,x_n\}_{_{x\in L}}\,$
generate $\hh_{t,c}(\ga)$.

Assume now that $t=0$. By Theorem \ref{proj_intro} and
Theorem \ref{wreath00} (Appendix E), we may
(and will) view  $\Spec\ZZ_{0,c}(\ga)$
as the moduli space of simple
$\hh_{0,c}(\ga)$-modules. \medskip

\noindent{\bf{The map: $\text{{\bf{Irreps}}}\bigl(\hh_{0,c}(\ga)\bigr) 
\;{\boldsymbol{\longrightarrow}}\; 
{\boldsymbol{\mathcal M}}_{_{\G,n,c}}$}}\,\,\quad
Let $S_{n-1}\subset S_n$
be the subgroup of permutations of the
set $\,\{2,\ldots,n\},\,$ acting trivially on the label $\,1\in
\{1,\ldots,n\}.\,$ We have an imbedding:
 $\gb=S_{n-1}\ltimes\G^{n-1}\into 
S_n\ltimes\G^{n}=\ga$.
The commutation relations in $\hh_{0,c}(\ga)$, cf. (\ref{rel}),
imply that, for any $x\in L$,
the element $x_1=(x,0,\ldots,0)\in L^{\oplus n}\subset\hh_{t,c}(\ga)
$ commutes with  $\gb\subset
\hh_{t,c}(\ga)$.

Let $E$ be a simple 
$\hh_{0,c}(\ga)$-module. Write $E^{^{\gb}}\subset E$ for the
subspace of $\gb$-fixed vectors in $E$. The elements
 $\,\{x_1\}_{x\in L}\,
$ commute with  $\gb$, hence
preserve the subspace $E^{^{\gb}}$, hence, give rise
to a $\G$-equivariant map: $L\to\End(E^{^{\gb}})\,,\,
x\mapsto x_1\big|_{_{E^{^{\gb}}}}.$

Theorem \ref{struc} says  that
$E|_{_{\C\ga}}\simeq \C\ga$,
 as a $\ga$-module. It follows that
$\dim E^{^{\gb}}$
$=n\cdot|\G|.$ Moreover,
the space $E^{^{\gb}}=(\C\ga)^{^{\gb}}$ may be identified with
$\C^n\otimes\C\G$, as a $\G$-module, as follows.
View $\C{\ga}$ as the vector space of $\C$-valued functions on 
$\ga$. Then $(\C{\ga})^{\gb}$ gets identified with
$\C[\gb\backslash{\ga}]$,
the vector space of $\C$-valued functions on
$\gb\backslash{\ga}$.
By definition,  we have a $\ga$-module isomorphism:
$\C[\gb\backslash{\ga}] ={\mathtt{Ind}}_{\gb}^{\ga}{\mathbf{1}}.\,$
On the other hand, the space $\C^n\otimes\C\G=
\C\G^{\oplus n}$ has  a natural $\ga$-module 
structure, and there is a canonical $\ga$-module isomorphism:
$\,\C^n\otimes\C\G\simeq {\mathtt{Ind}}_{\gb}^{\ga}{\mathbf{1}}.\,$
Thus, we get: $E^{^{\gb}}=(\C\ga)^{^{\gb}}=
\C[\gb\backslash{\ga}] ={\mathtt{Ind}}_{\gb}^{\ga}{\mathbf{1}}=
\C^n\otimes\C\G.\,$
Explicitly, the isomorphism between the RHS and LHS is given by the
formulas:
\begin{equation}\label{digamma}
\C^n\otimes \C\Gamma\iso (\C\ga)^{^{\gb}} \,\simeq\,E^{^{\gb}}
\quad,\quad
e_i\otimes \gamma  \mapsto 
\gamma _1\cdot s_{1i}\in (\C\ga)^{^{\gb}},\quad \gamma \in \C\Gamma,
\end{equation}
where
 $\lbrace{e_i\rbrace}$
denotes the standard basis of $\C^n$. 

Next, recall the matrix
$\bp$, see  (\ref{s_matrix}).
\begin{lemma}\label{irreps1} The endomorphim
$[x_1,y_1]\Big|_{E^{^{\gb}}}$ corresponds, under the
bijection (\ref{digamma}), to an endomorphism of $\C^n\otimes \C\Gamma$
given by the formula:
$\,[x_1,y_1]=k\cdot |\Gamma|\cdot  \bp\otimes {\bf{e}}_{_{\G}}+
\id_{_{\C^n}}\otimes c'.$ 
\end{lemma}

\Pf .
Let $x,y\in L$ and 
$x_1=(x,0,...,0)$, $y_1=(y,0,...,0)\in L^{\oplus n}$.
 Let $E$ be a representation of 
$\hh_{0,c}(\ga)$ and $E_1$ be the space of coinvariants under 
the subgroup $\gb=S_{n-1}\ltimes \Gamma^{n-1}$. Then 
$$
[x_1,y_1]\Big|_{E_1}=\sum\nolimits_{\gamma\in\G\,,\,j>1}\;
k\cdot \omega_{s_{1j}\cdot \gamma_1\cdot \gamma_j^{-1}}(x_1,y_1)\cdot 
s_{1j}\cdot \gamma_1\;+\;\omega_L(x,y)\cdot \sum\nolimits_{\gamma\in
\G\smallsetminus\{1\}}\;
c_\gamma'\cdot \gamma_1.
$$
It is easy to calculate that, independently of $\gamma$, one has:
$\omega_{s_{1j}\gamma_1\gamma_j^{-1}}(x_1,y_1)=\omega_L(x,y)$.

Now, a key point, verified by a direct calculation, is that
the action of the element   $\sum_{j>1}\,s_{1j}\,$ on
$\C[\gb\backslash{\ga}]$ goes, under the bijection
$\C[\gb\backslash{\ga}]\iso\C^n\otimes \C\Gamma,$
to the operator with matrix $\bp\otimes\id_{_{\C\G}}$
(in the  standard basis
$\,\{e_i\otimes \gamma\}_{i=1,\ldots,n,\,\gamma\in\G}\,$
of the vector space $\C^n\otimes \C\Gamma$).
The Lemma now follows since
 $\sum_{\gamma\in \Gamma}\gamma=|\Gamma|\cdot {\bf{e}}_{_{\G}}$.
 \quad\qed
\medskip

By definition of the variety 
$\mm_{_{\G\!,n,c}}=M_{_{\G\!,n,c}}/\PG_{_{\G\!,n}}$, 
we have: $\C[\mm_{_{\G\!,n,c}}]=
\C\big[M_{_{\G\!,n,c}}\big]^{\PG_{_{\G\!,n}}}.$
Here, $M_{_{\G\!,n,c}}$ is an affine subvariety in
the vector space $\dis\Hom_{_\G}\bigl(L\,,\,\End_{_\C}(\C^n\otimes
\C\Gamma)\bigr),$ see (\ref{M}).
Thus, $\C[\mm_{_{\G\!,n,c}}]$ is a quotient of
$\dis\C\big[\Hom_{_\G}\bigl(L\,,\,\End_{_\C}(\C^n\otimes
\C\Gamma)\bigr)\big]^{\PG_{_{\G\!,n}}},$
the algebra of invariant polynomials.
Therefore, the standard  increasing filtration on the 
polynomial algebra (by degree of polynomial) induces
an  increasing filtration on the algebra 
$\C[\mm_{_{\G\!,n,c}}]$.
On the other hand, the algebra $\hh_{0,c}(\ga)$
also comes equipped with the  canonical increasing filtration.

The main result of this section is the following explicit
description of the variety $\Spec\bigl(\ZZ_{0,c}(\ga)\bigr)$,
which is a more precise version of Theorem \ref{wreath_thm}.

\begin{theorem}\label{irreps} Let $c\in\CC(\ga)$ be generic. Then

\vi The morphism: $E\mapsto {}^E\nabla,$
 assigning to $E\in \irreps(\hh_{0,c}(\ga))$
 the $\GGG$-conjugacy class of the
map $\,\,{}^E\nabla: L \to \End_{_\C}(E^{^{\gb}})
\,\,,\,\,x\mapsto {}^E\nabla_x:=x_1\Big|_{E^\gb},\,$
induces an isomorphism
$\phi: \Spec\ZZ_{0,c}(\ga)\iso \mm_{_{\G\!,n,c}}\,
$ of Poisson
algebraic varieties.

\vii
The pull-back morphism $\phi^*: \C[\mm_{_{\G\!,n,c}}]\iso\ZZ_{0,c}(\ga)$
is an isomorphism of filtered algebras.
\end{theorem}

To prove the Theorem we need two lemmas.
\begin{lemma}\label{wang}
Let $A$ be a commutative algebra with unit, and $n\ge 2$ an integer.
Then $\bigl(A^{\otimes n}\bigr)^{S_n},$
the symmetric part of the commutative algebra $A^{\otimes n}$, is generated
by elements of the form:
$$\{a_1+\ldots+a_n\;\;\big|\;\; a\in A\}\;,
\quad\mbox{where}\quad a_k
:= 1^{\otimes (k-1)}\,\btimes\, a \btimes 1^{\otimes
(n-k)}\;,\;\; 1\leq k\leq n\,.$$
\end{lemma}

This lemma is a
generalisation of the classical
result of H. Weyl [We], who considered the case $A=\C[x,y]$.
Weyl's proof applies verbatim in the general case
(in the special case $A=\C[x,y]^\G$, to be exploited below,
 the 
lemma can be found in [Wan, Lemma 1]).

Recall the standard increasing filtration $F_\bullet\hh_{0,c}(\ga)$.
Given $a\in \hh_{0,c}(\ga)$, we write
${\mathtt{filt.\,deg}}(a)=m$ provided
$a\in F_m\hh_{0,c}(\ga)\sminus F_{m-1}\hh_{0,c}(\ga).$
Let $\rr_\chi$ denote the simple $\hh_{0,c}(\ga)$-module
corresponding to a smooth point $\chi\in\Spec \ZZ_{0,c}(\ga)$.

\begin{lemma}\label{ppp2}
 Let $a\in \hh_{0,c}(\ga)$ and $ z\in \ZZ_{0,c}(\ga)$ be elements 
such that, for any smooth point $\chi\in\Spec \ZZ_{0,c},$ one has 
$\bigl(\Tr\big|_{\rr_\chi}(a)\bigr)\cdot\id
=z\big|_{\rr_\chi}.$ Then: ${\mathtt{filt.\,deg}}(z)\leq
{\mathtt{filt.\,deg}}(a).$ 
\end{lemma}

\Pf . Assume the contrary, i.e. let ${\mathtt{filt.\,deg}}(a)=m,$ and
${\mathtt{filt.\,deg}}(z)=m+l,$ where $l>0.$ 
Let $z_0$ in $\grd\ZZ_{0,c}(\ga)\simeq
\ZZ_{0,0}(\ga)$ be the principal symbol of $z.$ 
Let $E$ be a generic irreducible module over $\hh_{0,0}(\ga)$ such that 
$z_0\big|_E$ is not zero. Let $E_\lambda$ be a 1-parameter deformation 
of $E,$ such that $E_\lambda$ is a representation of $\hh_{0,\lambda^2 c}(\ga)$ 
for each $\lambda.$ Let $f_\lambda: \hh_{0,c}(\ga)\iso \hh_{0,\lambda^2 c}(\ga)
\,,\,v\mapsto \lambda^{-1}v\,,\, v\in V\subset \hh_{0,c}(\ga),$
be the obvious isomorphism (identical on $\C\ga$).
Then we have 
$$
\Tr\big|_{E_\lambda}\bigl(f_\lambda(a)\bigr)\cdot\id_{_{E_\lambda}}
\;=\;f_\lambda(z)\big|_{E_\lambda}.
$$
Thus, 
$
\lambda^m\cdot\Tr\big|_{E_\lambda}\bigl(f_\lambda(a)\bigr)
\cdot\id_{_{E_\lambda}}=
\lambda^m\cdot f_\lambda(z)\big|_{E_\lambda}.\,$
But this is a contradiction: the LHS has a finite limit 
as $\lambda$ goes to $0$, since the element
$a$ has degree $m,$ while RHS
goes to infinity since $f_\lambda(z)\big|_{E_\lambda}$
is asymptotic to $\lambda^{-m-l}\cdot z_0\big|_E.\;\,$~\qed
\medskip

\noindent
{\bf Proof of Theorem \ref{irreps}:}\quad
We set $\hh_{0,c}:=\hh_{0,c}(\ga)$ and $\ZZ_{0,c}:=\ZZ_{0,c}(\ga)$.

Let $f\mapsto \hat f$ denote the standard symmetrisation map:
$\Sym L \to TL$, and $a\mapsto a_1+\ldots$
$+a_n$
the map: $TL\to (TL)^{\otimes n}$, 
where $a_k
:= 1^{\otimes (k-1)}\,\otimes\, a \otimes 1^{\otimes
(n-k)}$, as in Lemma \ref{wang}.
We form a  map: $f\mapsto a_f,$ given by the
 following composition:
$$
{
\diagram
\enspace \Sym L\;\rto^<>(.5){f\mapsto \hat f}&
\; (TL)\;\rrto^<>(.5){a\mapsto a_1+\ldots+a_n}
&&\;
(TL)^{\otimes n}\;=\;
T(\C^n\otimes L)\;\rrto^<>(.5){_{\sf projection}}&&
\;\hh_{0,c}(\ga)
\enddiagram}
$$
Restricting this map to
$\G$-equivariants we obtain
a map $\, \sigma: (\Sym L)^\G \to 
\e\hh_{0,c}\e\,,\,f\mapsto \sigma(f)=\e a_f\e$.
We equip  the algebra $(\Sym L)^\G$ with the filtration induced by
the grading.
Further, the algebra $\e\hh_{0,c}\e$ is filtered, and we have
a natural isomorphism: $\grd(\e\hh_{0,c}\e)\simeq
(\Sym V)^\ga=
\Bigl(\bigl((\Sym L)^\Gamma\bigr)^{\otimes n}\Bigr)^{S_n}$. 
It is clear, that the map $\sigma$ is filtration preserving,
and the associated graded map can be identified
with the map: $(\Sym L)^\G \to  
\Bigl(\bigl((\Sym L)^\Gamma\bigr)^{\otimes n}\Bigr)^{S_n}$
given by the formula: $(\grd^{\,}\sigma)(f)=
f_1+...+f_n$. Hence, the image of the map
$\grd^{\,}\sigma$ generates the algebra 
$\Bigl(\bigl((\Sym L)^\Gamma\bigr)^{\otimes n}\Bigr)^{S_n}$,
by  Lemma \ref{wang} applied to $A=(\Sym L)^\Gamma$.
It follows that the  image of the map
$\sigma$ generates the algebra  $\e\hh_{0,c}\e$.

Now choose a basis $\,\{x,y\}\,$ in $L$, and identify
$\Sym L$ with $\C[x,y],$ a polynomial 
algebra.
We write: $f(x,y) \mapsto \hat f(x,y)$ for the corresponding
symmetrisation map. Further,
compose  $\sigma$ with the Satake isomorphism
to get a map $\,\sigma^\sharp: \C[x,y]^\G=
(\Sym L)^\Gamma\too \e\hh_{0,c}\e\iso \ZZ_{0,c}$.

Let  $\e_{_{\mathsf{diag}}}: \C^n \onto
\C_{_{\mathsf{diag}}}\subset\C^n$
denote the projection on the principal diagonal, so that
$\e_{_{\mathsf{diag}}}\otimes\e_{_\G}$ 
is the projection to the 1-dimensional subspace
 $\C_{_{\mathsf{diag}}}\otimes\e\subset
 \C^n\otimes \C\Gamma$,
on which $[x_1,y_1]=
[\nabla_x,\nabla_y]=(n-1)\cdot k\cdot |\Gamma|$, 
and on which $S_n$ acts trivially.
It follows from the construction that
for the map $\phi^*: \C[\mm_{_{\G\!,n,c}}]\to
\ZZ_{0,c}$ we have:
\begin{equation}\label{construction}
\phi^*\left((\e_{_{\mathsf{diag}}}\otimes\e_{_\G})\ccirc
\hat f(\nabla_x,\nabla_y)\ccirc
(\e_{_{\mathsf{diag}}}\otimes\e_{_\G})\right)
\;=\;\mbox{$\frac{1}{n}$}\cdot\sigma^\sharp(f)\quad,\quad
\forall f\in \C[x,y]^\G\,.
\end{equation}
Since the image of the map
$\sigma$ generates the algebra  $\e\hh_{0,c}\e$,
it follows from this formula
 that the map $\phi^*$ is surjective.

We now show that $\phi^*$ is injective. 
Assume  $\Ker(\phi^*)\neq 0$, and
let $Y\subset \mm_{_{\G\!,n,c}}$ be
the zero variety of the ideal $\Ker(\phi^*)\subset\C[\mm_{_{\G\!,n,c}}].$
Thus, $Y$ an affine subscheme in $\mm_{_{\G\!,n,c}}$
and, moreover, $\dim Y<\dim \mm_{_{\G\!,n,c}}=
\dim V,$ since
$\mm_{_{\G\!,n,c}}$ is irreducible, by Corollary \ref{connected}.
But then we have:
$\dim Y< \dim \Spec\ZZ_{0,c} (=\dim V)$, hence the map
$\,\phi^*: \C[Y]= \C[\mm_{_{\G\!,n,c}}]/\Ker(\phi^*)
\too \ZZ_{0,c}$ cannot be surjective,
contrary to what we have proven above.

To prove part (ii) of the Theorem, we note
that the filtration 
on $\C[\mm_{_{\G\!,n,c}}]$ is such that
 $\deg(\nabla_x)=\deg(\nabla_y)=1$.
The filtration on $\ZZ_{0,c}$ is induced from that on $\hh_{0,c}.$
It follows directly 
from the proof of part (i)
that the map $(\phi^*)^{-1}: \ZZ_{0,c}\to 
\C[\mm_{_{\G\!,n,c}}]$ is filtration preserving.

Next, we claim that the map $\phi^*$ is also filtration preserving.
Indeed, any regular function $f$ on $\mm_{_{\G\!,n,c}}$ 
such that ${\mathtt{filt.\,deg}}(f)=d$
is representable as a sum of products of traces of monomials 
of $\nabla_x$ and $\nabla_y$, of total degree at most $d$. 
So $\phi^*(f)$ is an element of the form $\phi^*(f)(E)=\Tr_E(h_f)$, where 
$h_f\in \hh_{0,c}(\ga)$ is an element of  filtration degree $d$. 
Thus, it remains to check that 
if $h\in 
\hh_{0,c}(\ga)$ and $z\in \ZZ_{0,c}=\C[\irreps\hh_{0,c}(\ga)]$
is the element corresponding to the function: $E\to \Tr_E(h)$,
then ${\mathtt{filt.\,deg}}(z)\leq
{\mathtt{filt.\,deg}}(h)$. But this follows from Lemma
\ref{ppp2}.

To complete the proof, it remains 
to show that the map $\phi$ respects  the Poisson brackets.
First of all, the Poisson bracket $\{-,-\}_{_{\mm_{_{\G\!,n,c}}}}$
on the algebra $\C[\mm_{_{\G\!,n,c}}]$ is easily seen to
have filtration degree $\leq -2$, in the sense explained above
Lemma \ref{pois_filt}. Transporting this bracket
to $\ZZ_{0,c}$ via the isomorphism $\phi^*$,
we obtain, since the map $\phi^*$ is filtration compatible,
a  Poisson bracket $\phi^*\{-,-\}_{_{\mm_{_{\G\!,n,c}}}}$ on  $\ZZ_{0,c}$
of  filtration degree $\leq -2$ again. Hence,
Lemma \ref{pois_filt} implies that there exists
a  constant $\lambda\in \C$,
such that $\phi^*\{-,-\}_{\mm_{_{\G\!,n,c}}}=\lambda\cdot\{-,-\}_{_{\ZZ_{0,c}}}$.

We now, consider the line
$\C\cdot c\subset\CC(\ga)$. Hence, for any $r\in\C^*$ we get similarly:
$\phi^*\{-,-\}_{\mm_{_{\G,n,r\cdot c}}}=\lambda(r)\cdot
\{-,-\}_{_{\ZZ_{0,r\cdot c}}}$. Comparing homogeneity degrees
of both sides of this equation we conclude that the function
: $r\mapsto \lambda(r)$ must be homogeneous of degree zero,
hence a constant $\lambda(r)=\lambda$.
To compute this constant explicitly, we let $r\to 0$.
Thus, we are reduced to the case $c=0$, which is easy to 
calculate. 
Namely, on $\gln\oplus\gln=\gln\otimes L$
introduce the symplectic 2-form $\Tr\otimes \om_L$.
Then $\phi^*\om_{_\ZZ}=\frac{1}{|\Gamma|}\cdot(\Tr\otimes \om_L),$
where $\om_{_\ZZ}$ is the 2-form on $\Spec \ZZ_{0,c}$
induced by the (nondegenerate)
Poisson structure on $\ZZ_{0,c}$.
Thus, the theorem is proved.\,~\qed
\medskip

\noindent
{\bf Interpretation via Quiver varieties.}
\quad\,\,
Let $Q=Q(\G)$ be a quiver obtained by choosing an orientation
on the 
affine Dynkin graph (of ADE type) arising from the group 
$\Gamma\subset \SL_2$ via the McKay correspondence.
Thus, the set $\Vx$ of vertices of $Q$ is identified
with the set $\irreps(\G)$ of the isomorphism
classes of simple $\G$-modules,
so that the trivial representation corresponds to 
a special vertex $0$ 
(in a somewhat degenerate case: $\Gamma=1$,
the quiver $Q$ has only one (rather than 2) vertex and one edge-loop
at this vertex). Given $v\in \Vx$, we write $V_v\in \irreps(\G)$
for the corresponding $\G$-module.

Let $n$ be a positive integer,
and $\RQ_n$ the affine variety (a vector space)
 of representations of $Q$
in the vector space: $\bigoplus_{v\in \Vx}\;V_v^{\oplus n}
,\,$ that is
representations having vector space $\C^n\otimes V_v=
V_v^{\oplus n}$ at the vertex $v
\in \Vx$.
The group $\bigl(\prod_{v\in
\Vx}\;\GL(\C^n\otimes
V_v)\bigr)\big/\C^*$, which is isomorphic to
$\PG_{_{\Gamma\!,n}}$, 
acts naturally on  $\RQ_n$. 
Recall that $\Lie\PG_{_{\Gamma\!,n}}=
\pgg$ is the quotient of 
$\ggn=\bigoplus_{v\in\Vx}\;
\gl(\C^n\otimes V_v)$ by the 1-dimensional subalgebra
of scalar matrices, $\pgg=\ggn/\C$,
and $\pgg^*\simeq\sgg$.

Now, fix an element $c'=\sum_{\gamma\in \G\smallsetminus\{1\}}\;
c'_\gamma\cdot\gamma \in\ZZ\G$. For each vertex $v\in\Vx$,
we introduce the complex number
$\,C_v:=\sum_{\gamma\in \G\smallsetminus\{1\}}\;
c'_\gamma\cdot\Tr(\gamma\big|_{V_v})\in\C$.
In particular,  $C_0=\sum_{\gamma\in
\G\smallsetminus\{1\}}\;c'_\gamma.$

The $\PG_{_{\Gamma\!,n}}$-action on $\RQ_n$ induces
a Hamiltonian  $\PG_{_{\Gamma\!,n}}$-action  on 
$T^*(\RQ_n)$,
the cotangent bundle on $\RQ_n$.
Let $\mu_{_{\RQ_n}}:\, T^*(\RQ_n)\too
(\Lie\PG_{_{\Gamma\!,n}})^*=\pgg^*$ denote the corresponding
moment map.  By definition of the constants
$\,\{C_v\}_{v\in\Vx},\,$ the total trace of the
element $\oplus_{v\in\Vx}\; C_v\cdot \id_{_{\C^n\otimes V_v}}\in
\ggn$ equals 
$$
\sum_{v\in\Vx}\; C_v\cdot\dim(\C^n\otimes V_v)\;=\;
n\cdot\!\sum_{v\in\Vx}\sum_{\gamma\in \G\smallsetminus\{1\}}\;
c'_\gamma\cdot\Tr(\gamma\big|_{V_v})\cdot\dim V_v=
n\cdot\!\!\!\sum_{\gamma\in
\G\smallsetminus\{1\}}\;c'_\gamma\cdot\Tr(\gamma\big|_{\C\G}).
$$
Since the trace of any $\gamma\in
\G\smallsetminus\{1\}$ in the regular representation vanishes,
we see that 
 $\oplus_{v\in\Vx}\; C_v\cdot \id_{_{\C^n\otimes V_v}}\in
\sgg.$
 Therefore we may (and will) regard
the conjugacy class 
$$
\O_c:= \bigl(k\cdot|\G|\cdot\O+ C_0\cdot \id_{_{\C^n}}\bigr)
\;\;\bigoplus\;\;\bigl(\bplus_{\!\!_{v\in\Vx\smallsetminus\{0\}}}\;
C_v\cdot \id_{_{\C^n\otimes V_v}}\bigr)\quad\in\quad
\sgg\; \subset\; \bplus_{\!\!_{i\in\Vx}}\;
\gl(\C^n\otimes V_i)
$$
as a coadjoint orbit $\O_c\subset \pgg^*$, via the
trace duality. Then, Proposition \ref{ham_red} yields:
$\mm_{\Gamma,n,c}\simeq
\mu_{_{\RQ_n}}^{-1}(\O_c)/\PG_{_{\Gamma\!,n}}\,.$
Thus, Theorem \ref{irreps} says:
\begin{equation}\label{tra2}
\ZZ_{0,c}(\ga)\;\simeq\;
\C\big[\mu_{_{\RQ_n}}^{-1}(\O_c)/\PG_{_{\Gamma\!,n}}\big]
=\C\big[\mu_{_{\RQ_n}}^{-1}(\O_c)\big]^{\PG_{_{\Gamma\!,n}}}\,.
\end{equation}

\noindent
{\bf 
${\boldsymbol{\G}}$-analogue of the  Harish-Chandra homomorphism.}
\quad
We now discuss a quantization  of (\ref{tra2}).
Fix $c=(k,c')\in \CC(\ga)$, and identify $c'$ with an element 
$c'=\sum_{\gamma\in \G\smallsetminus\{1\}}\;
c'_\gamma\cdot\gamma \in\ZZ\G$. For each $v\in\Vx$,
we have defined
the complex number
$\,C_v:=\sum_{\gamma\in \G\smallsetminus\{1\}}\;
c'_\gamma\cdot\Tr(\gamma\big|_{V_v})$.

  Write $x\mapsto x_v$ for the projection
$\,\pr_v:\;\; \ggn=\mathop{\bplus}_{i\in\Vx}
\;
\gl(\C^n\otimes V_i)\;\;\onto\;\;
\gl(\C^n\otimes V_v)\,.
$
Let $\chi_{_{c'}}: \ggn \to \C$ be a character given
by the formula: $x\mapsto \sum_{v\in\Vx}
\;C_v\cdot \Tr(x_v\big|_{\C^n\otimes V_v}).\,$
The character $\chi_{_{c'}}$ vanishes on the scalar
subalgebra $\C\subset \ggn$, hence descends to a well-defined
homomorphism $\chi_{_{c'}}: \pgg=\ggn/\C\to \C.$

In the case of
the special vertex `0' corresponding to the
trivial representation, we get a Lie algebra  homomorphism 
$\pr_0: \ggn\onto \gln=\gl(\C^n\otimes\triv)$. We compose this homomorphism with the
$\gln$-representation $V_k$, see \S7 and Appendix B,
to obtain a representation $\,\varrho_k:
\ggn \to \End_{_\C}V_k$. Again, since the representation $V_k$ vanishes
 on the center of the Lie algebra $\gln$, the representation
$\varrho_k$ vanishes on the scalar
subalgebra $\C\subset \ggn$, hence descends to a well-defined
homomorphism $\,\varrho_k:
\pgg=\ggn/\C \to \End_{_\C}V_k$.
We 
let $\,\varrho_{_{k,c'}}=\varrho_{_k}\otimes\chi_{_{c'}}:
\;\pgg\to \bigl(\End_{_\C}V_k\bigr)\bigotimes\C$
be the tensor product representation,
and let $\Ann\varrho_{_{k,c'}}\subset \uga$ denote its annihilator
in the enveloping algebra of $\pgg$.

Let $\dd(\RQ_n)$ be the algebra of polynomial
differential operators on 
$\RQ_n$. The action of the group
$\PG_{_{\Gamma\!,n}}$ on $\RQ_n$ gives, by differentiation,
a Lie algebra map $\, \ad:\, \pgg \to \dd(\RQ_n)$.
We will also use a `shifted' map
$\,\, \ad -\chi_{_{c'}}:\, \pgg \too \dd(\RQ_n),\,$
which is  a  Lie algebra homomorphism again.
The Lie algebra map `$\ad$'
can be uniquely extended to an associative algebra
homomorphism $\, \ad: \uga\too\dd(\RQ_n)$.
Let $\ad(\Ann\varrho_{_{k,c'}}) \subset \dd(\RQ_n)$ denote
the image of $\Ann\varrho_{_{_{k,c'}}}\subset \uga$ under this homomorphism.
Similarly,
let $\,(\ad-\chi_{_{c'}})(\pggn) \subset \dd(\RQ_n)\,$  denote
the image of the Lie algebra $\pgg$ under the `shifted'
 homomorphism.

\begin{remark}
In a subsequent publication we plan to prove the following
$\G$-analogue of Theorem \ref{maintech}.

\begin{equation}\label{Gamma_2}
\,\,\bigl(\dd(\RQ_n)\cdot\ad(\Ann\varrho_{_{0,c'}})\bigr)^{\PG_{_{\Gamma\!,n}}}\;=\;
\bigl(\dd(\RQ_n)\cdot(\ad-\chi_{_{c'}})(\pggn)\bigr)^{\PG_{_{\Gamma\!,n}}}.
\end{equation}

The proof of this result is obtained by extending the argument
of \S10. Details will be given elsewhere. 

\end{remark}

Further, for any $c=(k,c')\in\CC(\ga)$,
 set 
$$I_{_{k,c'}}:=\,\bigl(\dd(\RQ_n)\cdot
\ad(\Ann\varrho_{_{k,c'}})\bigr)^{\PG_{_{\Gamma\!,n}}}=
\bigl(\dd(\RQ_n)\cdot\ad(\Ann\varrho_{_{k,c'}})\bigr) \,\cap\,
\dd(\RQ_n)^{^{\PG_{_{\Gamma\!,n}}}}.
$$
This is a two-sided ideal in $\dd(\RQ_n)^{^{\PG_{_{\Gamma\!,n}}}}$,
and we conjecture the following 
$\G$-analogue of the spherical Harish-Chandra isomorphism,
see Corollary 
\ref{isom}. 

\begin{conjecture}\label{Gamma_4}
For  any $c=(k,c')\in \CC(\ga)$,
 there is a filtration preserving algebra  isomorphism:
$\dis\dd(\RQ_n)^{^{\PG_{_{\Gamma\!,n}}}}/I_{_{k,c'}}\iso\e\hh_{1,c}\e$. 
\end{conjecture}
This isomorphism should be a quantization of  isomorphism 
(\ref{tra2}). In the special case: $n=1$ (and arbitrary group
$\G\subset \SL_2$) the conjecture is known to be true, due 
to a result of Holland [Ho].
\medskip

\noindent
{\bf Case: $\boldsymbol{\Gamma}\boldsymbol{=}\boldsymbol{\{1\}}$.}\quad
 From now on we will restrict ourselves to the special case of the trivial
group $\G$, so that: $\ga=S_n\,,\,
\ggn=\gln,$ and  $\PG_{_{\Gamma\!,n}}=\PGL_n$.
In this  case
the variety $\mm_{_{\G\!,n,c}}$ reduces to
\begin{equation}\label{mmn}
{\mathcal M}_{n}\;=\;\lbrace{(\nabla_1,\nabla_2)\in \gln\times\gln
\;\;\big|\;\;
[\nabla_1,\nabla_2]
+\id_{_{\C^n}}
=\;\text{rank } 1 \text{ matrix}\rbrace}/\text{Ad}\,\PGL_n
\,,
\end{equation}
the Calogero-Moser   space, introduced by
Kazhdan-Kostant-Sternberg \cite{KKS} and Wilson \cite{Wi}
in connection with  the  Calogero-Moser integrable system.
 We can also write:
 $\mm_n=\mu^{-1}(\O)/\PGL_n=\mu^{-1}(\bp)/\PGL_n(\bp),$
where $\mu: \gln\times\gln\to\sln\,,\,(X,Y)\mapsto [X,Y]$, cf.
 (\ref{M3}).
\medskip

 Consider the trivial
vector bundle: $\C^n\times \mu^{-1}(\bp)\to\mu^{-1}(\bp), $ on
$ \mu^{-1}(\bp)$. The fiber $\C^n$ over a point
$(X,Y) \in \mu^{-1}(\bp)$
carries the tautological action of the matrices $X$ and $Y$,
 such that $[X,Y]=\bp$. The imbedding: $\GL_{n-1}\into \GL_n\,,\,
A\mapsto A\oplus \id_{_{\Image(\bp+\id)}},\,$ makes 
$\C^n=\Ker(\bp+\id)\oplus\Image(\bp+\id)$  a $\GL_{n-1}$-module,
and  makes the vector bundle: $\C^n\times \mu^{-1}(\bp)\to\mu^{-1}(\bp) $
 a $\GL_{n-1}$-equivariant  bundle, with 
$\GL_{n-1}$ acting diagonally on  $\C^n\times \mu^{-1}(\bp)$.
 This $\GL_{n-1}$-equivariant bundle descends to a well-defined
algebraic vector bundle
 $\vv_n$ on $\mm_n= \mu^{-1}(\bp)/\GL_{n-1}$, 
to be called the {\it tautological} vector bundle.

Write $\hh_{0,c}=\hh_{0,c}(S_n)$ for the symplectic reflection algebra
attached to the Symmetric group $S_n$.
Recall that the family of 
simple $\hh_{0,c}$-modules forms, by Theorem \ref{struc},
a distinguished vector  bundle $\rr$ on  $\Spec \ZZ_{0,c}$.
Using the identification $\,\phi: \Spec \ZZ_{0,c} \iso \mm_n,$
provided by Theorem \ref{irreps} we transport $\rr$ to
 a vector bundle $\phi_*\rr$
on $\mm_n$. This  vector bundle carries a canonical $S_n$-action along the
fibers, and we proceed now to describing
some particular $S_n$-isotypic components of $\phi_*\rr$.

For any (not necessarily irreducible)
 $S_n$-module $\chi$, we
let $\rr_{(\chi)}=\Hom_{_{S_n}}(\chi,\rr)$
 denote the $\chi$-isotypic component of the
vector bundle $\rr$ on $\Spec\ZZ_{0,c}$.
In particular, take $\chi$ to be the reflection
representation in $\C^n$,
the sum of the irreducible $(n-1)$-dimensional,
and the trivial representation ${\mathbf{1}}$ of $S_n$.
Write $\rr_{_{\C^n}}$ for the corresponding
rank $n$ vector bundle on $\Spec\ZZ_{0,c}$.

\begin{proposition}\label{bundles}
Under   the isomorphism of
Theorem \ref{irreps}, the
tautological vector bundle $\vv_n$ on $\mm_n$ gets identified
 with the  vector bundle $\phi_*\rr_{_{\C^n}}$.
\end{proposition}

\Pf .
We use an $S_n$-module isomorphism:
$\,\C^n\simeq {\mathtt{Ind}}_{S_{n-1}}^{S_n}{\mathbf{1}}.\,$
By Frobenius reciprocity, this yields
$$
\rr_{_{\C^n}}\;=\;\Hom_{S_n}(\C^n\,,\,\rr)\;=\;
\Hom_{S_n}\bigl({\mathtt{Ind}}_{S_{n-1}}^{S_n}{\mathbf{1}}\,,\,
\rr\bigr)\;=\;\rr^{S_{n-1}}\,.
$$
Thus, to prove the Proposition we should identify $\rr^{S_{n-1}}$
with the vector bundle on\linebreak
 $\mm_n=\mu^{-1}(\bp)/\GL_{n-1},$
see (\ref{M3}),
 arising from the trivial
rank $n$ bundle on  $\mu^{-1}(\bp)$.

Next, recall the variety $\rep_{_{\C{W}}}(\hh_{0,c})$, see
(\ref{con_hom}).
In the special case $\G=S_n$, $c=1$, we get
the projection $\pi:  {\rep_{_{\C{S_n}}}(\hh_{0,c})}\onto \Spec \ZZ_{0,c}$.
According to Theorem \ref{components}(ii),
the trivial $\PAut_{_{S_n}}(\C{{S_n}})$-equivariant vector bundle
on ${\rep_{_{\C{S_n}}}(\hh_{0,c})}$ descends to a vector bundle on
${\rep_{_{\C{S_n}}}(\hh_{0,c})}/\PAut_{_{S_n}}(\C{S_n})= \Spec \ZZ_{0,c}$,
which is isomorphic to $\rr$.
Obviously, the pull-back $\pi^*(\rr^{S_{n-1}})$ 
may be identified with the trivial  vector
bundle on ${\rep_{_{\C{S_n}}}(\hh_{0,c})}$ 
with fiber $(\C{S_n})^{S_{n-1}}$.
Thus, the canonical isomorphism
$F: (\C{S_n})^{S_{n-1}} \iso \C^n$
provided by formula (\ref{digamma}) 
(applied in the special case $\G=\{1\}$)
and the map $\phi$ of Theorem \ref{irreps} give  a map
of rank $n$ 
 trivial vector
bundles
\begin{equation}\label{bunF}
F_{\mathtt{bun}}:\; \pi^*(\rr^{S_{n-1}})=
{\rep_{_{\C{S_n}}}(\hh_{0,c})}\times (\C{S}_n)^{S_{n-1}}\too
\mu^{-1}(\bp)\times \C^n\,.
\end{equation}

Further, 
the action of the
group $\Aut_{_{S_n}}(\C{S_n})$ on $(\C{S_n})^{S_{n-1}}$ 
may be identified
with an action of a subgroup of the group algebra $\C{S_n}$ acting 
on $\C[S_{n-1}\backslash{S_n}]$ by means of {\it right}
translations. The latter action gives, using the vector space isomorphism
$F$
and the imbedding: $\PAut_{_{S_n}}(\C{S_n})\into \Aut_{_{S_n}}(\C{S_n})$
constructed below Proposition \ref{criterion},
 a group homomorphism: $\PAut_{_{S_n}}(\C{S_n}) \to \GL_n$.
It is easy to see that the image of this homomorphism 
preserves the hyperplane: $\C^{n-1}=\{(x_1,\ldots,x_n)\in\C^n\;|\;
x_1+\ldots+x_n=0\}\,\subset\,\C^n,\,$ and acts
identically on the complementary line. Thus, one
obtains a (surjective) group homomorphism
$F_{\mathtt{group}}: \PAut_{_{S_n}}(\C{S_n}) \onto \GL_{n-1}$.

Now, 
the trivial bundle on the LHS of (\ref{bunF}) is
equivariant with respect to the group $\PAut_{_{S_n}}(\C{S_n})$, 
while the  bundle on the RHS
 is equivariant
 with respect to  group $\GL_{n-1}$
(using the imbedding (\ref{PG-split})).
It is easy to see that the vector  bundle morphism $F_{\mathtt{bun}}$
intertwines the
equivariant structures on the left and on the right via
the group homomorphism
$F_{\mathtt{group}}: \Aut_{_{S_n}}(\C{S_n}) 
\onto \GL_{n-1}$ described in the previous
paragraph.
It follows  that the morphism $F_{\mathtt{bun}}$
 descends to a vector  bundle isomorphism: $\rr^{S_{n-1}}\iso
\vv_n.\,$
\sq
\medskip

\noindent
{\bf Remark.}  In this remark
we let $n$ vary, i.e., $n=2,3,\ldots,$
and set $\rr_n:=\rr$.
Consider the projection $p_n: \mm_n \to \Sym^n(\C),$
$(X,Y)\mapsto \{\mbox{$n$-{\it tuple of eigenvalues of} $X$}.\,$
Set $Z_n= (p_n)^{-1}(0)$, and let $\rr^\circ_n:=
\rr_n\big|_{Z_n}$.
Wilson showed, see [Wi, Lemma 7.1], that for any partition
$n= n_1 + \ldots+ n_k,$ and any point $n_1\cdot
x_1 +\ldots+ n_k\cdot x_k \in
\Sym^n(\C),$ where $x_1,\ldots,x_k\in \C$ are pairwise distinct,
one has an isomorphism:
\begin{equation}\label{factor_iso}
(p_n)^{-1}(n_1\cdot x_1+\ldots+n_k\cdot x_k) \;\simeq\;
Z_{n_1} \times \ldots\times Z_{n_k}.
\end{equation}
In light of our results, one might expect that
for any partition $n= n_1 + \ldots+ n_k,$ there is a natural
$S_n$-equivariant vector bundle isomorphism (relative to
factorisation (\ref{factor_iso})):
\begin{equation}\label{factor_conj}
\!\!\!\begin{array}{ll}
\!\!\mbox{\bf Factorisation:}&
\rr_n\big|_{(p_n)^{-1}(n_1\cdot x_1 +\ldots n_k\cdot x_k)}\;
\simeq\;\Ind_{S_{n_1} \times\ldots\times S_{n_k}}^{S_n}(
\rr^\circ_{n_1}\boxtimes\ldots
\boxtimes \rr^\circ_{n_k})\,.
\end{array}
\end{equation}

\smallskip   \noindent
{\bf{Harish-Chandra homomorphism $\boldsymbol{\Phi_{0,c}}$.}}$\;$
Set $\g=\gln$, and let
 $\mreg\subset {\mathcal M}_n$ be the open subset formed by
conjugacy classes of pairs
$(X,Y)\in\g\times\g$ such that the matrix $X$ is semisimple.
We can represent points of $\mreg$ by pairs $(X,Y)$ such
 that $X$ is a diagonal matrix
and $([X,Y]+1)_{ij}=1,$ for all $i,j$,
see (\ref{s_matrix}).
Such a representation is unique
up to permuting diagonal elements.
It is easy to show,
see \cite{KKS}, that if $(X,Y)\in \mreg$ and $X=
\text{diag}(x_1,\ldots,x_n)$, then $X\in \hr$, i.e.,
$x_i\neq x_j\,,\,\forall i\neq j,$ and, moreover,
we have: $Y=\text{diag}(y_1,\ldots,y_n)+
\sum_{i\ne j}(x_i-x_j)^{-1}E_{ij}$. We regard
$x=(x_1,\ldots,x_n)$ as a point
in $\hreg\subset \h= \C^n$,
and $y= (y_1,\ldots,y_n)$ as  a point
in $\h^*=\C^n$. Then the assignment: $(X,Y)\mapsto (x,y)$
gives a well defined map
$\pi: \mreg\to (\h^*\times\hr)/S_n$, 
i.e., the image of a point of $\mreg$ depends only on the conjugacy class
of the pair $(X,Y)$.

Recall that  ${\mathcal M}_n$ has a natural symplectic structure.
Also  equip $\h^*\times\hr=T^*(\hreg)$ with the standard
 symplectic structure on the cotangent bundle. 

\begin{proposition}[\cite{Wi}] The set $\mreg$ is Zariski open dense in
$\mm_n$, and the
map $\pi: \mreg\to (\h^*\times\hr)/S_n$
 is a
symplectic isomorphism.\quad\sq
\end{proposition}
\begin{remark} The density claim in the  Proposition  follows
from the irreducibility of $\mm_n$, which is a special
case of Corollary \ref{connected}.$\quad\lozenge$
\end{remark}

The Proposition above implies that
 the inverse  map $\pi^{-1}:(\h^*\times\hr)/S_n 
\iso \mreg \into {\mathcal M}_n$
 induces an injective algebra homomorphism
 $\pi^\dag:\C[{\mathcal M}_n]\into
\C[\h^*\times\hr]^{S_n}$. Under this  homomorphism, the subalgebra of
$\C[\g\oplus \g]^\g$, formed by $\ad\g$-invariant functions in the
$X$-variable (first factor $\g$), 
goes to the subalgebra $\C[\h^*]^{S_n} \subset \C[\h^*\times\hr]^{S_n}$.
The function:
$(X,Y) \mapsto \Tr(Y^2)$ goes to $\,{\sf{L}}_{0,c}=\sum_i\, y_i^2-\sum_{i\ne
j}\,(x_i-x_j)^{-2},\,$ the Calogero-Moser hamiltonian,
cf. Definition \ref{CM_ham}. Therefore, the
 subalgebra of 
$\C[\g\oplus \g]^\g$, formed by $\ad\g$-invariant functions in the
$Y$-variable (second  factor $\g$), 
goes to $\cc_{0,c}$, the centraliser in $\C[\h^*\times\hr]^{S_n}_-$ of the
Calogero-Moser hamiltonian.

Recall the subvariety $M\subset \g\oplus \g$ defined in
(\ref{CMM}). Let $\I$ be the ideal in $ \C[\g\oplus \g]$
generated by quadratic functions of the form:
 $(X,Y)\mapsto \lambda\bigl(\bigwedge^2([X,Y]+\id)\bigr)$,
where $\lambda$ runs over $(\Lambda^2\g)^*$.
According to the  results of \cite{Wi} cited
above, $\I$ is a radical  ideal
and, moreover, coincides with the ideal of polynomials
 vanishing at $M$. We set $I_{0,c} =\I^\g$, an ideal in 
$ \C[\g\oplus \g]^\g$.
By definition, the kernel of the natural projection
$\pr: \C[\g\oplus \g]^\g\onto \C[{\mathcal M}_n]$
equals $I_{0,c}$.
Note further that this projection is a homomorphism of 
 Poisson algebras.
Composing 
the projection $\pr$
 with the homomorphism $\pi^\dag:\C[{\mathcal M}_n]\into
\C[\h^*\times\hr]^{S_n}$, one obtains
 a Poisson algebra homomorphism
$\,\Phi_{0,c}= \pi^\dag\ccirc\pr:\; 
\C[\g\oplus\g]^\g\to \C[\h^*\times\hr]^{S_n}$.

Recall that $\bcal_{0,c}$ denotes the  Poisson
subalgebra in $\C[\h^*\times\hr]^{S_n}$
generated, as a  Poisson algebra,
by $\C[\h^*]^{S_n}$ and $\cc_{0,c}\subset \C[\hr]^{S_n}$.
We claim that $\,\Image^{\,}\Phi_{0,c}=\bcal_{0,c},$ moreover,
we have

\begin{theorem}\label{class_isom}
\vi The map $\Phi_{0,c}$ gives a  Poisson  algebra
isomorphism 
$$\,\Phi_{0,c}:\,\,\C[\g\oplus \g]^\g/I_{0,c}\, =\,
\C[{\mathcal M}_n]
\,\iso\,\bcal_{0,c}\,.$$
\vii The  isomorphism $\,\phi^*:\C[{\mathcal M}_n]
\iso  \ZZ_{0,c}$ of Theorem \ref{irreps}$\mbox{(ii)}$
coincides with the following
composite map
$$
{
\xymatrix{
\Phi_{0,c}^{{\tt{spher}}}:\enspace
{\C[{\mathcal M}_n]\, =\, \C[\g\oplus \g]^\g/I_{0,c}\enspace}
\ar[r]^<>(.5){{\Phi_{0,c}}}
& \enspace\bcal_{0,c}\enspace
\ar[rr]^{\bigl(\Theta_{0,c}^{{\tt{spher}}}\bigr)^{-1}}
 &&
 \enspace\e\hh_{0,c}\e\enspace
\ar@{=}[rr]^{{\tt{Satake}}}
& &\enspace\ZZ_{0,c}\,.
}}
$$
\end{theorem}

\begin{remark}
The isomorphism of part (i) is a `commutative analogue'
of the  Harish-Chandra isomorphism
of Corollary \ref{isom}. The isomorphism  $\Phi_{0,c}^{{\tt{spher}}}$
of part (ii)
yields Theorem \ref{key_intro},
 formulated in the Introduction.\quad$\lozenge$
\end{remark}\medskip

\noindent
{\bf Proof of Theorem \ref{class_isom}}$\;\;$ 
View  the coordinate ring of
the variety ${\mathcal{M}}_n$ as a quotient of the
algebra $\C[\g\oplus\g]^\g$.
Thus, for each  $k\geq 0,$
we regard the function: $X,Y \mapsto \Tr(X^k)$
on $\g\oplus\g$ as an element $\Tr_{x,k}\in\C[\mm_n]$,
and the function: $X,Y \mapsto \Tr(Y^k)$,
as an element $\Tr_{y,k}\in\C[\mm_n]$
Also, to each  $k\geq 0,$
we attach the  elements: 
$\,z_{x,k}=\sum_i x_i^k\,,\,
z_{y,k}=\sum_i y_i^k\in\,\ZZ_{0,c}$.

We first claim  that
\begin{equation}\label{agree}
\Phi_{0,c}^{{\tt{spher}}}(\Tr_{y,2})=z_{y,2}
\;\quad\mbox{and}\;\quad
\Phi_{0,c}^{{\tt{spher}}}(\Tr_{x,k})=z_{x,k}\quad,\quad\forall 
k\geq 0\,.
\end{equation}
Here, the equations on the right are immediate. To
prove the equation on the left,
recall the Calogero-Moser Hamiltonian 
${\mathsf{L}}_{0,c}\in\C[\h\reg\times\h^*]$,
see Definition \ref{CM_ham}. We know that:
$\Theta^{{\tt{spher}}}_{0,c}(z_{y,2})={\mathsf{L}}_{0,c}$.
On the other hand, a result of [KKS] says:
$\Phi_{0,c}(\Tr_{y,2})={\mathsf{L}}_{0,c}.$ Thus,
$\Phi_{0,c}^{{\tt{spher}}}(\Tr_{y,2})=
\bigl(\Theta^{{\tt{spher}}}_{0,c}\bigr)^{-1}\ccirc\Phi_{0,c}(\Tr_{y,2})=
z_{y,2},$ and (\ref{agree}) is proved.

Next, given a simple $\hh_{0,c}$-module $E$, 
write  $E_{_{S_{n-1}}}$ for the space of
 $S_{n-1}$-{\it co}invariants.
We claim next that, for any  $ k\geq 0$, we have:
\begin{equation}\label{Claim}
(n-1)! \cdot\Tr\bigl((x_1)^k\big|_{_{E_{_{S_{n-1}}}}}\bigr)=
\Tr\big|_E(x_1)^k\quad,\quad
(n-1)! \cdot\Tr\bigl((y_1)^k\big|_{_{E_{_{S_{n-1}}}}}\bigr)=
\Tr\big|_E(y_1)^k\,.
\end{equation}
To prove this, observe that 
a generic simple $\hh_{0,c}$-module $E$
has a basis $\{v_s\,,\, s\in S_n\},$ such that,
for any $s\in S_n$ we have:
 $x_1(v_s)=
\lambda_{s^{-1}(1)}\cdot v_s,$ for certain constants:
$\lambda_1,...,\lambda_n\in \C.$
This yields the first equation in (\ref{Claim}).
The second is proved similarly by choosing a basis of $E$
adapted to $y_1$ instead of $x_1$.

Now, each of the central elements 
 $\,z_{x,k},z_{y,k}\in \ZZ_{0,c}$ acts in $E$ by a scalar,
say $\,z_{x,k,E},z_{y,k,E}.$
Using equation  (\ref{Claim}), for any $k\geq
0,$ we find
$$
\Tr\bigl((x_1)^k\big|_{_{E_{_{S_{n-1}}}}}\bigr)
\;=\;\frac{1}{(n-1)!}\cdot\Tr_{_E}(x_1)^k\;=\;
\frac{1}{n!}\Tr_{_E}\Bigl(\sum_{i=1}^n\; (x_i)^k\Bigr)
\;=\;
\frac{1}{n!}\Tr_{_E}(z_{x,k})\;=\;z_{x,k,E}.
$$
Thus, we have: $\phi^*(\Tr_{x,k})=z_{x,k}\in \ZZ_{0,c}$,
and similarly: $\phi^*(\Tr_{y,k})=
=z_{y,k}$,
for any $k\geq
0\,.$
Thus, 
 (\ref{agree}) yields:
\begin{equation}\label{agree2}
\phi^*(\Tr_{y,2})=
\Phi_{0,c}^{{\tt{spher}}}(\Tr_{y,2})
\quad\mbox{and}\quad
\phi^*(\Tr_{x,k})=
\Phi_{0,c}^{{\tt{spher}}}(\Tr_{x,k})
\quad,\quad\forall
k\geq
0\,.
\end{equation}

We claim next that
\begin{equation}\label{generate}
\mbox{\it The elements:
$\;\{\Tr_{y,2}\,;\,
\Tr_{x,k}\,,\,k\geq 0\}\;$ generate
 $\;\C[\mm_n]\;$ as a  Poisson algebra.}
\end{equation}
To prove this, consider the standard filtration on 
 $\C[\mm_n]$. Clearly, our claim reduces to a similar
statement for $\grd^{\,}\C[\mm_n]$. The Poisson algebra
$\grd^{\,}\C[\mm_n]$ is isomorphic
to $\grd^{\,}\ZZ_{0,c}\simeq\C[\h\oplus\h^*]^{S_n}$,
due to Theorem \ref{irreps}. Moreover, we have shown above that
the isomorphism  $\grd(\phi^*)$
sends $\grd(\Tr_{x,k})$ to $
\grd(z_{x,k})=\sum_i\,x_i^k\in \C[\h^*]^{S_n}$,
and sends $\grd(\Tr_{y,2})$ to $\grd(z_{y,2})=
\sum_i\,y_i^2\in \C[\h]^{S_n}$.
Thus, it suffices to show that these elements 
generate $\C[\h\oplus\h^*]^{S_n}$  as a  Poisson algebra.

To this end, observe first that the set of  elements
$\,\{\sum_i\,x_i^k\,,\,\sum_i\,y_i^k\}_{k\geq 0}\,$
 generates ${\C[\h\oplus\h^*]^{S_n}}$  as a  Poisson algebra,
by the theorem of H. Weyl cited earlier.
Thus, we have only to show that the set above generates
the same Lie algebra as the element $\sum_i\,y_i^2$
together with the set $\,\{\sum_i\,x_i^k\}_{k\geq 0}.\,$
But the two elements:
$\sum_i\,x_i^2\,,\,\sum_i\,y_i^2$
 generate a Lie subalgebra
in $\C[\h\oplus\h^*]^{S_n}$  isomorphic to
$\slt$  and, moreover,
this Lie algebra gives the standard
$\slt$-action on $\C[\h\oplus\h^*]^{S_n}$.
It is therefore clear  that applying the operator
$\,\{\sum_i\,y_i^2,-\}\,$ to $\sum_i\,x_i^k$ sufficiently
many times one gets $\sum_i\,y_i^k$, for any $k\ge 2$.
 This proves~(\ref{generate}).

We can now complete the proof of the Theorem.
Since $\Phi_{0,c}(\Tr_{y,2})= {\mathsf{L}}_{0,c}$,
it is immediate 
from (\ref{generate}) that $\Image^{\,}\Phi_{0,c}\subset\bcal_{0,c}.$
On the other hand, the discussion preceding Theorem \ref{class_isom}
implies that the map $\Phi_{0,c}$ is injective and contains
the algebra $\bcal_{0,c}$ in its image. This yields part (i) of
Theorem \ref{class_isom}.
Part (ii)  follows from equations (\ref{agree2}),
(\ref{generate}), and the fact that both
$\phi^*$ and $\Phi_{0,c}^{{\tt{spher}}}$ are Poisson algebra maps.
\quad\qed\medskip

\noindent
{\bf Remarks} \vi 
Write $\,\phi_+=\phi$
for the isomorphism: $\,\text{Irreps}(\hh_{0,c}) \iso {\mathcal M}_n\,$
of Theorem \ref{irreps}.
One can show that
replacing in the definition of $\phi_+$
the space $E^{^{S_{n-1}}}\subset E$ of $S_{n-1}$-invariants
 by the space of $S_{n-1}$-anti-invariants, one obtains a
similar map
$$
\phi_-:\;\text{Irreps}(\hh_{0,c}) \to
{\mathcal M}_n^-\;=\;\lbrace{(X,Y)\in \g\times\g\;\big|\;
[X,Y]-\id_{_{\C^n}}=\;\text{rank } 1 \text{ matrix}\rbrace}/\text{Ad}\,\PGL_n.
$$
Furthermore, one can verify by an explicit computation, that the maps
$\phi_+$ and $\phi_-$ are related to each other by the formula:
$\phi_-= \varphi\ccirc\phi_+$, where $\varphi: {\mathcal M}_n
\iso {\mathcal M}_n^-$ is an isomorphism of algebraic varieties given
by the formula $\,\varphi: (X,Y)\mapsto (X^t,Y^t).$

\vii Recall that 
Corollary \ref{cen} provides an explicit formula
for the Poisson brackets: $\lbrace{\psi_k,\phi_l\rbrace}$,
where the
"power sum" elements:
$\,\psi_\ell,\phi_\ell\in
\ZZ_{0,c}\,,\,\ell\geq 0,$ have been defined in
(\ref{power_sum}). For each $k,l\geq 0,$ put
$\,z_{k,l}:= \lbrace{\psi_k,\phi_l\rbrace}$.
It is easy to show that the set of elements:
$\,\lbrace{z_{k,l}\rbrace}_{k,l\geq 0}\,$
generates $\ZZ_{0,c}$ as an {\it associative} algebra.
Moreover, these elements
are images of the functions:
$(X,Y)\mapsto k\cdot l\cdot \Tr(X^kY^l)$ under the
isomorphism $\phi^*$.
This gives, in principle, a constructive way to identify
$\ZZ_{0,c}$ and $\C[\mm_n]$ as {\it associative} algebras.
$\quad\lozenge$
\medskip

Observe next that the $\GG$-action on $\hh_{t,c}$
arising from Theorem \ref{main_der} is trivial on $\C{W}\subset \hh_{0,c}$,
hence induces a $\GG$-action on
the $\hh_{0,c}$-module $\hh_{0,c}\e$. The latter may be viewed, by
Theorem \ref{key_intro},
as a $\GG$-action on the vector bundle $\rr$ on $\Spec\ZZ_{0,c}$.

Recall that Berest and Wilson have defined in [BW] a
$\GG$-action on the Calogero-Moser variety
${\mathcal M}_n$ by the formulas, cf. (\ref{ab}):
$$
a_{\tau,\ell}(X, Y)= (X+\tau Y^\ell\,,\, Y)\quad\text{and}\quad
b_{\tau,\ell}(X, Y)= (X\,,\, Y+\tau X^\ell)\quad,\enspace
(X,Y)\in\gln\times\gln\,.
$$
They showed further that the group $\GG$ acts {\it transitively}
on ${\mathcal M}_n$.

Repeating the proof of Proposition \ref{GG_phi}
in the `limiting' case $\ka=\infty$, we get

\begin{corollary}\label{action_onR} 
\vi The spherical Harish-Chandra homomorphism
$\,\Phi_\infty^{^{\tt{spher}}}:
\dd(\gln)^{\gln} \to \e\hh_{0,c}\e\,$
intertwines the $\GG$-actions, resp. the $\LL$-actions,
on $\dd(\gln)^{\gln}$ and on $\e\hh_{0,c}\e$.

\vii The (transitive) $\GG$-action on ${\mathcal M}_n$
lifts to a $\GG$-action on the vector bundle $\rr$ on ${\mathcal M}_n$.
A similar result holds for the Lie algebra $\LL$. \sq
\end{corollary}

\medskip \noindent
{\bf{Example: ${\mathbf{{\frak {sl}}_2}}$-case.}}\quad
The algebra $\hh_{0,c}({\mathbf{A_1}})$
has 3 generators: $s,x,y$ with defining relations
$$
s^2=1\enspace,\enspace\
sx=-xs\enspace,\enspace\ sy=-ys\enspace,\enspace\ [y,x]=-2s.
$$
In the following proposition, we collect the
results about $\hh_{0,c}({\mathbf{A_1}})$. The proofs are
straightforward and are left to the reader.

\begin{proposition}
\vi The center  of $\hh_{0,c}({\mathbf{A_1}})$ is generated by the 3 elements:
$A=x^2\,,\,B=y^2,\,$ and $C=(xy+yx)/2$, which satisfy the relation:
$\;AB=C^2-1.\;$
Thus, $\Spec\ZZ_{0,c}$ is a smooth quadric:
$\,\{(a,b,c)\in\C^3\;\big|\; ab=c^2-1\}.$

\vii Let $(a,b,c)\in \Spec\ZZ_{0,c}\subset
\C^3$. The unique irreducible representation
of $\hh_{0,c}({\mathbf{A_1}})$ corresponding to the point $(a,b,c)$
has a basis $\{v_+,v_-\}$ such that:
$$
s(v_\pm)=\pm v_\pm\enspace,\enspace
x(v_+)=\alpha\cdot v_-\enspace,\enspace x(v_-)=\beta\cdot
 v_+\enspace,\enspace
y(v_+)=\gamma\cdot  v_-\enspace,\enspace
y(v_-)=\delta\cdot  v_+,
$$
where $(\alpha,\beta,\gamma,\delta)$ is a solution
of the system of equations:
$$
\alpha\beta=a\quad,\quad\gamma\delta=b\quad,\quad\beta\gamma=c+1\quad,\quad
\alpha\delta=c-1,
$$
defined uniquely modulo  rescaling transformations:
$(\alpha,\beta,\gamma,\delta)\mapsto (t\cdot \alpha\,,\,
t^{-1}\cdot \beta\,,\,t\cdot \gamma\,,\,
t^{-1}\cdot \delta)$.

\viii The vector bundle $\rr$ is a direct sum of two line bundles:
$\rr=\e \rr\oplus {\mathbf{e_-}}\rr$, where $\e \rr=\oo$
is the trivial bundle, and
${\mathbf{e_-}}\rr$ is isomorphic to the pullback of $\oo(1)$
under the fibration $\pi:\Spec\ZZ_{0,c}\to \CP^1$, given by:
$\,\pi(a,b,c)=(\beta,\delta)$.

${\sf{(iv)}}\;$ The Poisson bracket on
${\mathcal M}_2$ is
induced by a rotation-invariant symplectic form.~\sq
\end{proposition}

\section{Appendix A: Almost commuting matrices.}
\setcounter{equation}{0}

In this appendix we  discuss a statement in
commutative algebra, which is a classical analogue
of Theorem \ref{maintech}. Throughout this Appendix
$\g=\gln$ and $G=\GL_n$.

Let $\,\II\subset\C[\g\oplus\g]=\C[X,Y]\,$
be the ideal generated by the matrix entries
of $[X,Y]\in \g=\End(\C^n)$, i.e. by all functions on
$\g\oplus\g$ of the form: $X,Y\mapsto \lambda([X,Y]),$
where $\lambda$ is a linear function on $\g$.
It is clear that the zero variety of $\II$ equals
${\mathcal Z}$, the commuting variety.
It is not known whether $\II=\sqrt{\II}$.

Recall next
that any linear map $A: \C^n\to\C^n$ induces
a linear map $\bigwedge^2\!A: \bigwedge^2\!\C^n\too\bigwedge^2\!\C^n$.
Let $\,\II_1\subset\C[\g\oplus\g]=\C[X,Y]\,$
 be the ideal generated by  the matrix entries
of $\bigwedge^2[X,Y]\in \End(\bigwedge^2\!\C^n)$, i.e. all functions on
$\g\oplus\g$ of the form: $X,Y\mapsto
\eta(\bigwedge^2[X,Y])$,
where $\eta$ is a linear function on $\End(\bigwedge^2\!\C^n).$
The  zero variety of $\II_1$ is
the variety  $\,\zz_1=$\linebreak
$\big\{(X,Y)\in \g\oplus \g\;\;\big|\;\;
\rk\bigl([X,Y]\bigr)\le 1\big\}\,.$
Again, it is not known (to us)
if $\II_1=\sqrt{\II_1}$.

We observe that $\II_1\subset \II$, in particular
$\zz\subset\zz_1$.
The classical analog of
theorem \ref{maintech} is
\begin{theorem} \label{schemeiso}
$\quad{\II}^G=(\II_1)^G$.
\end{theorem}

The proof of the Theorem requires several auxiliary lemmas.
We write $X=(x_{ij})$, $Y=(y_{ij})$ for
$n\times n$-matrices, where $x_{ij}, y_{ji}$
are indeterminates.

\begin{lemma} \label{Weyl1}
The space of $\g$-invariant polynomials on $\g\oplus \g$
with values in $\g$ is spanned by functions of the form:
$\dis
(X,Y)\,\mapsto\,
 \prod\nolimits_{i=1}^m \;\Tr\bigl(P_i(X,Y)\bigr)\cdot Q(X,Y),
\;$
where $P_1,\ldots,P_m,Q$  are noncommutative polynomials.
\end{lemma}

\Pf . Immediate  from the Weyl fundamental
theorem of invariant theory \cite{We}.
\sq\medskip

\begin{lemma} \label{iprime1} \quad
Any element of $\II$
is a linear combination of
functions of the form:
\linebreak
$\dis\;(X,Y)\,\mapsto\,
\prod\nolimits_i\;\Tr\bigl(P_i(X,Y)\bigr)\cdot
\Tr\bigl(Q(X,Y)\cdot [X,Y]\bigr),
\;$
where $P_i,Q$ are noncommutative
polynomials.
\end{lemma}

\Pf . Any element of $\II$ is of the form
$\Tr(Q(X,Y)\cdot [X,Y])$, where $Q$ is an invariant
polynomial on $\g\oplus \g$ with values in $\g$. Thus the lemma follows
 from Lemma \ref{Weyl1}.
\sq \medskip

By the definition, the ideal $\II_1$ is generated
by the matrix elements of $\bigwedge^2\bigl([X,Y]\bigr)$, i.e.
by the elements
$$
R_{jk}^{il}=[X,Y]_{ji}\cdot [X,Y]_{lk}-
[X,Y]_{li}\cdot [X,Y]_{jk}.
$$

We say that
$D\in \II$ has level $\le d$ if
it is a linear combination of
elements of the form
$\prod_i\;\Tr\bigr(P_i(X,Y)\bigr)\cdot 
\Tr\bigr(Q(X,Y)\cdot [X,Y]\bigr)$
with degree of $Q$ being $\le d$.
By Lemma \ref{iprime1}, level is defined for all
elements $D\in \II$.

To prove Theorem \ref{schemeiso},
we must show  that any element $D\in \II^G$ is zero in
$\C[\g\oplus\g]/\II_1$.
We will prove that $D=0$ in
$\C[\g\oplus\g]/\II_1$ by induction in the level of $D$.
The base of induction $(d=0)$ is obvious.
Assume that the level of $D$ is $d$ and for levels $\le d-1$
it has been proved that elements of $\II$ vanish in
$\C[\g\oplus\g]/\II_1$. By Lemma
\ref{iprime1}, we may assume that
$D=\Tr\bigl(Q(X,Y)\cdot [X,Y]\bigr)$, where the degree of
$Q$ is $d$.

 The following lemma plays the main role
in the proof of the induction step. It is here that we use
the ``rank 1'' condition. Write $\,Q_i:= Q_i(X,Y),$ for short.

\begin{lemma}\label{twocom1} Let $Q_1$, $Q_2$ be noncommutative
  polynomials in $X,Y$ of degrees $d_1,d_2$
such that $d_1+d_2\le d-4$. Then in $\C[\g\oplus\g]/\II_1$
one has:
$\;\dis
\Tr\bigl(Q_1\cdot [X,Y]\cdot Q_2\cdot [X,Y]\bigr)=0
\,.$
\end{lemma}

\Pf . We write the trace of $Q_1\cdot [X,Y]\cdot Q_2\cdot [X,Y]$
as the following sum:
$$
\sum\nolimits_{pqrs}\,(Q_1)_{pq}\cdot [X,Y]_{qr}\cdot (Q_2)_{rs}\cdot [X,Y]_{sp}=
\sum\nolimits_{pqrs}\,(Q_1)_{pq}\cdot (Q_2)_{rs}\cdot [X,Y]_{qr}\cdot [X,Y]_{sp}\,.
$$

Using the relations
of $\II_1$, we conclude that
in $\C[\g\oplus\g]/\II_1$ one has:
{\small {$$
\Tr\bigl(Q_1\cdot [X,Y]\cdot Q_2\cdot [X,Y]\bigr)=
\sum_{pqrs}(Q_1)_{pq}\cdot (Q_2)_{rs}\cdot 
[X,Y]_{sr}\cdot [X,Y]_{qp}
=\Tr\bigl(Q_1\cdot [X,Y]\bigr)\cdot\Tr\bigl(Q_2\cdot [X,Y]\bigr),
$$}}
which is zero in $\C[\g\oplus\g]/\II_1$
by the induction hypothesis.
\sq\medskip

\begin{lemma}\label{twocom21} Let $Q_1$, $Q_2$ be noncommutative
  polynomials in $X,Y$ of degrees $d_1,d_2$
such that $d_1+d_2\le d-4$. Then in $\C[\g\oplus\g]/\II_1$
one has:
$\,\dis
\Tr\bigl(Q_1\cdot X\cdot Y\cdot Q_2\cdot [X,Y]\bigr)=$\linebreak
$
\Tr\bigl(Q_1\cdot Y\cdot X\cdot Q_2\cdot [X,Y]\bigr).$
\end{lemma}

\Pf .
We have
$$
\Tr\bigl(Q_1\cdot X\cdot Y\cdot Q_2\cdot [X,Y]\bigr)=
\Tr\bigl(Q_1\cdot Y\cdot X\cdot Q_2\cdot [X,Y]\bigr)
+\Tr\bigl(Q_1\cdot [X,Y]\cdot Q_2\cdot [X,Y]\bigr),
$$
By Lemma \ref{twocom1},
$\;\dis
\Tr\bigl(Q_1\cdot [X,Y]\cdot Q_2\cdot [X,Y]\bigr)
\,$ is zero
in $\C[\g\oplus\g]/\II_1$.
The lemma is proved.
\sq\medskip

\noindent
{\bf {Proof of Theorem \ref{schemeiso}.}}\quad
The last lemma shows that in order to prove
the theorem it is enough to check the
statement of the theorem for
the elements $D_m=\Tr\bigl(Q_m\cdot [X,Y]\bigr)$,
$m=0,...,d$, where $Q_m=\sum_{k=0}^m X^{m-k}Y^{d-m}X^k$.
Thus, the theorem follows from the equation
\begin{align*}
D_m
&
=\sum\nolimits_{k=0}^m\;\Tr\bigl(X^{m-k}Y^{d-m}X^k[X,Y]\bigr)\\
&=
\sum\nolimits_{k=0}^m\;\Tr(Y^{d-m}X^k[X,Y]X^{m-k})=
\Tr\bigl([Y^{d-m}X^m,Y]\bigr)=0.\quad\square
\end{align*}

\begin{corollary} \label{varietyiso}
$(\sqrt{\II_1})^G=(\sqrt{\II})^G$.\sq
\end{corollary}

There is a more geometric proof of
Corollary \ref{varietyiso}, which does not
use theorem \ref{schemeiso}.
It is based on the following beautiful linear algebraic lemma,
see e.g., \cite{Gu}:

\begin{lemma} \label{Rud} Let $A,B$ be two $n\times n$-matrices whose
commutator, $[A,B]$ is a rank $1$ matrix. Then $A,B$ can be
simultaneously conjugated into upper
triangular matrices.
\end{lemma}

The proof of this lemma in \cite{Gu} is rather complicated,
so we would like to
give a simpler proof which we learned from A. Rudakov.

\Pf . Without loss of generality we can suppose that
the kernel of $A$ is non-zero; otherwise we can replace $A$ by
$A-\lambda\cdot\id$, where $\lambda$ is an eigenvalue of $A$.
It is enough to show that there is a proper
subspace which is both $A$- and $B$-stable.
Put $C:=[A,B]$.

If $\Ker^{\,}{A}$
 belongs $\Ker^{\,}{C}$ then $\Ker^{\,}{A}$ is $B$-invariant:
$A Bx=B Ax + Cx = 0$. Thus, $\Ker^{\,}  A$ is the needed subspace.

If $\Ker^{\,}  A$ is not in $\Ker^{\,}  C$
 then there exists $z$ such that $Az=0, Cz\ne 0$.
This means $ABz=Cz$ is not zero.
Hence $\Image^{\,}  C$, which is a 1-dimensional space,
is contained in $\Image^{\,}  A$.
Now $\Image^{\,}  A$ is $B$-invariant, because $B Ax= A Bx+ Cx$, and both
 summands are in
$\Image^{\,}  A$. Thus, $\Image^{\,}  A$
 is the needed subspace. $\square$\medskip

We observe now that Lemma \ref{Rud} implies Corollary \ref{varietyiso}.
In more detail, recall the variety
$\,\zz_1=\big\{(X,Y)\in \g\oplus \g\;\;\big|\;\;
\rk\bigl([X,Y]\bigr)\le 1\big\}\,.$
It follows from Lemma \ref{Rud}
 that the closure of any $G$-orbit
in ${\mathcal Z}_1$ contains a point of
${\mathcal Z}$ (since a pair of upper triangular
matrices can be conjugated arbitrarily  close
 to their diagonal parts).
Hence, any $\g$-invariant polynomial on $\g\oplus \g$
that vanishes on ${\mathcal Z}$,
 must vanish on the whole of ${\mathcal Z}_1$.\quad\sq
\medskip

\noindent
{\bf Remarks} \vi We have shown that Lemma \ref{Rud} implies
Corollary \ref{varietyiso}. In fact, the corollary also easily implies the
lemma. Indeed, it is easy to show that any two matrices
the closure of whose $GL_n$-orbit contains a pair of commuting matrices,
generate a solvable Lie algebra.

\vii There is a natural analogue of the claim of
Lemma \ref{Rud} for a general semisimple Lie algebra $\g$
saying
that any pair of elements $A,B\in \g$
such that $[A,B]$ is a root element, generates a solvable Lie algebra.
This claim is {\it false} in such generality;
 it is false for the Lie algebra $\g={\mathfrak{s}\mathfrak{o}}_8$,
for instance.

\section{Appendix B: Geometric construction of
 $V_k$}\setcounter{equation}{0}

In this appendix we explain the geometric construction of the representation
$V_k$ over $\g=\sln$ using the notion of twisted differential
operators, see e.g., \cite{BB}.

Let $\om$ be the canonical
(line) bundle on $\CP^{n-1}$, and let
$\LB$ denote its total space with  zero-section  removed.
Thus, $\LB$ is a principal $\C^*$-bundle on $\CP^{n-1}$.
Write $\eu$ for the  Euler vector field on $\LB$ generating the
$\C^*$-action, and $\dd(\LB)^{\eu}$  for the centralizer of the Euler field
in the algebra of global algebraic
 differential operators on $\LB$.
Note that $\eu$ is a central element in $\dd(\LB)^{\eu}$.
For any $k\in\C$, we set
$\,\dd_k(\LB):=\dd(\LB)^{\eu}/(\eu-k)\cdot\dd(\LB)^{\eu}\,.$
\medskip

\noindent
{\bf Remark.} The algebra $\dd_k(\LB)$ is known as a ring of
{\it twisted differential operators} on $\CP^{n-1}$.
In the special case when $k$ is an integer,
the algebra $\dd_k(\LB)$ turns out to be isomorphic to
$\dd(\CP^{n-1},\,\om^{\otimes(-k)})$, the algebra
of global algebraic differential operators acting on the sections of
 the line bundle $\om^{\otimes(-k)}$.$\quad\lozenge$\medskip

 The torus $T\subset SL_n$, formed by
 diagonal matrices, has a unique open dense orbit in $\CP^{n-1}$.
Let  $\LB\reg\subset \LB$ be the preimage of this
$T$-orbit in $\LB$, and let
$\,\dd_k(\LB\reg):=$
$\dd(\LB\reg)^{\eu}/(\eu-k)\cdot \dd(\LB\reg)^{\eu}\,$
be the corresponding ring of twisted differential operators.

For each $k\in \C$, we now define an $\Ug$-module $V_k$
as follows.
The tautological $SL_n$-action on
$\CP^{n-1}$ lifts naturally to an action on $\LB$, that commutes
with the $\C^*$-action. Differentiating
the $SL_n$-action
gives, by restriction to $\LB\reg$, a Lie algebra morphism
$\;\varrho_k: \g=\sln \too \C^*$-{\it equivariant
vector fields on} $\LB\reg$.
This Lie algebra morphism extends to an associative algebra
homomorphism: $\Ug \to\dd(\LB\reg)^{\eu}\,,$
and we form the composite map:
\begin{equation}\label{varrho}
\varrho_k:\; \Ug \too \dd(\LB\reg)^{\eu}\onto
\dd(\LB\reg)^{\eu}/(\eu-k)\cdot\dd(\LB\reg)^{\eu}\;=\;
\dd_k(\LB\reg)\,.
\end{equation}
Let $\dd_k(\LB\reg)\cdot\varrho_k(\h)\subset \dd_k(\LB\reg)$
be the left ideal generated by the image of Cartan subalgebra
$\h=\Lie T$. Thus,
$\,\dd_k(\LB\reg)/\dd_k(\LB\reg)\cdot\varrho_k(\h)\,$ is
a left $\dd_k(\LB\reg)$-module, hence, an $\Ug$-module via
(\ref{varrho}). The proof of the
following result is left to the reader.

\begin{proposition} The $\g$-module $\,\dd_k(\LB\reg)/\dd_k(\LB\reg)\cdot\varrho_k(\h)\,$
is isomorphic to $V_k$.\sq
\end{proposition}

\section{Appendix D: Small representations.}
\setcounter{equation}{0}

In this appendix we will consider the special case of the radial part
construction of \S6
and the Dunkl-Cherednik construction in \S4
 in the case when $V$ is a {\it small}
finite dimensional representation of $\g$.
This situation is of special interest because in the case of
small representations both constructions give the same answer.

Let $\g$ be a finite dimensional simple Lie algebra. 
Let $\varrho: \Ug\to\End V$ be a finite dimensional
(not necessarily irreducible) representation
such that all of its weights belong to the root lattice of $\g$.
Then, $\varrho$
 is called {\it small}
if $\varrho(e_\alpha)^2{v}=0$, for any root vector $e_\alpha\in \g$
and any $v\in V\!\langle 0\rangle$.
This condition is equivalent to saying that
$2\alpha$ is not a weight of $V$, for any root $\alpha$.
Note that $V=\g$,  the adjoint representation of any reductive Lie algebra,
is small.  We refer to the paper [Br] for more information about small
representations.

For $\g=\sln$, one obtains a complete understanding of
small representations as follows. Note that, for any
finite-dimensional $\g$-module $V$, the Weyl group $W$ acts naturally
on $V\!\langle 0\rangle$.
Write $\C^n$ for the tautological representation,
and $(\C^n)^{\otimes n}$ for its $n$-th tensor power.
We let $\Cat$ denote the (semisimple) abelian category
of finite dimensional $\sln$-representations 
$V$ such that all
simple constituents of $V$ are among the simple  constituents of 
$(\C^n)^{\otimes n}$.

\begin{proposition}\label{small_rep}
\vi The representation $(\C^n)^{\otimes n}$, and its dual, are both small.

\vii An  irreducible $\sln$-module is small if and
only if it is  a submodule of  either $(\C^n)^{\otimes n}$,
or of its dual. In particular, any object of $\Cat$ is a small 
 representation.

\viii The assignment: $V\,\rightsquigarrow
\, V\!\langle 0\rangle$ gives an equivalence
between
the category $\Cat$ and the
category of finite-dimensional $S_n$-modules.
This equivalence sends $(\C^n)^{\otimes n}$  to 
the regular representation of $S_n$.
\end{proposition}

\Pf . We leave the proof to the reader, and only comment on the last
statement. Let $e_1,\ldots,e_n$  the standard basis in $\C^n$.
It is clear that,
for $V=(\C^n)^{\otimes n}$,
 the element $e_{i_1}\otimes\ldots\otimes
e_{i_n}$ belongs to the zero weight subspace if and only
if the indices ${i_1},\ldots,{i_n}$ are pairwise distinct.
Hence
the
space $V\!\langle 0\rangle$ has a basis
formed by the elements $e_{\sigma(1)}\otimes\ldots\otimes
e_{\sigma(n)}\,,\, \sigma\in S_n.\,$
Therefore we may identify $V\!\langle 0\rangle$ with $\C{S_n}$
via the correspondence: $\,\sigma\,\longleftrightarrow\,
e_{\sigma(1)}\otimes\ldots\otimes
e_{\sigma(n)}\,.$\sq\medskip

Recall the Harish-Chandra homomorphism $\Psi_V:
\dd(\gr)^\g \too \dd\bigl(\hreg,\, \End_{_\C}V\!\langle
0\rangle\bigr)^W$,
see (\ref{psiv}).
On the other hand, consider
the rational Cherednik algebra
$\hh_{_{1,1}}$, with parameters: $t=1$ and $c=1$.
The {\it renormalized} Dunkl homomorphism
$\Theta_{_{1,1}}:\hh_{_{1,1}}\to \dd(\hreg,\C{W})$,
see above formula (\ref{theta_spher}), restricts to  a homomorphism
$\Theta_{_{1,1}}:\dis (\hh_{_{1,1}})^W\to \dd(\hreg,\C{W})^W$, 
where $(\hh_{_{1,1}})^W$ stands for the centralizer of $W$ in 
the algebra $\hh_{_{1,1}}$. Further, 
the action of $W$ on $V\!\langle 0\rangle$ gives
an algebra morphism
$\varpi_{_V}: \C{W}\to \End_{_\C} V\!\langle 0\rangle$.

The following result explains the relevance of
the notion of small representation for the study of
the Harish-Chandra homomorphism, cf. also [Br, Thm.1]

\begin{theorem}\label{small_conj}
If $t=1$ and $c=1$ then, for any small representation $V$, the following
diagram
 commutes
$$
{
\diagram
(S\g)^\g\;\;\rdouble^<>(.5){{\mbox{\footnotesize
{Chevalley}}}}_<>(.5){{\mbox{\footnotesize
{isomorphism}}}}
\dto_{{\mbox{\footnotesize
{inclusion}}}} &
\;\;S\h^W \rto^<>(.5){{\mbox{\footnotesize
{inclusion}}}} & (\hh_{_{1,1}})^W\dto^{\Theta_{_{1,1}}}\\
\dd(\g)^{\g} \rto^<>(.5){\Psi_V}
& \dd\bigl(\hreg\,,\,\End_{_\C} V\!\langle 0\rangle\bigr)^{W}&
\dd(\hreg\,,\,\C{W})^W\lto_<>(.5){\varpi_{_V}}\;.
\enddiagram}
$$
\end{theorem}

\noindent
{\bf Example:}\quad Let $\varrho: \g\to \End_{_\C}V$ be a small
representation.
Then a simple calculation with
1- and 3-dimensional ${\mathfrak{s}\mathfrak{l}}_2$-modules
shows that, for any $\alpha\in R$, the endomorphism
$\varrho(e_\alpha\cdot e_{-\alpha}):
V\langle 0\rangle\to V\!\langle 0\rangle$
is given by the formula: $\frac{(\alpha,\alpha)}{2}\cdot
(\id-s_\alpha),$
where $s_\alpha$ denotes the action on $V\langle 0\rangle$
of the simple reflection $s_\alpha\in W$. It follows from
(\ref{CM}) that
 for the  Harish-Chandra homomorphism $\,\Psi_V:
\dd(\g)^\g\to \dd\bigl(\hreg,\, \End_{_\C}V\!\langle
0\rangle\bigr)^W,$ associated to  $V$, one obtains:
\begin{equation}\label{psi_g}
 \Psi_{V}(\Delta_{_\g})\;=\;
\Delta_{_\h} -\sum\nolimits_{\alpha\in R}\;
\frac{(\alpha,\alpha)}{2}\cdot
\frac{(\id-s_\alpha)}{\alpha^2}\,.
\end{equation}\medskip

\noindent
{\bf Proof of Theorem \ref{small_conj}:}$\;$
Let $\Delta\in S\h^W\subset\hh_\ka$
be the quadratic $W$-invariant corresponding to 
the Laplacian $\Delta_{_\h}$
on $\h$.
In the special case $t=1$ and  $c=1,$ one finds:
\begin{equation}\label{psi_hh}
\Theta_{_{1,1}}(\Delta)\;=\;
\Delta_{_\h} -\sum\nolimits_{\alpha\in R}\;
\frac{(\alpha,\alpha)}{2}
\cdot \frac{(1-s_\alpha)}{\alpha^2}\;
\in \;\dd\bigl(\hreg,\C{W}\bigr)^W \,.
\end{equation}
Since $\Delta_\g\in (S\g)^\g$ goes under the Chevalley isomorphism to $\Delta\in S\h^W$,
comparing formulas  (\ref{psi_hh}) and (\ref{psi_g}) yields
an equality: 
$\varpi_V\ccirc\Theta_{_{1,1}}(\Delta_\h)= \Psi_V(\Delta_\g).\,$
Write $D$ for this element of  $\dd\bigl(\hreg,\, \End_{_\C}V\!\langle
0\rangle\bigr)^W$.

Now, fix $f\in (S\h)^W$, a homogeneous element of degree $m$.
First, view $f$ as an element of $(S\g)^\g$,
via the  Chevalley isomorphism, write $f_\g$ for the
corresponding $\ad \g$-invariant constant coefficient differential
operator
on $\g$, and set: $D_1=\Psi_V(f_\g)$.
Next, view $f$ as an element of $(\hh_{_{1,1}})^W$,
and set: $D_2= \varpi_V\ccirc\Theta_{_{1,1}}(f)$.
The algebra $(S\h)^W$, being commutative and all the maps
in the diagram of the Theorem being algebra homomorphisms,
we conclude that both $D_1$ and $D_2$ commute with $D$.
Further,  $D_1,D_2$ are differential operators of order $\leq m$, whose
principal symbols are both equal to $f$. Moreover, it is clear that
$D_1,D_2$  both have
homogeneity degree $-m$ with respect to the dilation action on $\hreg$.
Thus, $D_1-D_2$
is a differential operator on $\hreg$  of order $\leq m-1$,
that commutes with $D$ and has homogeneity degree $-m$.
The proof  of the Theorem is now completed by the following result.

\begin{proposition}\label{opzero} 
Let $M$ be a $\End_{_\C}V\langle 0\rangle$-valued differential operator  on
an open subset of $\h$  that commutes with
$D$, has order $\le m-1$, and homogeneity degree $-m$.
Then $M=0$.
\end{proposition}

The remainder of the section is devoted to the proof of
Proposition \ref{opzero}. We need an easy lemma from classical
representation theory. Write $E:p\mapsto E(p)=|p|^2$ for the
quadratic polynomial on $\h^*$ corresponding to
 $\Delta\in S^2\h$,
and view it as a function on $\h\oplus \h^*$ constant along the
first factor.

\begin{lemma}\label{symbols} 
Let $S: (x,p)\mapsto S(x,p)$ be a rational function on
$\h\oplus \h^*$, which is polynomial in $p\in \h^*$.
Suppose that the Poisson bracket $\lbrace E ,S\rbrace$
equals to zero. Then $S\in \C[\h\oplus \h^*]$.
\end{lemma}

{\it Proof of Lemma.}$\;$ We may assume without loss of generality that
the function $S$ is homogeneous in the $p$-variable, say of degree
 $N$.
By shifting $x$, we may assume further that
$S$ is regular at $x=0$. Let $S=\sum_{k=0}^\infty\; S_k$
be the  Taylor expansion of $S$ at the point $(0,0)\in \h\oplus \h^*$.
Thus, $S_k$ is a  homogeneous
polynomial on $\h\oplus \h^*$ of total degree $k$. 
Clearly, $S_k$ has degree  $N$  in the $p$-variable.
Further, separating homogeneous 
components in the equation   $\lbrace E,S\rbrace=0$
yields: $\lbrace E,S_k\rbrace=0$, for any $k\geq 0$.

The quadratic functions: $E(x,p)= |p|^2\,,\, H(x,p)=\langle x,p\rangle$, and
$F(x,p)=|x|^2,\,$
$(x,p)\in \hr\times \h^*$,  are well-known
to form an $\slt$-triple: $\,\{E,H,F\},\,$ with respect to the Poisson bracket
on $\C[\h\oplus~\h^*]$. Hence, taking the brackets with these functions
gives, for each $m\geq 0$, an $\slt$-action on $\C[\h\oplus \h^*]^{(m)}$,
the space of degree $m$ homogeneous
polynomials.
By representation theory of $\slt$, any
vector annihilated by $E$ (i.e. a highest weight vector) has to have
 a nonnegative weight. This means that the $p$-degree of such a polynomial
is greater than or equal to its $x$-degree. In particular, the $x$-degree of $S_k$
is $\le N$, which means that its total degree is $\le 2N$.
Therefore, for
$k>2N$, we have $S_k=0$.
\quad\sq\medskip

\noindent
{\bf Proof of Proposition \ref{opzero}:}$\;$
Assume the contrary. Let $S: (x,p)\mapsto S(x,p)$ be the principal symbol
of $M$, where $(x,p)\in \h\oplus \h^*$.  Then $S$ is a polynomial in $p$
and $\lbrace E, S\rbrace=0$. We know that $S$ has homogeneity degree
$-m$, in the grading where: $\deg(x)=1 =$
$ -\deg(p)$. Furthermore, $S$ has order
$\le m-1$ with respect to $p$. This implies that $S$ cannot be a
polynomial in $x$. That contradicts the lemma, and we are done.\sq

\section{Deformations, Poisson brackets and cohomology}\label{coh}
\setcounter{equation}{0}

\noindent
{\bf{Deformations and Poisson structures.}}\quad
Let $\BA$ be a flat formal (not necessarily commutative) 
deformation of an associative
{\it commutative} algebra $A$, in other words, we are given a
topologically free
$\C[[t]]$-algebra $\BA$ such that
$\BA/t\BA=A$.
Then, there is a canonical Poisson bracket on
$A$ defined as follows.
First, for each
${\widetilde{u}}, {\widetilde{v}}\in
\BA$  we have, since $A$ is commutative:
$[{\widetilde{u}},{\widetilde{v}}]\in t\BA,$
where $[{\widetilde{u}},{\widetilde{v}}]:=
{\widetilde{u}}{\widetilde{v}} -{\widetilde{v}}{\widetilde{u}}$.
Let $m({\widetilde{u}},{\widetilde{v}})\geq 1$ be the greatest
integer (possibly $= \infty$)
such that $[{\widetilde{u}},{\widetilde{v}}]\in
t^{m({\widetilde{u}},{\widetilde{v}})}\BA$.
Let ${\mathbf{m}} \geq 1$ be the minimum of the integers
$m({\widetilde{u}},{\widetilde{v}})$ over all pairs
${\widetilde{u}},{\widetilde{v}}$ as above.
If ${\mathbf{m}}=\infty$ we set the Poisson bracket on $A$ to be zero.
If ${\mathbf{m}}$ is finite then, given $u,v\in A$,
choose any representatives
${\widetilde{u}},{\widetilde{v}}\in \BA$,
so that $u={\widetilde{u}}\,\,\mbox{\sl mod}\;
(t\BA)$ and $v=\widetilde{v}\,\,\mbox{\sl mod}\;
(t\BA)$, and put
$\,\{u,v\}=
{(t^{-{\mathbf{m}}}[{\widetilde{u}},
{\widetilde{v}}])}
\,\,\mbox{\sl mod}\;
(t\BA).\,$
It is known that the assignment:
$u,v \mapsto \{u,v\}$ gives rise to a well-defined Poisson
bracket on $A$, independent of the choices involved.

 Observe that  the equality: ${\mathbf{m}}=\infty$
in the construction above implies $[{\widetilde{u}},{\widetilde{v}}]
\in t^\ell\cdot\BA,$
$\forall{\widetilde{u}},{\widetilde{v}}\in \BA,$
and $\forall \ell\geq 1$. Since
 $\,\bigcap_{\ell\geq 1}\;t^\ell\BA =0,\,$
we  see that the
vanishing of the Poisson
bracket on $A$ forces
the whole deformation to be
commutative,
i.e., we have
\begin{lemma}\label{commut} If $\,{\mathbf{m}}=\infty,\,$ then
$\BA$ is a commutative algebra.\quad\qed
\end{lemma}

\begin{remark} As has been observed by Hayashi [Ha],
a slight modification of the above construction still applies in the
case of a {\it non-commutative} algebra $A$, and gives
a Poisson bracket on $\ZZ(A)$, the center of $A$.
We note that, as opposed to the
commutative case considered above, the  Hayashi construction only
works if ${\mathbf{m}}=1$, in particular,
the bracket on $\ZZ(A)$ may turn out to be identically zero.
This is so, for instance, if $\ZZ(A)$ is a specialization
of the family $\,\{\ZZ(A_t)\}_{t\neq 0}.\,$
$\quad\lozenge$
\end{remark}\medskip

     From now until the end of this section we 
keep the following setup.
Let $A$ be a commutative  associative
$\C$-algebra, and $\BA$ a complete topological $\C[[t]]$-algebra
such that $\BA/t\cd\BA=A$.
Write $\KK=\C((t))$ for the field of formal Laurent series,
and set: $\BA_{_{\KK}}=\KK\,\mbox{$
\widehat{\otimes}$}_{\!_{\C[[t]]}}\BA$,
a $\KK$-algebra obtained by localizing at $t$
(here `$\,\widehat{\otimes}\,$' stands for the
completed tensor product).
Let $\HH^\bullet(\BA_{_{\KK}})$ denote the 
 Hochschild cohomology of  $\BA_{_{\KK}}$
(relative to the ground field $\KK$).
Further, let $H^\bullet(\Spec A)$ denote the ordinary
singular cohomology (with complex coefficients) of
$\Spec A$ viewed as a topological space.
For any vector space $E$, write $E[[t]]$ for the
space of formal power series with coefficients in $E$.

The deformation $\BA$ gives rise
to a Poisson structure $\,\{-,-\}\,$ on $A$,
making $\Spec A$ a Poisson variety.
The following result is based on the celebrated
{\it Formality theorem} of Kontsevich
[Kon].

\begin{proposition}\label{konts}
Assume that the canonical isomorphism: $\BA/t\cdot\BA\iso A$
can be lifted to 
a topological $\C[[t]]$-module isomorphism: 
\begin{equation}\label{ass}
\phi: \BA\iso A[[t]]\,.
\end{equation}
If $\Spec A$ is smooth and the Poisson structure  $\,\{-,-\}\,$
on $A$ is non-degenerate, i.e., makes  $\Spec A$ a symplectic manifold,
then there is a graded $\KK$-algebra isomorphism:
 $\HH^\bullet(\BA_{_{\KK}})$
$\simeq\KK\otimes_{_\C}H^\bullet(\Spec A).$
\end{proposition}

To prove the Proposition we transport the
multiplication map from $\BA$ to $A[[t]]$
via the bijection (\ref{ass}) to obtain
a $\C[[t]]$-bilinear 
 map
$\,\star: A[[t]] \otimes A[[t]] \too A[[t]]
\,,$
$\,
(P,Q)\mapsto P\star Q= \phi\bigl(\phi^{-1}(P)\cdot\phi^{-1}(Q)\bigr).\,$
Expanding into powers of $t$ one can write it
in the form of a {\sf star-product}
\begin{equation}\label{star}
P\star Q
\;=\;
P\cdot Q + t\cdot \mu_1(P,Q)
+ t^2\cdot \mu_2(P,Q) +\ldots \,.
\end{equation}
The algebra $(A[[t]], \star)$ thus obtained
is clearly isomorphic to 
$\BA$.
The following result is quite standard.

\begin{lemma}\label{triv}
Let $\BA$ be a deformation of a commutative
algebra
$A$, as in Proposition \ref{konts}.
 Then there exists a $\bigstar$-product on $A[[t]]$
equivalent to (\ref{star}), such that:

$\bullet\enspace$
 Each map $\mu_i: A\times A\to A$
in this new $\bigstar$-product
is given by
a certain regular bi-differential operator on
$\Spec A$;\newline
$\hphantom{x}\quad\bullet\enspace$
$\mu_1(P,Q)=\{P,Q\},\,$
where
$\,\{-,-\}\,$ stands for the Poisson bracket on
$A$ induced by the deformation $\BA$;\newline
$\hphantom{x}\quad\bullet\enspace$ 
The $\C[[t]]$-algebra $\BA$
is isomorphic to $\bigl(A[[t]],\bigstar\bigr).$
\end{lemma}

\noindent
{\bf Proof of Lemma \ref{triv}:}\quad
Put $\mm=\Spec A$. Then for any $k\ge 0$,
the Hochschild cohomology,
 $\HH^k(A)$,
are well-known  to be isomorphic to, $\Om^k(\mm)$, the space of
 regular differential $k$-forms on $\mm$. 
Let $C^\bullet(A)$ denote the  standard Hochschild cochain
complex, and $C^\bullet_{_{\sf{diff}}}(A)\subset C^\bullet(A)$
the subcomplex formed by poly-differential operators.
Differential forms on $\mm$ being local objects, it follows that the
imbedding:
$C^\bullet_{_{\sf{diff}}}(A)\into C^\bullet(A)$ is a 
quasi-isomorphism. Recall that
all the terms  $\mu_i$ in the
expansion (\ref{star}) are controlled by (second and third)
Hochschild cohomology groups.
Therefore, changing the star-product by an equivalent one,
one may achieve that:  $\mu_i\in C^2_{_{\sf{diff}}}(A)$,
for all $i\geq 1$. This proves the claim about
 bi-differential operators. The  claim
concerning the Poisson bracket is well-known.
The Lemma is proved.
\quad\qed\medskip

\noindent
{\bf Proof of Proposition \ref{konts}:}\quad
 Set $\mm:=\Spec A$, a smooth affine symplectic variety.
Hence, the Poisson structure on $\mm$
is given
by a regular bi-vector field: $\beta\in \bigwedge^2T\mm$.
Now, fix a star-product on $A[[t]]$ as in  Lemma
\ref{triv}.
By (an algebraic version of) the 
{\it Formality theorem} of Kontsevich
[Kon], this star-product is equivalent
to the one obtained  by the Kontsevich
quantization procedure applied to an appropriate
formal Poisson structure of the form:
$\widehat{\beta}=\beta_0+t\cdot\beta_1+t^2\cdot\beta_2+\ldots
\in (\bigwedge^2T\mm)[[t]],\,$
such that $\beta_0=\beta$.
Write  $H\!P^\bullet(-)$,
for the Poisson
 cohomology of a  Poisson algebra,
see Brylinski [Br]  and references therein.
Kontsevich's theorem on the
{\it cup-product on tangent cohomology}, see [Kon, \S8],
yields a graded $\KK$-algebra isomorphism:
\begin{equation}\label{tangent}
H\!H^\bullet(\BA_{_{\KK}})\, \simeq\, 
H\!P^\bullet\bigl(\KK[\mm]\,,\, \widehat{\beta}^{\,}
\bigr)\,.
\end{equation}
Here, the RHS stands for the   Poisson
 cohomology of the commutative $\KK$-algebra $\KK[\mm]
:=\KK\widehat{\otimes}_{_\C}\C[\mm]$ equipped
with the  $\KK$-linear 
Poisson structure given by $\widehat{\beta}$.

The cohomology on the RHS of (\ref{tangent}) 
are computed by means of $\KK\widehat{\otimes}_{_\C}(\bigwedge^\bullet T\mm)$,
the complex 
formed
by poly-vector fields,
with  differential being induced by the Schouten bracket
with the bivector $\widehat{\beta}$.
It follows, since the Poisson structure $\beta=\beta_0$ on $\mm$
is non-degenerate and $\C[[t]]$ is a local ring,
that the natural pairing with $\widehat{\beta}$
induces a $\C[[t]]$-module isomorphism: $(T^*\mm)[[t]]
\iso (T\mm)[[t]]$, hence, a graded isomorphism:\linebreak
$\KK\widehat{\otimes}_{_\C}(\bigwedge^\bullet T^*\mm)\iso
\KK\widehat{\otimes}_{_\C}(\bigwedge^\bullet T\mm)$.
Brylinski showed [Br] that the de Rham differential
goes under the latter isomorphism to the Poisson differential.
Thus, writing $H_{_{\!DR}}^\bullet(\mm)$ for 
the algebraic de Rham cohomology of $\mm$
 one obtains a natural
$\KK$-algebra  isomorphism: 
$\,H\!P^{j}(\KK[\mm]) \simeq \KK\otimes_{_\C}H_{_{\!DR}}^\bullet(\mm).\,$
But de Rham cohomology of the smooth affine variety $\mm$
is known (Grothendieck) to be isomorphic to $H^\bullet(\mm),$
the ordinary singular cohomology of $\mm$, viewed as a
topological space.
Combining these considerations with (\ref{tangent}), we find
\begin{equation*}
H\!H^\bullet(\BA_{_{\KK}})\, \simeq\, 
H\!P^\bullet\bigl(\KK[\mm]\,,\, \widehat{\beta}^{\,}
\bigr)\,\simeq \, \KK\otimes_{_\C}H_{_{\!DR}}^\bullet(\mm)
\,\simeq\, \KK\otimes_{_\C}H^\bullet(\mm)
\,.\qquad\qquad\square
\end{equation*}
\medskip

In applications,  Proposition \ref{konts} will be used
in the following situation.
Let  $\A$ be
a flat formal deformation
of $A$, as above. Assume given an increasing
filtration:
 $\C[[t]]=F_0\BA\subset F_1\BA\subset\ldots,$
 on 
$\BA$ such that:
$\; \cup_{i\geq 0}\;F_i\BA= \BA,\,$ and
$\;F_i\BA\cdot F_j\BA\subset F_{i+j}\BA,\quad\forall i,j\geq 0\,.
$
The filtration on $\BA$ induces a filtration on $A=\BA/t\cd\BA$,
and we write $\grd\BA$, resp. $\grd A$,
 for the corresponding associated graded
algebra. The following result is proved by a routine argument
involving, for each $i\geq 0$,
 an arbitrary choice of complementary subspace
to $F_iA$ in $F_{i+1}A$. 

\begin{lemma}\label{split}
Assume  
that: \vi $\grd A$ is a finitely generated (commutative)
$\C$-algebra, in particular,  $\,\dim F_iA < \infty,$
for any $i\geq 0$. 

\vii
The algebra $\A$ is complete, and the
isomorphism: $\grd\BA/t\cd\grd\BA\iso \grd A$
can be lifted to
a graded topological $\C[[t]]$-algebra isomorphism:
$\,
\grd\BA= (\grd A)[[t]]\,.$

Then (\ref{ass}) holds.\quad\qed
\end{lemma}
\medskip

Given  an algebra $A$  and
 an idempotent  $e\in A$, we say that the pair $(A,e)$ has 
Morita property if the map: $A\otimes A\to A$ given by 
$x\otimes y\to x\cdot e\cdot y$ is surjective. 
A similar definition applies 
for a complete topologically free  $\C[[t]]$-algebra $\A$ in 
which case the map: $\A\widehat{\otimes}_{_{\C[[t]]}}\A\to\A$ involves
the {\it completed} tensor product over $\C[[t]]$.

Let $\A$ be a complete topologically free  $\C[[t]]$-algebra,
and $A=\A/t\cd\A$.

\begin{proposition}\label{morita_e}
Let ${\mathbf{e}}\in\A$ be an idempotent in $\A$,
and $e\in A$ its reduction modulo $t$.
If $(A,e)$ has Morita property, then so does $(\A,{\mathbf{e}}).$
\end{proposition}

\Pf . This is an immediate consequence of the following general claim:
$\;${\it
Let ${\mathbf{f}}
:{\mathbf{V}}\to {\mathbf{W}}$  be a continuous  homomorphism of complete
topologically free C[[t]]-modules. If the induced
map  $f:{\mathbf{V}}/t\cd {\mathbf{V}}\to {\mathbf{W}}/t\cd 
{\mathbf{W}}$ is surjective,
then ${\mathbf{f}}$ is surjective as well.}

To prove the Claim, 
 we can assume that ${\mathbf{V}}=V[[t]], {\mathbf{W}}=W[[t]]$, 
and ${\mathbf{f}}=\sum_{j\ge 0}\,t^j\cdot f_j,$ where 
$f_j:V\to W.$ Suppose we want to solve the equation 
${\mathbf{f}}(\sum\, t^j\cdot v_j)=\sum\, t^i\cdot w_i,$ where $w_i$ are given. 
Suppose we've solved it up to order $n-1$, let's solve in order $n.$ 
We have $f_0(v_n)+ {\textsl{known things}}  =w_n.$ 
Since $f_0$ is surjective, we can always solve 
it. \quad\qed

\section{Some examples}
\setcounter{equation}{0}
Let $\Gamma$ be a finite subgroup of $\SL_2$. 
Consider the wreath product $\ga  =S_n\ltimes \Gamma^n$,
and use the notation as at the beginning of \S9.
We will consider the
symplectic reflection algebra $\hh_{0,k,c}(\ga)$,
 for generic $k,c$. 
Our first result in this section is 

\begin{theorem}\label{wreath00} Any simple $\hh_{0,k,c}$-module
 restricts to the regular representation of $\ga$. 
\end{theorem}

This theorem and proposition \ref{criterion} imply the following corollary:

\begin{corollary}\label{wreath1}
For wreath-products, the Calogero-Moser 
space $\Spec \ZZ_{0,k,c}$ is smooth.\qed
\end{corollary}

Before proving the theorem, we prove the following lemma.

\begin{lemma}\label{wreath0} Let $E$ be a finite dimensional representation 
of $\hh_{0,k,c}$. Then 
 as a $\ga$-module, $E$ has the form 
${\mathtt{Ind}}_{S_n}^\ga  L$, where $L$ is a unique $S_n$-module. 
\end{lemma}

\noindent
{\bf Proof of Lemma \ref{wreath0}.} 
To show that $E$ is induced, it suffices to show that $\text{tr}_E(g)=0$ 
if $g\in \ga$ is an element not conjugate to an element of $S_n$. 
We may assume that $g=(g_1, g_2,\ldots, g_m)$, where $g_i\in 
S_{n_i}\ltimes \Gamma^{n_i}$, with $\sum_{i=1}^m n_i=n$, and 
where $g_1=\sigma\cdot(\alpha,1,...,1)$, for some cycle $\sigma$,
and $\alpha\in \Gamma,$  a nontrivial element.
 
Let $a,b\in \C^2$, and consider the commutator  
$[(a,0,...,0),(b,b,\ldots,b, 0,...,0)]$ in $\hh_{0,k,c}$
(here and below,  the number of copies of $b$ in the element
$(b,b,\ldots,b, 0,...,0)\in V=(\C^2)^{\oplus n}\subset\hh_{0,k,c}$
equals $n_1=${\it length of the cycle} $\sigma$). According to
the main commutation relation, we have 
$$
[(a,0,...,0),(b,b,\ldots,b,0,...,0)]=(a,b)\cdot
\sum\nolimits_{\gamma\in \Gamma\smallsetminus 1}\;
c_\gamma\cdot (\gamma,1,...,1)+k\cdot f(a,b), 
$$
where $f(a,b)\in \C[\ga]$. 
Therefore, 
$$
[(\sigma,g_2,...,g_m)(a,0,...,0),(b,b,\ldots,b,0,...,0)]=
(a,b)\sum_{\gamma\in \Gamma\smallsetminus 1}
c_\gamma\cdot (\sigma\cdot(\gamma,1,...,1),g_2,...,g_n)+k\cdot\tilde f(a,b),
$$
where $\tilde f(a,b)\in \C[\ga]$. Taking the trace of both sides of this 
equation in $E$, and using the facts that $k$ and $c_\gamma$ are generic
and that traces of group elements are algebraic integers, 
we get 
$$
\text{tr}_{_E}(\sigma\cdot(\gamma,1,...,1),g_2,...,g_n)=0.
$$
The lemma is proved. 
\qed\medskip

\noindent
{\bf Proof of Theorem \ref{wreath00}.} 
In the previous lemma we showed that, for any $\gamma\in
\ga,$ the trace of $\gamma$ in $E$ vanishes
unless $\gamma$ is not conjugate to an element of $S_n\subset
\ga$.
Thus, it remains only
to show   that the trace of any nontrivial permutation $\sigma$ in $E$
vanishes.

Let us take $y=(b,b,...,b)$, and
$x=(0,...,0,a,0,...,0)$, $a,b\in \C^2$,
where $a$ stands in the $i$-th place. For any permutation 
$\sigma$ with $m$ cycles, let us consider 
the element $\sigma[x,y]=[\sigma x,y]$ of $\hh_{0,k,c}$. 
The trace of this element in $E$ is obviously zero. 
Thus, writing this trace explicitly using the main
commutation relation $\bigl(\ref{P-bracket}(iii)\bigr)$, 
and taking into account that $k,c$ are generic, we get
$$
\sum\nolimits_{\{j\;|\;j\ne i\}}\;\sum\nolimits_{\gamma\in
  \Gamma}\;\Tr(\sigma  \cdot s_{ij}\cdot \gamma_i\cdot \gamma_j^{-1},\,E)\cdot
\omega_{s_{ij}\cdot \gamma_i\cdot \gamma_j^{-1}}(x,y)=0
$$

The last sum consists of two parts:
the terms with $m+1$ cycles and with $m-1$ cycles. 
The first part is the sum over $j$
belonging to the same cycle of $\sigma  $ as $i$. In this sum, the
element $\sigma  \cdot s_{ij}$ has $i,j$ in different cycles, so 
the term containing 
the trace of $\sigma  \cdot s_{ij}\cdot \gamma_i\cdot \gamma_j^{-1}$ is
zero. 
Indeed, for $\gamma=1$ this term comes with a zero coefficient, and 
for $\gamma\ne 1$ the element has nontrivial cycle monodromies
(= conjugacy class of the product of the elements of 
$\Gamma$ in a cycle),
 hence, it cannot be conjugated into $S_n$. Thus, its trace is zero by 
Lemma \ref{wreath0}. Thus, only the second part of the sum, over 
$j$ belonging to the cycles of $\sigma  $ other than the one that contains
$i$, remains. 
In this case, $i$ and $j$ are in the same cycle of
$\sigma  \cdot s_{ij}$, so the trace is $\gamma$-independent, and 
after computing $\omega_{s_{ij}\cdot \gamma_i\cdot \gamma_j^{-1}}$ 
and simplifications we have 
$$
\sum\nolimits_{\{j\;|\;i\notin <\sigma  >j\}}\;\Tr(\sigma  \cdot s_{ij},\,E)=0
$$
The last equality involves only permutations with a fixed number
of cycles, namely\linebreak
 $m-1$. So the terms of the sum correspond to
partitions of $n$ into $m-1$ parts. 
In terms of partitions, the latter sum can be expressed in a very
simple way. Namely, let the length of the cycle of $\sigma  $
containing $i$ be $l$. Then the other cycles define a partition $\mu=
(\mu_1,\mu_2,...,\mu_{m-1})$, $\mu_1\ge \mu_2\ge \ldots$,  
of $n-l$. The last equation has the form 
$$
\sum\nolimits_{p=1}^{m-1}\; \Tr(\mu_1,...,\mu_p+l,...,\mu_{m-1})=0
$$
where trace of a partition is by definition the trace in $E$
of the
corresponding permutation (to view 
$\,(\mu_1,...,\mu_p+l,...,\mu_{m-1})\,$ as a partition one may 
rearrange the elements involved
in the decreasing order if this 
happens not to be so).

Now, let $\lambda$ be a partition of $n$ into $n-m-1$ parts,
$\lambda=(\lambda_1,\ldots,\lambda_{m-1})$,
$\lambda_1\ge \lambda_2\ge \ldots$. Let 
$p$ be the largest index for which $\lambda_p\ne 1$. Let 
$l=\lambda_p-1$, and $\mu$ be the partition
$(\lambda_1,...,\lambda_p-l=1,1,1,\ldots,1)$
of $n-l$. Then the above equation for this $\mu$ allows one to
express $\Tr(\lambda)$ through traces of partitions which 
are bigger with respect to the standard ordering than $\lambda$
(i.e. obtainable from $\lambda$ by increasing $\lambda_q$ and 
decreasing $\lambda_r$ with $q<r$). This implies that
$\Tr(\lambda)=0$, 
as desired. \quad
\qed\medskip

\begin{proposition}\label{smooth_A}\vi If $W$ is the Weyl group
of type ${\mathbf{A}}$ or of type 
${\mathbf{B}}={\mathbf{C}}$, then the variety
 $\Spec{\ZZ(\hh_{0,c})}$ is smooth, for generic $c$.

\vii If $W$ is the Weyl group
of type  ${\mathbf{G_2}}$, then the variety
 $\Spec{\ZZ(\hh_{0,c})}$ is never smooth.
\end{proposition}

\Pf . Part (i) is a special case of Corollary \ref{wreath1}
with $\G=\{1\}$ for type ${\mathbf{A}}$,
and $\G=\Z/2\Z$, for type ${\mathbf{B}}$.
To prove (ii), by
Proposition
\ref{criterion}(i), it suffices to produce,
for Weyl group $W$
of type ${\mathbf{G_2}}$, a simple $\hh_{0,c}$-module
of dimension less than $|W|$. 
We have $\dim\CC=2$,
hence, $c=(c_1,c_2),$ where $c_1\in \C$ and  $c_2\in \C$
correspond to long 
and short roots, respectively.

Recall the structure and representation theory of $W=W({\mathbf{G_2}})$. 
The group $W$ is the group of symmeries of the regular hexagon, 
so it contains rotations $\,\{\ro^j\}_{j=0,...,5},\,$ where $\ro$ is the 
rotation by $\pi/3$, and reflections: $s_0,...,s_5$, where 
$s_j$ is the reflection with respect to the line making angle 
$\pi j/6$ with the $x$-axis. The element $\ro^3$ 
(the central symmetry) is central in 
$W$, and $W=\Z/2\times S_3$, where 
$\Z/2=\lbrace{1,\ro^3\rbrace}$, 
and $S_3=\lbrace{1,\ro^2,\ro^4,s_0,s_2,s_4\rbrace}$.

Therefore, the group $W$ has 
 six irreducible representations:
$\,T_+,T_-,S_+,S_-,D_+,D_-,\,$ where $T,S,D$ are the 
trivial, sign, and 2-dimensional representations 
of $S_3$, and $\pm$ denotes the action of $\ro^3$. 
In particular, $\h=D_-$, and the sign representation of $W$ is
$S_+$. We note that $|W|=12$, and
 consider the following two {\it reducible} 6-dimensional
representations 
$E_+,E_-$ of $W$:
$E_+=T_+\oplus S_+\oplus 2D_-$, and 
$E_-=T_-\oplus S_-\oplus 2D_+$. 
Clearly, $E_+\oplus E_-=\C{W}$. 

 The proof of part (ii) of
Proposition \ref{smooth_A} is now completed by the following result.

\begin{proposition} Let $W$ be  of type  ${\mathbf{G_2}}$ and 
assume $c\in\CC$ is generic enough. Then

\vi Any finite dimensional representation of $\hh_{0,c}$
is a linear combination of $E_+$ and $E_-$ as a $W$-module. 

\vii Each of the $W$-modules $E_+$ and $E_-$ admits a unique extension 
to an $\hh_{0,c}$-module; in particular, $E_+\,,\,E_-$
are $\hh_{0,c}$-modules of dimension $\frac{1}{2}|W|$.
\end{proposition}

\Pf .
We first prove (i). Let $E$ be any finite dimensional $\hh_{0,c}$-module. 
It is clear from $W$-invariance of the
symplectic form $\om$ on $\h^*\oplus\h$,
that the forms
$\sum_{j=0}^2\omega_{s_{2j}}(x,y)$, and 
$\sum_{j=0}^2\omega_{s_{2j+1}}(x,y)$ on $\h^*\oplus \h$ are 
both nonzero multiples of the standard pairing.
Taking the trace in the $\hh_{0,c}$-module
$E$ of both sides of the main commutation
relation (\ref{P-bracket}iii) for the algebra $\hh_{0,c}$,
we get:
$$
[x,y]=c_1\cdot \sum\nolimits_{j=0}^2\omega_{s_{2j}}(x,y)\cdot s_{2j}+
 c_2\cdot \sum\nolimits_{j=0}^2\omega_{s_{2j+1}}(x,y)\cdot s_{2j+1}.
$$
Recall that $c_1,c_2$ are generic, all the
$s_{2j}$ are conjugate, and all the
$s_{2j+1}$ are conjugate. Thus, we get 
$
\Tr(s_i,E)=0.
$
Next, we consider the commutator 
$[s_0\cdot x, y(e_1)]$, where $y(e_1)\in\h$ corresponds
to the vector $e_1=(1,0)\in\C^2$. Then 
$\omega_{s_0}(x,e_1)=0$, and we get 
\begin{align*}
[s_0\cdot x,e_1]=s_0\cdot [x,e_1]&=c_1\cdot 
\sum\nolimits_{j=1}^2\omega_{s_{2j}}(x,e_1)\cdot s_0\cdot s_{2j}+
 c_2\cdot \sum\nolimits_{j=0}^2\omega_{s_{2j+1}}(x,e_1)\cdot s_0\cdot 
s_{2j+1}
\\
&=c_1\cdot \sum\nolimits_{j=1}^2\omega_{s_{2j}}(x,e_1)\cdot \ro^{2j}+
 c_2\cdot \sum\nolimits_{j=0}^2\omega_{s_{2j+1}}(x,e_1)\cdot \ro^{2j+1}.
\end{align*}
Let us now take the trace in $E$ of both sides. 
Recall that $\ro^2$ is conjugate to $\ro^4$ and $\ro$ to $\ro^5$. Thus,
using again that $c_1,c_2$ are generic, we get 
$$
\Tr(\ro^2,E)=0\quad,\quad \Tr(\ro^3,E)+\frac{1}{2}\Tr(\ro,E)=0. 
$$
It is easy to check that these equations are equivalent to saying that 
$E=n_+E_+\oplus n_-E_-$. Statement (i) is proved.

Now we prove (ii), and
introduce the structure of an $\hh_{0,c}$-module on $E_+$. 
We have $E_+=T_+\oplus S_+\oplus L\otimes D_-$, where $L$ is a 2-dimensional 
space. An action of $\hh_{0,c}$ is defined by  $W$-equivariant 
linear maps $f_1:\h\to \End(E_+)$ and $f_2: \h^*\to \End(E_+)$. 

Let $I:\h\to \h$ be the operator of rotation by $\pi/2$. 
Identifying $\h$ and $\h^*$ with $D_-$, and taking into account 
the $W$-equivariance of $f_j$, we find that 
$$
f_j(v)=
\begin{pmatrix} 0&0& a_j\otimes (v,*)\\ 0&0& b_j\otimes (Iv,*)\\
\alpha_j\otimes v& \beta_j\otimes Iv& 0
\end{pmatrix},
$$
where $a_j,b_j\in L^*$, $\alpha_j,\beta_j\in L$.
Now consider the relation 
$[f_j(v),f_j(w)]=0$. 
It is easy to calculate that 
this relation is equivalent to the conditions $a_j\otimes \alpha_j+
b_j\otimes \beta_j=0$, 
and $(a_j,\beta_j)=(b_j,\alpha_j)=0$. It is not hard to show that 
$\alpha_j,\beta_j,a_j,b_j$ must be nonzero. Thus, by changing scale in 
$S_+$ we can make sure that $b_j=\sqrt{-1}\cdot
a_j\,,\,\beta_j=\sqrt{-1}\cdot\alpha_j$, and we
have: 
$(a_j,\alpha_j)=0$.  

Next consider the main relation for $[f_1(v),f_2(w)]$. 
After a straightforward calculation, 
this relation turns out to be equivalent 
to an additional requirement that 
$a_1\alpha_2+a_2\alpha_1=0$, and 
$a_1\alpha_2-a_2\alpha_1=k$ (where $k$ is a constant multiple of $c_1+c_2$). 
There is obviously a unique solution to these equations, up 
to linear transformations of $L$:
namely, $a_1,a_2$ is the dual basis to 
$\,\{2k^{-1}\alpha_2\,,\,-2k^{-1}\alpha_1\}$.

The case of $E_-$ is treated similarly. In this case, $c_1+c_2$ is replaced 
with $c_1-c_2$. Part (ii) is proved.\quad\qed\medskip 

\noindent
{\bf Remarks}\quad\vi 
It is easy to check that in $E_\pm$ the right hand side of the 
main relation depends only on $c_1\pm c_2$ since the sums which 
are multiplied by $c_1,c_2$ are the same in $E_+$ and differ by sign in $E_-$. 
Therefore, if ${{\sf{f}}}_+,{{\sf{f}}}_-
\in\C{W}\,$ are the idempotents corresponding to 
$E_+,E_-$ then one has:
${{\sf{f}}}_\pm\cdot \hh_{0,c_1,c_2}({{\mathbf{G_2}}})\cdot  {{\sf{f}}}_\pm=
\C[\Z/2]\ltimes \bigl({{\sf{f}}}_\pm \cdot 
\hh_{0,c_1\pm c_2}({{\mathbf{A_2}}})\cdot {{\sf{f}}}_\pm\bigr)$. 
Thus in effect we are taking a 6-dimensional representation 
of $\hh_{0,l}({{\mathbf{A_2}}})$ (regular as $S_3$-module, as it must be), 
which is compatible with the action of the central element of 
$W$. The claim we proved is that such a module exists 
and is unique. 

\vii Similar result holds more generally for $W=W_n$ being
the Coxeter group of the
symmetries of a regular $n$-gon, whenever $n\geq 5$
(with the ${{\mathbf{G_2}}}$-case corresponding to $n=6$).
Specifically, the algebra $\hh_{0,c}(W_n)$ always has
a simple module whose restriction to $\C{W_n}$ is isomorphic
to a direct sum of the trivial representation, the
sign representation, and of two copies of the 
2-dimensional reflection-representation. Thus,
$\hh_{0,c}(W_n)$ has a 6-dimensional module,
hence the variety $\Spec\ZZ_{0,c}(W_n)$ is not
smooth, if $n\geq 5$.

\viii Note  that our result  on 
singularity of $\Spec \ZZ_{0,c}({{\mathbf{G_2}}})$
agrees perfectly with
a similar result for the  ${{\mathbf{G_2}}}$-type
Hilbert scheme, obtained in [Ka2].$\quad\lozenge$

\section{Open questions}
\setcounter{equation}{0}

Throughout this section $(V,\om)$ is a finite dimensional
symplectic vector space over $\C$, and  $\G\subset Sp(V)$
is a finite group generated by symplectic
reflections
and, moreover, the triple 
$(V,\om,\G)$ is assumed to be indecomposable.
 We will identify $V^*$ with $V$ without further notice.
Recall that $\ZZ_{0,c}:=\ZZ(\hh_{0,c})$. We have, cf. \S16:
\medskip

\noindent
{\bf Question 17.1}\quad {\it
For which groups $\G$
is the algebraic variety $\Spec \ZZ_{0,c}$  smooth, at least
for generic values of $c\in \CC$ ?
}\medskip

Recall that if $\Spec \ZZ_{0,c}$ is smooth then the Poisson structure
on $\ZZ_{0,c}$ is induced by a symplectic form, by Lemma \ref{normal}(ii).
\medskip

\noindent
{\bf Question 17.2}\quad {\it Assume that 
$c\in \CC$ is such that $\Spec \ZZ_{0,c}$ is smooth. Is there
a natural hyper-K\"ahler structure on $\Spec \ZZ_{0,c}$
compatible with the holomorphic symplectic 2-form ?
}\medskip

\noindent
{\bf Question 17.3}\quad {\it For which finite groups
$\G\subset Sp(V)$ the dimensions of the
Poisson cohomology groups of the Poisson algebra $\C[V]^\G$
are given by:
$\dim H\!P^{2i}(\C[V]^\G) = \bn(i),\,$ and
$\,\,\dim H\!P^{^{{\sf{odd}}}}(\C[V]^\G) = 0$} ?
\medskip

An argument analogous to our proof of Theorem \ref{deform_t}
yields the following
\medskip

\noindent
{\bf Proposition 17.4}\quad {\it  If the answer to
question 17.3 is affirmative for $\Gamma$, then
the
family $\,\{\ZZ_{0,c}\}_{c\in\CC}\,$ gives a universal
deformation of $(SV)^\G$ in the class of Poisson algebras.
}\quad\qed\medskip

\noindent
{\bf Problem 17.5}\quad {\it Study  representation theory of the
algebra $\hh_{t,c}$. Specifically,
describe the set of those $(t,c) \in \C\oplus\CC$ such that:

\vi The algebra $\hh_{t,c}$ is {\sf not simple}, and/or has
modules whose Gelfand-Kirillov dimension is 
less than $\frac{1}{2}\dim V$.

\vii The algebra $\hh_{t,c}$ has finite-dimensional modules.}
\medskip

\noindent
For rational Cherednik algebras, problem 17.5(ii) has been studied 
by Cherednik (see e.g. [ChM] for the ${\mathfrak{s}\mathfrak{l}}_2$-case).
\medskip

Let ${\sf{Q}}$ denote the skew field of fractions
of the Weyl algebra $\Weyl_1(V)$.
The group $\G$ acts naturally
on ${\sf{Q}}$ by algebra automorphisms.
We have the following analogue of the Gelfand-Kirillov
conjecture 
\medskip

\noindent
{\bf Conjecture 17.6}\quad {\it For any $(t,c)
\in \C^*\oplus\CC$, there is an isomorphism between  ${\sf{Q}}\#\G$
and the Goldie ring
of fractions of the algebra
$\hh_{t,c}$, which is identical on the subalgebra $\C\G$.}
\medskip

This would imply the following result

\noindent
{\bf Proposition 17.6*}\quad {\it If Conjecture 17.6
holds then, for any $(t,c)
\in \C^*\oplus\CC$, the  skew field
of fractions of the algebra
$\e\hh_{t,c}\e$ is isomorphic to ${\sf{Q}}^\G$.}
\medskip

A `quasi-classical' analogue of Proposition 17.6* is
\medskip

\noindent
{\bf Conjecture 17.7}\quad {\it 
 For any $c\in \CC$, the variety
$\Spec \ZZ_{0,c}$ is birationally isomorphic
to $V/\G$ as a Poisson variety.}
\medskip

We remark that the existence of a (faithful)
Dunkl-type representation,
see Proposition \ref{theta2}, for the rational
Cherednik algebra $\hh_{0,c}$ associated with any complex
reflection group $\G$ yields (without assuming Conjecture 17.6):
\medskip

\noindent
{\bf Proposition 17.7*} \quad {\it Proposition
17.6* and Conjecture 17.7
 hold true in the case where $\G$ is a complex
reflection group in $\h$, and $V=\h\oplus\h^*$.}\quad\qed
\medskip

\noindent
{\bf Question  17.8}\quad {\it Does the
isomorphism:
$\HH^{\bullet}(\hh_{t,c})\simeq \grd^F_\bullet(\ZZ\G)$
of Theorem {\ref{mckay}\vi}
 hold for all $(t,c)\in \C^*\times\CC$
without exception ?} 
\medskip 

\noindent
{\bf Problem 17.9}\quad {\it Prove that the equality in 
Theorem \ref{inject}\vii   holds for all
$k\in\C$ without exception.} 
\medskip

Recall the setup of \S5, esp. Definition \ref{LLL}.
Let $P\in \C\langle x,y\rangle
$ be a noncommutative polynomial in two variables.
Call $P$ $\,${\it central}$\,$ if, for any $n$, the element $\sum_i P(x_i,y_i)$
is central in $\hh_{0,c}$,
where the algebra $\hh_{0,c}$ is associated with   $\g=\gln$.
\medskip

\noindent
{\bf Question 17.10}\quad {\it 
Find all central polynomials.}
\medskip 

The above arguments show that the space of central polynomials
is an $\LL$-submodule of $\C\langle x,y\rangle$ which contains $x^k$ and $y^k$
for all $k$. So denoting by $E$ the $\LL$-module generated by
$x^k$ and $y^k$, we might ask if any central polynomial belongs to $E$.
In particular, any element of $E$ is palindromic (i.e. invariant under
the antiinvolution of $\C\langle x,y\rangle$ which fixes $x$ and $y$),
 so one might
ask whether every central polynomial is palindromic.
\medskip 

Recall that according to Proposition \ref{k_theory}, we have a
canonical isomorphism: $K(\G)$
$\simeq K(\Spec\ZZ_{0,c}).$
Further, Theorem \ref{mckay} says: 
$H^{\bullet}(\Spec\ZZ_{0,c}) \simeq \grd^F_\bullet(\ZZ\G)$.
Therefore, the Chern character map
$\,ch: \C\otimes_{_\Z}K(\Spec\ZZ_{0,c})\iso 
H^{\bullet}(\Spec\ZZ_{0,c})\,$
may be viewed as an isomorphism
$\,ch: \C\otimes_{_\Z}K(\G)\iso \grd^F_\bullet(\ZZ\G)$.
On the other hand, identifying $\ZZ\G$ with the space of class functions
on $\G$, and associating to any $\G$-module $E$ its character
$\Tr_{_E}$, gives a vector space isomorphism:
$\C\otimes_{_\Z}K(\G)\iso \ZZ\G\,,\,E\mapsto \Tr_{_E}$.
\medskip

\noindent
{\bf Problem 17.11}\;\;\; {\it Compute the composite
map} $\,ch\ccirc\Tr^{-1}: \ZZ\G \iso \C\otimes_{_\Z}K(\G)\to
\grd^F_\bullet(\ZZ\G)$.
\medskip 

Let $\hh_{0,c}$ be the Cherednik algebra associated to 
a root system $R$ in $\h^*$. Recall the finite map $\Upsilon:
\Spec\ZZ_{0,c}\too \h^*/W\times\h/W$, see  Proposition \ref{nilp_rep}.
Let
$\chi\in \Upsilon^{-1}(0)$ be a point in the zero-fiber,
 viewed as a maximal ideal $\chi\subset\ZZ_{0,c},$ 
and $(\chi)\subset \hh_{0,c}$ the two-sided ideal
in $\hh_{0,c}$ generated by $\chi$.
 We assume that $\Spec\ZZ_{0,c}$
is smooth, so that one has an algebra
isomorphism: $\hh_{0,c}/(\chi) \simeq \End_{_\C}\C\G.$
It follows readily from this isomorphism that the
space $\,\e\cdot\bigl(\hh_{0,c}/(\chi)\bigr)\cdot
\e_{_{\boldsymbol{-}}}\,\subset\,\hh_{0,c}/(\chi)\,$ is 1-dimensional.
Let ${\boldsymbol{\varepsilon}}_\chi$
denote a base vector in $\e\cdot\bigl(\hh_{0,c}/(\chi)\bigr)\cdot
\e_{_{\boldsymbol{-}}}$.

Now, the action of two commutative subalgebras:
$\C[\h^*]\,,\, \C[\h]\subset \hh_{0,c}$
generate two subspaces: $\C[\h^*]\cdot
{\boldsymbol{\varepsilon}}_\chi\subset \hh_{0,c}/(\chi),$
and $
{\boldsymbol{\varepsilon}}_\chi\cdot\C[\h]\subset \hh_{0,c}/(\chi)$,
which are stable under left, resp. right, $\C{W}$-action
on $\hh_{0,c}/(\chi)$ by multiplication.
\medskip

\noindent
{\bf Conjecture 17.12}\quad  {\it Assume $R$ is the root system of type
${\mathbf{A_{n-1}}}$. Then,
for each $\chi\in \Upsilon^{-1}(0):$

\vi There is a unique (up to $W$-conjugacy) pair of subsets
 $S_1, S_2\subset R$, such that
one has $W$-module isomorphisms:
$$\C[\h^*]\cdot
{\boldsymbol{\varepsilon}}_\chi\simeq \Ind_{W(S_1)}^W{\mathbf{1}}
\quad,\quad
{\boldsymbol{\varepsilon}}_\chi\cdot\C[\h]
\simeq \Ind_{W(S_2)}^W{\mathsf{sign}},
$$
where $W(S_i)\subset W$ denotes the
Weyl subgroup generated by $S_i$. Moreover, each of the sets $S_i$
(but not $S_1$ and $S_2$ simultaneously) is a subset of a 
basis of simple roots in $R$.

\vii The two $W$-modules in \vi have a single non-zero irreducible
constituent, $L_\chi\in$
$\irreps(W)$, in common} (see [Gi, 4.16, 4.17]).
\medskip

Further, for $\g=\sln$, let $\,\phi: \Spec\ZZ_{0,c}\iso \mm_n$ be the
isomorphism
of $\Spec\ZZ_{0,c}$ with the Calogero-Moser space. Let
$\OO_\chi\subset \sln\times\sln$ be the $\PGL_n$-conjugacy class
corresponding to a point $\chi\in\Spec\ZZ_{0,c}$ under that
isomorphism. Set $\OO_\chi^{(i)}\,,\, i=1,2,$ for
the image of $\OO_\chi$ under the corresponding projection
$\pr_i:\sln\times\sln\to \sln$. It is clear that,
for any $\chi\in \Upsilon^{-1}(0)$, the set $\OO_\chi^{(i)}$
is a nilpotent conjugacy class in $\sln$.
\medskip

\noindent
{\bf Conjecture 17.13} \quad 
{\it The simple $W$-module  $L_\chi$ introduced in Conjecture 17.12
corresponds to the nilpotent orbit $\OO_\chi^{(1)}\subset \sln$
via the Springer correspondence. Similarly,
the module ${\mathsf{sign}}\otimes L_\chi$ 
corresponds to the nilpotent orbit $\OO_\chi^{(2)}\subset \sln$
via the Springer correspondence.
\medskip

\noindent
{\bf Conjecture 17.14}\quad 
{\it The scheme-theoretic
multiplicity of $\Upsilon^{-1}(0)$ at the point $\chi$
equals $(\dim L_\chi)^2$.}
\medskip

\noindent
{\bf Problem 17.15} \quad 
{\it For any  Weyl group $W$, and $c\in\CC$
 such that $\Spec\ZZ_{0,c}$ is smooth,
find a natural bijection:} $\Upsilon^{-1}(0)
\;\stackrel{_\sim}\longleftrightarrow \; {\sf{Irreps}}(W)\,,\,
\chi \mapsto L_\chi$, such that the scheme-theoretic
multiplicity of $\Upsilon^{-1}(0)$ at any point $\chi\in
\Upsilon^{-1}(0)$
equals $(\dim L_\chi)^2$.}
\medskip

In the setup of Theorem \ref{Gamma_2},
we expect the following  $\G$-analogue of 
 \cite{LS1} to hold
\medskip

\noindent
{\bf Conjecture 17.16}\quad 
{\it 
Let $u\in \dd(\RQ_n)^{^{\PG_{_{\Gamma\!,n}}}}$ be a differential 
operator  which annihilates all
$\PG_{_{\Gamma\!,n}}$-invariant rational functions on  $\RQ_n$. Then 
$u\in  \bigl(\dd(\RQ_n)\cdot
(\ad-\chi_{_{c'}})(\pggn)\bigr)^{\PG_{_{\Gamma\!,n}}}$.
}

\section{Appendix G: Results of Kostant and shift operators}
\setcounter{equation}{0}
The notion of  shift operator is quite important in the theory
of Calogero-Moser systems.
The goal of this Appendix is to place this notion
 in the framework of Lie theory.

Fix $\g$, a semisimple Lie algebra
with
a triangular decomposition $\g=\frakn_-\oplus\h\oplus\frakn$.
Let $\prh:
S\g \onto S\g/S\g\cd(\frakn+\frakn_-)=S\h$ be the corresponding
projection. Set  ${\mathfrak{b}}=\h+\frakn$,
the fixed Borel subalgebra in $\g$.
Write $\rho$
for the half-sum of positive roots, and
 ${\check{\delta}}\in S\h$ for the product of positive coroots.
Let $L_\mu$ denote   the irreducible
$\g$-module with highest weight~$\mu$.

Below $\g=\sln$. Let $\hh_{1,k}$ be the rational Cherednik algebra associated with
the root system of $\sln$.
We view the family $\,\{\hh_{1,k}\}_{k\in\C}\,$ as a single
flat
$\C[k]$-algebra.
Also, we
view ${\check{\delta}}$ as an element ${\check{\delta}}=\delta_y\in
S\h\subset \hh_{1,k}$.

Recall the symmetriser idempotent $\e$ and
the "anti-symmetriser idempotent'':
$\e_{\pmb{-}}=\sum_{w\in W}\; (-1)^{\ell(w)}\cdot w.\,$
The Dunkl representation $\Theta_k:
\hh_k \to \bigl(\dd(\hreg)\#{W}\bigr)[k],\,$ see \S4,
restricts to a map
$$\Theta_k:\;\,\e_{_{\pmb{-}}}\cd \hh_k\cd \e \too 
\e_{_{\pmb{-}}}\cd \bigl(\dd(\hreg)\#{W}\bigr)\cd \e \otimes_{_\C} \C[k]\,\simeq\,
\dd(\hreg)^{^{\sf{sign}}}[k]\,,
$$
where $\dd(\hreg)^{^{\sf{sign}}}$ stands for the subspace
of $\dd(\hreg)$ formed by the $W$-anti-symmetric
 differential
operators.
 The operator $\,\SS:=\Theta_k(\e_{_{\pmb{-}}}\cd\delta_y\cd\e)
\in \dd(\hreg)^{^{\sf{sign}}}[k]\,$
is called the shift operator (this definition differs from the
original one, and is due to Heckman).
It is easy to check that $\SS$ is a differential
operator of order $\,m=\mbox{\it number of positive roots},\,$
with the principal symbol $\sigma_m(\SS)={\check{\delta}}$.
Furthermore, it is known
from  work of Opdam [O] that,
for $\g={\frak {sl}}_n$, we have:
$\SS\ccirc \CM_k=\CM_{k+1}\ccirc\SS$. Moreover,
there exists an isomorphism $\tau_k:\cc_k\iso \cc_{k+1}$
such that for all $L\in \cc_k$ one has $\SS\ccirc L=\tau_k(L)\ccirc\SS$.
\medskip

\noindent
{\bf{Example.}}\quad
For $n=2$, we have: $\SS=\partial -\frac{k+1}{x}$.
\medskip

\noindent
{\bf Remark.} In fact, Opdam considers the trigonometric
Calogero-Moser system, while we are considering its rational
degeneration. Let us point out that the results of Opdam
are valid in the rational limit, while our construction below
easily generalizes to the trigonometric case by replacing
the Lie algebra with the corresponding group.$\quad\lozenge$ \medskip

We  are going to give a Lie theoretic construction of the shift operator.

Write $\C^\circ\![x_1^{\pm 1},
\ldots, x_n^{\pm 1}]\subset \C[x_1^{\pm 1},
\ldots, x_n^{\pm 1}]$ for the subspace of Laurent polynomials of
total degree zero.
Recall the family
$\,V_k=(x_1\cdot\ldots\cdot x_n)^k\cdot\C^\circ\![x_1^{\pm 1},
\ldots, x_n^{\pm 1}]\,,\,k\in\C,\,$  of $\sln$-representations,
see the beginning of \S7.
Let $V$ be the
nonzero finite dimensional subrepresentation of the representation
$V_1$, i.e., $V=S^n(\C^n)\,=\, \mbox{\it degree $n$ polynomials}\,.$
The  natural `multiplication' map:
$$
f\,\otimes\, (x_1\cdot\ldots\cdot x_n)^k\cd P
\;\mapsto\; (x_1\cdot\ldots\cdot x_n)^{k}\cd (f P)
\;=\;
(x_1\cdot\ldots\cdot x_n)^{k+1}\cdot
\bigl(\frac{f}{x_1\cdot\ldots\cdot x_n}P\bigr)
$$
gives a $\g$-module
homomorphism: $ V\otimes V_k\to V_{k+1}$.
This morphism induces a $\g$-equivariant pairing
$\,\nu: \dd(\gr,V)^\g \,\otimes\, \C[ \gr,V_k]^\g
\too   \C[ \gr,V_{k+1}]^\g.\,$
Let $\{v_i\}$ be a basis of $V$, and $\{v^i\}$ the dual basis of
$V^*$. Recall that since
$V_k\!\langle 0\rangle$ is 1-dimensional, the
restriction to $\hreg$ gives an isomorphism:
$\,
\C[ \gr,V_k]^\g\iso \C[\hreg,
V_k\!\langle 0\rangle]^W=\C[\hreg]^W.\,$

 View elements  of $S\g$ as
constant coefficient differential
operators on $\g$.

\begin{theorem}\label{EK}
\vi  For $\g=\sln$,
there exists a $\g$-module embedding  $\,j :V^*
\into S^{n(n-1)/2}\g,$
such that  $\prh\bigl(j(V^*)\bigr)$
is the line in $S\h$ spanned by ${\check{\delta}}$.

\vii   The assignment: $f\mapsto \sum j(v^i)f\otimes v_i$ defines
a $V$-valued invariant
differential operator
${\sf{\widehat{S}}}\in\dd(\g,V)^\g=\bigl(\dd(\g)\otimes V\bigr)^\g,\,$
such that,
 for any
$\g$-invariant function $f: \gr\to V_k$, one has:
$$
(\nu({\sf{\widehat{S}}}\otimes f))\big|_{_{\hreg}}=
\SS(f\big|_{_\hreg})\,.
$$
\end{theorem}

Before entering the proof of the Theorem, we are going to put it
into a more general framework of Kostant's results.
It is worth noticing
that Kostant's motivation had nothing to do with
 the theory of Calogero-Moser systems
and shift operators.
\medskip

\noindent
{\bf{Kostant results.}}\quad
We keep the notation as at the beginning of this
section. In particular,
let $\g=\frakn_-\oplus\h\oplus\frakn$
 be an arbitrary semisimple Lie algebra with a fixed
triangular decomposition, and ${\mathfrak{b}}=\h+\frakn$.
Set $\,m=\mbox{\it number of positive
roots},\,$
 so that $\,\dim\g = \rk^{\,}\g +2m$.
Write $\Lambda^\bullet\g$ and
$\Lambda^\bullet\g^*$ for the exterior
algebras over $\g$ and
$\g^*$, respectively.

Let ${\sf{Ab}}$ denote the set
parametrising all  the abelian  ideals  $\afrak\subset\bfrak$
 having dimension $\rk^{\,}\g$. For
each $\afrak\in {\sf{Ab}}$, the element $\Lambda^{\rk^{\,}\g}\afrak
\in \Lambda^{\rk^{\,}\g}\g$ is clearly an
$\ad^{\,}\bfrak$-weight vector, hence a highest weight vector
of a simple $\g$-submodule in $\Lambda^{\rk^{\,}\g}\g$.
The corresponding highest weight, $\mu_\afrak$,
equals the sum of all the
 positive   roots whose root
vectors form a basis of $\afrak$.
Further, the Hodge `star-operator' relative to an invariant form
on $\g$ gives a $\g$-equivariant isomorphism
$\,\star:\,\Lambda^{\rk^{\,}\g}\g\iso
\Lambda^{\dim\g-\rk^{\,}\g}\g=\Lambda^{2m}\g.\,$
We let $L_\afrak\subset \Lambda^{2m}\g\,$
denote the simple $\g$-submodule corresponding
to the $\g$-submodule in $\Lambda^{\rk^{\,}\g}\g$
generated by  $\Lambda^{\rk^{\,}\g}\afrak$.
Thus, $L_\afrak \simeq
L_{\mu_\afrak}$; moreover, all the modules $L_\afrak\,,
\,a\in {\sf{Ab}}$,
are known to be pairwise nonisomorphic.

Write
$\Lambda^{^{\!\sf{ev}\!}}\g^*\subset \Lambda^\bullet\g^*$
for the even part, which is
a commutative algebra.
Let $d: \g^*\too \Lambda^2\!\g^*$ denote the map
dual to the commutator map: $\g\wedge\g\stackrel{[,]}{\too}\g$.
The map $d$ can be uniquely extended
to a commutative algebra
homomorphism $s: S\g^* \to \Lambda^{^{\!\sf{ev}\!}}\g^*.$
Let $A(\g) =  \bigoplus_i\;  A^{2i}\g,$ be the image of $s$,
 naturally graded by even integers.

\begin{theorem}[\cite{Ko2},\cite{Ko3}]
\label{k_triv}\label{k_A}
\vi The maximal  $i$ such that $A^{2i}\g\neq 0$
 equals $m$.

\vii $\;A^{2m}\g\;=\; \bigoplus_{\afrak \in {\sf{Ab}}} \;L_\afrak.\;$

\viii  A simple $\g$-module occurs in  $L_\rho\otimes L_\rho$
with non-zero multiplicity
if and only if it occurs in $\Lambda^\bullet\g$. Thus,
$L_\afrak$ occurs in  $L_\rho\otimes L_\rho$
with non-zero multiplicity, for any
$\afrak \in {\sf{Ab}}$. \sq
\end{theorem}

\noindent
{\bf{Example: $\sln$-case.}}\quad In this case $m=n(n-1)/2.$
The root vectors in  any Young diagram occupying  the upper
right hand
corner of the Borel subalgebra of
the upper triangular matrices  spans an ideal in the
Borel subalgebra,
 where the  Young diagram progresses downward and to the
left of $E_{1n}$. If there are exactly $n-1$
 boxes then the ideal is abelian.
Given such a Young diagram,
 the sum of the roots corresponding to the boxes
of the diagram
form the corresponding dominant weight $\mu_\afrak.$
$\enspace\lozenge$\medskip

Following [Ko2], consider
 the  map
$\, s^\dag: \Lambda^{^{\!\sf{ev}\!}}\g \too S\g,$ transposed
to $s$ (this map is  related to the
Amitsur-Levitski standard identity, as explained in [Ko3]).
 Put: $\,Q(\g)=s^\dag(\Lambda^{2m}\g).$
By self-duality of $A(\g)$, one gets  a $\g$-module isomorphism:
$\,Q(\g)\simeq A^{2m}(\g)= \bigoplus_{\afrak \in {\sf{Ab}}}
\;L_\afrak.\;$
     From now on we will think of the  $L_\afrak$'s as
sitting inside $Q(\g)$ via this decomposition. Recall the notion
of harmonic elements in  $S\g$ introduced in [Ko1].

\pagebreak[3]
\begin{theorem}[Kostant \cite{Ko4}]\label{k_B}
\vi Any element of $\,Q(\g)\subset S^m\g\,$ is harmonic.

\vii  For any $\afrak\in {\sf{Ab}}$,
the space $\prh(L_\afrak)\subset S\h$
 is the 1-dimensional space\footnote{this result
has been first obtained (via a case by case argument)
by Ranee Brylinski in her PhD Thesis.} spanned by
${\check{\delta}}$.

\viii For any $\afrak\in {\sf{Ab}},$ and any integer $k=1,2,\ldots,$
there exists a $\g$-equivariant imbedding
$\,j: L_{k\cdot\mu(a)} \into S^{k\cdot m}\g$ such that
$\prh\bigl(j(L_{k\cdot\mu(a)})\bigr)
\subset S^{k\cdot m}\h$
 is the one dimensional space spanned by  ${\check{\delta}}^k,\,$
the $k$-th power of the
product of the positive roots.
\end{theorem}

\noindent
{\bf Proof (after Kostant).} $\;$
Recall the standard Koszul differential $d$ on
$\Lambda^\bullet\g^*$.
 Note that the algebra $A(\g)$ is generated by
exact 2-forms, relative to $d$. Hence, any element of
$A(\g)$ is exact. It follows that, for any
$\ad^{\,}\g$-invariant $u\in S\g^*$,
the element $s(u)$ is $\ad^{\,}\g$-invariant
and exact. But it is well-known that
any non-zero $\ad^{\,}\g$-invariant
in $\Lambda^\bullet\g^*$ represents a non-zero
cohomology class in $H^*(\g)$, so that
any exact $\ad^{\,}\g$-invariant must be equal to zero.
Thus,
$s$
 must vanish on the ideal $I \subset S\g$
 generated by the $\ad^{\,}\g$-invariants
of positive degree. Part (i) follows from this,
by duality.

  Let $\omega\in\Lambda^{2m}\g$
 be defined as the
wedge product of all the root vectors,
with respect to the triangular decomposition.
It is a zero weight vector which
transforms according to the sign representation of the Weyl group
since any simple reflection clearly takes $\omega$ to $-\omega$.
Identify $\g^*$ with $\g$ and $\h^*$ with $\h$.
It follows that the element
$\prh\bigl(s^\dag(\omega)\bigr)$ is an element
of $S^{m}\h$ that transforms,
under the Weyl group action, according to the
sign character of $W$.
But any such element is known to be proportional to ${\check{\delta}}$.
Thus, $\prh\bigl(s^\dag(\omega)\bigr)\in \C\cdot{\check{\delta}}.$

To prove (ii), take the orthocomplement $D$
in $\Lambda^{2m}\g$ of the one dimensional
space spanned by $\omega$ . Since $s(x^m)$ is a non-zero
 multiple of $\omega$, for any $x\in
\hreg$, it follows that $s^\dag(D)$, viewed as
an element of $S\g^*$, vanishes on
the elements of the form $x^m\,,\,x\in\hreg$.
Hence,
by the previous paragraph,
 $\prh\bigl(s^\dag(\Lambda^{2m}\g))$,
is contained in the one dimensional space spanned by
the product of positive coroots.
 But $L_\afrak$
is contained in $\Lambda^{2m}\g$.
Hence, for any  $a\in {\sf{Ab}}$, either
$\prh\bigl(s^\dag(L_\afrak)\bigr)$
vanishes or is   contained in
$\prh\bigl(s^\dag(L_\afrak)\bigr)\subset \C\cdot{\check{\delta}}.$
  We claim it cannot vanish.
 Indeed if $\prh\bigl(s^\dag(L_\afrak)\bigr)$
vanished
then since it is stable under the $\Ad^{\,}G$-action
 it follows that $
\Ad^{\,}g\Bigl(\prh\bigl(s^\dag(L_\afrak)\bigr)\Bigr)\,$
 vanishes, for all $g\in G$.
Identifying $S\g$ with $\C[\g]$ via an invariant form,
that would mean that elements of
$\prh\bigl(s^\dag(L_\afrak)\bigr)$
vanish  on all adjoint transforms of $\h$.
 But that would imply that they  vanish
identically on $\g$, which is impossible.

To prove (iii), view $L_\afrak$ as a vector subspace in
$S^m\g$. Clearly,  elements of the form $\,\{v^k\,,\,v\in
L_\afrak\subset S^m\g\}\,$ span an $\ad^{\,}\g$-stable subspace
in $S^{k\cdot m}\g$. Furthermore, let $v_\mu$ be a highest weight
vector for $L_\afrak$. Then $(v_\mu)^k$ is a highest weight
vector,
and the $\g$-module generated by this vector has
highest weight $k\cdot \mu_\afrak$. This gives an imbedding
$j_k: L_{k\cdot \mu_\afrak}\into S^{k\cdot m}\g$.
It remains to observe that, for any $v\in S^m\g$,
one has: $\prh(v^k) = \bigl(\prh(v)\bigr)^k.$
It follows, that $\prh(L_{k\cdot \mu_\afrak})$ is contained
in the linear span of $\,\{\prh(v^k)\,,\,v\in
L_\afrak\},\,$ which has been shown to be the line
spanned by ${\check{\delta}}^k$. Finally,
$\prh(L_{k\cdot \mu_\afrak})\neq 0$ since
$\prh\bigl((v_\mu)^k\bigr)\neq 0.$
\sq\medskip

Fix $a \in  {\sf{Ab}}$.
Let $\{v_i\}$ be a basis of $L_\afrak$, and $\{v^i\}$ the dual basis of
$L_\afrak^*$. Since, $Q(\g)$ is a self-dual
module, the above results hold for
$L_\afrak^*$ as well as for $L_\afrak$. In particular,
we have the canonical $\g$-module imbedding
$j: L_\afrak^*\into S^m\g$, as a harmonic subspace.
Put $S_\afrak:=\sum j(v^i)\otimes v_i$. This
is
a $\g$-invariant  $L_\afrak$-valued differential operator
on $\g$ with constant coefficients, of order $m$.
Moreover, part (ii) of Theorem \ref{k_B} insures
that for $\sigma_m(S_\afrak)$,
the principal symbol of $S_\afrak$, one has:
$\prh\bigl(\sigma_m(S_\afrak)\bigr)
={\check{\delta}}$.

For any integer $k\ge 0$,  one has a canonical
 $\g$-module projection:
$L_\afrak \otimes L_{k\cdot \mu_\afrak}
\onto L_{(k+1)\cdot \mu_\afrak}$. It induces
a natural
 pairing
$\,\dis
\nu: \dd(\gr, L_\afrak)^\g \otimes \C[ \gr,L_{k\cdot \mu_\afrak}]^\g
\too   \C[ \gr, L_{(k+1)\cdot \mu_\afrak}]^\g.\,$
Recall the  isomorphism:
$\,
\C[ \gr, L_{k\cdot \mu_\afrak}]^\g\iso
\C[\hreg, L_{k\cdot \mu_\afrak}\!\langle 0\rangle]^W,\,$
induced by
restriction to $\hreg$. Thus,
  the radial part of $S_\afrak$ gives a differential operator
\begin{equation}\label{S_a}
\Psi(S_\afrak) \in \dd\bigl(\hreg\,,\,
 \Hom_{_\C}(L_{k\cdot \mu_\afrak}\!\langle 0\rangle,
L_{(k+1)\cdot \mu_\afrak}\!\langle 0\rangle)\bigr)^W\,.
\end{equation}
\medskip

\noindent
{\bf{Proof of Theorem \ref{EK}.}}\quad
If the Young diagram
corresponding to  $a \in  {\sf{Ab}}$ is the top
row or the last column we get $
L_\afrak=
S^n(\C^n)(=V)$, or its dual. Thus, Theorem  \ref{k_triv}(iii),(vi)
 yields a canonical
imbedding $j: V\into S^{n(n-1)/2}\g$.
Part (i) of Theorem \ref{EK} now follows
from Theorem
\ref{k_B}(ii).

Further, since $V_k\!\langle 0\rangle$ is a 1-dimensional
space, the operator (\ref{S_a}) may be viewed in this case
as a scalar valued differential operator on $\hreg$.
The principal symbol of this differential operator
equals ${\check{\delta}}$, by part (i) of the theorem.

Moreover, we have:
$\Psi(S_\afrak)\ccirc \CM_k= \CM_{k+1}\ccirc \Psi(S_\afrak),\,$
because $S_\afrak$ transforms $G$-invariant functions
with values in $V_k$ into $G$-invariant functions with values
in $V_{k+1}$.

 Thus, the operator
$\Psi(S_\afrak)$  satisfies all
the properties of the shift operator.
Therefore, it must coincide with the shift operator, by the
uniqueness results, due to Opdam.
\sq\medskip

\noindent
{\bf{Remark.}}\quad An alternative proof of Theorem
\ref{EK} can be deduced from the results of
 \cite{EK}. Specifically, first it was proven in \cite{EK}
 that there exists a nonzero
homomorphism $\Xi: L_\rho\to L_\rho\otimes V^*$.
For any $X\in \g$ the element $e^X\in G$ gives a well-defined
automorphism of the vector space $L_\rho$, and one considers
the composite map
 $\,\Xi\ccirc e^X: \,L_\rho\stackrel{e^X}{\too}
L_\rho\stackrel{\Xi}{\too}L_\rho\otimes V^*\,$
as an element of
$\bigl(\End_{_\C}L_\rho\bigr) \otimes V^*$.
Thus, taking the trace: $\End_{_\C}(L_\rho)\to \C$,
on obtains a $V^*$-valued holomorphic function:
$\dis\,X\,\mapsto\,\text{Tr}\big|_{L_\rho}(\Xi\ccirc e^{X}),\,$
 on $\g$. It has been shown in [EK] that,
for any $X\in \h$, one has the identity:
$\text{Tr}\big|_{L_\rho}(\Xi\ccirc e^{X})=\prod_{\alpha>0}\;\bigl(
e^{\langle\alpha,X\rangle/2}
-e^{-\langle\alpha,X\rangle/2}\bigr)=$
{\it Weyl denominator},  where
we have identified $V^*\langle 0\rangle$ with $\C$. This implies that
$$
F:\;\;X\;\mapsto\;\frac{1}{\frac{n(n-1)}{2}!}\cdot
\text{Tr}\big|_{L_\rho}\bigl(\Xi\ccirc X^{n(n-1)/2}
\bigr)\quad,\quad X\in \g,
$$
is a $V^*$-valued function on $\g$ which defines the required
embedding, with property (i).

We note that this remark is in good agreement with
Theorem \ref{k_triv}(viii), since
$L_\rho^*=L_\rho$.

\nopagebreak
{\footnotesize{

}}

\footnotesize{
{\bf P.E.}: Department of Mathematics, Rm 2-165, MIT,
77 Mass. Ave, Cambridge, MA 02139\\
\hphantom{x}\quad\, {\tt etingof@math.mit.edu}

{\bf V.G.}: Department of Mathematics, University of Chicago,
Chicago, IL
60637, USA;\\
\hphantom{x}\quad\, {\tt ginzburg@math.uchicago.edu}}

\end{document}